\documentclass{book}
\usepackage{amsmath,amsfonts,amssymb}
\usepackage[utopia]{mathdesign}

\begin{document}
\title{Introduction to Combinatorial Topology}
\author{Kurt Reidemeister\\
\small{Translated by John Stillwell, with the assistance of Warren Dicks}}

\maketitle

\pagenumbering{roman}
\chapter*{Foreword}
\addcontentsline{toc}{chapter}{Foreword}

The first three chapters of this book deal with infinite groups; the last four
with line and surface complexes and especially with 2-dimensional manifolds.
This choice of material may be justified from several standpoints.

In the first place I was concerned to work out the profound connections
between groups and complexes. The close connection between these two
fields has been known since the basic work of \textsc{Henri Poincar\'e}.
If it has not been plainly evident in the further development of
combinatorial topology, then this is due to the problems where topology
and group theory meet: it seems unfruitful to pursue connections which
primarily permit only the translation of unsolved topological questions into
unsolved group-theoretic questions. Today such thoughts are no longer
justifiable. Since generators and defining relations of subgroups of groups
presented by generators and relations may be determined, group theory 
provides a profitable instrument of computation for topology, with which
several previously inaccessible questions become subject to systematic
investigation. Conversely, complexes perform a valuable service in making
group-theoretic theorems more intuitive and making geometric examination
fruitful for groups; e.g., planar complexes give information about the 
structure of planar discontinuous groups.

Accordingly, I have developed the theory of groups presented by generators
and relations as fully as possible, and have favored those fields of topology
that best demonstrate the connection between groups and complexes, and
which permit new group-theoretic results to be obtained. If, as a result, the
topology of 3-dimensional manifolds is not explicitly mentioned, nevertheless
all methods necessary to attack problems in this area are presented.

I hope that the reader who wants to acquire a few tools for further work
and to learn positive geometric results in a polished logical structure will
also be content with the choice of material. In the last four chapters there
are numerous results immediately accessible to the intuition and which are
derived in a logically transparent way from a few simple axioms about
points, line segments, and surface simplexes. I believe that Chapters 4
to 6 in particular will present no difficulties in comprehension. Perhaps
it is advisable for the reader to begin with these chapters and, where
necessary, to refer back to the prerequisites on groups. But Chapters  1
and 2 will also present little difficulty for most readers. The sections on
groups with operators, as well as Chapters 3 and 7 (Chapter 7 is 
written somewhat tersely), may be omitted the first time.

Apart from a small sketch in \textsc{F. Levi's} ``Geometrische
Konfigurationen'' (Hirzel 1929), combinatorial topology has not previously
been given a comprehensive presentation from the viewpoint adopted
here. So this book relies on the original literature, the results of which
must often be newly derived in the interests of a uniform structure, and
a few lectures of the author. I have received the most fruitful
encouragement from the conversations and work of my late friend
\textsc{Otto Schreier}. Chapters 4 and 6 represent the execution of a
program drawn up by \textsc{M. Dehn} in an address to the
Naturforscherversammlung in Leipzig, 1922. Unfortunately, this address
has not appeared in print. I am pleased to be able to present the
Freiheitssatz formulated by \textsc{Dehn} and recently proved by
\textsc{Magnus}. Chapter 7 contains various previously unpublished
theorems on branched coverings.

In completing the manuscript and overseeing corrections, 
\textsc{L. Goeritz}, \textsc{W. Magnus}, \textsc{E. Podehl} and
\textsc{G. Szeg\"o} have rendered me valuable assistance.
\bigskip

K\"onigsberg, May 1932
\bigskip

\hfill\textbf{Kurt Reidemeister}

\tableofcontents
\newpage
\pagenumbering{arabic}
\chapter*{Introduction}
\addcontentsline{toc}{chapter}{Introduction}

If one wants to investigate surfaces globally, it is often appropriate to
divide them into finitely many surface pieces, bounded by finitely many
curve pieces, and identified along these curve pieces in a certain way.
This is appropriate for the intuition, which can better control pieces in
their totality than the idea of a complicated surface; appropriate
from the standpoint of differential geometry, the methods of which
initially give access only to elementary pieces of surfaces and curves;
and appropriate finally for questions in which the local properties of the
surface play no decisive role, e.g., when it is to be decided whether two
surfaces may be mapped onto each other one-to-one and continuously.
If we call the boundary relations between the surface pieces, curve pieces,
and their endpoints the structure of the decomposition then it may be 
easily shown, e.g., that surfaces possessing decompositions of the same
structure may be mapped onto each other one-to-one and continuously,
and it is possible to prove the even more plausible converse theorem that
surfaces that can be mapped one-to-one can continuously onto each other
possess decomositions with the same structure.

However, as soon as the decomposition of a surface is introduced as a tool,
it becomes an unavoidable question which properties of the surface are
expressed in the structure of a decomposition, or how the structures
of different decompositions of the same surface are related. These
questions are the starting point of combinatorial topology.

In order to answer them, one first asks about the properties of the
elementary surface and curve pieces that bring about the structure. There
are a few simple and intuitive facts; say, the fact that a curve piece always
has two boundary points, and that a surface piece always has a boundary
curve determined by finitely many curve pieces. By formulating these
facts in axioms it is then possible to delimit an area of geometry---the
combinatorial topology of line segments and surface complexes---the
foundation of which is as logically clear as it is intuitively satisfactory.
The results of this theory present a noteworthy contrast to the axioms,
in that they lead very quickly to questions that are as difficult to answer
as they are easy to ask. The classical and popular example of this is the
four color problem. No situation can show more clearly that
mathematics does not live by logic and intuition alone, and that any theory
requires not only unobjectionable axioms but also fruitful ideas.

Here we must mention the group concept, which very soon proves to be 
an Ariadne thread through the labyrinth of complexes, and which we
therefore discuss first, as far as it relevant to combinatorial topology.
In this connection, the methods of determining groups in terms of
generators and relations, and deriving properties of the group from these
relations, stand in the foreground. This field already has a combinatorial
character; the visualization of the group with generators by a line 
segment complex and the visualization of relations in such a group by
elementary surface pieces is just as clear as the attempt to reverse this 
process by converting arbitrary complexes into ``group diagrams.''
This connection is also important for topology. Topological manifolds,
e.g., are in the final analysis none other than representations of groups
in which a few order relations compatible with the group structure are
established between the generators and relations.

It should be pointed out, however, that the sections on line and surface 
complexes in Chapters 4 and 5 may be understood without knowledge of
the first three chapters on groups.

\chapter{Groups}

\section{Definition of the Group Concept}

We begin our exposition with an explanation of the group concept and a
few simple theorems on groups, subgroups, factor groups, and 
isomorphisms of groups.

A class $\mathfrak{F}$ of elements is called a group when each ordered 
pair\footnote{Notice that Reidemeister denotes an ordered pair simply by
juxtaposing its elements. The more usual ordered pair notation $(a,b)$
is used by him as the notation for the \emph{greatest common divisor}
of integers $a,b$ in Section 1.3. (Translator's note.)} 
$F_1 F_2$ is associated with
a certain element $F_{12}$ of $\mathfrak{F}$, in symbols
\[
F_1 F_2 =F_{12}\qquad\text{($F_1$ times $F_2$ equals $F_{12}$),}
\]
and this linking, or multiplication, satisfies the following rules:
\begin{enumerate}
\item[A. 1.]
\emph{If $F_1,F_2,F_3$ are any three elements of $\mathfrak{F}$ and if}
\[
F_1 F_2 = F_{12},\quad F_2 F_3 = F_{23}
\]
\emph{then}
\[
F_{12}F_3=F_1 F_{23}.
\]
We will write this
\begin{equation}
(F_1 F_2)F_3=F_1 (F_2 F_3) \tag{1}
\end{equation}
for short.
\end{enumerate}

A. 1. is called the associative law, and a multiplication that satisfies it is
called associative.

Because of (1) we can write the product of an ordered triple of factors
$F_1,F_2,F_3$ as
\[
F_1 F_2 F_3=(F_1 F_2)F_3=F_1(F_2 F_3)
\]
and, as may be seen by induction, the product of $n$ factors is expressible
analogously as
\[
F_1 F_2 \cdots F_n,
\]
independent of bracketing.

\begin{enumerate}
\item[A. 2.]
\emph{There is an element $E$ in $\mathfrak{F}$ for which}
\[
EF=FE=F
\]
\emph{for any $F$ in $\mathfrak{F}$.} Such an element is known
as an identity element.
\end{enumerate}

There can be only one such element; for, if $E^{*}$ were a second
element of this kind, we must have $EE^{*}=E$ on the one hand,
and $EE^{*}=E^{*}$ on the other, so that $E=E^{*}$.

\begin{enumerate}
\item[A. 3.]
\emph{For each element $F$ of $\mathfrak{F}$ there is an element
$X$ for which}
\begin{equation}
FX=E. \tag{2}
\end{equation}
$X$ is called an element inverse to $F$.
\end{enumerate}

For an element $X$ inverse to $F$ we also have
\begin{equation}
XF=E. \tag{3}
\end{equation}
This is because there is an element $Y$ for which
\begin{equation}
XY=E. \tag{4}
\end{equation}
Then multiplication by $F$ gives, on the one hand, that
\[
F(XY)=FE=F.
\]
While, on the other hand, it follows from (1) and (2) that
\[
F(XY)=(FX)Y=EY=Y
\]
and hence $Y=F$, so the asserted equation (3) follows from (4).
It follows further that there is only one inverse element. On the one hand
it follows from $FX_1=E$ and $FX_2=E$ that
\[
X_2(FX_1)=X_2 E=X_2.
\]
And since $X_2 F=E$ by (3) we have, on the other hand,
\[
X_2(FX_1)=(X_2 F)X_1=EX_1=X_1,
\]
and thus $X_1=X_2$.

If $X$ is the element inverse to $F$, then $F$ is the element inverse to $X$.
We denote the element inverse to $F$ by $F^{-1}$, so
\[
(F^{-1})^{-1}=F.
\]

This symbolism may be extended by the following convention. By $F^{1}$
we mean $F$ itself. $F^n$ is defined for positive integers $n>1$ by
induction as
\[
F^n=(F^{n-1})F.
\]
By $F^0$ we mean the identity element $E$, by $F^{-n}$ ($n>0$) we mean
$(F^{-1})^n$; it then follows from the associative law and the properties of 
the inverse that
\[
F^n F^m=F^{m+n},\quad (F^n)^m=F^{nm}
\]
for arbitrary $m,n$.

$F^{-n}$ is therefore the element inverse to $F^n$. We call $F^n$ the
$n$th power of the element $F$. If $F_1$ and $F_2$ are two different
elements of $\mathfrak{F}$ then we can construct the elements
\begin{equation}
F^{n_{11}}_1 F^{n_{12}}_2 F^{n_{21}}_1 F^{n_{22}}_2\cdots
F^{n_{r1}}_1 F^{n_{r2}}_2 \tag{5}
\end{equation}
by iterated multiplication, which may be called \emph{power products}
of $F_1$ and $F_2$. We call
\begin{equation}
F^{-n_{r2}}_2 F^{-n_{r2}}_1\cdots  F^{-n_{22}}_2 F^{-n_{21}}_1
F^{-n_{12}}_2 F^{-n_{11}}_1 \tag{6}
\end{equation}
the power product formally inverse to (5), because one computes that the
product of (5) and (6) is equal to $E$. The power products of elements
$F_1,F_2,\ldots,F_k$ may be constructed analogously.

If $F_1 F_2=F_2 F_1$ for any two elements $F_1$ and $F_2$ of
$\mathfrak{F}$ then the group $\mathfrak{F}$ is called \emph{commutative}.
If each element of $\mathfrak{F}$ may be written as a power of a fixed
element $F$ then $\mathfrak{F}$ is called a \emph{cyclic} group with
generator $F$. Since
\[
F^n F^m = F^{n+m} = F^{m+n} = F^m F^n,
\]
a cyclic group is commutative.

With a view towards our objectives we will always assume that the groups 
under consideration have only a denumerable number of elements. The
number of elements of a group is called its \emph{order}.

\section{Cyclic Groups}

Cylic groups are easy to exhibit. E.g., the positive and negative integers
and zero constitute such a group when one takes addition as the group
operation: all numbers may then be regarded as power products of
$+1$ and $-1$.

The residue classes with respect to a modulus, under addition,
constitute another example. If $m$ is any positive integer we call two
integers $n_1$ and $n_2$ congruent modulo $m$, denoted
$n_1\equiv n_2$ (mod $m$), if the difference $n_1-n_2$ is
divisible by $m$. If
\[
n_1\equiv n_2\text{ (mod $m$)}\quad\text{and}\quad
n_2\equiv n_3\text{ (mod $m$)}
\]
then also $n_1\equiv n_3$ (mod $m$).

We now understand the residue class $[n]$ to be all numbers congruent
to $n$ (mod $m$). Obviously, the class $[n]$ is identical with the
class $[n']$,
\[
[n]=[n'],\quad\text{if}\quad n\equiv n'\text{ (mod $m$)}
\]
and conversely. For each residue class there is exactly one representative
$r$ satisfying the inequality
\[
0\le r<m.
\]
Thus there are $m$ different residue classes.

We define an operation on these residue classes, denoted by the symbol 
$+$ and called addition, by
\[
[n_1]+[n_2]=[n_1+n_2].
\]
This definition is not contradictory. Namely, if $[n'_i]=[n_i]$ then
$n'_i\equiv n_i$ (mod $m$)  ($i=1,2$) and hence
\[
n'_1+n'_2\equiv n_1+n_2 \text{ (mod $m$)},
\]
as one easily verifies, so
\[
[n_1+n_2]=[n'_1+n'_2].
\]
This operation satisfies the group axioms. It is associative because
addition of whole numbers is:
\begin{align*}
([n_1]+[n_2])+[n_3]&=[(n_1+n_2)+n_3]\\
                              &=[n_1+(n_2+n_3)]\\
                              &=[n_1]+([n_2]+[n_3]).
\end{align*}
$[0]$ is the identity element and $[-n]$ is the element inverse to $[n]$.
The group is cyclic as well, because all $m$ residue classes result from
iterated addition of the residue class $[1]$.

The multiplication property of cyclic groups is easily seen to lead to the
two examples given. If all elements of $\mathfrak{F}$ are powers 
$F^n$ (where $F$ is the generator of $\mathfrak{F}$) then either
$F^n$ is different from $F^{n'}$ as long as $n$ is different from 
$n'$---in which case multiplication of elements of $\mathfrak{F}$
reduces to addition of whole numbers---or else there are two different
exponents $n>n'$ which yield equal elements of $\mathfrak{F}$. Then
\[
F^n F^{-n'}=F^{n-n'}=E
\]
and there is a smallest positive exponent $f$ for which
\[
F^f=E.
\]
In this case $F^n=F^{n'}$ if and only if $n\equiv n'$ (mod $f$).
In fact
\[
F^{kf}=(F^f)^k=E^k=E
\]
and so $F^{n+kf}=F^n$. Conversely, if $F^n=F^{n'}$ ($n>n'$)
then $F^{n-n'}=E$ Here we must have $n-n'\ge f$, and if
\[
n-n'=kf+r\quad(0\le r<f)
\]
then $F^r=1$, so $r=0$. In this case the cyclic group has order $f$.

\section{Multiplication of Residue Classes}

We can construct groups from the residue classes mod $m$ in another way
by taking the group operation to be multiplication of residue classes.
The product of $[n_1]$ and $[n_2]$ is defined by
\[
[n_1][n_2]=[n_1 n_2].
\]
The product residue class is uniquely determined, because it follows from
$[n'_i]=[n_i]$, or $n'_i\equiv n_i$ (mod $m$) ($i=1,2$), that
$n'_1 n'_1\equiv n_1 n_2$ (mod $m$) also, and hence
$[n'_1 n'_2]=[n_1 n_2]$. This multiplication is associative and commutative,
because multiplication of integers is. [1] is the identity element because
\[
[n][1]=[1][n]=[n].
\]
On the other hand, there is not a multiplicative inverse for each residue
class. Namely, if $n$ is a number with which the modulus $m$ has a
greatest common divisor $d=(n,m)\ne 1$ then all numbers of the 
residue class $[n]$, the numbers $n+km$, have the same greatest common 
divisor in common with $m$:
\[
(n+km,m)=(n,m)=d.
\]
Obviously all numbers in the residue class $[n][n']=[nn']$ then have a
greatest divisor in common with $m$ that is divisible by $d$ and hence
$\ne 1$. So it cannot be the case that $[n][n']=[1]$.

However, it may be shown that \emph{the residue classes $[n]$ for which
the greatest common divisor}
\[
d=(n,m)=1
\]
\emph{constitute a group under multiplication.} The product of two
residue classes relatively prime to $m$ is again relatively prime to $m$.
Further, if
\[
[r_1],\quad [r_2],\quad \ldots, \quad [r_{\overline m}]
\]
is the totality of these residue classes, and $[r]$ is any one of them, then
\[
[rr_1],\quad [rr_2],\quad \ldots, \quad [rr_{\overline m}]
\]
is again the totality of residue classes relatively prime to $m$. For if
$[rr_i]=[rr_k]$ then $rr_i\equiv rr_k$ (mod $m$), so $r(r_i-r_k)$
must be divisible by $m$, or $r_i\equiv r_k$ (mod $m$), and
\[
i=k.
\]
The $\overline{m}$ residue classes $[rr_i]$ are therefore all different,
and hence they exhaust the residue classes relatively prime to $m$.

Consequently, the residue class [1] appears among the $[rr_i]$, and
if
\[
[rr']=[1]
\]
then $[r']$ is $[r]^{-1}$, the residue class inverse to $[r]$. The assertion
now follows easily.

With a prime number $p$ as modulus,
\[
m=p,
\]
the residue classes $[1], [2], [3],\ldots,[p-1]$ are relatively prime to $p$
and so they constitute a group under the operation of residue class
multiplication.

An interesting result on cyclic groups follows easily from our development.
\emph{If $\mathfrak{S}$ is a cyclic group with generator $S$ and finite 
order $m$, then each element $S^k$ with}
\[
(k,m)=1
\]
\emph{is a generator of the group}; on the other hand, an element $S^l$
with $(l,m)\ne 1$ does not generate the group. This is because the
$S^{ki}$ ($i=0,1,\ldots,m-1$) are all different from each other, for it
follows from $S^{ki_1}=S^{ki_2}$ that $ki_1\equiv ki_2$ (mod $m$),
so $k(i_1-i_2)$ must be divisible by $m$ and hence $i_1=i_2$. The
assertion about $S^k$ then follows. On the other hand, the elements
$S^{li}$ correspond only to residue classes that have a common
divisor with $m$ greater than 1, and hence certainly not to all elements
of the group.

The existence of the inverse residue class $[r]^{-1}$ can also be
expressed as follows: if the greatest common divisor of $n$ and $m$,
$(n,m)=1$, then the congruence
\[
nx\equiv 1\text{ (mod $m$)}
\]
always has solutions, and indeed the set of solutions constitutes a
residue class; namely, the residue class $[n]^{-1}$. It follows that a
congruence
\begin{equation}
n_1 x\equiv n_2 \text{ (mod $m$)} \tag{1}
\end{equation}
is always uniquely solvable if $(n_1,m)=1$, because an equivalent to (1)
is
\[
[n_1][x]=[n_2],
\]
and from this we have $[x]=[n_1]^{-1}[n_2]$. I.e., there are 
solutions of the congruence (1), and the set of solutions is the
congruence class $[n_1]^{-1}[n_2]$.

\section{Groups of Transformations}

One can take, as the elements of a group, the transformations of any
domain $\mathfrak{X}$ of objects, with the composition of these
transformations as the group operation. Let
\[
\overline{x}=F(x)
\]
be a one-to-one onto transformation of $\mathfrak{X}$; i.e., each object
$x$ corresponds to a well-defined object $\overline{x}$; all objects
in $\mathfrak{X}$ appear among the $\overline{x}$; if $x_1$ and
$x_2$ are different then so are $\overline{x_1}$ and
$\overline{x_2}$. The inverse transformation $F^{-1}$ is defined
by
\[
F^{-1}(\overline{x})=x,
\]
and it is obviously also one-to-one and onto. If $F_1$ and $F_2$ are two 
one-to-one  onto transformations, and if
\[
\overline{x}=F_1(x),\quad \overline{\overline{x}}=F_2(\overline{x}),
\]
then the correspondence $\overline{\overline{x}}=F_{21}(x)$
is again one-to-one. We call $F_{21}$ the product of $F_2$ and
$F_1$ and write $F_{21}=F_2 F_1$. This operation satisfies the
associative law. namely, if
\[
\overline{x}=F_1(x),\quad
x'=F_2(\overline{x}),\quad
x^{*}=F_3(x')
\]
then
\[
x'=F_{21}(x),\quad x^{*}=F_{32}(\overline{x})
\]
and hence
\[
x^{*}=F_3(F_{21}(x))=F_{32}(F_1(x))
\]
is the same transformation of $x$.

Now \emph{if $\mathfrak{F}$ is a family\footnote{Strictly, a
\emph{nonempty} family, but Reidemeister always assumes nonempty
sets. (Translator's note.)} of such one-to-one 
transformations and if $\mathfrak{F}$ contains, along with each member
$F$, the inverse $F^{-1}$ and, along with each pair $F_1,F_2$, their 
product $F_{21}$, then $\mathfrak{F}$ is obviously a group.}

If the domain $\mathfrak{X}$ consists of finitely many objects
\[
x_1,\quad x_2,\quad \ldots,\quad x_m
\]
then a one-to-one onto transformation
\[
x_{n_i}=F(x_i)\qquad (i=1,2,\ldots,m)
\]
is called a permutation of the objects $x_i$. 

If suitable transformations in a group $\mathfrak{F}$ will send any $x$
to any other, then the transformation group is called \emph{transitive}.

As an example of a transformation group we introduce the
\emph{modular} group. The domain $\mathfrak{X}$ consists of the
complex numbers
\[
x=\xi_1+i\xi_2\quad\text{with}\quad \xi_2>0
\]
and the transformations are
\begin{equation}
x'=\frac{ax+b}{cx+d}, \tag{1}
\end{equation}
where $a,b,c,d$ are integers with determinant
\[
ad-bc=1.
\]
If
\begin{equation}
x''=\frac{a'x'+b'}{c'x'+d'} \tag{2}
\end{equation}
is a second such transformation, then the composite transformation
\begin{equation}
x''=\frac{a''x+b''}{c''x+d''} \tag{3}
\end{equation}
has 
\begin{align*}
&a''=a'a+b'c,\qquad c''=c'a+d'c,\\
&b''=a'b+b'd,\qquad d''=c'b+d'd. \tag{4}
\end{align*}
The determinant $a''d''-b''c''$ is equal to the product
\[
(a'd'-b'c')(ad-bc)=1.
\]
\begin{equation}
x=\frac{dx'-b}{-cx'+a} \tag{5}
\end{equation}
is the transformation inverse to (1). If (1) is the identity transformation
then we must have
\[
cx^2-(d-a)x-b=0
\]
for all $x$. It follows that $b=c=d-a=0$ and, because $ad=1$, either
$a=d=+1$ or $a=d=-1$. One sees from this that two transformations
defined by the formula (1) are identical if and only if their coefficients
$a,b,c,d$ are respectively equal or else respectively of the same magnitude
but oppositely signed.

\section{Subgroups}

In order to penetrate more deeply into the structure of a group
$\mathfrak{F}$ one considers its subgroups, i.e., groups $\mathfrak{f}$
whose elements all belong to $\mathfrak{F}$. In this context the group
operation for the elements of $\mathfrak{f}$ is the same as that for
$\mathfrak{F}$. Thus $\mathfrak{F}$ itself is a subgroup of
$\mathfrak{F}$. Each subgroup different from $\mathfrak{F}$ itself
is called a \emph{proper} subgroup. One can also characterize subgroups
as follows: a collection $\mathfrak{f}$ of elements of $\mathfrak{F}$ is 
called a subgroup when
\[
F_1 F_2=F_{12}
\]
belongs to $\mathfrak{f}$ along with $F_1$ and $F_2$ and, along with
each element $F$, its inverse $F^{-1}$ also belongs to $\mathfrak{f}$.
Obviously the identity element $E=FF^{-1}$ then belongs to $\mathfrak{F}$
and, since the product of elements in $\mathfrak{f}$ is naturally
associative, $\mathfrak{f}$ is in fact a group and hence a subgroup of
$\mathfrak{F}$ according to the first definition.

The elements representable as powers of an element $F$ constitute a
subgroup, because the product of two powers of $F$ and the inverse 
of a power are again powers of $F$. One calls the order of this subgroup
the \emph{order of the element $F$}. In the example of the whole 
numbers these groups consist of all the elements divisible by a given
number.

One concludes similarly that the power products (5) in Section 1.1, of
two or an arbitrary finite or infinite set of elements constitute a group.
Thus if $\mathfrak{m}$ is any set of elements of $\mathfrak{F}$ one
may speak of the subgroup of $\mathfrak{F}$ determined or
\emph{generated} by $\mathfrak{m}$. It is just the set of all power
products of elements of $\mathfrak{m}$.

An important subgroup defined in this way is the \emph{commutator
group} $\mathfrak{K}_1$ of $\mathfrak{F}$. By the commutator of
$F_1$ and $F_2$ we mean the element
\[
K=F_1 F_2 F^{-1}_1 F^{-1}_2.
\]
Now if $\mathfrak{k}_1$ is the set of all commutator elements of
$\mathfrak{F}$, then $\mathfrak{K}_1$ is the group generated by
$\mathfrak{k}_1$. By the commutators of second order, $\mathfrak{k}_2$,
we mean all commutators of an element of $\mathfrak{k}_1$ with
an element of $\mathfrak{F}$, and by $\mathfrak{K}_2$ the group
so generated, the \emph{second commutator group}. Commutator
groups of higher order may be defined by induction.

The elements $F$ that commute with a fixed element $F_0$, i.e.,
those for which $F_0 F=FF_0$, constitute a group. For if
\[
F_0 F_1=F_1 F_0\quad\text{and}\quad F_0 F_2=F_2 F_0
\]
then also
\[
F_0 F_1 F_2 = F_1 F_0 F_2 = F_1 F_2 F_0
\]
and
\[
F_0 F^{-1}_1=F^{-1}_1 F_1 F_0 F^{-1}_1
=F^{-1}_1 F_0 F_1 F^{-1}_1=F^{-1}_1 F_0.
\]
Similarly, one concludes that the set $\mathfrak{Z}$ of those elements
of $\mathfrak{F}$ that commute with all the elements of $\mathfrak{F}$
constitute a subgroup. It is called the \emph{center} of $\mathfrak{F}$.

It is also easy to construct subgroups of a group of transformations.
The set of all transformations that leave a given element $x_0$ fixed,
e.g., constitute a subgroup. For if $F_1(x_0)=x_0$ and $F_2(x_0)=x_0$
then also
\[
F_2(F_1(x_0))=F_{21}(x_0)=x_0
\]
and
\[
F^{-1}_1(x_0)=x_0.
\]
Similarly, the transformations that leave several points $x$ fixed
constitute a subgroup.

If $\mathfrak{f}_1$ and $\mathfrak{f}_2$ are subgroups of
$\mathfrak{F}$, then so is the collection $\mathfrak{f}_{12}$ of all
elements that belong to both $\mathfrak{f}_1$ and $\mathfrak{f}_2$.
Namely, if $F_1$ and $F_2$ are elements that belong to both
$\mathfrak{f}_1$ and $\mathfrak{f}_2$, then $F_1 F_2=F_{12}$
and $F^{-1}_1$ also belong to both $\mathfrak{f}_1$ and 
$\mathfrak{f}_2$. One concludes similarly that the intersection of 
arbitrarily many subgroups is also a subgroup.

\section{Conjugate Subgroups}

If $\mathfrak{f}$ is a subgroup of $\mathfrak{F}$, $\overline{F}$
runs through all the elements of $\mathfrak{f}$, and $F_0$ is a
fixed element of $\mathfrak{F}$, then the elements
\[
F_0\overline{F}F^{-1}_0=\overline{F}'
\]
run through a collection $\mathfrak{f}'$ of elements that also
constitute a group. For if
\[
\overline{F}_1\overline{F}_2=\overline{F}_{12}
\]
then
\[
\overline{F}'_1\overline{F}'_2
=F_0 \overline{F}_1 F^{-1}_0 F_0 \overline{F}_2 F^{-1}_0
=F_0\overline{F}_{12}F^{-1}_0
=\overline{F}'_{12}
\]
and $F_0\overline{F}^{-1} F^{-1}_0$ is the element inverse to
$F_0\overline{F} F^{-1}_0$. We also write 
$F_0\mathfrak{f}F^{-1}_0$ for $\mathfrak{f}'$ and call it a
\emph{subgroup conjugate to $\mathfrak{f}$}.

If $\mathfrak{f}'$ is conjugate to $\mathfrak{f}$, then
$\mathfrak{f}$ is also conjugate to $\mathfrak{f}'$. For in fact
$F^{-1}_0\mathfrak{f}' F_0$ is identical with $\mathfrak{f}$.
If $\mathfrak{f}_u$ is a proper subgroup of $\mathfrak{f}$,
then $F_0\mathfrak{f}_u F^{-1}_0$ is a proper subgroup of
$F_0\mathfrak{f}F^{-1}_0$.

If $\mathfrak{f}'$ and $\mathfrak{f}''$ are two subgroups
conjugate to $\mathfrak{f}$, and if
\[
\mathfrak{f}'=F_0\mathfrak{f} F^{-1}_0\quad\text{and}\quad
\mathfrak{f}''=F^{-1}_0\mathfrak{f} F_0,
\]
and if both $\mathfrak{f}'$ and $\mathfrak{f}''$ are contained in
$\mathfrak{f}$, then $\mathfrak{f},\mathfrak{f}',\mathfrak{f}''$
are identical. E.g., if $\mathfrak{f}'$ were a proper subgroup of
$\mathfrak{f}$ then we would also have
\[
F^{-1}_0\mathfrak{f}' F_0=\mathfrak{f}
\]
a proper subgroup of
\[
F^{-1}_0 \mathfrak{f} F_0=\mathfrak{f}'';
\]
thus $\mathfrak{f}=\mathfrak{f}''$ and consequently
$\mathfrak{f}'=\mathfrak{f}$ also.

If $F$ runs through all elements of $\mathfrak{F}$, then
$F\mathfrak{f}F^{-1}$ runs through a \emph{class of conjugate
subgroups}. Such a class is determined by any one of its elements.

Naturally, it can happen that formally different conjugate subgroups
are identical with each other. E.g., $F_0\mathfrak{f} F^{-1}_0$ is
identical with $\mathfrak{f}$ when $F_0$ belongs to $\mathfrak{f}$.
But the groups can also coincide when $F_0$ does not come from
$\mathfrak{f}$. In particular, all conjugate subgroups can be
identical with each other. Then $\mathfrak{f}$ is called an
\emph{invariant} subgroup.\footnote{Of course, this is what we now
call a \emph{normal} subgroup. However, the word ``invariant'' is
reasonable (if understood to mean ``invariant under conjugation'')
and the word ``normal'' is overused in mathematics. So I have allowed
``invariant'' to stand. (Translator's note.)} 
Transformation groups yield examples
of conjugate subgroups. If $\mathfrak{f}_x$ is the group of
transformations that leave the element $x$ fixed, and if
$F_0(x)=\overline{x}$, then
\[
F_0 \mathfrak{f} F^{-1}_0=\mathfrak{f}_{\overline{x}}
\]
is the group of transformations that leave $\overline{x}$ fixed.
In fact, each transformation in $\mathfrak{f}_{\overline{x}}$
carries the object $\overline{x}$ back to itself. Conversely, if
$F(\overline{x})=\overline{x}$ then
\[
F^{-1}_0 F F_0(x)=x
\]
is a transformation $F'$ in $\mathfrak{f}_{x}$, and then
\[
F=F_0 F' F^{-1}_0
\]
belongs to $\mathfrak{f}_{\overline{x}}$. If $\mathfrak{F}$ is a
transitive group, then the class of subgroups that leave some $x$
fixed constitute a class of conjugate subgroups of $\mathfrak{F}$.
If, for some $x$, $\mathfrak{f}_x$ consists only of the identity,
then this is the case for all $x$, and the transformation group is
called \emph{simply transitive}.

Now for a few examples of invariant subgroups. The center
$\mathfrak{Z}$ of the group $\mathfrak{F}$ is obviously an
invariant subgroup; for if $F$ is an arbitrary element of
$\mathfrak{F}$ and $Z$ is any element of the center, then always
\[
FZF^{-1}=Z.
\]

The commutator group $\mathfrak{K}_1$ is also an invariant
subgroup. Namely, if $K$ is the commutator
\[
K=F_1 F_2 F^{-1}_1 F^{-1}_2
\]
then
\[
FKF^{-1}=F F_1 F^{-1}
                F F_2 F^{-1}
                F F^{-1}_1 F^{-1}
                F F^{-1}_2 F^{-1}
\]
is also a commutator. Consequently, for each product of
commutators $K_1 K_2 \cdots K_r$ we also have
\[
FK_1 K_2\cdots K_r F^{-1}
=FK_1 F^{-1} FK_2 F^{-1}\cdots FK_r F^{-1}
\]
belonging to $\mathfrak{K}_1$. The higher commutator groups 
are also invariant subgroups.

The intersection $\mathfrak{D}$ of the subgroups in a class of
conjugate subgroups is an invariant subgroup. Namely, if $F^{*}$
is an element that appears in all groups $F\mathfrak{f}F^{-1}$
then the element $F_0 F^{*} F^{-1}_0$ appears in all groups
$F_0(F\mathfrak{f} F^{-1})F^{-1}_0$ and, since $F_0 F=F'$
likewise runs through all elements of $\mathfrak{F}$ when $F$
does, this means that $F_0 F^{*} F^{-1}_0$ also appears in
all groups $F\mathfrak{f} F^{-1}$.

\section{Congruence Subgroups of the Modular Group}

In the case of the modular group defined in Section 1.4, subgroups
may be easily defined in terms of number-theoretic properties of
the coefficients. E.g., the modular transformations with 
$c\equiv 0$ (mod $n$) constitute a a group $\mathfrak{U}_n$.
Namely, if
\[
x'=\frac{ax+b}{cx+d}\quad\text{and}\quad
x''=\frac{a'x'+b'}{c'x'+d'}
\]
are two transformations with $c\equiv c'\equiv 0$ (mod $n$),
then for the composite transformation
\[
x''=\frac{a''x+b''}{c''x+d''} 
\]
we have $c''=c'a+d'c$, which is obviously divisible by $n$, so
$c''\equiv 0$ (mod $n$) as well. If
\[
x=\frac{a'x'+b'}{c'x'+d'}
\]
is the transformation inverse to $x'=\frac{ax+b}{cx+d}$, then
by (5) of Section 1.4, $c'=-c$, so likewise $c'\equiv 0$ (mod $n$).

Another \emph{``congruence subgroup''} $\mathfrak{U}'_n$
consists of all the modular substitutions with $b\equiv 0$ (mod $n$).
One can either verify the group property of $\mathfrak{U}'_n$
directly, or else confirm that $\mathfrak{U}'_n$ is identical with a
group conjugate to $\mathfrak{U}_n$. Namely, if one denotes the
transformation
\[
x'=-\frac{1}{x}
\]
by $S$, then the group
\[
S\mathfrak{U}_n S^{-1}=S\mathfrak{U}_n S
\]
consists of the transformations
\begin{equation}
x'=-\frac{c\left(-\frac{1}{x} \right)+d}{a\left(-\frac{1}{x} \right)+b}
=\frac{dx-c}{-bx+a} \tag{1}
\end{equation}
with $c\equiv 0$ (mod $n$). These transformations obviously all belong to
$\mathfrak{U}'_n$, and since one confirms similarly that
$S\mathfrak{U}'_n S^{-1}$ is contained in $\mathfrak{U}_n$, we have
$\mathfrak{U}'_n=S\mathfrak{U}_n S$.

The modular transformations with
\[
a\equiv d\equiv 1\text{ (mod $n$)},\qquad
 b\equiv c\equiv 0\text{ (mod $n$)}
\]
likwise constitute a group $\mathfrak{I}_n$. For, along with each
transformation, the inverse also belongs to $\mathfrak{I}_n$, by
(5) in Section 1.4, and if (1) and (2) in Section 1.4 are two
transformations in $\mathfrak{I}_n$, then the composite transformation
(3) in Section 1.4 satisfies:
\[
\begin{array}{ll}
a''=a'a+b'c\equiv a'a\equiv 1  &\text{(mod $n$)}\\
d''=c'b+d'd\equiv d'd\equiv 1 &\text{(mod $n$)}\\
b''=a'b+b'd\equiv 0               &\text{(mod $n$)}\\
c''=c'a+d'c\equiv 0                &\text{(mod $n$)}.
\end{array}
\]
In contrast to the $\mathfrak{U}_n$, the $\mathfrak{I}_n$ are
invariant subgroups.
\[
\mathfrak{I}'_n=S\mathfrak{I}_n S^{-1}=S^{-1}\mathfrak{I}_n S
\]
is contained in $\mathfrak{I}_n$ by formula (1) and the assumptions
\[
a\equiv d\equiv 1\text{ (mod $n$)},\qquad
 b\equiv c\equiv 0\text{ (mod $n$)}.
\]
In addition, we construct
\[
T^{\eta}\mathfrak{I}_n T^{-\eta}\qquad (\eta=\pm 1),
\]
where $T$ denotes the transformation
\[
x'=x+1.
\]
This group consists of the transformations
\[x'=\frac{(a+c\eta)x+(d-a)\eta+b-c}{cx-c\eta+d}
\]
But
\begin{align*}
a+c\eta\equiv -c\eta+d\equiv 1 &\quad\text{(mod $n$)}\\
(d-a)\eta+b-c\equiv c\equiv 0 &\quad\text{(mod $n$)},
\end{align*}
so these transformations belong to $\mathfrak{I}_n$. Consequently,
the groups 
\[
T^{\eta}\mathfrak{I}_n T^{-\eta}\quad\text{and}\quad
S\mathfrak{I}_n S
\]
are identical with $\mathfrak{I}_n$ by the third paragraph of Section 1.6.
Since we will see in Section 2.9 that all modular substitutions may be
written as power products of $S$ and $T$, the invariance of
$\mathfrak{I}_n$ follows.

\section{Residue Classes modulo Subgroups}

The method of constructing residue classes in the domain of integers
may be extended to arbitrary groups. If $\mathfrak{F}$ is a group
and $\mathfrak{f}$ is a subgroup, we call two elements $F_1$ and 
$F_2$ right congruent modulo $\mathfrak{f}$, written
\[
F_1 \equiv_r F_2\quad\text{(mod $\mathfrak{f}$)},
\]
if there is an element $F_f$ of $\mathfrak{f}$ such that $F_1=F_2 F_f$.
If
\[
F_1\equiv_r F_2\quad\text{(mod $\mathfrak{f}$)}
\qquad\text{then also}\qquad
F_2\equiv_r F_1\quad\text{(mod $\mathfrak{f}$)},
\]
because in fact $F_2=F_1 F^{-1}_f$. And if
\[
F_1\equiv_r F_2\quad\text{(mod $\mathfrak{f}$)}
\qquad\text{and}\qquad
F_2\equiv_r F_3\quad\text{(mod $\mathfrak{f}$)}
\]
then also
\[
F_1 \equiv_r F_3\quad\text{(mod $\mathfrak{f}$)},
\]
because if $F_1=F_2 F_f$ and $F_2=F_3 F_{f'}$ then
$F_1=F_3 F_{f'} F_f=F_3 F_{f''}$. We understand the right-sided
residue class\footnote{This residue class if of course what we now 
call a \emph{coset} of the subgroup $\mathfrak{f}$. However, I have
thought it best to retain the term ``residue class '' (in German,
Restklasse) to reflect Reidemeister's
view of cosets as generalizations of residue classes in number theory.
(Translator's note.)}
determined by $F_1$ to be the collection $F_1\mathfrak{f}$
of elements right congruent to $F_1$ modulo $\mathfrak{f}$. If
$F_1\equiv_r F_2$ (mod $\mathfrak{f}$) then
\[
F_1\mathfrak{f}=F_2\mathfrak{f},
\]
and conversely. The residue class is determined by any one of its elements,
or in other words: two right-sided residue classes modulo $\mathfrak{f}$
that have a common element are identical. Thus the elements of
$\mathfrak{F}$ are partitioned by a subgroup into disjoint residue classes.
A \emph{system of representatives} of these residue classes is a set
$\mathfrak{r}$ of elements $R$ of the following kind: if $F$ is any
element of $\mathfrak{F}$ then there is an element $R$ of
$\mathfrak{r}$ such that
\[
F\equiv_r R\quad\text{(mod $\mathfrak{f}$)},
\]
and if $R_1$ and $R_2$ are two elements of $\mathfrak{r}$ then $R_1$
is not right congruent to $R_2$. The classes $R\mathfrak{f}$ then yield
all the residue classes, without repetitions, as $R$ runs through
$\mathfrak{r}$.

We can define left congruence analogously to right congruence. We say
\[
F_1\equiv_l F_2\quad\text{(mod $\mathfrak{f}$)}
\]
if $F_1=F_f F_2$, where $F_f$ belongs to $\mathfrak{f}$. We
analogously define the left-sided residue classes and a system of
representatives $\mathfrak{r}'$ of them. We will show: \emph{if
$\mathfrak{r}$ is a full system of representatives for the
right-sided residue classes, then $\mathfrak{r}^{-1}$, the collection
of inverses to the $R$ in $\mathfrak{r}$, is a left-sided system of
representatives}. Namely, if $F$ is an arbitrary element and
$F^{-1}=RF_f$, then $F=F^{-1}_f R^{-1}$, so
\[
F\equiv_l R^{-1}\quad\text{(mod $\mathfrak{f}$)}.
\]
Further, if
\[
R^{-1}_1\equiv_l R^{-1}_2,\quad\text{so}\quad
R^{-1}_1=F_f R^{-1}_2,
\]
then $R_2\equiv_r R_1$ (mod $\mathfrak{f}$).

If the number of right-sided residue classes modulo $\mathfrak{f}$
is a finite number $n$, then this shows that the number of left-sided
residue classes is also $n$.

$R\mathfrak{f}R^{-1}$ runs through all subgroups conjugate to
$\mathfrak{f}$ as $R$ runs through the set $\mathfrak{r}$. For if
$F_0$ is an element of $\mathfrak{F}$ and $F_0=R_0 F_f$ then
\[
F_0\mathfrak{f} F^{-1}_0 = 
R_0 F_f \mathfrak{f} F^{-1}_f R^{-1}_0=
R_0 \mathfrak{f} R^{-1}_0.
\]

\section{Residue Classes modulo Congruence Subgroups of the
Modular Group}

As an example, we determine the left-sided residue classes modulo
the subgroup $\mathfrak{U}_n$ of the modular group defined in
Section 1.7, for $n$ equal to a prime number $p$.

If we denote the transformation
\[
x'=\frac{-1}{x+k}
\]
by $G_k$, then the identity transformation $E$ and the $G_k$ 
($k=0,1,\ldots,p-1$) form a full system of representatives for the
$\mathfrak{U}_p M$, where $M$ denotes an arbitrary modular
substitution. Namely, for each $M$ that does not belong to
$\mathfrak{U}_p$ (that is, $c\not\equiv 0$ (mod $p$)) we may
determine a substitution $U$ in $\mathfrak{U}_p$
\[
x'=\frac{a'x+b'}{c'x+d'},
\]
so $c'\equiv 0$ (mod $p$), and a $G_k$ so that
\[
UG_k=M.
\]
If the substitution $MG^{-1}_k$ is given by
\[
x''=\frac{a''x+b''}{c''x+d''},
\]
then by (4) of Section 1.4
\[
c''=kc-d.
\]
Thus $k$ must satisfy the congruence
\[
kc-d\equiv 0\quad\text{(mod $p$)}
\]
and then by Section 1.3 the residue class $[k]$ is uniquely
determined, because $c$ is assumed to be relatively prime to $p$
and hence $k$ is also, because
\[
0\le k<p.
\]
But then the coefficients of $U$ are likewise determined by
\[
U=MG^{-1}_k.
\]

\section{Factor Groups}

The definition of addition of residue classes can also be generalized
under the hypothesis that the modulus $\mathfrak{f}$ is an invariant
subgroup of $\mathfrak{F}$. One sees first of all that: \emph{if
$\mathfrak{f}$ is an invariant subgroup and if}
\[
F_1\equiv_r F_2\quad\text{(mod $\mathfrak{f}$)}
\]
\emph{then also}
\[
F_1\equiv_l F_2\quad\text{(mod $\mathfrak{f}$)}
\]
\emph{and conversely}. For if
\[
F_1=F_2 F_f=F_2 F_f F^{-1}_2 F_2
\]
and $\mathfrak{f}$ is invariant, then $F_1=F'_f F_2$, because
$F_2 \mathfrak{f} F^{-1}_2=\mathfrak{f}$ and so 
$F_2 F_f F^{-1}_2$ is itself an element of $\mathfrak{f}$. Thus
we can simply speak of residue classes modulo $\mathfrak{f}$.
Now we can show further that if
\[
F_1\equiv F_2\quad\text{(mod $\mathfrak{f}$)}\qquad
\text{and}\qquad
F'_1\equiv F'_2\quad\text{(mod $\mathfrak{f}$)}
\]
then $F_1 F'_1\equiv F_2 F'_2$ (mod $\mathfrak{f}$) as well.
For if $F_1=F_2 F_f$ and $F'_1=F'_2 F'_f$ then
\[
F_1 F'_1=F_2 F'_2 F'^{-1}_2 F_f F'_2 F'_f=F_2 F'_2 F''_f F'_f.
\]

If now $F_1\mathfrak{f}$ and $F_2\mathfrak{f}$ are two residue
classes, we define their product to be the residue class 
$F_1 F_2 \mathfrak{f}$. By what has just been proved, this
multiplication is independent of the choice of $F_i$ from
$F_i \mathfrak{f}$.

The residue classes form a group under the product defined in this
way. The product is associative, because that for the $F$ is; 
$F^{-1}\mathfrak{f}$ is the element inverse to $F\mathfrak{f}$, and
$E\mathfrak{f}=\mathfrak{f}$ is the identity element. This group
is called the \emph{factor group} of $\mathfrak{F}$ by
$\mathfrak{f}$. We denote if by $\mathfrak{F}/\mathfrak{f}$.

We want to construct the factor group by the commutator group
$\mathfrak{K}_1=\mathfrak{K}$ and show that it is commutative.
We have to show that
\[
F_1\mathfrak{K}\cdot F_2\mathfrak{K}=
F_2 \mathfrak{K}\cdot F_1\mathfrak{K}.
\]
But now
\[
F_1\mathfrak{K}\cdot F_2\mathfrak{K}\cdot
F^{-1}_1\mathfrak{K}\cdot F^{-1}_2\mathfrak{K}=
F_1 F_2 F^{-1}_1 F^{-1}_2 \mathfrak{K}=
\mathfrak{K},
\]
because $F_1 F_2 F^{-1}_1 F^{-1}_2$ is the commutator of
$F_1$ and $F_2$.

\section{Isomorphisms}

Two different groups may possess the same product structure.
One captures this relation more precisely by an isomorphic
correspondence between the two groups. \emph{If $\mathfrak{F}$
and $\mathfrak{F}'$ are two groups, and if 
$F'=\mbox{\boldmath$I$}(F)$ is a
one-to-one correspondence between the elements 
of $\mathfrak{F}$ and $\mathfrak{F}'$ such that}
\[
\mbox{\boldmath$I$}(F_1)\mbox{\boldmath$I$}(F_2)=
\mbox{\boldmath$I$}(F_1 F_2),
\]
\emph{then the
group $\mathfrak{F}'$ is called isomorphic to $\mathfrak{F}$ and
the mapping $\mbox{\boldmath$I$}$ is called an isomorphism}.
If we drop the condition that $\boldsymbol{I}$ sends different elements 
of $\mathfrak{F}$
to different elements of $\mathfrak{F}'$, then
$\mathfrak{F}'$ is called \emph{homomorphic} to $\mathfrak{F}'$
and the mapping $\boldsymbol{I}$ is called a \emph{homomorphism}.

A factor group $\mathfrak{F}/\mathfrak{f}$ is homomorphic to
$\mathfrak{F}$. Namely, if we set 
$\mbox{\boldmath$I$}(F)=F\mathfrak{f}$ this is in fact a
homomorphism, not one-to-one unless $\mathfrak{f}$ consists
of a single element, the identity element. In this case the elements
$F$ that form the identity element of $\mathfrak{F}/\mathfrak{f}$
are just those that make up $\mathfrak{f}$.

Conversely, if $\mbox{\boldmath$I$}(F)=F'$ is a homomorphism,
then the identity element of $\mathfrak{F}$ must be associated
with the identity of $\mathfrak{F}'$, i.e.,
 $\mbox{\boldmath$I$}(E)=E'$, because indeed
\[
\mbox{\boldmath$I$}(E)\mbox{\boldmath$I$}(F)=
\mbox{\boldmath$I$}(F)\mbox{\boldmath$I$}(E)=
\mbox{\boldmath$I$}(F).
\]
Consequently, $\mbox{\boldmath$I$}(F)$ is inverse to
$\mbox{\boldmath$I$}(F^{-1})$, because
\[
\mbox{\boldmath$I$}(F)\mbox{\boldmath$I$}(F^{-1})=
\mbox{\boldmath$I$}(E)=E'.
\]
We now take $\mathfrak{f}$ to be the collection of elements $F$
of $\mathfrak{F}$ for which $\mbox{\boldmath$I$}(F)=E'$.
They form a group, because if
\[
\mbox{\boldmath$I$}(F_1)=\mbox{\boldmath$I$}(F_2)=E'
\]
then also
\[
\mbox{\boldmath$I$}(F_1 F_2)=E'.
\]
And since
\[
\mbox{\boldmath$I$}(F_1)\mbox{\boldmath$I$}(F^{-1}_1)=
E'\mbox{\boldmath$I$}(F^{-1}_1)=
\mbox{\boldmath$I$}(F^{-1}_1)
\]
and, on the other hand,
\[
\mbox{\boldmath$I$}(F_1)\mbox{\boldmath$I$}(F^{-1}_1)=
\mbox{\boldmath$I$}(F_1 F^{-1}_1)=E',
\]
$F^{-1}_1$ belongs to this collection along with $F_1$. Finally,
$\mathfrak{f}$ is an invariant subgroup of $\mathfrak{F}$.
Indeed,
\[
\mbox{\boldmath$I$}(F F_1 F^{-1})=
\mbox{\boldmath$I$}(F)E'\mbox{\boldmath$I$}(F^{-1})=E'.
\]

It follows easily from this that $\mbox{\boldmath$I$}$ associates
all elements of a residue class $F\mathfrak{f}$ with the same $F'$
and hence it realizes an isomorphism between
$\mathfrak{F}/\mathfrak{f}$ and $\mathfrak{F}'$.

From Section 1.2 it follows that cyclic groups of the same order are isomorphic.

\section{Automorphisms}

A one-to-one onto transformation of a group $\mathfrak{F}$ into
itself,
\[
F'=\mbox{\boldmath$I$}(F),
\]
is called an autoisomorphism, or simply \emph{automorphism},
when
\[
\mbox{\boldmath$I$}(F_1)\mbox{\boldmath$I$}(F_2)=
\mbox{\boldmath$I$}(F_1 F_2).
\]
\emph{The automorphisms of a group $\mathfrak{F}$ constitute a
group}. For the identity mapping is obviously an automorphism,
likewise the inverse $\mbox{\boldmath$I$}^{-1}$ of an
automorphism. Namely, if $F'_1 F'_2=F'_3$ and 
$F'_i=\mbox{\boldmath$I$}(F_i)$, so that 
$\mbox{\boldmath$I$}^{-1}(F'_i)=F_i$ and
$\mbox{\boldmath$I$}(F_1)\mbox{\boldmath$I$}(F_2)=
\mbox{\boldmath$I$}(F_3)$ then
\[
\mbox{\boldmath$I$}(F_1 F_2)=\mbox{\boldmath$I$}(F_3)
\]
because $\mbox{\boldmath$I$}$ is an automorphism. Thus
$F_1 F_2=F_3$ because $\mbox{\boldmath$I$}$ is one-to-one,
and since $F_i=\mbox{\boldmath$I$}^{-1}(F'_i)$ it follows that
\[
\mbox{\boldmath$I$}^{-1}(F'_1)\mbox{\boldmath$I$}^{-1}(F'_2)
=\mbox{\boldmath$I$}^{-1}(F'_3).
\]
Further, if $\mbox{\boldmath$I$}_1(F)=F'$ and
$\mbox{\boldmath$I$}_2(F)=F''$ are automorphisms then
\[
F''=\mbox{\boldmath$I$}_2(\mbox{\boldmath$I$}_1(F))=
\mbox{\boldmath$I$}_{21}(F)
\]
is likewise an automorphism, because
\begin{align*}
\mbox{\boldmath$I$}_{21}(F_1)\mbox{\boldmath$I$}_{21}(F_2)&=
\mbox{\boldmath$I$}_2
(\mbox{\boldmath$I$}_1(F_1)\mbox{\boldmath$I$}_1(F_2))\\&=
\mbox{\boldmath$I$}_2(\mbox{\boldmath$I$}_1(F_1 F_2))\\&=
\mbox{\boldmath$I$}_{21}(F_1 F_2).
\end{align*}
Finally, the product is associative because it is the composition of
functions.

It is easy to exhibit particular automorphisms. If $F_0$ is a fixed
element of $\mathfrak{F}$ and $F$ runs through all the elements
of $\mathfrak{F}$, then
\[
F'=F_0 F F^{-1}_0
\]
is a one-to-one onto transformation of group elements because
\[
F^{-1}_0 F'F_0=F
\]
is the inverse mapping, and also
\[
F_0 F_1 F^{-1}_0\cdot F_0 F_2 F^{-1}_0=F_0 F_1 F_2 F^{-1}_0.
\]
Such a mapping is called an \emph{inner automorphism}. If $F_0$
runs through all elements of $\mathfrak{F}$ we obtain the totality of
inner automorphisms of $\mathfrak{F}$. They constitute a group
that is homomorphic to $\mathfrak{F}$ itself. The product of two
inner automorphisms
\[
F'=F_1 F F^{-1}_1\quad\text{and}\quad F''=F_2 F' F^{-1}_2
\]
is another, namely,
\[
F''= F_2 F_1 F F^{-1}_1 F^{-1}_2=F_{21} F F^{-1}_{21}.
\]
In order to ascertain whether the homomorphism from
$\mathfrak{F}$ to the group of its inner automorphisms is
one-to-one we must establish which inner automorphisms correspond
to the identity mapping---but these are just those defined by elements
belonging to the center of $\mathfrak{F}$. The group of inner
automorphisms is therefore isomorphic to the factor group
$\mathfrak{F}/\mathfrak{Z}$.

\emph{The inner automorphisms are an invariant subgroup of all the
automorphisms}. Namely, if
\[
\mbox{\boldmath$A$}(F)=F^{*}
\]
is an arbitrary automorphism, and
\[
\mbox{\boldmath$I$}(F)=F'=F_0 F F^{-1}_0
\]
is an inner automorphism, then
\[
\mbox{\boldmath$A$}(\mbox{\boldmath$I$}(F))=
\mbox{\boldmath$A$}(F_0)\mbox{\boldmath$A$}(F)
\mbox{\boldmath$A$}(F^{-1}_0)=
F^{*}_0 F^{*}{ F^{*}_0}^{-1}.
\]
Thus if we set
\[
\mbox{\boldmath$I$}'(F^{*})=\overline{F}=
F^{*}_0 F^{*}{ F^{*}_0}^{-1}=
F^{*}_0 \mbox{\boldmath$A$}(F) {F^{*}_0}^{-1}
\]
then
\[
\mbox{\boldmath$A$}(\mbox{\boldmath$I$}(F))=
\mbox{\boldmath$I$}'(\mbox{\boldmath$A$}(F))
\]
so
\[
\mbox{\boldmath$A$}\mbox{\boldmath$I$}=
\mbox{\boldmath$I$}'\mbox{\boldmath$A$}
\quad\text{or}\quad
\mbox{\boldmath$A$}\mbox{\boldmath$I$}\mbox{\boldmath$A$}^{-1}
=\mbox{\boldmath$I$}'.
\]

If $\mathfrak{f}$ is an invariant subgroup of $\mathfrak{F}$ and $F_0$
is any element of $\mathfrak{F}$, then the mapping
\[
F_0 F F^{-1}_0 =F'
\]
is an automorphism of $\mathfrak{f}$. If $F_0$ itself belongs to
$\mathfrak{f}$, then it is an inner automorphism. The totality of
automorphisms induced by elements $F_0$ of $\mathfrak{F}$
constitute a subgroup of all the automorphisms of $\mathfrak{f}$.
The elements of a residue class $F_0\mathfrak{f}$ modulo
$\mathfrak{f}$ correspond to automorphisms resulting from
multiplication by inner automorphisms.

One can see from these remarks that any group can be embedded as
an invariant subgroup of a larger group. We want to formulate the
situation as follows: given a group $\mathfrak{f}$ and its product
operation, together with a system $\mathfrak{r}$ of representatives
of residue classes of $\mathfrak{F}$ modulo $\mathfrak{f}$, one
then knows that each element of $\mathfrak{F}$ may be written as
a product $RF$, where $R$ is from $\mathfrak{r}$ and $F$ is from
$\mathfrak{f}$. In order to extend the group product to all of
$\mathfrak{F}$, i.e, to know the value of the product
\[
R_1 F_1 R_2 F_2 = R_1 R_2 R^{-1}_1 F_1 R_2 F_2,
\]
we must first know the automorphisms of $\mathfrak{f}$
corresponding to the elements $R$ and also, for any two elements
$R_1,R_2$, the product $R_{12} F_{12}$. Then the group
$\mathfrak{F}$ itself will be known.

\section{Groups with Operators}

When a group $\mathfrak{F}$ with a cyclic group of automorphisms
$\mbox{\boldmath$A$}^n$ is given\footnote{That is, consisting of the 
powers of an automorphism $\boldsymbol{A}$. (Translator's note.)} 
we can make the structure
connecting these two domains of elements clearer by means of a new
symbolism, which is particularly convenient in the case of a commutative
group $\mathfrak{F}$. So we will assume that $\mathfrak{F}$ is
commutative. By $F^x$ we will mean the element 
$\mbox{\boldmath$A$}(F)$ and by $F^{x^n}$ the element
$\mbox{\boldmath$A$}^n(F)$ ($n=0,\pm 1, \pm 2, \ldots)$.
For any integer $a_n$, $F^{a_n x^n}$ means $(F^{a_n})^{x^n}$.

If
\[
f(x)=a_n x^n+a_{n+1}x^{n+1}+\cdots+a_{n+m}x^{n+m}
\]
is an ``$L$-polynomial,''\footnote{The $L$ presumably stands for
``Laurent,'' since these polynomials can have terms with negative
exponent. (Translator's note.)} with integral coefficients $a_i$, then by
$F^{f(x)}$ we mean the element
\[
F^{f(x)}=F^{a_n x^n} F^{a_{n+1}x^{n+1}}
\cdots F^{a_{n+m}x^{n+m}}
\]

One can compute in this extended domain of exponents as in the
original domain of integers. We call two polynomials equal if they
are convertible into each other by deletion or insertion of terms
$a_i x^i$ with $a_i=0$. If $f(x)$ and $g(x)$ are two polynomials
and $n$ is the lowest, and $n+m$ the highest exponent of an 
$x^i$ appearing in $f$ and $g$ along with $a_i$ or $b_i\ne 0$, 
then
\begin{align*}
f(x)&=a_n x^n+a_{n+1} x^{n+1}+\cdots+a_{n+m}x^{n+m}\\
g(x)&=b_n x^n+b_{n+1} x^{n+1}+\cdots+b_{n+m} x^{n+m}.
\end{align*}
As  usual, we understand the sum of $f(x)$ and $g(x)$ to be the
polynomial
\[
f(x)+g(x)=\sum^{n+m}_{i=n} (a_i+b_i)x^i.
\]
This addition operation satisfies the laws of a commutative group,
because the integers under addition are such a group. The polynomial
$f=0$ plays the role of the identity element.

We understand the product of $f(x)$ and $b x^l$ to be the polynomial
\[
a_n b x^{n+l}+a_{n+1} b x^{n+1+l}+\cdots+a_{n+m} b x^{n+m+l}
\]
and as usual we understand the product of $f(x)$ and $g(x)$ 
to be the polynomial
\[
fg=f(x) b_n x^n+f(x) b_{n+1}x^{n+1}+\cdots+f(x)b_{n+m}x^{n+m}.
\]
This multiplication is associative and commutative; $f(x)=1$ is the
identity element. However, an inverse element does not exist in general;
e.g., the polynomial $f(x)=a$ has no inverse when $a\ne\pm 1$,
because the coefficients of all products $a\cdot f(x)$ are divisible
by $a$. The multiplication and addition are further related by the
distributive law:
\[
(f_1(x)+f_2(x))g(x)=f_1(x)g(x)+f_2(x)g(x).
\]
If $f(x)$ and $g(x)$ are two polynomials, both nonzero, and if
$a_n x^n$ and $b_m x^m$ are the lowest-order terms appearing in
$f(x)$ and $g(x)$ with $a_n\ne 0$ and $b_m\ne 0$, then the
lowest-order term appearing in $f(x)g(x)$ is $a_n b_m x^{n+m}$.
It follows from this that if $f(x)g(x)=0$ then at least one of the factors
$f(x)$ or $g(x)$ equals zero. Thus the polynomials constitute an
integral domain.\footnote{See a textbook of algebra, e.g.,
\textsc{H. Hasse} \emph{H\"ohere Algebra}, Band 1, Sammlung
G\"oschen.}

If the smallest exponent in a polynomial $f(x)$ is greater than or
equal to zero then $f(x)$ is an ordinary\footnote{Reidemeister calls
such a polynomial \emph{entire}, following the terminology of complex
analysis. But it seems harmless, and clearer, to call such polynomials
``ordinary.'' (Translator's note.)}
polynomial in $x$. For any ordinary polynomial
\[
f(x)=a_n x^n+a_{n+1} x^{n+1}+\cdots+a_{n+m}x^{n+m}
\]
with $a_n\ne 0$
\[
x^{-n}f(x)=a_n+a_{n+1}x+\cdots+a_{n+m}x^m
\]
is an ordinary polynomial with nonzero constant term.

If we now consider $f(x)$ and $g(x)$ as exponents of group
elements, it turns out that
\[
F^{f(x)} F^{g(x)}=F^{f(x)+g(x)}.
\]
This follows easily from the commutativity of the group and the
definition of $F^{f(x)}$ and $f(x)+g(x)$. Further,
\[
(F^{f(x)})^{g(x)}=F^{f(x)g(x)},
\]
because
\[
(F^{f(x)})^{g(x)}=
(F^{f(x)})^{b_n x^n} (F^{f(x)})^{b_{n+1} x^{n+1}}\cdots
(F^{f(x)})^{b_{n+m} x^{n+m}}.
\]
Then, on the one hand,
\[
(F^{f(x)})^{b_i}=F^{b_i f(x)}
\]
since
\[
(F^{f(x)})^{b_i}=F^{f(x)} F^{f(x)}\cdots F^{f(x)}\quad
\text{with $b_i$ factors}
\]
for positive $b_i$, and for negative $b_i$
\[
(F^{f(x)})^{b_i}=\left((F^{f(x)})^{-1}\right)^{-b_i}
                         =\left(F^{-f(x)}\right)^{-b_i}.
\]
While, on the other hand,
\begin{align*}
(F^{f(x)})^{x^i}
&=\left(F^{a_n x^n} F^{a_{n+1}x^{n+1}}
\cdots F^{a_{n+m}x^{n+m}}\right)^{x^i}\\
&=\left(F^{a_n x^n}\right)^{x^i}
    \left(F^{a_{n+1} x^{n+1}}\right)^{x^i}\cdots
    \left(F^{a_{n+m} x^{n+m}}\right)^{x^i}\\
&=(F^{a_n})^{x^{n+i}}
    (F^{a_{n+1}})^{x^{n+1+i}}\cdots
    (F^{a_{n+m}})^{x^{n+m+i}}\\
&=F^{x^i f(x)}.
\end{align*}
Consequently, $(F^{f(x)})^{b_i x^i}=F^{b_i x^i f(x)}$
and hence
\begin{align*}
(F^{f(x)})^{g(x)}
&=F^{b_n x^n f(x)} F^{b_{n+1} x^{n+1} f(x)}\cdots
    F^{b_{n+m} x^{n+m} f(x)}\\
&=F^{b_n x^n f(x)+b_{n+1}x^{n+1}+\cdots+
          b_{n+m} x^{n+m} f(x)}\\
&=F^{f(x)g(x)}.
\end{align*}

Thus one can compute with the formally introduced
$L$-polynomials as exponents just as with integral exponents.

\section{Groups and Transformation Groups}

We call a transformation group $\mathfrak{T}$ that is
homomorphic to an arbitrary group $\mathfrak{F}$ a
representation of $\mathfrak{F}$, and we further examine the
different representations of a group.

We take as our domain $\mathfrak{X}$ of objects the right-sided 
residue classes modulo a subgroup $\mathfrak{f}$, so
$x=R\mathfrak{f}$, and define
\[
F(x)=FR\mathfrak{f}=x'
\]
to be the transformation of this domain corresponding to the
group element $F$.

This mapping is one-to-one and onto, because $F^{-1}$ yields the 
inverse mapping. The transformations that carry the element 
$x=\mathfrak{f}$ into itself are exactly those that correspond to
elements of $\mathfrak{f}$. The transformations that carry
$R\mathfrak{f}$ to itself  correspond to the elements of the 
group $R\mathfrak{f}R^{-1}$ conjugate to $\mathfrak{f}$.
We can now easily give a criterion for isomorphism between
$\mathfrak{F}$ and the group just defined.

Those elements that correspond to transformations leaving all
$x$ fixed must therefore belong to the intersection
$\mathfrak{D}$ of the groups $R\mathfrak{f}R^{-1}$ conjugate 
to $\mathfrak{F}$. The group $\mathfrak{T}$ is therefore 
isomorphic to the factor group $\mathfrak{F}/\mathfrak{D}$.

If $\mathfrak{f}$ is an invariant subgroup, then the transformations
that correspond to $\mathfrak{f}$, and thus carry the element
$x=\mathfrak{f}$ to itself, also carry all the remaining $x$ to
themselves, because $\mathfrak{D}$ in this case is equal to
$\mathfrak{f}$. The transformation group is then simply transitive.

Conversely, given any simply transitive group of transformations
isomorphic to $\mathfrak{F}$, an arbitrary element $x_0$, and
$\mathfrak{f}_{x_0}$ the subgroup of transformations that leave
$x_0$ fixed, the transformations that carry $x_0$ to $x$
constitute a residue class modulo $\mathfrak{f}_{x_0}$. Namely,
if $R_x$ carries the element $x_0$ to $x$, so also do the
transformations in $R_x\mathfrak{f}_{x_0}$, and if $R'$ is any
transformation that carries $x_0$ to $x$, then $R^{-1}_x R'$
carries the element $x_0$ to itself and it therefore belongs to
$\mathfrak{f}_{x_0}$.

If we associate with each group element
\[
F(x)=x'
\]
the transformation
\[
F(R_x \mathfrak{f}_{x_0})=FR_x\mathfrak{f}_{x_0}
\]
in the domain of residue classes, then
\[
FR_x\mathfrak{f}_{x_0}=R_{x'}\mathfrak{f}_{x_0}.
\]
Thus the new transformation group is simply the original one with
renaming of the objects transformed. A representation of a group
$\mathfrak{F}$ by a transitive transformation group is an isomorphism,
by the remarks above, if and only if the domain of objects can be
viewed as a system of right-sided residue classes modulo a subgroup
$\mathfrak{f}$, where the intersection of $\mathfrak{f}$ with its
conjugate subgroups is the identity.

\section{The Groupoid}

For many topological questions a generalization of the group
concept, the groupoid,\footnote{\textsc{H. Brandt}, Math. Ann. 
\textbf{96}, 360.} is a useful auxiliary.

A collection $\mathfrak{G}$ of elements $G$ with a product
$G_1 G_2=G_3$ is called a \emph{groupoid} when the following
conditions are satisfied.
\begin{enumerate}
\item[A. 1.]
\emph{If a relation $G_1 G_2=G_3$ holds between three elements
$G_1,G_2,G_3$, then each of them is uniquely determined by the 
other two.}
\item[A. 2.]
\emph{If $G_1G_2$ and $G_2 G_3$ exist then $(G_1 G_2)G_3$
and $G_1(G_2 G_3)$ also exist; if $G_1 G_2$ and $(G_1 G_2)G_3$
exist then $G_2 G_3$ and $G_1(G_2 G_3)$ also exist; if $G_2 G_3$ 
and $G_1(G_2 G_3)$ exist then $G_1 G_2$ and $(G_1 G_2)G_3$ also 
exist; and in each case $(G_1 G_2)G_3=G_1(G_2 G_3)$, so that it
can also be written $G_1 G_2 G_3$.}
\item[A. 3.]
\emph{For each element $G$ the following elements exist: the right
identity $E$, the left identity $E'$, and the inverse $G^{-1}$, for
which the following relations hold:}
\[
GE=G,\quad E'G=G,\quad G^{-1}G=E.
\]
\item[A. 4.]
\emph{For any two identities $E,E'$ there is an element $G$ for
which $E$ is the right identity and $E'$ is the left identity.}
\end{enumerate}

One sees that the generalization consists in relinquishing general
applicability of the product and admitting several identities. 
Groupoids with a single identity are groups.

Just as in Section 1.1, one can prove the  appropriate analogue of
the associative law for the product of arbitrarily many elements.

Two elements $G_1$ and $G_2$ are composable in that order
if and only if the right identity of $G_1$ is identical with the left
identity of $G_2$. The subclass $\mathfrak{G}_i$ of elements
$G$ for which the right and left identities both equal $E_i$
constitute a group. The groups $\mathfrak{G}_i$ associated with
different identities are isomorphic.

An example of a groupoid may be constructed from a group
$\mathfrak{T}$ of transformations $T$ of objects $x$, which 
carry the object $x_1$ in particular into finitely many objects
\[
x_1,\quad x_2,\quad \ldots,\quad x_n.
\]
The groupoid then has $n$ corresponding identities $E_1,E_2,\ldots,E_n$;
also, for each element $G$ of the groupoid there is an associated element
$\mbox{\boldmath{$A$}}(G)=T$, where $T$ is a transformation
carrying $x_a$ to $x_b$ if $G$ has left identity $E_a$ and right identity
$E_b$ and, conversely, for each such element $T$ there is an element
$G$ with the corresponding identities. One observes that an element
$T$ corresponds to $n$ different $G$, because $T$ permutes the
$x_i$ among themselves. Also, it it the case that
\[
\mbox{\boldmath{$A$}}(G_1 G_2)=
\mbox{\boldmath{$A$}}(G_1)\mbox{\boldmath{$A$}}(G_2).
\]

A groupoid is uniquely determined by these conditions. The groups
corresponding to the identities are isomorphic to the subgroups of
$\mathfrak{T}$ that leave the $x_i$ fixed.
\chapter{Free Groups and their Factor Groups}

\section{Generators and Defining Relations}

The groups that appear in combinatorial topology are defined in a
way that itself has a combinatorial character. The peculiar difficulties of
topology can be better appreciated when one has at hand the analogous 
problems of group theory, which we are about to present.

If $\mathfrak{F}$ is any group and $\mathfrak{m}$ is a class of
elements from which all elements of $\mathfrak{F}$ may be constructed
as power products, then $\mathfrak{m}$ is called a \emph{system of
generators} for the group $\mathfrak{F}$. Thus a system of generators
for the integers is just 1, and for an additive residue class group it is
the class [1]. These examples already draw attention to the fact that
formally different power products can yield the same group element.
If one wants to be able to derive the product of elements of $\mathfrak{F}$
from the product of power products, then one must be able to decide
which power products represent equal group elements. This reduces to
the question of which products represent the identity 
element.\footnote{The so-called \emph{word problem} for the group
$\mathfrak{F}$. (Translator's note.)}

We call each product $R(\mathfrak{m})$ of elements of $\mathfrak{m}$
that equals the identity a relation, and call the totality of relations
$\mathfrak{R}$. Now if $P$ is any power product and $R$ is any
relation, then obviously $P$ and $PR$ are the same group element
\[
P=PR.
\]
Conversely, if $P_1$ and $P_2$ are two power products that denote the
same group element, and if $P^{-1}_1$ is the product formally inverse
to $P_1$, then $P^{-1}_1 P_2$ is a relation $R'$, and the power product
$P_2$ results from the product $P_1 R'$ by deletion of
adjacent factors $FF^{-1}$. We thus obtain all representations of the element
$P$ in the form $PR$ when $R$ runs through the class $\mathfrak{R}$,
if we also include those products that result from deletion of formally
inverse adjacent factors of $PR$.

The power products of $\mathfrak{R}$ have the following properties.

If $R=P_1 P_2$ belongs to $\mathfrak{R}$, so does $R'=P_1 FF^{-1} P_2$,
and conversely, if $R'$ belongs to $\mathfrak{R}$ so does $R$. If $R$
belongs to $\mathfrak{R}$, so does the formal inverse $R^{-1}$. If $P$
is an arbitrary power product, $P^{-1}$ its formal inverse, and if $R$
belongs to $\mathfrak{R}$, then $PRP^{-1}$ also belongs to 
$\mathfrak{R}$. If $R_1$ and $R_2$ belong to $\mathfrak{R}$, then
the product $R_1 R_2$ also belongs to $\mathfrak{R}$. By means of
these four processes, ``consequence relations'' may be derived from
relations originally given. We call a class $\mathfrak{r}$ of defining relations
from which all relations in $\mathfrak{R}$ may be derived by the four
processes a \emph{system of defining relations}. With a class
$\mathfrak{m}$ of generators and a class $\mathfrak{r}$ of defining
relations the product law is obviously defined for all elements of
$\mathfrak{F}$, hence the name ``defining relations.'' Establishing
generators and defining relations for groups given in other ways is a
far from trivial problem.\footnote{Cf. Sections 2.9 and 3.1 and, e.g.,
\textsc{J. Nielsen}, Kgl. Dan. Vid. Selsk., Math. fys. Med. 
\textbf{V}, 12 (1924). }

Just as for groups, one can speak of \emph{generators} for a
\emph{groupoid}. We will assume that one can find, from generators
$S_i$ ($i=1,2,\ldots,m$) of a groupoid $\mathfrak{G}$ with
identities $E_i$ ($i=1,2,\ldots, n$), a system of generators $T_i$ of
the group $\mathfrak{G}_0$ of elements doubly associated with the identity
$E_0$.

Let $A_i$ ($i=1,2\ldots,n$) be a system of elements with left-sided
identity $E_0$ and right-sided identities including all $E_i$
($i=1,2\ldots,n$). Further, let $A_0=E_0$. Now if $S_i$ has the left
identity $E_{l_i}$ and the right identity $E_{r_i}$ then the element
\begin{equation}
T_i=A_{l_i}S_i A^{-1}_{r_i} \tag{1}
\end{equation}
may be called the generator of $\mathfrak{G}_0$ associated with $S_i$.
In fact, the $T_i$ ($i=1,2\ldots,n$) constitute a system of generators
for $\mathfrak{G}_0$. Namely, if
\[
S^{\varepsilon_1}_{\alpha_1}
S^{\varepsilon_2}_{\alpha_2} 
\cdots
S^{\varepsilon_a}_{\alpha_a} \tag{2}
\]
is any element of $\mathfrak{G}_0$, then the left identity of
$S^{\varepsilon_1}_{\alpha_1}$ and the right identity of
$S^{\varepsilon_a}_{\alpha_a}$ is the identity $E_0$ and, further,
the right identity of $S^{\varepsilon_i}_{\alpha_i}$ is identical with the
left identity of $S^{\varepsilon_{i+1}}_{\alpha_{i+1}}$. The product
\[
T^{\varepsilon_1}_{\alpha_1}
T^{\varepsilon_2}_{\alpha_2} 
\cdots
T^{\varepsilon_a}_{\alpha_a}
\]
that results from (2) when $S_i$ is replaced by $T_i$ may be converted
into (2) by means of equation (1) and cancellation of formally inverse
factors $S_i$.

\section{Free Groups}

Instead of starting with a group and constructing generators and
defining relations for it, we will now proceed from a class $\mathfrak{m}$
of symbols, define the power products of these symbols, take an
arbitrary system $\mathfrak{r}$ of these power products, and show that 
there is a group $\mathfrak{F}$ that has the symbols in $\mathfrak{m}$
as generators and the products in $\mathfrak{r}$ as defining relations.
For this purpose we first explain relation-free groups, or simply
\emph{free groups with $n$ generators}. Let
\[
S^{+1}_1,S^{+1}_2,\ldots,S^{+1}_n,S^{-1}_1,S^{-1}_2,\ldots,
S^{-1}_n
\]
be letters, which we combine into ``words''
\begin{equation}
W=
S^{\varepsilon_1}_{\alpha_1}
S^{\varepsilon_2}_{\alpha_2}
\cdots
S^{\varepsilon_m}_{\alpha_m}\quad
(\alpha_i=1,2,\ldots,n;\; \varepsilon_i=\pm 1) .
\tag{1}
\end{equation}
Let $W_0$ be the ``empty'' word, which contains no letters $S^{\pm 1}_i$.
The word
\[
W^{-1}=
S^{-\varepsilon_m}_{\alpha_m}
\cdots
S^{-\varepsilon_2}_{\alpha_2}
S^{-\varepsilon_1}_{\alpha_1}
\]
is called the formal inverse of $W$. If
\[
W_1=
S^{\varepsilon_1}_{\alpha_1}
S^{\varepsilon_2}_{\alpha_2}
\cdots
S^{\varepsilon_m}_{\alpha_m}
\]
and
\[
W_2=
S^{\eta_1}_{\beta_1}
S^{\eta_2}_{\beta_2}
\cdots
S^{\eta_{m'}}_{\beta_{m'}}
\]
are two such words we set
\[
W_1 W_2=
S^{\varepsilon_1}_{\alpha_1}
S^{\varepsilon_2}_{\alpha_2}
\cdots
S^{\varepsilon_m}_{\alpha_m}
S^{\eta_1}_{\beta_1}
S^{\eta_2}_{\beta_2}
\cdots
S^{\eta_{m'}}_{\beta_{m'}}
\]
and
\[
W_0 W_1=W_1 W_0=W_1.
\]

By an \emph{elementary transformation} of a word we mean the cancellation
or insertion of two symbols
\[
S^{\varepsilon_i}_{\alpha_i} S^{\varepsilon_{i+1}}_{\alpha_{i+1}}
\quad\text{when}\quad \alpha_i=\alpha_{i+1}\quad\text{and}\quad
\varepsilon_i+\varepsilon_{i+1}=0.
\]
An elementary transformation in $W_1$ is always one in $W_1 W_2$,
but the converse does not hold. Two words $W_1$ and $W_n$
[not the same as the $n$ above] are called \emph{equivalent},
denoted
\[
W_1\equiv W_n,
\]
if there is a chain of words $W_1,W_2,\ldots,W_n$ of which any 
successive two are convertible into each other by an elementary
transformation. If $W_1\equiv W_2$ and $W_2\equiv W_3$ then also
$W_1\equiv W_3$, and $W_1\equiv W_2$ implies $W_2\equiv W_1$.
We can therefore speak of the class $[W]$ of words equivalent to $W$,
and
\[
[W_1]=[W_2]
\]
if and only if $W_1\equiv W_2$.

We define a product for the classes $W$ which, as we will show, satisfies
the group axioms. We set
\[
[W_1][W_2]=[W_1 W_2].
\]
This product is well-defined. Namely, if $W'_1\equiv W_1$ and
$W'_2\equiv W_2$ then
\[
W'_1 W'_2\equiv W_1 W_2
\]
because elementary transformations of the $W_1$ ($i=1,2$) are also
elementary transformations of $W_1 W_2$. The product is associative.
$[W_0]$ is the identity element and $[W^{-1}]$ is the element inverse to
$[W]$.

The group so defined is called the free group with $n$ generators.
For the
\[
[S^{+1}_i]\quad (i=1,2,\ldots,n)
\]
obviously constitute a set of generators for this group. E.g., $[W]$ in (1)
is equal to
\[
[S^{+1}_{\alpha_1}]^{\varepsilon_1}
[S^{+1}_{\alpha_2}]^{\varepsilon_2}
\cdots
[S^{+1}_{\alpha_m}]^{\varepsilon_m}.
\]
We can now regard the word $W$ as a sign for the element $[W]$ and as a
power product of the elements $[S^{+1}_i]$. We also write $S_i$ for
$[S_i]$ and use $S^n_i$ in the way explained in Section 1.1.

One can define the free group with denumerably many generators in quite
an analogous way.\footnote{Another foundation for free groups is due to
\textsc{O. Schreier}, Hamb. Abh. \textbf{5} (1927).}

\section{The Word Problem for Free Groups}

It is easy to survey the representations of the identity element by words
$W$, i.e., the totality of relations $\mathfrak{R}$ in the generators $S_i$,
and to solve the word problem in general, i.e.,  decide directly when two
words $W_1$ and $W_2$ are equivalent. For this purpose we define the
concept of a reduced word and show that there is only one reduced word
$|W|$ in a class $[W]$. A word $W$ is called \emph{reduced} if no two
letters $S^{\varepsilon_i}_{\alpha_i} S^{\varepsilon_{i+1}}_{\alpha_{i+1}}$
with $\alpha_i=\alpha_{i+1}$ and $\varepsilon_i+\varepsilon_{i+1}=0$
appear in $W$.

In order to prove our theorem, we give a unique reduction process for the
word $W$ in Section 2.1 (1). Let
\[
W_1=S^{\varepsilon_1}_{\alpha_1},\quad
W_2=S^{\varepsilon_1}_{\alpha_1} S^{\varepsilon_2}_{\alpha_2},\quad
\ldots,\quad
W_m=W.
\]
Then $|W_1|=S^{\varepsilon_1}_{\alpha_1}$, $|W_2|$ is $W_0$ if
$\alpha_1=\alpha_2$ and $\varepsilon_1+\varepsilon_2=0$, otherwise
$|W_2|=W_2$. $|W_i|$ is defined inductively: if $|W_{i-1}|=W_0$, then
$|W_i|=S^{\varepsilon_i}_{\alpha_i}$; if $|W_{i-1}|\ne W_0$, 
$S^{\varepsilon}_\beta$ is the last letter in $|W_{i-1}|$ and $\beta=\alpha_i$,
$\varepsilon+\varepsilon_i=0$, then $|W_i|$ is the word that results from
cancelling the $S^{\varepsilon}_\beta$ off the end of $|W_{i-1}|$; if
$\beta=\alpha_i$ and $\varepsilon+\varepsilon_i=0$ do \emph{not}
both hold then $|W_i|$ is the word $|W_{i-1}|S^{\varepsilon_i}_{\alpha_i}$.
Obviously, all words $|W_i|$ are reduced, and so $|W_m|$ is a reduced
word equivalent to $W$.

Now let $W'$ be the word that results from $W$ by insertion of
$S^{\varepsilon}_\alpha S^{-\varepsilon}_\alpha$ between
$S^{\varepsilon_k}_{\alpha_k}$ and 
$S^{\varepsilon_{k+1}}_{\alpha_{k+1}}$. We show that our process,
applied to $W'$, leads to the same reduced word $|W'|=|W|$. We set
\begin{align*}
W'_i=&W_i\quad (i=1,2,\ldots, k)\\
W'_{k+1}=W_k S^{\varepsilon}_\alpha,\quad
&W'_{k+2}=W_k S^{\varepsilon}_\alpha S^{-\varepsilon}_\alpha,
\quad\ldots,\quad W'_{m+2}=W'.
\end{align*}
Then
\[
|W'_i|=|W_i|\quad (i=1,2,\ldots, k).
\]
$|W'_k|$ ends either with $S^{-\varepsilon}_\alpha$, in which case
$|W'_{k+1}|$ equals the word resulting from $|W'_k|$ by cancellation of
$S^{-\varepsilon}_\alpha$, so that $|W'_{k+2}|=|W'_k|=|W_k|$ and
in general $|W'_{k+l+2}|=|W_{k+l}|$; or else $|W'_k|$ does not end in
$S^{-\varepsilon}_\alpha$, in which case
\[
|W'_{k+1}|=|W'_k|S^{\varepsilon}_\alpha=|W_k|S^{\varepsilon}_\alpha,
\]
so $|W'_{k+2}|=|W_k|$ and again in general
\[
|W'_{k+l+2}|=|W_{k+l}|\quad (l=1,2,\ldots,m-k).
\]

Now if $W$ and $W^{*}$ are any words equivalent to each other they may
be embedded in a chain of words of which each is convertible to its neighbor
by an elementary transformation, and thus our reduction process must lead
each of $W,W^{*}$ to the same reduced word $|W|=|W^{*}|$. Since
$W=|W|$ for a reduced word, \emph{reduced words are equivalent only if 
they are identical}.

\section{The Transformation Problem in Free Groups}

Closely related to the word problem is the more general question of the
transformation problem.\footnote{Today, this is called
the \emph{conjugacy} problem, and words $W_1$ and $W_2$ such that
$W_2=W_3 W_1 W^{-1}_3$ are called \emph{conjugate} rather than
``transforms'' of each other. Note that Reidemeister has spoken about
conjugate subgroups already in Section 1.5. (Translator's note.)}
Given two words $W_1$ and $W_2$, one has to decide whether there
is a third word $W_3$ such that
\[
W_2=W_3 W_1 W^{-1}_3.
\]
In this case $W_2$ is called a ``transform'' of $W_1$; it results from
$W_1$ by ``transformation by $W_3$.''

We first define a special class of words, the \emph{short words}. These
are reduced words
\[
W=S^{\varepsilon_1}_{\alpha_1}
     S^{\varepsilon_2}_{\alpha_2}
     \cdots
    S^{\varepsilon_m}_{\alpha_m}
\]
in which 
\[
\alpha_1=\alpha_m\quad\text{and}\quad
\varepsilon_1+\varepsilon_m=0
\]
do not both hold.\footnote{Today such words are called \emph{cyclically
reduced}. (Translator's note.)}

If one permutes the letters of such a word cyclically,
\[
W'=S^{\varepsilon_2}_{\alpha_2}
     \cdots
    S^{\varepsilon_m}_{\alpha_m}
    S^{\varepsilon_1}_{\alpha_1},
\]
then $W'$ is also a short word. By $\{W\}$ we mean the class of short
words that result from $W$ by cyclic permutations. Since 
\[
W'=S^{-\varepsilon_1}_{\alpha_1}
     S^{\varepsilon_1}_{\alpha_1}
     S^{\varepsilon_2}_{\alpha_2}
     \cdots
    S^{\varepsilon_m}_{\alpha_m}
    S^{\varepsilon_1}_{\alpha_1},
\]
all elements of $\{W\}$ are transforms of $W$.

If $W$ is any reduced word that is not a short word, then
\[
W= S^{\varepsilon_1}_{\alpha_1}
     \left(S^{\varepsilon_2}_{\alpha_2}
           \cdots
           S^{\varepsilon_{m-1}}_{\alpha_{m-1}}\right)
    S^{-\varepsilon_1}_{\alpha_1},
\]
so by continuing in this way we eventually obtain
\[
W=W_1 \overline{W} W^{-1}_1,
\]
where $\overline{W}$ is a short word. $\overline{W}$ may be called
the kernel of $W$. By $\{\{W\}\}$ we will mean all those words that
have a kernel in $\{W\}$. All words in $\{\{W\}\}$ obviously
correspond to elements that are transforms of each other.

One now sees that, along with $W^{*}$, the word that results from
$S^{\varepsilon}_i W^{*} S^{-\varepsilon}_i$ by reduction also
belongs to $\{\{W\}\}$, and from this it follows that the elements
belonging to $\{\{W\}\}$ are all the transforms of this element.

\section{Groups with Arbitrary Relations}

We now construct a \emph{group with generators}
\[
S_1,\quad S_2,\quad \ldots,\quad S_n
\]
\emph{and defining relations}
\[
R_1(S),\quad R_2(S),\quad \ldots, \quad R_m(S),
\]
where the $R_i$ are any words in the $S_i$. We first construct the free
group $\mathfrak{S}$ determined by the $S_i$. The $R_i$ are extended
by adjoining all $LR_i L^{-1}$, where $L$ is an arbitrary element of
$\mathfrak{S}$, and we construct the subgroup $\mathfrak{R}$ of
$\mathfrak{S}$ consisting of all power products of the $R_i$ and their
transforms $LR_i L^{-1}$. This is obviously an invariant subgroup of
$\mathfrak{S}$. Thus we can construct the factor group
$\mathfrak{F}=\mathfrak{S}/\mathfrak{R}$ of $\mathfrak{S}$ by
$\mathfrak{R}$ by Section 1.10. We claim that the residue classes
\[
S_1\mathfrak{R},\quad
S_2\mathfrak{R},\quad
\ldots,\quad
S_n\mathfrak{R}
\]
generate this group and that the $R_i$ yield a system of defining relations
for $\mathfrak{F}$ in the generators $S\mathfrak{R}$ when $S_i$ is
replaced by $S_i\mathfrak{R}$. The products $R_i(S\mathfrak{R})$
that result in this way are certainly relations, for it follows from
\[
S^{\pm 1}_i\mathfrak{R} S^{\pm 1}_k\mathfrak{R}=
S^{\pm 1}_i  S^{\pm 1}_k\mathfrak{R}
\]
that
\[
R_i(S\mathfrak{R})=R_i(S)\mathfrak{R}=\mathfrak{R}.
\]
Conversely, if $R(S\mathfrak{R})$ is any relation in the group
$\mathfrak{S}/\mathfrak{R}$ then $R(S)\mathfrak{R}=\mathfrak{R}$,
so $R(S)$ must belong to $\mathfrak{R}$; i.e., $R(S)$ may be written
as a power product of the $R_i(S)$ and their transforms $LR_i L^{-1}$.
Thus the $R_i(S\mathfrak{R})$ really are a system of defining relations
for $\mathfrak{S}/\mathfrak{R}$. 

Since each word in the $S_i$ corresponds to a well-defined element of
the group $\mathfrak{S}/\mathfrak{R}$, we can regard it as a notation
for this element and, e.g., speak of the element $S_i$ of the group
$\mathfrak{S}/\mathfrak{R}=\mathfrak{F}$ and hence call $\mathfrak{F}$
the group with generators $S_i$ ($i=1,2,\ldots, n$) and defining
relations $R_k$ ($k=1,2,\ldots,m$).\footnote{\textsc{O. Schreier},
Hamb. Abhdl. \textbf{5} (1927) 161.} On the other hand, \emph{If
$\mathfrak{F}'$ is a group with the generators $S'_i$ ($i=1,2,\ldots,n$)
and defining relations $R'_k(S')$ ($k=1,2,\ldots,m$), then
$\mathfrak{F}'$ is isomorphic to a factor group of the free group with
$n$ free generators.}

If $L$ and $M$ are arbitrary power products from $\mathfrak{S}$, and
$R$ is a power product from $\mathfrak{R}$, then the element $LRM$
is equal to $LM$ in $\mathfrak{F}$. For $LRM$ is in fact equal to
$LM\cdot M^{-1}RM$.

If $A$ is any element of $\mathfrak{F}$ then the power products of $A$
and its transforms constitute an invariant subgroup $\mathfrak{A}$ of
$\mathfrak{F}$. If we now construct $\mathfrak{F}/\mathfrak{A}=
\mathfrak{F}'$, then each power product $F(S)$ of the $S_i$ also
represents a certain element of $\mathfrak{F}'$, and in fact it represents
the identity of $\mathfrak{F}'$ if and only if it represents an element $A'$
of $\mathfrak{A}$ in $\mathfrak{F}$, i.e., if $F(S)\equiv A' R$ in
$\mathfrak{S}$, where $R$ is a consequence relation of the $R_i$.
Thus one sees that $F(S)$ can be regarded as a consequence of the
$R_i$ ($i=1,2,\ldots,m$) and the relation $A=R_{m+1}$.

One shows analogously: if $A_1,A_2,\ldots,A_l$ are elements, which
together with their transforms generate an invariant subgroup 
$\mathfrak{A}$ of $\mathfrak{F}$, then the relations $R_i$
($i=1,2,\ldots,m$) of $\mathfrak{F}$ and the relations $R_{m+i}=A_i$
($1=1,2\ldots,l$) constitute a system of defining relations for
$\mathfrak{F}'=\mathfrak{F}/\mathfrak{A}$.

\section{The general word problem}

The peculiar difficulties of combinatorial problems show themselves for 
the first time when one tries to solve the word problem for a group with
arbitrary defining relations, i.e., to decide when two products of the 
generators $S_i$ denote the same element of the group $\mathfrak{F}$.
We are far from a general solution of the problem and we have reached 
the goal in only a few cases.

That this is in the nature of things is shown by the following 
remark:\footnote{\textsc{W. Hurewicz}, Hamb. Abhdl. \textbf{8} (1931) 
307.}

We assume that in a group $\mathfrak{G}$ with generators $S_i$
($i=1,2\ldots,n$) there is a certain power product $P_G$ in the
generators $S_i$ for each element $G$, satisfying the condition
\begin{equation}
P_{G_1 G_2}=P_{G_1} P_{G_2}, \tag{1}
\end{equation}
i.e., the condition that both sides of (1) are identical in the free
group $\mathfrak{S}$ generated by the $S_i$. Then the group
$\mathfrak{G}$ is a free group. Namely, $\mathfrak{G}$ is
isomorphic to the subgroup of the free group $\mathfrak{S}$
generated by $P_G$ and, as we will show in Sections 3.9, 4.17,
4.20, and 7.12, the subgroups of free groups are free.

As an example\footnote{\textsc{O. Schreier}, Hamb. Abhdl.
\textbf{3} (1924) 167.} of a word problem we consider the \emph{group
$\mathfrak{F}$ with generators $S_1$ and $S_2$ and the defining
relations}
\begin{equation}
R_1=S^{a_1}_1,\quad R_2=S^{a_2}_2\quad (a_1,a_2>1). \tag{2}
\end{equation}
Here
\[
S^{m_i}_i\equiv S^{n_i}_i\quad\text{when}\quad
m_i\equiv n_i\text{ (mod $a_i$)}\quad (i=1,2).
\]
We call a product
\begin{equation}
S^{r_{11}}_1 S^{r_{21}}_2 S^{r_{12}}_1 S^{r_{22}}_2 \cdots
S^{r_{1l}}_1 S^{r_{2l}}_2 \tag{3}
\end{equation}
``reduced in $\mathfrak{F}$'' when $0\le r_{ik}<a_i$ ($k=1,2,\ldots,l;
i=1,2$) and all $r_{ik}$ except possibly $r_{11}$ and $r_{2l}$ are
nonzero. The reduced product in $\mathfrak{F}$
\[
S'^{r'_{2l}}_2 S'^{r'_{1l}}_1\cdots S'^{r'_{21}}_2 S'^{r'_{11}}_1
\]
is the element inverse to (3) when $r_{ik}+r'_{ik}=a_i$. We will show
that each element of our group is representable in only one way as a
reduced product.

First we give a process that associates with each reduced word $W$ of
the form (1) from Section 2.2 a unique reduced word $|W|$, ``equivalent''
in $\mathfrak{F}$. We set $W_1=S^{\varepsilon_1}_{\alpha_1}$ and
$|W_1|=S^{r_1}_{\alpha_1}$, where $\varepsilon_1\equiv r_1$
(mod $a_{\alpha_1}$), $0\le r_1<a_{\alpha_1}$,
$W_2=S^{\varepsilon_1}_{\alpha_1} S^{\varepsilon_2}_{\alpha_2}$
and $|W_2|=S^{r_2}_{\alpha_1}$ where 
$\varepsilon_1+\varepsilon_2\equiv r_2$ (mod $a_{\alpha_1}$),
$0\le r_2<a_{\alpha_1}$, when $\alpha_1=\alpha_2$,
and $|W_2|=S^{r_1}_{\alpha_1} S^{r_2}_{\alpha_2}$ where
$\varepsilon_2\equiv r_2$ (mod $a_{\alpha_2}$), $0\le r_2<a_{\alpha_2}$
when $\alpha_1\ne \alpha_2$.

In general, let 
$W_i=S^{\varepsilon_1}_{\alpha_1} S^{\varepsilon_2}_{\alpha_2}
\cdots S^{\varepsilon_i}_{\alpha_i}$ and $|W_i|=W'_i S^{r_i}_\beta$.
If $\alpha_{i+1}=\beta$, let $|W_{i+1}|=W'_i S^{r'}_\beta$, where
$r_i+\varepsilon_{i+1}\equiv r'$ (mod $a_\beta$), $0\le r'<a_\beta$.
If $\alpha_{i+1}\ne \beta$, let 
$|W_{i+1}|=|W_i|S^{r_{i+1}}_{\alpha_{i+1}}$, where
$r_{i+1}\equiv\varepsilon_{i+1}$ (mod $a_{\alpha_{i+1}}$),
$0\le r_{i+1}<a_{\alpha_{i+1}}$. $|W|$ equals $|W_m|$.

Now if $W=W'W''$ and 
$W^{*}=W'S^{\varepsilon}_\alpha S^{-\varepsilon}_{\alpha}W''$
one sees that $|W'|=|W'S^{\varepsilon}S^{-\varepsilon}|$ and hence
$|W|=|W^{*}|$. Further, if
\[
W=W'W''
\]
and 
\[
W^{*}=W' 
S^\varepsilon_{\alpha}S^\varepsilon_{\alpha}\cdots S^\varepsilon_{\alpha}
W'',
\]
where $a_\alpha$ factors $S^\varepsilon_{\alpha}$ are inserted, then
likewise
\[
|W'|=|W' S^\varepsilon_{\alpha}S^\varepsilon_{\alpha}\cdots S^\varepsilon_{\alpha}|
\]
and hence also
\[
|W|=|W^{*}|.
\]
It follows that each word representing an element of the group 
$\mathfrak{R}$ generated by the $R_i$ ($i=1,2$) and their transforms is
converted into the empty word by reduction in $\mathfrak{F}$. Because
each such word results from a word
\[
R=L_1 R^{\varepsilon_1}_{\alpha_1} L^{-1}_1
    L_2 R^{\varepsilon_2}_{\alpha_2} L^{-1}_2 \cdots
    L_m R^{\varepsilon_m}_{\alpha_m} L^{-1}_m
\quad (\alpha_i=1,2;\;\varepsilon_i=\pm 1)
\] 
by reduction in the free group on $S_1,S_2$. However, these $R$ result
from $W_0$ by successive elementary transformations and insertions of
factors $(S^{\varepsilon}_\alpha)^{a_\alpha}$. Further, if $W$ and $W'$
are two words that denote the same word in $\mathfrak{F}$, then $W'$
may be converted to the form $WR$ by elementary transformations in
$\mathfrak{S}$, where $R$ belongs to $\mathfrak{R}$, so $W$ and $W'$
go to the same word by reduction in $\mathfrak{F}$.

Groups with generators
\[
S_1,\quad S_2,\quad \ldots,\quad S_n
\]
and defining relations
\begin{equation}
R_i=S^{a_i}_i\quad (i=1,2,\ldots,n)
\tag{4}
\end{equation}
can be handled quite analogously. It follows easily from the solution of 
the word problem that if $S$ is an element of $\mathfrak{F}$ of finite
order, then
\[
S=LS^{s}_i L^{-1}.
\]

The word problem in the group $\mathfrak{F}'$ with two generators
$S_1$ and $S_2$  and a single relation
\begin{equation}
R=S^{a_1}_1 S^{a_2}_2 \tag{5}
\end{equation}
may be easily reduced to the case treated above.\footnote{\textsc{M. Dehn},
Math. Ann. \textbf{75} (1915) 402 and \textsc{O. Schreier} loc.cit.}
Here $S^{a_1}_1=S^{-a_2}_2$, from which it follows that the element
$S^{a_1}$ commutes with all elements of $\mathfrak{F}'$, because
\[
S_1 S^{a_1}_1=S^{a_1}_1 S_1
\]
and
\[
S_2 S^{a_1}_1=S_2 S^{-a_2}_2=S^{-a_2}_2 S_2=S^{a_1}_1 S_2,
\]
so $S^{a_1}$ commutes with all power products of the $S_i$. Each element
may then be converted into a reduced word of the form
\[
S^{r_{11}}_1 S^{r_{21}}_2 \cdots S^{r_{1m}}_1 S^{r_{2m}}_2
S^{ka_1}_1\quad\text{with}\quad 0\le r_{il}<a_i.
\]
One proves quite analogously as for the groups with defining relations (2)
that each word is representable in only one way as a reduced word. From
the solution of the word problem\footnote{Further solutions of word problems
are found in Section 7.14. See also \textsc{W. Magnus}, Math. Ann. 
\textbf{105} (1931) 52 and \textbf{106} (1932) 295; \textsc{E. Artin},
Hamb. Abhdl. \textbf{4} (1925) 47; \textsc{K. Reidemeister}, ibid.
\textbf{6} (1928) 56; \textsc{M. Dehn}, Math. Ann. \textbf{72} (1912) 41.} 
one easily obtains that the subgroup of $\mathfrak{F}'$ generated by 
$S^{a_1}$ is the center of $\mathfrak{F}'$.

\section{The free product of groups}

The methods of Section 2.2 may be extended without difficulty to the
so-called free product\footnote{\textsc{O. Schreier}, Hamb. Abhdl.
\textbf{5} (1927) 16.} of groups. Let $\mathfrak{G}_1$ and
$\mathfrak{G}_2$ be two groups with the elements $G_{1i}$ and
$G_{2i}$ respectively. From these elements we construct words
\[
W=G_1 G_2 \cdots G_n,
\]
where the $G_i$ are any elements of $\mathfrak{G}_1$ or
$\mathfrak{G}_2$ different from the identity. Thus
\[
G_i=G_{k_i l_i}\quad (k_i=1\text{ or }2).
\]
By an elementary expansion of this word $W$ we mean the insertion
of a word $G_{i1}G_{i2}$ that equals the identity when regarded as a
product in $\mathfrak{G}_i$, or replacement of a letter $G_{k_i l_i}$
by two, $G'_{k_i l_i}G''_{k_i l_i}$, the product of which equals 
$G_{k_i l_i}$ in $\mathfrak{G}_{k_i}$. By a reduction we mean the 
reverse process.

Again the words $W$ may be divided into equivalence classes $[W]$
and the product defined as in Section 2.2 by
\[
[W_1][W_2]=[W_1 W_2].
\]
The resulting group $\mathfrak{G}=\mathfrak{G}_1 * \mathfrak{G}_2$
is called the \emph{free product\footnote{Reidemeister uses the notation
$\mathfrak{G}_1 \times \mathfrak{G}_2$, which I have dropped because of
its potential for confusion with the direct product. (Translator's note.)} 
of $\mathfrak{G}_1$ and
$\mathfrak{G}_2$}. One can define the free product of any number of
groups by iteration.

The free group with $n$ free generators $S_i$ is the free product of the
$n$ infinite cyclic groups generated by the $S_i$. The groups (4) of
Section 2.6 are the free products of $n$ finite cyclic groups generated
by the $S_i$ with $S^{a_i}=1$. This construction is important for the 
word problem, because one can obviously solve the problem in a free
product $\mathfrak{G}$ as soon as it is solved in the original groups
$\mathfrak{G}_i$. This is because the reduced word $|W|$ may be 
defined analogously as in Section 2.3---a word is called reduced when
any two neighboring factors $G_i,G_{i+1}$ do not belong to the same
group---and it is then demonstrable that each class of reduced words
contains only one in reduced form.

\emph{If}
\[
S_{1k}\quad(k=1,2\ldots,n_1),\qquad S_{2k}\quad(k=1,2,\ldots,n_2)
\]
\emph{are systems of generators for the groups $\mathfrak{G}_1$
and $\mathfrak{G}_2$, and}
\[
R_{1l}(S_{1k})\quad(l=1,2,\ldots,m_1),\qquad
R_{2l}(S_{2k})\quad(k=1,2,\ldots,m_2)
\]
\emph{are the respective sets of defining relations of $\mathfrak{G}_1$,
$\mathfrak{G}_2$, then all the $S_{ik}$ and all the $R_{il}$
together constitute a system of generators and defining relations for
the free product $\mathfrak{G}$}. It is clear that the $R_{il}$ are
satisfied in $\mathfrak{G}$. On the other hand, one can carry out 
expansion and reduction of the word $W$, where $G_i$ is now
viewed as a power product of the $S_{ik}$, on the basis of the
relations $R_{il}$, because these operations take place only between
elements of the same group $\mathfrak{G}_i$. Hence the
$R_{il}(S_{ik})$ are in fact the defining relations of $\mathfrak{G}$.

It follows conversely that, given a group with generators $S_{1k}$
and $S_{2k}$ and a system of defining relations that can be
divided into two classes $R_{1l}$ and $R_{2l}$, in which only the
$S_{il}$ appear in the $R_{ik}$, then the group in question is the free
product of the subgroups generated by the $S_{1k}$ and the $S_{2k}$.

The concept of the free product may be extended in the following way.
The group $\mathfrak{G}_1$ may possess a subgroup $\mathfrak{U}_1$
that is isomorphic to a subgroup $\mathfrak{U}_2$ of $\mathfrak{G}_2$.
Let $\mbox{\boldmath$I$}(\mathfrak{U}_1)=\mathfrak{U}_2$ be a
specific isomorphism between the $\mathfrak{U}_i$. Under these
assumptions we add the following process to expansion and reduction
of words (1): if $G_i$ is an element of $\mathfrak{U}_1$ or
$\mathfrak{U}_2$ then $G_i$ may be replaced by 
$\mbox{\boldmath$I$}(G_i)$ from $\mathfrak{U}_2$ or
$\mbox{\boldmath$I$}^{-1}(G_i)$ from $\mathfrak{U}_1$.

Classification of words can again be carried out and it leads, again with
the help of equation (2), to the definition of a group $\mathfrak{G}$,
which may be called the \emph{free product of $\mathfrak{G}_1$ and
$\mathfrak{G}_2$ with the subgroups $\mathfrak{U}_1$ and
$\mathfrak{U}_2$ amalgamated}.

A uniquely determined normal form may now be produced as follows:
in the groups $\mathfrak{G}_i$ we choose a system of representatives
for the residue classes modulo $\mathfrak{U}_1$, say
$\mathfrak{U}_1 N_{1k}$, and modulo $\mathfrak{U}_2$, say
$\mathfrak{U}_2 N_{2k}$,  and then one can put each word $W$ in
the form
\[
UN_1 N_2 \cdots N_n,
\]
where $U$ belongs to $\mathfrak{U}_1$, the $N_l$ are certain
representatives $N_{ik}$, and two neighboring $N_i,N_{i+1}$ do not
belong to the same group $\mathfrak{G}_l$.

From this one can solve the word problem in $\mathfrak{G}$ if one
can give each element $G_i$ in $\mathfrak{G}_l$the representation 
$U_i N_{ik}$ ($i=1,2$). If $S_{1k}$ ($k=1,2,\ldots,u$) are the generators of
$\mathfrak{U}_1$, $S_{1k}$ ($k=1,2,\ldots,n_1$) the generators of
$\mathfrak{G}_1$, and analogously if $S_{2k}$ ($k=1,2,\ldots,u$)
are the generators of $\mathfrak{U}_2$, and
$S_{2k}$ ($k=1,2,\ldots, n_2$)
are those of $\mathfrak{G}_2$; and if $R_{1l}(S_{1k})$
($l=1,2,\ldots,m_1$) and $R_{2l}(S_{2k})$ ($l=1,2,\ldots,m_2$) are
the defining relations of $\mathfrak{G}_1$ and $\mathfrak{G}_2$
respectively; and if finally the mapping
\[
\mbox{\boldmath$I$}(S_{1k})=S_{2k}\quad (k=1,2,\ldots,u)
\]
is an isomorphism between $\mathfrak{U}_1$ and $\mathfrak{U}_2$;
then the $S_{ik}$ ($k=1,2,\ldots,n_i;\; i=1,2$), together with the
relations $R_{il}(S)$, ($l=1,2,\ldots,n_i;\: i=1,2$) and $S_{1k}=S_{2k}$,
($k=1,2,\ldots,u$)
are generators and relations for the free product with amalgamated
subgroup, as one may prove analogously with the theorem on the 
free product itself.

\section{A transformation problem}

The groups treated in Section 2.6 admit an easy solution of the 
transformation problem. However, for what follows we will need only
the special case of the relations\footnote{\textsc{K. Reidemeister},
Hamb. Abhdl. \textbf{8} (1930), 187.}
\[
R_1=S^3_1,\quad R_2=S^2_2.
\]
We alter the normal form of Section 2.6 by always writing $S^{-1}_1$
in place of $S^2_1$. If $\varepsilon_i=\pm 1$ then each element
different from the identity can be brought into one of the following 
reduced forms
\begin{equation}
W=S^{\varepsilon_1}_1 S_2 S^{\varepsilon_2}_1 S_2 \cdots 
S^{\varepsilon_m}_1;\quad
WS_2;\quad S_2 W;\quad S_2 W S_2;\quad S_2 \tag{1}
\end{equation}
By $W^{-1}$ we mean the power product formally inverse to $W$,
\[
W^{-1}=S^{-\varepsilon_m}_1 \cdots S_2 S^{-\varepsilon_2}_1
S_2 S^{-\varepsilon_1}_1,
\]
and similarly for the other reduced products. Now for the solution of
the transformation problem we remark that the first and last factors of
a product (1) are either a) formally inverse to each other or b) not.

In the first case a) we can put the product in the form
\[
H=LH'L^{-1}
\]
where $L$ and $L^{-1}$ are formally inverse to each other and where
the kernel $H'$ of the product begins and ends with factors that are not
formally inverse to each other. The kernel $H'$ has the form $W$ of (1)
with $\varepsilon_1=\varepsilon_m$ when $L$ ends with $S_2$, but it
contains only one factor $S_2$ when $L$ ends with $S^\varepsilon_1$.
In the second case b) the product $H$ has one of the forms $S_2$,
$W$ with $\varepsilon_1=\varepsilon_m$, $WS_2$, or $S_2 W$.

We will call the products $S^\varepsilon_1$, $S_2$, $WS_2$, and $S_2 W$
short words of the first kind. The products $S^\varepsilon_1$, $S_2$,
and $W$ with $\varepsilon_1=\varepsilon_m$ will be called short
words of the second kind. Each element has a transformed product
that is a short word of the first kind; this is because it is either a short word 
of the first or second kind or else it has a kernel that is a short word of the
second kind, and a short word of the second kind becomes a short word 
of the first kind by transformation with an $S^\varepsilon_1$ and 
reduction.

We now let $K$ denote a short word of the first kind and let $\{K\}_1$
denote the collection of products that result from $K$ by cyclic 
interchange of factors. By $\{K\}_2$ we mean the collection of short 
words of the second kind that result from a word
$WS_2=S^{\varepsilon_1}_1 S_2\cdots S_2$ out of $\{K\}_1$ by the
process
\[
S^{\varepsilon_1}_1 WS_2 S^{-\varepsilon_1}=
S^{-\varepsilon_1}_1 S_2 S^{\varepsilon_2}_1 S_2 \cdots 
S^{\varepsilon_m}_1 S_2 S^{-\varepsilon_1}_1,
\]
as well as those short words of the first kind from $\{K\}_1$ that are
also of the second kind. By $\{K\}_3$ we mean all those words $H$
that have a kernel $H'$ in $\{K\}_2$. Finally, let $\{K\}$ denote the
totality of elements from the classes $\{K\}_i$ ($i=1,2,3$).

Each element obviously belongs to exactly one class $\{K\}$. Further,
it is clear on the one hand that any two products in $\{K\}$ are
convertible into each other by transformation and reduction, and hence
they denote transforms of each other in our group, while on the other
hand, if $H$ is any word in $\{K\}$ then 
$S^\varepsilon_1 H S^{-\varepsilon}_1$ and $S_2 H S_2$ yield other
words in $\{K\}$ by reduction. One verifies this by considering the
cases where $H$ lies in $\{K\}_1$, $\{K\}_2$, or $\{K\}_3$. It follows
in general that $MHM^{-1}$ yields a word, by reduction, that lies in the
same class $\{K\}$ as $H$.

Now, on the one hand, we can decide whether two reduced products
belong to the same class $\{K\}$, and on the other hand each element
of our group corresponds to a unique reduced product, so the
transformation problem is solved.

One more remark about the powers of an element $H$. If  $H$
belongs to a class $\{K\}_i$, then each power $H^k$ belongs to a class
$\{\overline{K}\}_i$ with the same index $i$.

\section{Generators and relations for the modular group}

The modular group defined in Section 1.4 is isomorphic to the group
discussed in the previous section. Thus we have solved the transformation
problem for the modular group.

One can of course also solve the transformation problem by proceeding
from the arithmetic representation of the substitutions and asking what
conditions the coefficients
\[
a,b,c,d\quad\text{and}\quad a',b',c',d'
\]
must satisfy for the associated substitutions to be transformable into each
other in the modular group. However, this way is much more difficult. It is
connected with the question of when two binary quadratic forms
\[
Ax^2+Bxy+Cy^2\quad\text{and}\quad A'x'^2+B'x'y'+C'y'^2
\]
are equivalent, i.e., when there there are integers
\[
a,b,c,d\quad\text{with}\quad ad-bc=1
\]
such that the unprimed form goes to the primed form when $x,y$ are
replaced by
\[
x=ax'+by',\quad y=cx'+dy'.
\]

We now apply ourselves to the proof that \emph{the modular group is
generated by two elements $S_1$ and $S_2$ which satisfy the relations}
\[
R_1=S^{3}_1,\quad R_2=S^{2}_2
\]
\emph{and no others independent of them}.

We let $T$ be the substitution
\[
x'=x+1,
\]
so that $T^n$ is the substitution
\[
x'=x+n.
\]
By $S$ we mean
\[
x'=-\frac{1}{x}.
\]
If $A$ is the substitution
\[
x'=\frac{ax+b}{cx+d}\quad\text{with}\quad |b|\ge|d|>0,
\]
and
\[
A'=T^n A
\]
corresponds to the substitution
\[
x''=\frac{a'x+b'}{c'x+d'},
\]
then
\[
b'=b+nd
\]
and hence by suitable choice of $n$ one can obtain
\[
|b'|<|d|\le |b|.
\]
If
\[
0<|b|<|d|
\]
then the substitution $SA$ or
\[
x'=\frac{cx+d}{-ax-b}
\]
satisfies the previous condition. Hence it follows by induction that: for each
transformation $A$ there is a power product
\[
M=S^\varepsilon T^{n_1} S T^{n_2} S \cdots T^{n_m} S^\eta
\]
($\varepsilon$ and $\eta$ equal 0 or 1) such that, in the transformation
\[
x'=\frac{ax+b}{cx+d}
\]
corresponding to $MA$, we must have $d=0$. It must then be that $-bc=1$,
i.e.,
\[
x'=-\frac{1}{x}+a,
\]
and this is $T^a S$. Consequently, $S$ and $T$ are generators of the
modular group.

Now we set $S_1=TS$, $S_2=S$ and confirm that $S^{2}_2=1$.
Further, $S_1$ corresponds to the substitution
\[
x'=-\frac{1}{x}+1=\frac{x-1}{x}.
\]
$S^{2}_1$ corresponds to
\[
x''=-\frac{1}{-\frac{1}{x}+1}+1=\frac{x}{x-1}+1=\frac{1}{-x+1},
\]
so that $S^{3}_1=1$. Since $T=S_1 S^{-1}_2$, $S_1$ and $S_2$ are
also generators of the modular group. They satisfy the two given relations,
and it remains only to prove that they satisfy no other relations apart from
consequences of $R_1$ and $R_2$. We will show that, if one computes the
substitution
\[
x'=\frac{ax+b}{cx+d}
\]
for a reduced word (1) from Section 2.8 in which one replaces the $S_i$ by
the corresponding modular substitutions, then it is never the identity substitution.
It suffices to prove this for words of the form $WS_2$ since, by Section 2.8,
each element of the group may be converted into a word $WS_2$ by
transformation with $S^\varepsilon_1$ or $S_2$.

For the proof we convert $WS_2$ back to a certain power product of $S$ and 
$T$. Namely, we combine all neighboring elements $S_1 S_2$ into powers
$(S_1 S_2)^{\delta_i}$ and likewise the elements $S^{-1}_1S_2$ into
powers $(S^{-1}_1 S_2)^{\delta_k}$ and then set
\[
(S_1 S_2)^{\delta_i}=T^{\delta_i},\quad
(S^{-1}_1 S_2)^{\delta_k}=ST^{-\delta_k}S,\quad
(\delta_i,\delta_k>0).
\]
One sees that this gives a product in $S$ and $T$ in which the exponents
have alternating signs. But it is easy to see that such an element is never the
identity substitution by computing the coefficients of the corresponding
modular substitution.\footnote{Cf. \textsc{Dirichlet-Dedekind}, \emph{
Vorlesungen \"uber Zahlentheorie}, 2nd edition, 1871, \S81.}

Another method of detemining generators and defining relations for the
modular group consists in the construction of its fundamental domain in the
complex number plane.\footnote{Cf. a textbook on function theory, e.g., that
of \textsc{Bieberbach}, vol. II.}

\section{A theorem of \textsc{Tietze}}

It is clear that a group may be defined in various ways by generators and
relations. If
\[
S_1,\quad S_2,\quad \ldots,\quad S_m
\]
is a system of generators for a group $\mathfrak{F}$ and the set
$\mathfrak{r}$ of products
\[
R_1(S),\quad R_2(S),\quad \ldots,\quad R_r(S)
\]
in the $S_i$ is a system of defining relations, and if $R_{r+1}(S)$ is any
consequence of these relations, then, e.g., the set that results from
$\mathfrak{r}$ by addition of $R_{r+1}$ is also a system of defining
relations. If, on the other hand, $R_r(S)$ is a consequence of
$R_1(S)_,R_2(S)\ldots$, $R_{r-1}(S)$ then the latter set is also a system 
of defining relations for $\mathfrak{F}$.

Further, if $T$ is a letter denoting any power product of the $S_i$,
\[
T=T(S),
\]
then
\[
R_{r+1}=T(S)T^{-1}
\]
is a relation, and 
\[
S_1,\quad S_2,\quad\ldots,\quad S_m,\quad T
\]
is a system of generators and, as we will show,
\[
R_1(S),\quad R_2(S),\quad \ldots,\quad R_r(S),\quad R_{r+1}(S,T)
\]
is a system of defining relations. This is because each relation containing
only the $S$ is a consequence of the $R_i$ ($i=1,2,\ldots,r$) and,
using the relation $R_{r+1}$, each power product containing a factor
$T$ may be converted into one in the $S$ alone. Namely, if
\[
F=A(S) T B(S,T)
\]
then
\begin{align*}
F&=A(S)TT^{-1}(S)T(S)B(S,T)\\
 &=A(S)R^{-1}_{r+1}T(S)B(S,T)\\
 &=A(S)T(S)B(S,T)
\end{align*}
and the latter product contains one $T$ factor fewer than $F$ does.
In this way the factors $T^\varepsilon$ ($\varepsilon=\pm 1$)
may be removed successively.

On the other hand, if $S_m$ is representable as a power product of
$S_1,S_2,\ldots,S_{m-1}$ then $S_1,S_2,\ldots,S_{m-1}$ obviously
constitute a system of generators. One can successively eliminate $S_m$
from all power products. Further, if the defining relations $R_1,R_2,$
$\ldots,
R_{r-1}$ contain only the generators $S_1,S_2,\ldots,S_{m-1}$ and if
\[
R_r=S_m(S_1,S_2,\ldots,S_{m-1})S^{-1}_m
\]
then the $R_i$ ($i=1,2,\ldots,r-1$) constitute a system of defining relations
in the generators
\[
S_1,\quad S_2,\quad \ldots,\quad S_{m-1}.
\]
This is because the group defined by the
\[
S_i\; (i=1,2,\ldots,m-1),\quad R_k\; (k=1,2,\ldots,r-1)
\]
is, as we saw above, identical with that defined by
\[
S_i\; (i=1,2,\ldots,m),\quad R_k\; (k=1,2,\ldots,r).
\]

We now have an important theorem (of 
\textsc{Tietze}\footnote{\textsc{H. Tietze}, Mon. f. Math. u. Phys. 
\textbf{19}, p. 1.}) that \emph{any two systems of generators and
defining relations for the same group are always convertible to each
other by successive applications of the transformations above.}

Let 
\begin{align*}
S_1,S_2,\ldots,S_m;&\quad R_1(S),R_2(S),\ldots,R_r(S) \tag{1}\\
S'_1,S'_2,\ldots,S'_{m'};&\quad R'_1(S'),R'_2(S'),\ldots,R'_{r'}(S') \tag{2}
\end{align*}
be two systems of generators and defining relations for the same group
$\mathfrak{F}$. The $S'_k$ must then be expressible in terms of the
$S_i$ and, conversely, the $S_i$ in terms of the $S'_k$. If
\[
S'_k=S'_k(S);\quad S_i=S_i(S')
\]
then we set
\[
U_k(S,S')=S'_k(S)S'^{-1}_k;\quad
V_i(S,S')=S_i(S')S^{-1}_i.
\]
Obviously the $S_i,S'_k$ are a system of generators and the relations
\begin{equation}
R_l(S),\quad U_k(S,S') \tag{3}
\end{equation}
on the one hand, as well as the relations
\[
R'_l(S'),\quad V_i(S,S') \tag{4}
\]
on the other, are systems of defining relations for $\mathfrak{F}$
that result from (1) and (2) respectively by successive addition of the
respective generators $S'_k$ and $S_i$ with the respective relations
$U_k$ and $V_i$.

But now the relations (4) must be consequences of (3), because the
relations are indeed relations in the $S_i,S'_k$. Similarly, the
relations (3) are consequences of (4). Hence by addition of
consequence relations we can extend both systems, (3) and (4), to
the same system
\begin{equation}
S_i,S'_k;\quad R_l(S'),R'_l(S),U_k(S,S'),V_i(S,S'), \tag{5}
\end{equation}
and hence convert the system (1) to the system (2) by a sequence of
the transformations described. 

One can apply this theorem to a purely \emph{combinatorial
characterization of the properties of a group} given by generators $S_i$
and defining relations $R_k$. Each property of a system of generators
$S_i$ and relations $R_k$ that is invariant under the above 
transformations of the $S_i$ and the $R_k$ is a property of the group
$\mathfrak{F}$ defined by $S_i,R_k$. This is because such a property
holds for all presentations of the group by generators and relations and
hence it is a property of the group itself. Despite this simple connection
between different presentations of the same group it is in general not
possible to decide whether two groups presented by generators and
relations are isomorphic to each other.\footnote{This remarkable claim
was first made by Tietze (1908) in the paper cited above. At the time when
Reidemeister wrote, a precise concept of algorithm---formalizing what
it means to ``decide''---was still a few years away from being 
published. It first appeared in publications of Church, Post, and most
convincingly by Turing in 1936. Another two decades elapsed before
Adyan and Rabin proved that the isomorphism problem is 
algorithmically unsolvable, in 1958. Their work also established the 
unsolvability of the problems next mentioned by Reidemeister: deciding  
whether a given finitely-presented group is free, or trivial. It may be
worth mentioning that Reidemeister could have had some intimation of
the coming wave of unsolvability results, because he organized
the conference in K\"onigsberg in 1930 at which G\"odel first
announced his famous (and related) result on the incompleteness
of formal systems. (Translator's note.)} One also cannot decide whether
such a group is a free group on ``non-free'' generators, or whether it
follows from the relations $R_k(S)$ that all the $S_i$ equal the identity
$E$.

We make a simple application of the transformation rules to the
modular group presentation by the generators $S_1,S_2$
and relations
\[
R_1=S^{3}_1\equiv 1;\quad R_2=S^{2}_2\equiv 1.
\]
As we have seen, the operations $S$ and $T$ defined in Section 2.9 
also generate the modular group. We now ask what are the defining
relations in the group generated by $T=S_1S^{2}_2$ and $S=S_2$.
For this purpose we take $T$ as a generator in addition to $S_1$ and
$S_2$ and add
\[
R_3=S_1 S^{-1}_2 T^{-1}
\]
as a third relation.

With the help of this equation we now eliminate $S_1$ from $R$, by first 
constructing $R'_1=R^{-1}_3 R_1=TS_2 S^{2}_1$ and then deriving
$R_1$ as a consequence of $R'_1$ and $R_3$. Then we replace $R'_1$
by $R''_1=S_1 T S_2 S_1$, and this in turn by
\[
R'''_1=R^{-1}_3 R''_1 S^{-1}_1 R^{-1}_3 S_1=(TS_2)^2.
\]
In this way we obtain the defining relations of the modular group in the
generators $S_2=S$ and $T$ as
\[
R_1=(TS)^3\equiv 1\quad\text{and}\quad R_2=S^2\equiv 1.
\]

\section{Commutative groups}

We will use the theorem of \textsc{Tietze} to characterize the
commutative or ``abelian'' group $\mathfrak{F}$ with finitely many
generators and relations through properties of these relations. In a
commutative group with generators $S_i$ ($i=1,2\ldots,n$) each of
the relations
\begin{equation}
R_{ik}(S)=S_i S_k S^{-1}_i S^{-1}_k \tag{1}
\end{equation}
holds, since this says that $S_i$ and $S_k$ commute with each other.
It follows that all power products of the $S_i$ commute with each other.
Hence each relation $R(S)$ may be brought into the form
\[
R(S)=S^{r_1}_1 S^{r_2}_2 \cdots S^{r_n}_n.
\]
Thus we can take the system of defining relations to be in the form
\begin{equation}
R_i(S)=S^{r_{i1}}_1 S^{r_{i2}}_2 \cdots S^{r_{in}}_n
\quad (i=1,2,\ldots m). \tag{2}
\end{equation}
The characteristic properties of a particular commutative group must
then reside in the relations (2), because the relations (1) are satisfied
in any commutative group. We now construct the matrix
\[
\rho=(r_{ik})\quad (i=1,2,\ldots,m;\; k=1,2,\ldots,n)
\]
and show that $\mathfrak{F}$ has certain characteristic numbers that
may be derived from $\rho$, the so-called elementary divisors of $\rho$.

By $\delta^{(k)}_i$ ($i=1,2,\ldots,k\le m,n$) we mean the collection
of $k$-rowed subdeterminants obtainable from $\rho$ by striking out
$m-k$ rows and $n-k$ columns. If all $\delta^{(s+1)}_i=0$ while there is
a $\delta^{(s)}_i\ne 0$ then $s$ is called the rank of $\rho$. By
$\delta^{(k)}>0$ we mean the greatest common divisor of all the
$\delta^{(k)}_i$ for $k\le s$. Then $\delta^{(k)}$ is always divisible by
$\delta^{(k-1)}$, because all the $k$-rowed determinants are linear 
combinations of $(k-1)$-rowed determinants. We now set
\[
d_1=\delta^{(1)};\quad \delta^{(k)}=d_k \delta^{(k-1)}\quad
(k=2,3,\ldots,s)
\]
and call $d_k$ the $k$th elementary divisor of $\rho$. We claim

Theorem 1. \emph{The $d_k\ne 1$ and $n-s$ are the same for all
relation systems for $\mathfrak{F}$}.

Theorem 2. \emph{New generators}
\[
T_1,\quad T_2,\quad \ldots,\quad T_n
\]
\emph{may be introduced, for which the defining relations take the form}
\[
R_i(T)=T^{d_i}_i\quad (i=1,2,\ldots,s).
\]

By leaving out the generators $T_j$ for which $d_j=1$ one obtains a
unique normal form for $\mathfrak{F}$; the number of relation-free
generators is $n-s$. On the basis of Theorems 1 and 2 we then have:
\emph{$\mathfrak{F}$ is characterized by the elementary divisors of
the matrix $\rho$ different from $1$ and the difference $n-s$ between
the number of generators and the rank of $\rho$}.

\section{A theorem on matrices}

To prove the theorems of the last section we first define an equivalence
of matrices with respect to the following transformations. The matrix
$\rho=(r_{ik})$ is called equivalent to the matrix $\rho'=(r'_{ik})$
\begin{enumerate}
\item
if $\rho'$ results from $\rho$ by an exchange of rows or columns,
\item
if $\rho'$ results from $\rho$ when the elements $r_{1i}$ of the first
row are replaced by $r'_{1i}=r_{1i}+ar_{2i}$, or when the elements
$r_{i1}$ of the first column are replaced by$r'_{i1}=r_{i1}+ar_{i2}$
($a$ an arbitrary integer), while all the remaining rows or columns
remain unaltered,
\item
if all elements in some row or column have their signs reversed,
\item
if there is a chain of matrices 
$\rho_1=\rho,\rho_2,\rho_3,\ldots,\rho_n=\rho'$
in which $\rho_{i+1}$ results from $\rho_i$ by one of the elementary transformations 1, 2, or 3. [Again, this $n$ does not denote the number of
generators.] E.g. it is a permissible transformation to
add $k$ times a row, or column, to any other row or column, 
respectively.
\end{enumerate}

Theorem 1. \emph{Equivalent matrices have the same rank and the same
elementary divisors}.

This is clear for matrices convertible into each other by the transformation 1.
If $\rho'$ results from $\rho$ by a row transformation 2, then any of its
determinants $\delta'^{(k)}_i$ results from $\delta^{(k)}_i$ when
$r_{ik}$ is replaced by $r'_{ik}$, hence it equals $\delta^{(k)}_i$ if the
first row does not contribute any elements to $\delta^{(k)}_i$.
Otherwise, we expand $\delta'^{(k)}_i$ along the first row and obtain
\[
\delta'^{(k)}_i=\delta^{(k)}_i\quad\text{or}\quad
\delta'^{(k)}_i=\delta^{(k)}_i+a\delta^{(k)}_j
\]
according as the second row appears in $\delta'^{(k)}_i$ or not. It follows
that the rank $s'$ of $\rho'$ satisfies $s'\le s$ and that the elements
$\delta'^{(k)}$ of $\rho'$ are divisible by $\delta^{(k)}$. But since $\rho$
also results from $\rho'$ by a row transformation 2, because
$r_{1i}=r'_{1i}-ar'_{2i}$, it follows that $s=s'$ and 
$\delta^{(k)}=\delta'^{(k)}$. Hence we have Theorem 1 for arbitrary
equivalent matrices.

To clarify the meaning of the $d_k$ we now assert:

Theorem 2. \emph{The matrix $\rho$ is equivalent to the matrix
$\delta=(d_{ik})$, where $d_{ik}=0$ if $i\ne k$, $d_{ii}=d_i$
($i=1,2,\ldots,s$) and $d_{ii}=0$ for $i>s$}.

We first prove the following Lemma 1: \emph{if $|r_{ik}|>d_1$ for all
nonzero $r_{ik}$, then there is a matrix equivalent to $\rho$ that contains
a nonzero $r'_{i_1 k_1}$ smaller than all $|r_{ik}|$}.

Namely, let $r_{i_1 k_1}$ be a term of $\rho$ of smallest absolute value:
\[
|r_{ik}|\ge|r_{i_1 k_1}|.
\]
Now suppose there is either an element $r_{i_1 k_2}$ of the $i_1$th row, 
or an element $r_{i_2 k_1}$ of the $k_1$th column, which is nonzero and 
not divisible by $r_{i_1 k_1}$. Then by subtraction of a suitable multiple,
either of the $k_1$th column from the $k_2$th column, or of the $i_1$th
row from the $i_2$th row, we obtain a matrix with the property claimed.

If, on the other hand, all the $r_{i_1 k}$ and $r_{ik_1}$ are divisible by
$r_{i_1 k_1}$, then one can construct, by elementary transformations,
an equivalent matrix $\rho'$ in which all elements $r'_{i_1 k}=r'_{ik_1}=0$,
except for $r'_{i_1 k_1}=r_{i_1 k_1}$. Then if there is an $r'_{ik}$ with
\[
|r'_{ik}|<|r_{i_1 k_1}|
\]
there is nothing more to prove. If all
\[
|r'_{ik}|\ge|r_{i_1 k_1}|
\]
then certainly not all $r'_{ik}$ are divisible by $r_{i_1 k_1}$, otherwise
\[
d_1=d'_1<|r_{i_1 k_1}|.
\]
If $r'_{i_2 k_2}$ is not divisible by $r_{i_1 k_1}$ then I construct $\rho''$
by adding the $k_2$th column to the $k_1$th column, whence
\[
r''_{ik_1}=r'_{ik_2}\;(i\ne i_1);\quad r''_{i_1 k_1}=r_{i_1 k_1},
\]
and the second case is reduced to the first.

From this we get Lemma 2: \emph{For each matrix $\rho$ there is an equivalent $\rho'=(r'_{ik})$ with}
\[
r'_{11}=d_1;\quad r'_{1i}=r'_{i1}=0\; (i\ne 1).
\]

Firstly, by Lemma 1 there is an equivalent matrix containing an element
equal to $\pm d_1$. I can bring this element into the first row and first 
column, and then make all other elements of the first row and column zero 
by subtraction of suitable multiples of the first row and column. This is the
desired matrix $\rho'$.

By $\rho^{*}$ we mean the matrix that results from $\rho'$ by striking
out the first row and first column. The rank of $\rho^{*}$ equals $s-1$,
essentially because any $l$-rowed nonzero determinant from $\rho^{*}$
can be used to construct an $(l+1)$-rowed determinant of $\rho'$ that
is likewise nonzero.

Now Theorem 2 comes about as follows:

Let $d^{*}_2=F_{22}$ be the greatest common divisor of the nonzero
$r^{*}_{ik}$. Then $d_1$ is a divisor of $d^{*}_2$, since $d_1$ is a
divisor of all $r^{*}_{ik}$, and by Lemma 2 there is a matrix
${\rho^{*}}'=(r_{ik}^{*'})$ equivalent to $\rho^{*}$ with
\[
r_{11}^{*'}=d^{*}_2;\quad r_{1i}^{*'}=r_{k1}^{*'}=0\quad
(i,k\ne 1)
\]
But then $\rho$ itself is equivalent to the matrix $\rho''=(r''_{ik})$
with
\begin{align*}
&r''_{11}=d_1;\quad r''_{22}=d^{*}_{22}\\
&r''_{1i}=r''_{k1}=0\quad (i,k\ne 1);\quad r''_{2i}=r''_{k2}=0
\quad (i,k\ne 2)\\
&r''_{i+2,k+2}=r^{*}_{ik}\quad (i,k>1).
\end{align*}
By iteration of this process we find that $\rho$ is equivalent to a
matrix $\overline{\rho}=(F_{ik})$ with
\[
F_{ik}=0\quad (i\ne k);\quad \overline{r}_{ii}>0\quad
(i=1,2,\ldots,s),\quad\overline{r}_{ii}=0 \quad (i>s)
\]
and $\overline{r}_{ii}$ is a divisor of $F_{i+1,i+1}$.

But the $\overline{r}_{ii}$ are the elementary divisors of 
$\overline{\rho}$ and hence also of $\rho$, because all $k$-rowed
subdeterminants from $\overline{\rho}$ that are nonzero have a value
\[
\overline{r}_{i_1,i_1}
\overline{r}_{i_2,i_2}
\cdots
\overline{r}_{i_k,i_k}
\]
where all the $i_l$ ($l=1,2,\ldots,k$) are different. Such a product is
divisible by 
$\overline{r}_{11}\overline{r}_{22}\cdots\overline{r}_{kk}$, so
\[
\overline{\delta}^{(k)}=
\overline{r}_{11}\overline{r}_{22}\cdots\overline{r}_{kk}
\]
and hence
\[
\overline{r}_{kk}=\overline{d}_k=d_k.
\]
One more remark: if $d_i\ne 1$ then also $d_{i+1}\ne 1$, because
$d_i$ is a divisor of $d_{i+1}$.

\section{Characterization of commutative groups}

We now return to the commutative group $\mathfrak{F}$ and see how the
matrix $\rho=(r_{ik})$ of the exponents $r_{ik}$ in the defining relations
(2) of Section 2.11 are altered when we transform the generators and
defining relations as in Section 2.10.

If $R$ is a consequence relation of the $R_i$ then, by means of the
relations (1) of Section 2.11 for exchange of factors, $R$ may be written
on the one hand as a power product 
$R^{p_1}_1 R^{p_2}_2 \cdots R^{p_m}_m$ and on the other hand it 
may be brought into the form
\[
S^{r_1}_1 S^{r_2}_2 \cdots S^{r_n}_n.
\]
We therefore must have
\[
r_i=\sum^{m}_{k=1} p_k r_{ki}.
\]
Thus if we extend the defining relations $R_i$ by addition of $R_{m+1}=R$
and construct the matrix of coefficients $\rho'=(r'_{ik})$ for the new system,
then
\[
r'_{ik}=r_{ik};\quad
r_{m+1,k}=\sum^{m}_{i=1} p_i r_{ik};\quad
(i=1,2,\ldots,m;\; k=1,2,\ldots,n).
\]

One can now replace $\rho'$ by an equivalent matrix $\rho''$ that contains
only zeros in the $(m+1)$th row by successively subtracting $p_i$ times
the $i$th row from the last row. Since $\rho$ results from $\rho'$ by
omitting the last row, the elementary divisors and rank of $\rho$ and
$\rho''$, and hence also of $\rho$ and $\rho'$, are identical.

Now let $T$ be any power product
\[
T=S^{q_1}_1 S^{q_2}_2 \cdots S^{q_n}_n.
\]
Take $T$ as a new generator and
\[
S^{q_1}_1 S^{q_2}_2 \cdots S^{q_n}_n T^{-1}
\]
as a new relation. The new coefficient matrix is then $\rho'=(r'_{ik})$,
where $r'_{ik}=r_{ik}$ ($i=1,2,\ldots,m;\; k=1,2,\ldots,n$); $r_{i,n+1}=0$
($i\ne m+1$); $r_{m+1,k}=q_k$ ($k\ne n+1$); and $r_{m+1,n+1}=-1$.
We can convert $\rho'$ into a matrix $\rho''$ by successively adding $q_k$
times the last column to the $k$th column. In $\rho''$ all elements of the
$(m+1)$th row and the $(n+1)$th column apart from $r''_{m+1,n+1}=-1$
are zero. Now $d''_1$ is certainly equal to 1, because $r''_{m+1,n+1}=-1$.
Further, $d''_i=d_{i-1}$, because all $i$-rowed nonzero determinants from
$\rho''$ are either determinants from $\rho$ or else they contain the
element $r''_{m+1,n+1}$ and hence are equal to  an $(i-1)$-rowed
determinant from $\rho$. Conversely, from each $(i-1)$-rowed
determinant of $\rho$ we can construct an $i$-rowed determinant of
$\rho''$ with the same absolute value by taking suitable elements from the
$(m+1)$th row and the $(n+1)$th column. Since all $\delta^{(k)}_i$ 
are divisible by $\delta^{(k-1)}$ we have
\[
\delta''^{(k)}=\delta^{(k-1)}.
\]
Therefore the elementary divisors $d'_i$ of $\rho'$ are equal to $d_{i-1}$
for $i>1$, and $d'_1=1$. From the connection between the determinants of
$\rho$ and $\rho''$ it also follows that the rank $s''$ of $\rho''$ is equal to
$s+1$. Consequently, the rank $s'$ of $\rho'$ is also  $s+1$. Theorem 1
of Section 2.11 then follows.

To prove Theorem 2 in Section 2.11 we show that \emph{the transformations
defined in Section $2.12$ may be accomplished for the matrix $\rho$ of
exponents $r_{ik}$ by alteration of the generators and defining relations.}
Transformation 1 may be accomplished by changing the numbering of
generators and relations, and transformation 3 by changing to the inverse of 
a generator or relation. We accomplish transformation 2 by first taking the
consequence relation $R_1 R^{k}_2$,
\[
S^{r_{11}+kr_{21}}_1
S^{r_{12}+kr_{22}}_1
\cdots
S^{r_{1n}+kr_{2n}}_1=R'_1.
\]
But then $R_1=R'_1 R^{-k}_2$ is a consequence relation of
$R'_1,R_2,\ldots,R_m$ and hence may be omitted. We accomplish the
column transformation 2 by taking the new generator $S'_2$ and relation
$R_{m+1}=S'^{-1}_2 S^{-k}_1 S_2$. Then $S_2=S^{k}_1 S'_2$ and
if we now replace $S_2$ in all relations by $S^{k}_1 S'_2$, using 
$R_{m+1}$, then
\[
R'_i=S^{r_{i1}+kr_{i2}}_1 S'^{r_{i2}}_2 \cdots S^{r_{in}}_n.
\]
The $R'_i$ are consequence relations of the $R_i$ ($i=1,2,\ldots,m+1$).
But conversely, the $R_i$ are also consequence relations of the $R'_i$
($i=1,2,\ldots,m$) and $R_{m+1}$, since indeed 
$R_i=R'_i R^{-r_{i2}}_{m+1}$. Consequently, the $R'_i$ and $R_{m+1}$ 
form a system of defining relations, and hence so do the $R'_i$ alone,
when $S_2$ and $R_{m+1}$ are both omitted. Theorem 2 in Section 2.11
now follows from Theorem 2 in Section 2.12.

The \emph{word problem} may be simply solved for a commutative group
in the normal form given by Theorem 2 of Section 2.11. All representations 
of the identity are comprised by
\[
R^{k_1}_1 R^{k_2}_2 \cdots R^{k_s}_s=
T^{k_1 d_1}_1 T^{k_2 d_2}_2 \cdots T^{k_s d_s}_s.
\]
The $d_i$ are zero for $i>s$. If 
\[
T^{n_1}_1 T^{n_2}_2 \cdots T^{n_n}_n\quad\text{and}\quad
T^{n'_1}_1 T^{n'_2}_2 \cdots T^{n'_n}_n
\]
are two words in the $T$, then they are the same element if and only if
\[
n_i\equiv n'_i\quad\text{(mod $d_i$)}.
\]

If all $d_i=0$ the group is called a \emph{free commutative} or a
\emph{free} \textsc{Abelian} \emph{group}. Free \textsc{Abelian} groups
are characterised by the number of their generators.

\section{Commutative groups with operators}

Using the coefficients defined in Section 1.13 for a commutative group with
operator $x$, the concepts of ``generator,'' ``relation,'' and ``defining
relation'' may be extended as follows.\footnote{\textsc{J. W. Alexander},
Trans. Amer. Math. Soc. \textbf{30} (1928), 275.} \emph{The elements}
\begin{equation}
S_1,\quad S_2,\quad \ldots,\quad S_n \tag{1}
\end{equation}
\emph{are called the generators of the commutative group $\mathfrak{F}_x$
with operator when each element of $\mathfrak{F}$ may be written as a
power product}
\begin{equation}
\prod^{n}_{i=1}S^{f_i(x)}_i. \tag{2}
\end{equation}
Such a product is called a relation when it is equal to the identity element of the
group. \emph{The relations}
\[
R_1,\quad R_2,\quad \ldots,\quad R_m
\]
\emph{are called defining relations of $\mathfrak{F}_x$ in the generators
$S$ if each relation $R(S)$ may be derived, by rearrangement of terms,
from a product}
\begin{equation}
\prod^{n}_{i=1}R^{g_i(x)}_i. \tag{3}
\end{equation}

Conversely, \emph{given any system of generators}
\[
S_1,\quad S_2,\quad \ldots,\quad S_n
\]
\emph{and a system of relations}
\begin{equation}
R_i(S)=S^{r_{i1}(x)}_1 S^{r_{i2}(x)}_2 \cdots S^{r_{in}(x)}_n
\tag{$i=1,2,\ldots,m$}
\end{equation}
\emph{there is always a commutative group with operator defined by this
system.} To prove this we introduce new symbols
\begin{equation}
S^{x^k}_i=S_{i,k} \tag{$i=1,2,\ldots,n;\; k=0,\pm 1,\pm 2,\ldots$}
\end{equation}
and set
\[
S^{a_n x^n+a_{n+1}x^{n+1}+\cdots+a_{n+m}x^{n+m}}_i=
S^{a_n}_{i,n} S^{a_{n+1}}_{i,n+1} \cdots S^{a_{n+m}}_{i,n+m}.
\]
The relations $R_l(S_i)$ may be transcribed as relations $R_i(S_{i,k})$
in the $S_{i,k}$. We include all relations $(R_l)^{x^p}=R_{l,p}(S_{i,k})$
$(p=0,\pm 1, \pm 2,\ldots)$ expressed in the $S_{i,k}$. Then there is a commutative group $\mathfrak{F}$ generated by the $S_{i,k}$ and
defined by the relations
\[
R^{x^p}_l=R_{l,p}(S_{i,k}).
\]

In this group the mapping defined on power products $F$ of the $S_{i,k}$
by
\begin{align*}
\mbox{\boldmath{$A$}}(S_{i,k})&=S_{i,k+1}\\
\mbox{\boldmath{$A$}}(F_1 F_2)&=
\mbox{\boldmath{$A$}}(F_1)\mbox{\boldmath{$A$}}(F_2)
\end{align*}
is an automorphism, because it
sends each power product $R_{l,k}$ to the power product $R_{l,k+1}$
and hence each relation goes to another relation. Thus if $F_1$ and $F_2$
are two different power products in the $S_{i,k}$ which denote the same
element of $\mathfrak{F}$, so
\[
F_1=F_2 R
\]
(i.e., $F_1$ is convertible to $F_2 R$ by rearranging and applying
the relation $S^{a}_{i,k} S^{b}_{i,k}=S^{a+b}_{i,k}$), then
\[
\mbox{\boldmath{$A$}}(F_1)=
\mbox{\boldmath{$A$}}(F_2)\mbox{\boldmath{$A$}}(R),
\]
and hence also
\[
\mbox{\boldmath{$A$}}(F_1)\equiv \mbox{\boldmath{$A$}}(F_2).
\]
Likewise, one concludes from
\[
F_1 F_2 \equiv  F_{12}
\]
that 
\[
\mbox{\boldmath{$A$}}(F_1)\mbox{\boldmath{$A$}}(F_2)\equiv
\mbox{\boldmath{$A$}}(F_{12}).
\]
Furthermore, the mapping $\mbox{\boldmath{$A$}}$ is invertible, and
hence it is an automorphism of $\mathfrak{F}$.

If we now introduce the exponent $x$ into the group by setting
\[
\mbox{\boldmath{$A$}}(F)=F^x
\]
then one sees that
\[
S_{i,k}=S^{x^k}_{i,0}=S^{x^k}_i
\]
and that we have a system of generators $S_i$ and a system of defining 
relations $R_l(S_k)$ for the group when the exponents $f(x)$ are admitted.

\section{Characterization of groups with operators}

As for ordinary groups, one can ask how the various ways of defining
a group with operators by generators and relations are connected to
each other. It is clear that one can add any consequence to the
defining relations, or omit any relation $R_m$ when it is a consequence
of the others. Likewise, it is permissible to introduce a new generator
$S_{n+1}$ defined as a power product of the $S_1,S_2,\ldots,S_n$
with the help of a new relation, or to eliminate a generator $S_n$ that 
may be expressed in terms of the others. One can then prove, by
considerations quite similar to those in Section 2.10, that any two
systems of generators and defining relations may be converted to each
other by such steps.

As a result, the properties of the defining relations characteristic of the
group $\mathfrak{F}$ with operator $x$ can be given purely formally
as matrix properties. If
\[
R_i(S_k)=\prod^{n}_{i=1}S^{r_{ik}(x)}_k
\]
are the defining relations of $\mathfrak{F}_x$ and
\[
\rho=(r_{ik}(x))
\]
is the matrix of exponents $r_{ik}(x)$, and if
\[
R_{m+1}=\prod^{m}_{k=1}R^{p_k(x)}_k
=\prod^{n}_{i=1}S^{r_{m+1,i}(x)}_i
\]
is a consequence relation, then
\[
r_{m+1,i}=\sum^{m}_{k=1} p_k r_{ki}.
\]
If we add $R_{m+1}$ to the others as a defining relation, then the
exponent matrix of the new system will be denoted by 
$\rho'=(r'_{ik})$. The passage from $\rho$ to $\rho'$, as well as from
$\rho'$ to $\rho$, will be called a type I rearrangement  of matrices.
If
\[
\prod^{n}_{i=1}S^{r_{m+1,i}}_i
\]
is any power product, $S_{n+1}$ is a new generator and
$r_{m+1,n+1}=-1$, and if we add $S_{n+1}$ as a new generator and
\[
R_{m+1}=\prod^{n+1}_{i=1}S^{r_{m+1,i}}_i
\]
as a new relation, then the matrix corresponding to the new system will
be denoted by $\rho''$. The passage from $\rho$ to $\rho''$ and
conversely will be called a type II matrix rearrangement.

\emph{The properties of exponent matrices invariant under
rearrangements of the first and second kind characterize the group
$\mathfrak{F}$.}

It now remains to show that the elementary divisors of $\rho$ may
also be defined in this case, that the elementary divisors $\ne x^n$
are the same for all presentations of $\mathfrak{F}_x$, but that they
do \emph{not} chacterize $\mathfrak{F}_x$.

For this purpose we introduce the \emph{concept of divisibility and
greatest common divisor for $L$-polynomials with integral
coefficients}. We call $f(x)$ divisible by $g(x)$ when there is a
polynomial $h(x)$ for which
\[
f(x)=g(x)h(x).
\]
By the greatest common divisor $d(x)$ of the polynomials $f_i(x)$
$(i=1,2,\ldots,r)$,
\[
d(x)=(f_1(x),f_2(x),\ldots,f_r(x))
\]
we mean a polynomial which is a divisor of all the $f_i(x)$ and divisible
by all their common divisors. We show that greatest common divisors 
always exist, and if $d_1(x)$ and $d_2(x)$ are both greatest common 
divisors of the $f_i(x)$ then
\[
d_2(x)=\pm x^n d_1(x).
\]

Further: if $d(x)$ is a greatest common divisor of the $f_i(x)$
$(i=1,2,\ldots,r)$, and if
\[
f_{r+1}(x)=\sum n_i(x)f_i(x)
\]
is a linear combination of the $f_i(x)$, then $d(x)$ is likewise the
greatest common divisor of
\[
f_i(x) \qquad (i=1,2,\ldots,r+1).
\]
Elementary divisors of a matrix $\rho(x)$ may now be defined exactly
as the elementary divisors of the matrix $\rho$ were in Section 2.11.
Moreover, equivalence of matrices $\rho(x)$ may be defined as in
Section 2.12, except that integers are replaced by arbitrary 
$L$-polynomials in the definition of the transformation 2. It then follows
that \emph{equivalent matrices have the same elementary divisors},
and it follows in turn, by considerations like those in Sections 2.12 and
2.13, that the elementary divisors $\ne x^n$ of an exponent matrix are
invariant under the matrix rearrangements of types I and II. The proofs of 
these theorems for integral matrices may be carried over directly, since they
use only properties of the greatest common divisor $d(x)$.

An example later will illustrate that matrices with the same elementary
divisors need not be equivalent, and that groups $\mathfrak{F}_x$
are therefore not characterized by the elementary divisors of the
exponent matrix of their defining relations.

\section{Divisibility properties of $L$-polynomials}

We will reduce the divisibility relations between $L$-polynomials with
integral coefficients to those between ordinary polynomials with integral 
coefficients. \emph{We
call the integral domain of integral $L$-polynomials $\mathfrak{I}$,
and that of the integral ordinary polynomials, $\mathfrak{I}_g$}.
The ordinary polynomial $f(x)$ is said to be divisible by the
ordinary polynomial $g(x)$ in $\mathfrak{I}_g$ if there is an ordinary
polynomial $h(x)$ such that
\[
f(x)=g(x)h(x).
\]

If
\[
f(x)=a_n x^n + a_{n+1}x^{n+1} + \cdots + a_{n+m}x^{n+m},\quad
a_n\ne 0,
\]
then by $|f(x)|$ we mean the ordinary polynomial $x^{-n}f(x)$. We
now claim: \emph{if $f(x)$ is divisible by $g(x)$, then $|f(x)|$ is divisible
by $|g(x)|$ in $\mathfrak{I}_g$, and conversely, if $|f(x)|$ is divisible
by $|g(x)|$ in $\mathfrak{I}_g$ then $f(x)$ is divisible by $g(x)$}.

Namely, if
\[
f(x)=g(x)h(x)
\]
and 
\[
|f(x)|=x^{-n}f(x),\quad |g(x)|=x^{-m}g(x)
\]
then
\[
|f(x)|=|g(x)|x^{m-n}h(x).
\]
The polynomial $x^{m-n}h(x)$ must now be equal to
\[
c_0+c_1 x+\cdots+c_l x^l
\]
with $c_0\ne 0$, and thus
\[
|x^{m-n}h(x)|=x^{m-n}h(x).
\]
So $|f(x)|$ is divisible by $|g(x)|$ in $\mathfrak{I}_g$. The converse 
is trivial.

To find the divisors of $f(x)$ we therefore need only to find the divisors
$t(x)$ of $|f(x)|$ in $\mathfrak{I}_g$; $x^n t(x)$ then yields all the
divisors of $f(x)$ when $n$ runs through all the integers and $t(x)$
runs through all the divisors of $|f(x)|$ in $\mathfrak{I}_g$.

Now to find the divisors of a polynomial in $\mathfrak{I}_g$ we must
introduce yet another domain of polynomials and a new concept of
divisibility. By $\mathfrak{I}_r$ we mean the collection of ordinary
polynomials in one variable with rational coefficients. We denote
polynomials from $\mathfrak{I}_g$ with a subscript $g$, polynomials
from $\mathfrak{I}_r$ with a subscript $r$. $\mathfrak{I}_g$ is
contained in $\mathfrak{I}_r$. Addition and multiplication of
polynomials in $\mathfrak{I}_r$ are defined as in Section 1.13. One
verifies that $\mathfrak{I}_r$ is also an integral domain. A polynomial
$f_r(x)$ is said to be divisible by $g_r(x)$ in $\mathfrak{I}_r$ if there 
is a polynomial $h_r(x)$ such that
\[
f_r(x)=g_r(x)h_r(x).
\]

In order to describe the relation between divisibility in $\mathfrak{I}_g$
and $\mathfrak{I}_r$ we call a polynomial
\[
f_g(x)=a_0+a_1 x+\cdots+a_n x^n
\]
from $\mathfrak{I}_g$ \emph{primitive} when the greatest common
divisor
\[
(a_0,a_1,\ldots,a_n)=a
\]
is equal to 1. If
\[
f_r(x)=b_0+b_1 x+\cdots+b_m x^m
\]
is a polynomial from $\mathfrak{I}_r$ then we take the rational
numbers $b_i$ to their least common denominator $n\ge 1$ 
[not the same $n$] as
\[
b_i=\frac{b'_i}{n},
\]
denote by $b$ the greatest common divisor of the $b'_i$, and set
\[
b'_i=bb''_i\qquad (i=0,1,\ldots,m).
\]
The polynomial $\sum^{n}_{i=0}b''_i x^i$ is then a primitive
polynomial uniquely determined by $f_r(x)$, which we may denote by
$||f_r(x)||$, and
\[
f_r(x)=\frac{b}{n}||f_r(x)||.
\]
For the polynomial $f_g(x)$ we have
\[
f_g(x)=a||f_g(x)||,
\]
where $a$, as above, is the greatest common divisor of the coefficients
of $f_g(x)$.

The connection between the divisors $t_g$ of $f_g$ in $\mathfrak{I}_g$
and the divisors $t_r$ of $f_g$ in $\mathfrak{I}_r$ is now the
following:

\emph{If $t_r$ is a divisor of $f_g$ in $\mathfrak{I}_r$ and 
$f_g=a||f_g||$, then $t\cdot ||t_r||$ is a divisor of $f_g$ in $\mathfrak{I}_g$ when $t$ is a divisor of $a$. If, on the other hand, $t_g$ is a divisor of
$f_g$ in $\mathfrak{I}_g$, then $t_g$ is also a divisor of $f_g$ in
$\mathfrak{I}_r$, and if $t_g=t||t_r||$ then $t$ is a divisor of $a$.}
Thus the $t||t_r||$ comprise all divisors of $f_g$ in $\mathfrak{I}_g$.
The proof depends as usual on the theorem\footnote{This theorem is
commonly known as \emph{Gauss's lemma}, because of its appearance in
Article 42 of Gauss's \emph{Disquisitiones arithmeticae}. (Translator's note.)} on primitive polynomials: the product of primitive polynomials is again 
a primitive polynomial.

If $t_r(x)$ is a divisor of $f_g(x)$ in $\mathfrak{I}_r$, then $t_r(x)$ is
also a divisor of $||f_g(x)||$. Now let
\begin{align*}
||f_g(x)||&=t_r(x)\cdot h_r(x),\\
t_r(x)=\frac{c_1}{n_1}||t_r(x)||,&\quad h_r(x)=\frac{c_2}{n_2}||h_r(x)||.
\end{align*}
Then
\[
n_1 n_2 ||f_g(x)||=c_1 c_2 ||t_r(x)||\cdot ||h_r(x)||.
\]
Since $||t_r||\cdot ||h_r||$ is a primitive polynomial, the greatest common
divisor of the coefficients on the right hand side equals $c_1c_2$, and that of those on the left is $n_1 n_2$, so
\[
c_1 c_2=n_1 n_2.
\]
Thus
\[
||f_f(x)||=||t_r(x)||\cdot ||h_r(x)||,
\]
so $||t_r||$ is a divisor of $||f_g||$ in $\mathfrak{I}_g$, and hence also 
a divisor of $f_g$ in $\mathfrak{I}_g$. But then $t||t_r||$ is also a divisor
of $f_g(x)$ when $t$ is a divisor of $a$.

If, on other hand, $t_g(x)$ is a divisor of $f_g(x)$ in $\mathfrak{I}_g$,
so that
\[
f_g(x)=t_g(x)h_g(x),
\]
and if
\[
f_g(x)=a||f_g(x)||,\quad
t_g(x)=t||t_g(x)||,\quad
h_g(x)=h||h_g(x)||,
\]
then $a=th$ and the converse follows.

\section{Greatest common divisor}

In the domains $\mathfrak{I}_g$ and $\mathfrak{I}_r$ we define the
greatest common divisors of $n$ polynomials
\begin{align*}
&f_{g1}(x),\quad f_{g2}(x),\quad \ldots,\quad f_{gn}(x),\\
\text{respectively},\quad
&f_{r1}(x),\quad f_{r2}(x),\quad \ldots,\quad f_{rn}(x),
\end{align*}
in the usual way as follows: $d_g(x)$, respectively $d_r(x)$, is called the
greatest common divisor of the $f_{gi}$, respectively $f_{ri}$, if $d_g(x)$, respectively $d_r(x)$, is a common divisor of all the  $f_{gi}$, respectively 
$f_{ri}$, and each common divisor of all the  $f_{gi}$, respectively $f_{ri}$,
is also a divisor of $d_g(x)$, respectively $d_r(x)$.

In algebra one shows that a greatest common divisor of $n$ polynomials in
$\mathfrak{I}_r$ exists, and may be determined, and that two different
greatest common divisors of the same polynomials, say $d_r(x)$ and
$d'_r(x)$, differ only by a constant factor:
\[
d'_r(x)=cd_r(x).
\]
Thus $\pm ||d_r(x)||$ is uniquely determined. We now claim:

\emph{If $d_r(x)$ is a greatest common divisor of the polynomials}
\[
f_{g1}(x),\quad f_{g2}(x),\quad \ldots,\quad f_{gn}(x)
\]
\emph{in $\mathfrak{I}_r$, and $a$ is the greatest common divisor of the
coefficients of the $f_{gi}(x)$, then $a||d_r(x)||$ is the greatest common
divisor of the polynomials $f_{gi}$ in $\mathfrak{I}_g$}. Namely,
$a||d_r(x)||$ is a divisor of $f_{gi}(x)$, because if
\[
f_{gi}(x)=a^{(i)}||f_{gi}(x)||
\]
then $a$ is a divisor of $a^{(i)}$ and $d_r(x)$ is a divisor of $f_{gi}(x)$
in $\mathfrak{I}_r$, hence $||d_r(x)||$ is a divisor of $f_{gi}(x)$ in
$\mathfrak{I}_g$. Conversely, if $t_g(x)$ is a divisor of all the $f_{gi}(x)$
in $\mathfrak{I}_g$ and $t_g=t||t_g(x)||$ then $t$ must divide all of the
$a^{(i)}$ and hence also $a$. And $||t_g(x)||$ is a divisor of all the
$f_{gi}(x)$ in $\mathfrak{I}_r$, so $||t_g(x)||$ is likewise a divisor of
$d_r(x)$ in $\mathfrak{I}_r$, hence $||t_g(x)||$ is also a divisor of
$||d_r(x)||$ in $\mathfrak{I}_g$. Thus $a||d_r(x)||=d_g(x)$ is the
greatest common divisor of the polynomials $f_{gi}(x)$ in $\mathfrak{I}_g$.
If $d'_g(x)$ is another greatest common divisor in $\mathfrak{I}_g$ then
$d'_g(x)$ is also a greatest common divisor in $\mathfrak{I}_r$ and we
must therefore have $d'_g(x)=c||d_g(x)||$. Further, $c$ must divide all
$a^{(i)}$ and likewise $a$ must divide $c$, so $c=\pm a$ and hence
$d_g(x)=\pm d'_g(x)$.

Finally we return to the original integral domain $\mathfrak{I}$ of the 
polynomials $f(x)$. If $f_i(x)$ $(i=1,2,\ldots,m)$ are polynomials of
$\mathfrak{I}$ and $d_g(x)$ is the greatest common divisor of the
$|f_i(x)|$, then $d_g(x)$ is also the greatest common divisor of the
$f_i(x)$. Because $d_g(x)$ is a divisor of the $|f_i(x)|$ and hence also 
of the $f_i(x)$, and if $t(x)$ is a common divisor of the $f_i(x)$ then
$|t(x)|$ is also a common divisor of the $|f_i(x)|$ in $\mathfrak{I}_g$,
and thus a divisor of $d_g(x)$.

If $f(x)$ is any greatest common divisor of the $f_i(x)$ then $|d(x)|$ is 
a divisor of the $|f_i(x)|$ in $\mathfrak{I}_g$, so $|d(x)|$ is a divisor
$d_g(x)$ in $\mathfrak{I}_g$, i.e., $d_g(x)=|d(x)|h_g(x)$.
Conversely, $d_g(x)$ is a common divisor of the $|f_i(x)|$, hence also
of the $f_i(x)$, so $d_g(x)$ is among the divisors of $d(x)$ in
$\mathfrak{I}$, i.e.,
\[
d(x)=d_g(x)h(x),
\]
which implies
\[
|d(x)|=|d_g(x)||h(x)|.
\]
Since also
\[
|d_g(x)|=|d(x)||h_g(x)|,
\]
it follows that
\[
|d(x)|=|d(x)||h_g(x)||h(x)|,
\]
whence
\[
|h_g(x)||h(x)|=1
\]
and so
\[
|h_g(x)|=\pm 1,\quad |h(x)|=\pm 1.
\]
It follows that
\[
d(x)=\pm d_g(x) x^n.
\]

Finally, let
\[
f_{m+1}(x)=\sum^{m}_{i=1}n_i(x)f_i(x).
\]
We claim: \emph{a greatest common divisor in $\mathfrak{I}$ of the
$f_i(x)$ $(i=1,2,\ldots,m)$ is also a greatest common divisor of
$f_i(x)$ and $f_{m+1}(x)$, and conversely}. Namely, each common 
divisor of the $f_i(x)$ is also one of $f_i(x)$ and $f_{m+1}(x)$, and
conversely.

\section{An example}

A simple example of a group $\mathfrak{F}_x$ which is not characterized
by the elementary divisors of its exponent matrix is the group with the
relations
\begin{equation}
R_1=S^{x^2+1}\equiv 1,\quad R_2=S^2\equiv 1. \tag{1}
\end{equation}
The greatest common divisor $d_g(x)$ of $x^2+1$ and 2 is obviously 1.
Another matrix with the same elementary divisors is given by
\[
R'_1=S\equiv 1.
\]

We will show that the element $S$ does not equal the identity in the 
group defined by (1). If, on the contrary, $S$ were a consequence relation,
then, by Section 2.14,
\[
S=S^{(x^2+1)n_1(x)+2n_2(x)}
\]
would hold for suitable $n_i(x)$. That is,
\begin{equation}
1=(x^2+1)n_1(x)+2n_2(x). \tag{2}
\end{equation}
This polynomial relation must hold for all values of $x$, but if we set
$x=1$ we see that the right hand side is divisible by 2 and hence cannot
equal 1.

If we set $S_i=S^{x^i}$ and substitute the $S_i$ into $R_1$ it follows
that
\[
S_{i+2}=S^{-1}_i,\quad \text{so}\quad S_{i+4}=S_i.
\]
If we take $S_1$, $S_2$ as generators then the relation
\[
S^{x^{n}(x^2+1)}\equiv 1
\]
may be used to express the remaining $S_i$ in terms of $S_1$ and $S_2$.
$R_2$ then says
\[
S^{2}_1\equiv 1,\quad S^{2}_2\equiv 1
\]
and we have
\[
S^{x}_1=S^{-1}_2=S_2,\quad S^{x}_2=S^{-1}_1=S_1.
\]

The unsatisfiability of equation (2) also shows a property of the greatest
common divisor $d_g(x)$ of two polynomials $f_{g1}(x)$ and
$f_{g2}(x)$ in $\mathfrak{I}_g$. Namely, $d_g(x)$ is not in general
expressible as a linear combination of $f_{g1}(x)$ and $f_{g2}(x)$.

\section{Factor groups with respect to commutator groups}

A commutative group may be constructed from each group with
generators $S_1,S_2$, $\ldots$, $S_n$ and relations $R_1(S),R_2(S),\ldots,
R_m(S)$ by adding the relations
\begin{equation}
R_{ik}=S_i S_k S^{-1}_i S^{-1}_k. \tag{1}
\end{equation}
The group that results from $\mathfrak{F}$ in this way is called 
$\mathfrak{F}'$. We claim that  
$\mathfrak{F}'=\mathfrak{F}/\mathfrak{K}_1$, \emph{the factor group of
$\mathfrak{F}$ by the commutator group $\mathfrak{K}_1$}. Certainly
$\mathfrak{K}_1$ contains all the $R_{ik}$, and hence all transforms of the $R_{ik}$ and their power products, hence all consequence relations of the
$R_{ik}$. However, if $A$ and $B$ are any power products of the $S_i$
$(i=1,2,\ldots,n)$ then it follows from (1) that $A$ commutes with $B$,
so
\[
ABA^{-1}B^{-1}
\]
is a consequence relation of the $R_{ik}$. Therefore, since each
commutator element of $\mathfrak{F}$ may be represented as a power
product of the $R_{ik}$ and their transforms, the $R_{ik}$ and their
transforms generate the group $\mathfrak{K}_1$, so $\mathfrak{F}'$ is
indeed equal to $\mathfrak{F}/\mathfrak{K}_1$.

Similarly, one can also construct the factor group of $\mathfrak{F}$ by
the second commutator group $\mathfrak{K}_2$,
$\mathfrak{F}''=\mathfrak{F}/\mathfrak{K}_2$, by taking, in place of (1),
the relations
\begin{equation}
R_{ikl}=S_l R_{ik} S^{-1}_l R^{-1}_{ik}\qquad (i,k,l=1,2,\ldots,n).
\tag{2}
\end{equation}
This is because these relations have the consequence
\[
FR_{ik}F^{-1}\equiv R_{ik}
\]
for each element $F$ of $\mathfrak{F}$, so $S_l$ commutes with
$FR_{ik}F^{-1}$ and hence with all elements of $\mathfrak{K}_1$. The
equations (2) now similarly permit a simple representation of all elements by
power products, as with commutative groups.

The power product $F(S)$ of the $S_l$ is equivalent in $\mathfrak{F}''$ to
a product
\[
S^{r_1}_1 S^{r_2}_2 \cdots S^{r_n}_n K,
\]
where $K$ belongs to the commutator group $\mathfrak{K}_1$, and hence 
is a product of the $R_{ik}$ and their transforms. Now if we introduce new generators $T_{ik}$ by the relations
\[
T^{-1}_{ik} R_{ik}\equiv 1\qquad (i<k)
\]
then, in consequence of the relations (2),
\[
FT_{ik} F^{-1}\equiv T_{ik}
\]
and
\[
T_{ik}T_{lm}\equiv T_{lm}T_{ik},\quad T_{ik}\equiv R^{-1}_{ki},
\]
thus
\[
K\equiv T^{r_{12}}_{12}\cdots T^{r_{1n}}_{1n}
T^{r_{23}}_{23}\cdots T^{r_{n-1,n}}_{n-1,n}.
\]
The products $R_{ikl}$ go to the empty word under this conversion.

In particular, one can express the relations $R_i(S)$ in the given form, so we
have
\[
R_i(S)\equiv S^{r_{i1}}_1 S^{r_{i2}}_2 \cdots S^{r_{in}}_n
           T^{r_{i,12}}_{12} \cdots T^{r_{i,n-1,n}}_{n-1,n}.
\]
Striking out the $T_{ik}$, one obtains the relations of 
$\mathfrak{F}'=\mathfrak{F}/\mathfrak{K}_1$.

We will now deal in more detail with the case where all $r_{ik}=0$, so
that the group $\mathfrak{F}'=\mathfrak{F}/\mathfrak{K}_1$ is a free
\textsc{Abelian} group with $n$ generators. Then the 
\[
R_i\equiv T^{r_{i,12}}_{12} \cdots T^{r_{i,n-1,n}}_{n-1,n}
\]
are also the defining relations of the subgroup of $\mathfrak{F}''$ 
generated by the $T_{ik}$, the commutator group $\mathfrak{K}''_1$
of $\mathfrak{F}''$. This is a commutative group that is completely
characterized by the elementary divisors of the matrix $(r_{i,lk})$.
Since $\mathfrak{K}''_1$ is also the factor group of $\mathfrak{K}_1$
by $\mathfrak{K}_2$, $\mathfrak{K}_1/\mathfrak{K}_2$, in the 
elementary divisors of $(r_{i,lk})$ we have numbers determined by the 
group $\mathfrak{F}$ itself, and not dependent on the presentation of
$\mathfrak{F}$ by generators and 
relations.\footnote{Cf. \textsc{K. Reidemeister}, Hamb. Abhdl. \textbf{5},
33 and \textsc{H. Adelsberger}, J. f. reine u. angew. Math. \textbf{163}
(1930) 103.}

\chapter{Determination of Subgroups}

\section{Generators of Subgroups}

Many deeper insights into the structure of a group are obtained by a
process for determining generators and defining relations for
subgroups.\footnote{For the following sections cf. 
\textsc{K. Reidemeister}, Hamb. Abhandl. \textbf{5}, (1926), 8 and
\textsc{O. Schreier}, Hamb. Abhandl. \textbf{5}, (1926), 161 }
This process allows, e.g., the commutator subgroup of a group to be
constructed. A geometric interpretation of the following considerations
will be found in Sections 4.20 and 6.14.

\emph{Let $\mathfrak{m}$ be a set of generators $S_1,S_2,\ldots,S_n$
of a group $\mathfrak{F}$, $\mathfrak{U}$ a subgroup of $\mathfrak{F}$
and $\mathfrak{g}$ a system of representatives}
\[
G_1,G_2,,\ldots
\]
\emph{of the left-sided residue classes $\mathfrak{U}G$ of
$\mathfrak{U}$ in $\mathfrak{F}$.} The residue class $\mathfrak{U}$
itself may be represented by the identity element $E=1$, the remaining
$G_i$ are fixed power products of the $S$. If $F$ is any element of
$\mathfrak{F}$ and if $F$ belongs to the residue class $\mathfrak{U}G$,
then we define $\overline{F}$ by
\[
\overline{F}=G.
\]
If $U$ belongs to $\mathfrak{U}$ then $\overline{UF}=\overline{F}$ and
in particular $\overline{U}=1$. Under these assumptions we claim that
\emph{the elements}
\[
U_{G,S}=GS{\overline{GS}\,}^{-1}
\]
\emph{constitute a system $\mathfrak{U}$ of generators for the
subgroup $\mathfrak{U}$ when $G$ and $S$ run through the classes
$\mathfrak{g}$ and $\mathfrak{m}$ independently of each other.}

To prove this we first remark that the $U_{G,S}$ themselves belong to
the subgroup $\mathfrak{U}$. Namely, since $GS$ and $\overline{GS}$
belong to the same residue class of $\mathfrak{U}$ we have
\[
GS\equiv U\overline{GS},
\]
where $U$ is a suitable element of $\mathfrak{U}$. Consequently
\[
U_{G,S}\equiv U\overline{GS}\, {\overline{GS}\,}^{-1}=U.
\]
If we now note that 
\[
GS^{-1}{\overline{GS^{-1}}\,}^{-1}
\]
is formally inverse to $\overline{GS^{-1}} SG^{-1}$, and that the
latter is an element $U_{G,S}$ (namely, the element $U_{G',S}$ with
$G'=\overline{GS^{-1}}$, since 
$\overline{G'S}=\overline{\overline{GS^{-1}}S}=G$) then we can
easily see that each power product of the $S_i$ ($i=1,2,\ldots,n$)
\begin{equation}
S^{\varepsilon_1}_{\alpha_1}
S^{\varepsilon_2}_{\alpha_2}\cdots
S^{\varepsilon_m}_{\alpha_m}
\qquad(\varepsilon_k=\pm 1, k=1,2,\ldots,m)
\tag{1}
\end{equation}
which yields an element of $\mathfrak{U}$ may also be written as a
power product in the $U_{G,S}$. Namely, we set
\[
W_0=1,\quad
W_1=S^{\varepsilon_1}_{\alpha_1},\quad
W_2=S^{\varepsilon_1}_{\alpha_1}S^{\varepsilon_2}_{\alpha_2},\quad
\ldots,\quad
W_m=S^{\varepsilon_1}_{\alpha_1}
S^{\varepsilon_2}_{\alpha_2}\cdots
S^{\varepsilon_m}_{\alpha_m}
\]
and construct
\[
\overline{W}_0 S^{\varepsilon_1}_{\alpha_1} \overline{W}^{\,-1}_1
\overline{W}_1 S^{\varepsilon_2}_{\alpha_2} \overline{W}^{\,-1}_2
\overline{W}_2\cdots\overline{W}^{\,-1}_{\alpha_{m-1}}
\overline{W}_{\alpha_{m-1}} S^{\varepsilon_m}_{\alpha_m} \overline{W}^{\,-1}_m.
\]
Since $\overline{W}_0=\overline{1}=1$ and $W_m=1$ likewise,
because $W_m$ belongs to $\mathfrak{U}$, this power product
becomes (1) by elementary computations in the free group on the
$S$. But each of the factors
\[
\overline{W}_{i-1}S^{\varepsilon_i}_{\alpha_i}\overline{W}^{\,-1}_i
\]
is either a $U_{G,S}$ or else the inverse of such an element, so the
$U_{G,S}$ are a system of generators for $\mathfrak{U}$. One sees
that it is essential to have set $\overline{U}=1$.

Since there are different representative systems $\mathfrak{g}$ for
the residue classes $\mathfrak{U}G_i$ modulo $\mathfrak{U}$
which satisfy the latter requirement there are also different systems
of generators for $\mathfrak{U}$. We shall make use of this in
Section 3.6 in order to bring the defining relations of $\mathfrak{U}$
into a clearly arranged form. But first we must carry out the
determination of these relations.

\section{Generators of the subgroup as special generators of the
group}

It is quite simple to give relations which the generators $U_{G,S}$
must satisfy: any relation $R$ in the $S_i$ certainly yields an element
which also belongs to $\mathfrak{U}$, and hence may be expressed
in terms of the $U_{G,S}$; the resulting representation of the identity
element in $\mathfrak{U}$ is then a relation in the $U_{G,S}$. It is
more difficult to clarify how one obtains all relations in the $U_{G,S}$.
Naturally we cannot just say that each relation in the $U_{G,S}$
results from substitution in a relation in the $S_i$. Indeed we shall
see that this is in general false.

We prepare for the solution of this problem by introducing the
$U_{G,S}$ as new generators of the group $\mathfrak{F}$
subject to the relations
\[
U^{-1}_{GS} GS \,\overline{GS}^{\,-1}\equiv 1
\]
and give a somewhat modified procedure for expressing power
products of in the $S_i$ and $U_{G,S}$ which yield elements of
$\mathfrak{U}$ in terms of the $U_{G,S}$ alone. As a more convenient
way of writing we denote the new generators of $\mathfrak{F}$ by
$T_1,T_2,\ldots$ and their totality by $\mathfrak{t}$. A $T_i$ is
therefore a certain $S_k$ or a certain $U_{G,S}$. Now if
\begin{equation}
F=T^{\varepsilon_1}_{\alpha_1}
    T^{\varepsilon_2}_{\alpha_2}\cdots
    T^{\varepsilon_m}_{\alpha_m}
    \tag{1}
\end{equation}
is an arbitrary product in the $T$, we construct the subproducts
\[
W_0=1,\quad
W_1=T^{\varepsilon_1}_{\alpha_1},\quad
W_2=T^{\varepsilon_1}_{\alpha_1} T^{\varepsilon_2}_{\alpha_2},\quad
\ldots,\quad
W_m=T^{\varepsilon_1}_{\alpha_1}
    T^{\varepsilon_2}_{\alpha_2}\cdots
    T^{\varepsilon_m}_{\alpha_m}
\]
and, as in the previous section, taking $\overline{W}$ to be a certain
power product $G$ of the $S_i$ which represents the residue class to
which $W$ belongs, we set
\begin{equation}
F'=\overline{W}_0 T^{\varepsilon_1}_{\alpha_1} {\overline{W}\,}^{-1}_1
     \overline{W}_1 T^{\varepsilon_2}_{\alpha_2} {\overline{W}\,}^{-1}_2
     \cdots {\overline{W}\,}^{-1}_{m - 1}
\overline{W}_{m - 1} T^{\varepsilon_m}_{\alpha_m} {\overline{W}\,}^{-1}_m.
\tag{2}
\end{equation}
Again $\overline{W}_0=1$ and, if $F$ belongs to $\mathfrak{U}$, also
$\overline{W}_m=1$, and thus $F'$ is convertible into $F$ by elementary
manipulations. If we also express the terms
$\overline{W}_{i-1} T^{\varepsilon_i}_{\alpha_i} {\overline{W}\,}^{-1}_i$,
and in general the terms
\begin{equation}
GT^\varepsilon {\overline{GT^\varepsilon}\,}^{-1}\qquad
\varepsilon=\pm 1, \tag{3}
\end{equation}
where $G$ and $T$ run independently through the classes
$\mathfrak{g}$ and $\mathfrak{t}$, in terms of the $U_{G,S}$ in a
specific way, then we obtain a new rule for representing each $F$
that belongs to $\mathfrak{U}$ in terms of the $U_{G,S}$.

If we now make the replacement
\[
GT^{-1}{\overline{GT^{-1}}\,}^{-1}=
\left(G'T\overline{G'T^{-1}}\right)^{-1}
\]
where $G'=\overline{GT^{-1}}$, then in case $T_i=S_k$
\begin{equation}
GT_i{\overline{GT_i}\,}^{-1}   \tag{4}
\end{equation}
is replaced by $U_{G,S_k}$, and we express this by writing
\[
|GT_i{\overline{GT_i}\,}^{-1}|_{\mathfrak{U}}
=U_{G,S_k}.
\]
When the $T$ in expression (4) corresponds to a $U_{G',S}$ and
$G=1$, so that ${\overline{GT_i}\,}^{-1}=1$, we write $U_{G',S}$
for (4). If $G\ne 1$ and $T_i=U_{G',S}$ we suppose a rule
is given, which we need not specify precisely, to replace the 
$U_{G'',S}$ in a
power product $|GT_i{\overline{GT_i}\,}^{-1}|_\mathfrak{U}$ from 
expression (4). Since $G$ and ${\overline{GT_i}\,}^{-1}$
contain only $S$ generators we can proceed, e.g., by replacing
$T_i=U_{G',S}$ by the product $G'S{\overline{G'S}\,}^{-1}$ in the $S$,
and the resulting product $GU_{G',S}{\overline{GU_{G',S}}\,}^{-1}$ can
be expressed in terms of the $U_{G,S}$ as in the previous section.
We also denote this product by 
$|GT_i{\overline{GT_i}\,}^{-1}|_\mathfrak{U}$. If $F$ is an arbitrary
power product in the $T$ which yields an element of $\mathfrak{U}$,
then by $F_\mathfrak{U}$ we mean that product of the $U_{G,S}$
which consists of $F'=F$ in which the factors
$\overline{W}_{i-1} T^{\varepsilon_i}_{\alpha_i}\overline{W}^{\,-1}_i$
are replaced by
$|\overline{W}_{i-1} T^{\varepsilon_i}_{\alpha_i}\overline{W}^{\,-1}_i|_\mathfrak{U}$, 
respectively 
$(|\overline{W}_{i-1} T^{\varepsilon_i}_{\alpha_i}{\overline{W}\,}^{-1}_i|_\mathfrak{U})^{-1}$,
in the way just described.

To avoid misunderstanding, we emphasise that the symbol 
$|F_\mathfrak{U}|$ is defined only for the special elements (4) and
that $|F_\mathfrak{U}|$ is in general different from $F_\mathfrak{U}$.
As an abbreviation we set 
\[
|F^{-1}|_\mathfrak{U}=(|F|_\mathfrak{U})^{-1}.
\]

\section{Properties of the replacement process}

We illuminate the connection between the products $F$ and
$F_\mathfrak{U}$ by the following theorems.
\medskip

Theorem 1. \emph{If $F$ is a product of the $U_{G,S}$ alone, then
$F_\mathfrak{U}$ is identical with $F$.}

Namely, if we construct the product $F'$ in Section 3.2 (2) for such an
element $F$ in 3.2 (1) then $W_i=1$ ($i=1,2,\ldots,m$) because all
the factors of $F$ belong to $\mathfrak{U}$.
\medskip

Theorem 2. \emph{If $F_1 F_2 = F_{12}$ is the product of $F_1$ and
$F_2$ written in juxtaposition, and if $F_1$ and $F_2$ belong to
$\mathfrak{U}$, then}
\[
{F_1}_\mathfrak{U} {F_2}_\mathfrak{U}={F_{12}}_\mathfrak{U}.
\]

Namely, if
\[
F_k=\prod^{m_k}_{l=1} T^{\varepsilon_{l,k}}_{\alpha_{l,k}}
\quad(k=1,2),\qquad
W_{i,k}=\prod^i_{l=1} T^{\varepsilon_{l,k}}_{\alpha_{l,k}}
\quad(i=1,2,\ldots,m_k;k=1,2)
\]
then
\[
F_{12}=\prod^{m_1}_{l=1} T^{\varepsilon_{l,1}}_{\alpha_{l,1}}
            \prod^{m_2}_{l=1} T^{\varepsilon_{l,2}}_{\alpha_{l,2}}.
\]
The subproducts are
\[
W_{i,12}=\prod^i_{l=1}T^{\varepsilon_{l,1}}_{\alpha_{l,1}} = W_{i,1}
\quad(i\le m_1)
\]
and for the subproducts
\[
W_{i,12}=
W_{m_1,12}\prod^{i-m_1}_{l=1}T^{\varepsilon_{l,2}}_{\alpha_{l,2}}
\quad(i>m_1)
\]
we have
\[
\overline{W}_{i,12}=\overline{W}_{i-m_1,2}
\]
because $W_{m_1,12}$ belongs to $\mathfrak{U}$ and 
$\overline{UF}=\overline{F}$. Thus $F'_{12}=F'_1 F'_2$ and the
assertion follows.
\medskip

Theorem 3. \emph{If $F$ is a power product with}
\[
\alpha_{i-1}=\alpha_i,\quad \varepsilon_{i-1}+\varepsilon_i=0
\]
\emph{and $F^{*}$ is the product that results from $F$ by striking
out $T^{\varepsilon_{i-1}}_{\alpha_{i-1}}T^{\varepsilon_i}_{\alpha_i}$,
then $F_\mathfrak{U}$ results from $F^{*}_\mathfrak{U}$ by
elementary manipulations in the domain of the $U_{G,S}$.}

Because if
\[
F=\prod^m_{i=1}T^{\varepsilon_i}_{\alpha_i},\quad
F^{*}=\prod^{m-2}_{i=1}T^{\eta_i}_{\beta_i},\quad
W_j=\prod^j_{l=1}T^{\varepsilon_l}_{\alpha_l},\quad
W^{*}_j=\prod^{j}_{l=1}T^{\eta_l}_{\beta_l},
\]
then $W_j=W^{*}_j$ for $j\le i-2$. For $j>i-2$, $W^{*}_j$ results
from $W_{j+2}$ by an elementary manipulation, hence
$\overline{W^{*}_j}=\overline{W}_j$ for $j\le i-2$ and 
$\overline{W^{*}_j}=\overline{W}_{j+2}$ for $j>i-2$.
If we now construct $F'$ and $F^{*}$ then
\[
F'=
\prod^{i-2}_{l=1} \overline{W}_{l-1} T^{\varepsilon_l}_{\alpha_l}
{\overline{W}\,}^{-1}_l\cdot
\overline{W}_{i-2} T^{\varepsilon_i-1}_{\alpha_i-1}
{\overline{W}\,}^{-1}_{i-1}\cdot
\overline{W}_{i-1} T^{\varepsilon_i}_{\alpha_i}
{\overline{W}\,}^{-1}_{i}\cdot
\prod^{m}_{l=i+1} \overline{W}_{l-1} T^{\varepsilon_l}_{\alpha_l}
{\overline{W}\,}^{-1}_l
\]
and by the identity just established
\[
{F^{*}}'=
\prod^{i-2}_{l=1} \overline{W}_{l-1} T^{\varepsilon_l}_{\alpha_l}
{\overline{W}\,}^{-1}_l\cdot
\prod^{m}_{l=i+1} \overline{W}_{l-1} T^{\varepsilon_l}_{\alpha_l}
{\overline{W}\,}^{-1}_l.
\]
If now $\varepsilon_{i-1}$ is, say, $+1$ then $\varepsilon_i$ is
$-1$, hence
\[
\overline{W}_{i-1} T^{\varepsilon_i}_{\alpha_i}
{\overline{W}\,}^{-1}_i =
(\overline{W}_{i-2} T_{\alpha_{i-1}}
{\overline{W}\,}^{-1}_{i-1} )^{-1}
\]
and if we now express the factors
\[
\overline{W}_{l-1} T^{\varepsilon_l}_{\alpha_l}
{\overline{W}\,}^{-1}_l
\]
in $F'$ and ${F^{*}}'$ by the $U_{G,S}$ as prescribed, then the two 
factors for $l=i-1$ and $l=i$ in $F_\mathfrak{U}$ yield actual
formally inverse components in the $U_{G,S}$.

From this we also have Theorem 4: \emph{If $F$ and $F^{-1}$ are
two formally inverse power products in the $T$ which represent
elements of $\mathfrak{U}$, then the corresponding products
$F_\mathfrak{U}$ and $(F^{-1})_\mathfrak{U}$ are likewise formally
inverse to each other.} For, since $FF^{-1}$ may be reduced to the
identity by cancellation in the free group generated by the $T$, this
also holds for 
$(FF^{-1})_\mathfrak{U}=F_\mathfrak{U} (F^{-1})_\mathfrak{U}$
in the free group generated by the $U_{G,S}$.

\section{Defining relations}

We have now completed the preparations needed to give the
defining relations of $\mathfrak{U}$ in the generators $U_{G,S}$.

\emph{If}
\[
R_1(T),\quad
R_2(T),\quad
\ldots,\quad
R_r(T)
\]
\emph{is a system $\mathfrak{r}$ of defining relations of the group
$\mathfrak{F}$ in the generators $T$, then we obtain a system of
defining relations for $\mathfrak{U}$ in the generators $U_{G,S}$
by expressing the power products}
\[
GRG^{-1}
\]
\emph{in terms of the $U_{G,S}$. Here $R$ runs through all relations
in $\mathfrak{r}$ and $G$ through a complete system of
representatives of the left-sided residue classes $\mathfrak{U}G$
of $\mathfrak{U}$ in $\mathfrak{F}$.}

To see why, let $R$ be any relation in the $U_{G,S}$. We can also
regard it as a relation in the $T$, and hence it is a consequence of 
the relations $R_i$ in $\mathfrak{r}$, so $R(U_{G,S})$ is, after
renaming the $U_{G,S}$ in the corresponding $T$,
\[
R(T)=\prod_i L_i(T) R^{\varepsilon_i}_{\alpha_i}(T) L^{-1}_i(T).
\]
That is, we obtain a product convertible into a product of transforms
of the $R$ from $\mathfrak{r}$ and their inverses in the free group
generated by the $T$. By further manipulations in the free group on
the $T$ we obtain
\begin{align*}
R(T)&=\prod L_i(T) R^{\varepsilon_i}_{\alpha_i}(T) L^{-1}_i(T)\\
&=\prod (L_i(T){\overline{L}\,}^{-1}_i)
     (\overline{L}_i R^{\varepsilon_i}_{\alpha_i}{\overline{L}\,}^{-1}_i)
     (\overline{L}_i L^{-1}_i(T)).
\end{align*}
Here the three bracketed elements belong to $\mathfrak{U}$, and
the products $L_i{\overline{L}\,}^{-1}_i$ and $\overline{L}_i L^{-1}_i$ 
are formally inverse to each other. Now, on the one hand,
\[
\left(\prod L_i(T){\overline{L}\,}^{-1}_i
     \overline{L}_i R^{\varepsilon_i}_{\alpha_i}{\overline{L}\,}^{-1}_i
     \overline{L}_i L^{-1}_i(T)\right)_\mathfrak{U}
\]
results from $R(U_{G,S})$ by elementary manipulations in the
domain of the $U_{G,S}$, by Theorems 1 and 3 of Section 3.3.
And, on the other hand, this product equals
\[
\prod (L_i{\overline{L}\,}^{-1}_i)_\mathfrak{U}
(\overline{L}_i R^{\varepsilon_i}_{\alpha_i} {\overline{L}\,}^{-1}_i)_\mathfrak{U}
     (\overline{L}_i L^{-1}_i)_\mathfrak{U}
\]
by Theorem 2 of 3.3. Also, the $(L_i{\overline{L}\,}^{-1}_i)_\mathfrak{U}$ 
and $(\overline{L}_i L^{-1}_i)_\mathfrak{U}$ are formally inverse
to each other by Theorem 4 of 3.3, and hence $R(U_{G,S})$ is a
consequence relation of the $(GRG^{-1})_\mathfrak{U}$.

\section{\textsc{Schreier's} normalized replacement process}

The defining relations of $\mathfrak{F}$ in the generators $T$ fall
naturally into two classes: the class of relations that define the
generators $U_{G,S}$ and the class of relations that originate from
the defining relations of $\mathfrak{F}$ in the generators $S$.
Correspondingly,  we can also divide the defining relations of
$\mathfrak{U}$ into two classes: \emph{relations of the first kind}
\begin{equation}
(GR(S)G^{-1})_\mathfrak{U}=1, \tag{1}
\end{equation}
and \emph{relations of the second kind}
\begin{equation}
(G' U_{G,S} [GS\overline{GS}^{\,-1}]^{-1} G'^{-1})_\mathfrak{U}=1.
\tag{2}
\end{equation}
Following \textsc{Schreier}, we now show that, by skillful
\emph{use of the freedom which we still have in the definition of the
process} $F_\mathfrak{U}$, we can eliminate the relations of the
second kind. We can still decide how we shall express
\[
G' U_{G,S} {\overline{G'U}\,}^{-1}_{GS}
\]
in terms of the $U_{G',S}$ for $G'\ne 1$, i.e., how we shall define
\[
|G' U_{G,S} {\overline{G'U}\,}^{-1}_{GS}|_\mathfrak{U},
\]
and we can choose the representatives of the residue classes $G_i$
in many ways. We first determine the operation 
$|F|_\mathfrak{U}$ more exactly. Let
\[
G' U_{G,S} [GS\overline{GS}^{\,-1}]^{-1} G'^{-1}=
\prod^{m}_{i=1} T^{\varepsilon_i}_{\alpha_i},
\]
let $W_i$ be the subproducts of this expression, and in particular
let
\[
W_k=G',\quad\text{so}\quad
T^{\varepsilon_k + 1}_{\alpha_k + 1} = U_{G,S}.
\]
If we construct
\[
\left( \prod^{m}_{i=1}T^{\varepsilon_i}_{\alpha_i}\right)'=
\prod^{m}_{i=1}\overline{W}_{i-1}T^{\varepsilon_i}_{\alpha_i}
{\overline{W}\,}^{-1}_i,
\]
then all $T_{\alpha_i}$ ($i\ne k+1$) correspond to generators $S$,
because the representatives $G$ of the residue classes are power
products of the $S$ alone. Hence
\[
|\overline{W}_{i-1}T^{\varepsilon_i}_{\alpha_i}
{\overline{W}\,}^{-1}_i|_\mathfrak{U}\quad (i\ne k+1)
\]
is a fixed power product of the $U_{G,S}$. We now set
\[
|\overline{W}_{k}T^{\varepsilon_{k+1}}_{\alpha_{k+1}}
{\overline{W}\,}^{-1}_{k+1}|_\mathfrak{U},
\quad\text{respectively}\quad
|G' U_{G,S}\overline{G'U}^{\,-1}_{GS}|_\mathfrak{U},
\]
equal to
\[
\left(\prod^{k}_{i=1} |\overline{W}_{i-1}T^{\varepsilon_i}_{\alpha_i}
{\overline{W}\,}^{-1}_i|_\mathfrak{U}\right)^{-1}
\left(\prod^{m}_{i=k+2} |\overline{W}_{i-1}T^{\varepsilon_i}_{\alpha_i}
{\overline{W}\,}^{-1}_i|_\mathfrak{U}\right)^{-1}.
\]
This is permissible, because the product of the $U_{G,S}$ on the right
hand side really represents the element
\[
G' U_{G,S} {\overline{G'S}\,}^{-1}_{GS}
\]
on the basis of the relation. This arrangement ensures that the relations
\[
(G' U_{G,S} [GS\overline{GS}^{-1}]^{-1} G'^{-1})_\mathfrak{U}
\]
are identically satisfied, as long as $G'\ne 1$.

\section{\textsc{Schreier's} choice of representatives $G$}

Now the only remaining relations of the second kind are the
\[
(U_{G,S} [GS\overline{GS}^{\,-1}]^{-1})_\mathfrak{U}\equiv 1.
\]
We simplify them by a suitable choice of the representatives $G$.
Namely, we impose the following condition on the $G$
(\textsc{Schreier's} condition):

($\Sigma$) \emph{Whenever}
\[
G=\prod^{r}_{i=1} S^{\varepsilon_i}_{\alpha_i}
\]
\emph{is the representative of its residue class, then the initial
segment products}
\[
\prod^{j}_{i=1} S^{\varepsilon_i}_{\alpha_i}\quad (j=1,2,\ldots,r-1)
\]
\emph{are also representatives of their residue classes.}

This condition is always satisfiable:

In each residue class there is at least one power product
\begin{equation}
\prod^{l}_{i=1} S^{\varepsilon_i}_{\alpha_i}\quad
(\eta_i=\pm 1) \tag{1}
\end{equation}
with the smallest possible number of factors $l$. We call $l$ the
length of the residue class determined by (1). $\mathfrak{U}$ is
the unique residue class of length zero; for we have in it the
empty power product, 1, which was chosen earlier as the
representative of $\mathfrak{U}$. Now suppose that we have succeeded
in choosing a representative, the expression for which satisfies
condition ($\Sigma$) and has length equal to that of the residue 
class, for each residue class with length $<l$ ($l>0$).

We will show that our condition can also be satisfied for residue
classes of length $l$. Let $\mathfrak{U}F$ be a residue class of
length $l$ and let (1) be a power product of length $l$ chosen
from $\mathfrak{U}F$. Then we construct the residue classes
\begin{equation}
\mathfrak{U}\prod^{l-1}_{i=1} S^{\eta_i}_{\alpha_i}. \tag{2}
\end{equation}
Their length is $l-1$. For (2) contains an expression of only
$l-1$ factors, but if it contains a shorter expression
\[
F=U\prod^{l-1}_{i=1} S^{\eta_i}_{\alpha_i}
\]
then
\[
F S^{\eta_l}_{\alpha_i}=U\prod^{l}_{i=1} S^{\eta_i}_{\alpha_i}
\]
belongs to the residue class defined by (1) and has less than $l$
factors, contrary to our assumption about $\mathfrak{U}F$. By
the induction hypothesis, (2) therefore contains a product
\[
\prod^{l-1}_{i=1} S^{\varepsilon_i}_{\beta_i}
\]
which satisfies ($\Sigma$). Now we take
\[
\prod^{l-1}_{i=1} S^{\varepsilon_i}_{\beta_i}\cdot S^{\eta_l}_{\alpha_i},
\]
which must have length $l$, as representative of $\mathfrak{U}F$,
and proceed similarly with all residue classes of length $l$. In this
way our assertion is proved.

\section{The relations of the second kind}

We now assume that the representatives $G$ satisfy ($\Sigma$).
The elements $U_{G,S}$ are divided into two classes: $U_{G,S}$ is
said to be of the first kind, or a member of $\mathfrak{u}_1$, if
\[
\overline{GS}=GS
\]
in the free group on the $S$. The remaining $U_{G,S}$ are said to
be of the second kind and they comprise the class $\mathfrak{u}_2$.
Now if $U_{G,S}$ is of the first kind then, because
\[
GS{\overline{GS}\,}^{-1}=1
\]
in the free group on the $S$, it also follows that 
\[
(GS{\overline{GS}\,}^{-1})_\mathfrak{U}=1
\]
in the free group on the $U_{G,S}$. That is, equation (2) of Section
3.5 can be replaced by
\[
U_{G,S}=1.
\]
In other words: \emph{the generators $U_{G,S}$ of class 
$\mathfrak{u}_1$ can be struck out, since they are equal to the
identity.}

We now investigate the elements of the class $\mathfrak{u}_2$.
Let
\[
G=\prod^{r}_{i=1} S^{\varepsilon_i}_{\alpha_i},\quad
\overline{GS}=\prod^{s}_{k=1} S^{\eta_k}_{\beta_k}.
\]
According to the rule for computing 
$(GS{\overline{GS}\,}^{-1})_\mathfrak{u}$ we have to introduce
representatives of those residue classes which are determined by
initial segments of the product
\[
\prod^{r}_{i=1} S^{\varepsilon_i}_{\alpha_i} S
\prod^{1}_{k=s} S^{-\eta_k}_{\beta_k}.
\]
By condition ($\Sigma$) we have the following series of expressions:
\begin{align*}
\overline{W}_0=\overline{W}_{r+s+1},\quad
&\overline{W}_j=\prod^{j}_{i=1} S^{\varepsilon_i}_{\alpha_i}\quad
(j=1,2,\ldots,r),\\
\overline{W}_{r+1}=\overline{GS}=\prod^{s}_{k=1} S^{\eta_k}_{\beta_k},
\quad
&\overline{W}_{r+i}=\prod^{s-i+1}_{k=1} S^{\eta_k}_{\beta_k}\quad
(i=2,3,\ldots,s).
\end{align*}
The following relations hold. For $\varepsilon_j=+1$,
\[
\overline{W}_{j-1} S^{\varepsilon_j}_{\alpha_j}=
\overline{\overline{W}_{j-1}S_{\alpha_j}}=\overline{W}_j
\quad (j\le r)
\]
and for $\eta_{s-i+1}=-1$,
\[
\overline{W}_{r+1} S^{-\eta_{s-i+1}}_{\beta_{s-i+1}}=
\overline{\overline{W}_{r+1}S_{\beta_{s-i+1}}}=
\overline{W}_{r+i+1}.
\]
For $\varepsilon_j=-1$,
\[
\overline{W}_j S^{-\varepsilon_j}_{\alpha_j}=
\overline{\overline{W}_j S_{\alpha_j}}=
\overline{W}_{j-1}\quad (j\le r)
\]
and for $\eta_{s-i+1}=+1$,
\[
\overline{W}_{r+i+1} S^{\eta_{s-i+1}}_{\beta_{s-i+1}}=
\overline{\overline{W}_{r+i+1} S_{\beta_{s-i+1}}}=
\overline{W}_{r+i},
\]
possibly after manipulations in the free group of the $S$.
Consequently, the elements
\[
|\overline{W}_{j-1} 
S^{\varepsilon_j}_{\alpha_j}{\overline{W}\,}^{-1}_j|_\mathfrak{U},\quad
|\overline{W}_{r+1}
S^{\eta_{s-i+1}}_{\beta_{s-i+1}}{\overline{W}\,}^{-1}_{r+i-1}|_\mathfrak{U}
\]
are all equal to $U_{G,S}$ from the class $\mathfrak{u}_1$, or inverses
of such elements, and thus equal to 1, so
\[
(GS{\overline{GS}\,}^{-1})_\mathfrak{U}=
|GS{\overline{GS}\,}^{-1}|_\mathfrak{U}
\]
and the equation in question reduces to the identity
\[
U_{G,S} U^{-1}_{G,S}=1.
\]

To summarize, we have arrived at the following: by suitable choice of
the expressions $|F|_\mathfrak{U}$ in Section 3.5, which are 
independent of the others, the equations (2) in 3.5 for $G\ne 1$ are
eliminated, so that only the equations for $G=1$ remain. If, in addition,
the representatives of the residue classes are chosen so as to satisfy
condition ($\Sigma$) in 3.6 then the result is that the relations of the
second kind may be replaced by those which say that the
\[
U_{G,S}\quad\text{for which}\quad\overline{GS}=\overline{G}S
\]
are equal to the identity.

Another remark on the computation of relations: if
\[
G=S^{\varepsilon_1}_{\alpha_1}S^{\varepsilon_2}_{\alpha_2}\cdots
    S^{\varepsilon_m}_{\alpha_m}
\]
and
\[
R=S^{\eta_1}_{\beta_1}S^{\eta_1}_{\beta_1}\cdots S^{\eta_n}_{\beta_n}
\]
then the subproducts $W_i$ of $GRG^{-1}$ are
\begin{align*}
W_j&=\prod^{j}_{i=1}S^{\varepsilon_i}_{\alpha_i}\qquad(j=1,2,\ldots,m)\\
W_{m+k}&=G\prod^{k}_{i=1}S^{\eta_i}_{\beta_i}\qquad(k=1,2,\ldots,n)\\
W_{m+n+l+1}&=GR\prod^{m-l}_{r=m}
                          S^{-\varepsilon_r}_{\alpha_r}\qquad(l=0,1,\ldots,m-1).
\end{align*}
Since
\[
\overline{GRG'}=\overline{GG'},
\]
because $RF$ and $F$ denote the same element in $\mathfrak{F}$,
then
\[
\overline{W}_l=\overline{W}_{m+n+m-l}
\]
and hence
\[
|\overline{W}_l S^{\varepsilon_{l+1}}_{\alpha_{l+1}}
{\overline{W}\,}^{-1}_{l+1}|_\mathfrak{U}=
|\overline{W}_{m+n+m-l} S^{\varepsilon_{l+1}}_{\alpha_{l+1}}
{\overline{W}\,}^{-1}_{m+n+m-l-1}|_\mathfrak{U}.
\]
Therefore, $(GRG^{-1})_\mathfrak{U}$ is converted by suitable
transformations into
\[
\prod^{n}_{i=1}|\overline{W}_{m+i-1} S^{\eta_{i}}_{\beta_{i}}
{\overline{W}\,}^{-1}_{m+i}|_\mathfrak{U}.
\]
If $G$ satisfies the conditions ($\Sigma$) then the collection of all
the elements $U_{G,S}$ determined by the initial segments of $G$
is contained in the class $\mathfrak{u}_1$, and can therefore be
left out.

The relationships obtained between the generators and defining
relations of a subgroup may be derived quite quickly in another
way.\footnote{\textsc{W. Hurewicz}, Hamb. Abhandl. \textbf{5}, 
1930, 307.} Let $\mathfrak{T}$ be the free group on the generators
$T_i$ of $\mathfrak{F}$, $\mathfrak{R}$ the invariant subgroup
determined by the defining relations $R_i(T)$, so $\mathfrak{F}$
equals $\mathfrak{T}/\mathfrak{R}$; let $\mathfrak{U}$ be a
subgroup of $\mathfrak{F}$, and $\mathfrak{U}_\mathfrak{T}$
that subgroup of $\mathfrak{T}$ the power products of which
yield elements of $\mathfrak{U}$. $\mathfrak{U}_\mathfrak{T}$
contains $\mathfrak{R}$, since all elements of $\mathfrak{R}$
are representable by the identity element of $\mathfrak{F}$ and
hence of $\mathfrak{U}$. Thus the residue classes $\mathfrak{U}G$
of $\mathfrak{U}$ in $\mathfrak{F}$ consist exactly of power
products of the residue classes $\mathfrak{U}_\mathfrak{T} L$
modulo $\mathfrak{U}_\mathfrak{T}$ in $\mathfrak{T}$, and a
full system of representatives $N_i(T)$ ($i=1,2,\ldots$) for the
$\mathfrak{U}G$ is at the same time a full system of
representatives for the $\mathfrak{U}_\mathfrak{T}L$.

Thus we have: the generators $U_{G,S}$ of $\mathfrak{U}$ are at
the same time a system of generators for the subgroup
$\mathfrak{U}_\mathfrak{T}$ and in fact they yield a system of free
generators of $\mathfrak{U}_\mathfrak{T}$ if the $N_i$ satisfy the
\textsc{Schreier} condition and those $U_{G,S}$ that equal 1
are left out.

We now ask more about the determination of $\mathfrak{R}$ as a
subgroup of $\mathfrak{U}_\mathfrak{T}$. The elements $R$ of
$\mathfrak{R}$ are power products of the elements
$L(T)R_k(T)L^{-1}(T)$, where $L$ is an arbitrary element of
$\mathfrak{T}$. Now if $L=U\overline{L}$, where $\overline{L}$ is
the representative of the residue class $\mathfrak{U}_\mathfrak{T}L$
and $U$ belongs to $\mathfrak{U}_\mathfrak{T}$, then
$LRL^{-1}=U\overline{L}R{\overline{L}\,}^{-1}U^{-1}$; thus the
elements of $\mathfrak{R}$ may be composed from the elements
$N_i R_k N^{-1}_i$ and their transforms in $\mathfrak{U}_\mathfrak{T}$.
Consequently, these $N_i R_k N^{-1}_i$, expressed in terms of the
$U_{G,S}$, yield a system of defining relations of $\mathfrak{U}$.

\section{Invariant subgroups}

If the subgroup $\mathfrak{U}$ is an invariant subgroup of
$\mathfrak{F}$ then each $F$ from $\mathfrak{U}$ yields an
automorphism of $\mathfrak{U}$ by the transformation
\[
FUF^{-1}=U',
\]
as we saw in Section 1.12. We now show that we can determine these
automorphisms. 

Obviously an automorphism
\[
\boldsymbol{A}(F)=F'
\]
of a group $\mathfrak{F}$ is determined when
\[
\boldsymbol{A}(S)=S'
\]
is given for all generators $S$. Because, since
\[
\boldsymbol{A}(F_1)\boldsymbol{A}(F_2)=
\boldsymbol{A}(F_1 F_2),
\]
we have
\[
\boldsymbol{A}
(S^{\varepsilon_1}_{\alpha_1}S^{\varepsilon_2}_{\alpha_2}\cdots
 S^{\varepsilon_m}_{\alpha_m})=
\boldsymbol{A}(S_{\alpha_1})^{\varepsilon_1}
\boldsymbol{A}(S_{\alpha_2})^{\varepsilon_2}\cdots
\boldsymbol{A}(S_{\alpha_m})^{\varepsilon_m}
\]
for each power product in the $S$. Applying this to the group
$\mathfrak{U}$ with the generators $U_{G,S}$, we need only
determine the elements
\[
F U_{G,S} F^{-1},
\]
i.e., we need only express the products
\[
F U_{G,S} F^{-1}
\]
in terms of the $U_{G,S}$ for all $G$ in $\mathfrak{g}$ and $S$ in
$\mathfrak{m}$. But this is possible immediately from the process
in Section 3.2. As $F$ runs through all elements of $\mathfrak{F}$
we obtain the totality of automorphisms of $\mathfrak{U}$, which
is a group. The automorphisms
\[
SUS^{-1}=U',
\]
as $S$ runs through the class $\mathfrak{m}$, can obviously be
chosen as generators of this group.

\section{Subgroups of special groups}

We now go to a few applications. If $\mathfrak{F}$ is a free group,
then each subgroup $\mathfrak{U}$ of $\mathfrak{F}$ is a free
group, as follows immediately from Section 3.7, generated by the
$U_{G,S}$ in the class $\mathfrak{u}_2$.\footnote{Other proofs
may be found in Sections 4.17, 4.20, and 7.12. For the literature
see the work of \textsc{O. Schreier} cited on p.~57 and 
\textsc{F. Levi}, Math. Zeit. \textbf{32}, (1930), 315.} It is in fact the case
that \emph{if $\mathfrak{F}$ is a group with generators}
\[
S_1,\quad S_2,\quad \ldots,\quad S_n
\]
\emph{and the relations}
\[
S^{a_i}_i\equiv 1,\qquad a_i\ge 0,
\]
\emph{then each subgroup $\mathfrak{U}$ of $\mathfrak{F}$ 
has a presentation with
the generators $U_{G,S}$ and defining relations of the
form}
\[
U^{u_{G,S}}_{G,S}\equiv 1.
\]

We prove this for the case where the $a_i\ne 0$ are all prime numbers.
The residue classes of
\[
GS^{r}_i,\quad 0\le r<a_i
\]
are either disjoint or all identical. Namely, if
\[
\mathfrak{U}GS^{r_i}_i=\mathfrak{U}GS^{r_2}_i,\quad r_1<r_2,
\]
then also
\[
\mathfrak{U}G=\mathfrak{U}GS^{(r_2-r_1)}_i.
\]
Further,
\[
\mathfrak{U}GS^{(r_2-r_1)}_i=\mathfrak{U}GS^{2(r_2-r_1)}_i,
\]
so
\[
\mathfrak{U}G=\mathfrak{U}GS^{2(r_2-r_1)}_i,
\]
and hence in general
\[
\mathfrak{U}G=\mathfrak{U}GS^{k(r_2-r_1)}_i.
\]
But, by Section 1.3, $k(r_2-r_1)$ runs through all residue classes
(mod $a_i$) and hence $\mathfrak{U}GS^{k(r_2-r_1)}_i$ runs
through all the residue classes $\mathfrak{U}GS^{r}_i$.

Further, if
\[
\mathfrak{U}GS^{r_i}_i=\mathfrak{U}G'S^{r'_i}_i
\]
then the residue classes
\[
\mathfrak{U}G'\text{ and }\mathfrak{U}GS^{r_i-r'_i}_i
\]
coincide, and in general so do the residue classes
\[
\mathfrak{U}G'S^{r'}_i\quad\text{and}\quad\mathfrak{U}GS^{r_i-r'_i+r'}_i.
\]

We now consider the consequence relations of the $S^{a_i}_i$ in
the $U_{G,S}$. They read
\begin{equation}
\prod^{a_i-1}_{r=0}\left|\overline{GS^{r}_i}S_i 
{\overline{GS^{r+1}_i}\,}^{-1}\right|_\mathfrak{U}=1. \tag{1}
\end{equation}

Now if the same generators appear in two such relations then for
two different representatives $G$ and $G'$ we must have
\[
\left| \overline{GS^{r_i}_i}S_i 
{\overline{GS^{r_i+1}_i}\,}^{-1}\right|_\mathfrak{U}=
\left| \overline{G'S^{r'_i}_i}S_i 
{\overline{G'S^{r'_i+1}_i}\,}^{-1}\right|_\mathfrak{U}
\]
and hence also
\begin{align*}
\mathfrak{U}G'S^{r'_i}_i&=\mathfrak{U}GS^{r_i}_i\\
\mathfrak{U}G'S^{r}_i&=\mathfrak{U}GS^{r_i-r'_i+r}_i
\end{align*}
so that the product
\[
\prod^{a_i-1}_{r=0}\left| \overline{G'S^{r}_i}S_i 
{\overline{G'S^{r+1}_i}\,}^{-1}\right|_\mathfrak{U}
\]
must result from (1) by a cyclic interchange, and hence is a 
consequence of (1). By omitting these superfluous relations we can
reach a stage where each generator $U_{G,S}$ apears in at most one
defining relation of $\mathfrak{U}$. But in such a relation either all
the $U_{G,S}$ are the same, in which case the relation has the form
\[
U^{u_{G,S}}_{G,S}=1,
\]
or else the $U_{G,S}$ appearing in the relation are formally different.
In the latter case either these $U_{G,S}$ all belong to the class
$\mathfrak{u}_1$ and hence equal the identity, or else they appear
as different elements of the class $\mathfrak{u}_2$, in which case
all but one of these generators can be eliminated by expressing
them in terms of the others. Thus our assertion is proved.

\section{Generators and defining relations of the congruence
subgroup $\mathfrak{U}_p$}

It follows from the result of the previous section that, e.g., all
subgroups of the modular group have relations of the form
\[
S^{a_i}_i\equiv 1\quad\text{with}\quad a_i=2\text{ or }3
\]
for suitable choice of generators $S_i$. We check this for the
congruence subgroups $\mathfrak{U}_p$ and establish the number
of these relations.\footnote{For the next two sections cf.
\textsc{H. Rademacher}, Hamb. Abhandl. \textbf{7}, (1930), 134.}

As generators of the modular group we take $S$ and $T$ with the 
defining relations
\[
S^2\equiv 1,\quad (ST)^3\equiv 1.
\]
The representatives of the residue classes $\mathfrak{U}_p M$
introduced in Section 1.9 are $G_k=ST^k$ ($k=0,1,\ldots,p-1$)
These products satisfy the condition ($\Sigma$) in 3.6. We now
construct the generators $U$ as in 3.1.
\begin{align*}
U_{G_k,T}&=G_k T \,{\overline{G_k T}\,}^{-1}=ST^k T(ST^{k+1})^{-1}
\qquad(k=0,1,\ldots,p-2)\\
U_{E,S}&=S{\overline{S}\,}^{-1}=SS^{-1}
\end{align*}
belong to the class $\mathfrak{u}_1$, the others to the class
$\mathfrak{u}_2$. We have
\begin{align*}
U_{E,T}=T;\quad&U_{G_{p-1},T}
=G_{p-1}T\,{\overline{G_{p-1}T}\,}^{-1}
=ST^{p}S^{-1};\\
U_{G_0,S}=S^2;\quad&U_{G_k,S}
=G_k S\,{\overline{G_k S}\,}^{-1}
=ST^k S\,{\overline{ST^k S}\,}^{-1}\qquad(k=1,2,\ldots,p-1)
\end{align*}
In order to determine the representatives
\[
\overline{G_k S}=\overline{ST^k S}=G_{k^{*}}
\]
we compute the substitution  corresponding to
\[
ST^{k}S=G_k S,
\]
namely
\[
x'=\frac{-1}{-\frac{1}{x}+k}=\frac{-x}{kx-1}.
\]
By Section 1.9
\[
kk^{*}\equiv -1\text{ (mod $p$)},\quad 0<k^{*}<p
\]
and so
\[
U_{G_k,S}=ST^k ST^{-k^*}S^{-1}.
\]
Now for the relations! They read:
\begin{align*}
(S^2)_\mathfrak{U}\equiv 1;\quad
\left((ST)^3\right)_\mathfrak{U}\equiv 1;&\quad
(G_k S^2 G^{-1}_k)_\mathfrak{U}\equiv 1;\quad
\left(G_k(TS)^3 G^{-1}_k\right)_\mathfrak{U}\equiv 1\\
&(k=0,1,\ldots,p-1).
\end{align*}
We begin with the consequences of the first relation,
\begin{align*}
(S^2)_\mathfrak{U}&
=|SG^{-1}_0|_\mathfrak{U}\;|G_0 S|_\mathfrak{U}
=U_{E,S}U_{G_0,S}=U_{G_0,S}\equiv 1,\\
(G_0 S^2 G^{-1}_0)_\mathfrak{U}&
=|G_0 S|_\mathfrak{U}\;|SG_0|_\mathfrak{U}
=U_{G_0,S}U_{E,S}=U_{G_0,S}\equiv 1,\\
(G_k S^2 G^{-1}_k)_\mathfrak{U}&
=|G_k S\,{\overline{G_k S}\,}^{-1}|_\mathfrak{U}\;
   |\overline{G_k S}\,SG^{-1}_k|_\mathfrak{U}
=U_{G_k,S}U_{G_{k^*},S}\equiv 1 \tag{1}
\end{align*}
The relations for $G_k$ and $G_{k^*}$ are convertible into each
other by cyclic interchange. By omission of superfluous relations
one can thus reach the stage where each generator appears in
only one relation.

Now for the consequences of the second relation:
\begin{align*}
(STSTST)_\mathfrak{U}=&\,
|SG^{-1}_0|_\mathfrak{U}\;
|G_0 TG^{-1}_1|_\mathfrak{U}\;
|G_1 S\,{\overline{G_1 S}\,}^{-1}|_\mathfrak{U}\;
|\overline{G_1 S}\, T\, {\overline{G_1 ST}\,}^{-1}|_\mathfrak{U}\\
&\cdot
|\overline{G_1 ST}\,S\,{\overline{G_1 STS}\,}^{-1}|_\mathfrak{U}\;
|\overline{G_1 STS}\,T\,{\overline{G_1 STST}\,}^{-1}|_\mathfrak{U}.
\end{align*}
Noting that $\overline{G_1 S}=G_{p-1}$, and hence 
$\overline{G_1 ST}=G_0$, it follows by omitting the generators
already known to be equivalent to 1 that
\begin{equation}
U_{G_1,S}U_{G_{p-1},T}U_{E,T}\equiv 1. \tag{2}
\end{equation}

Further
\begin{align*}
&\left(G_k(ST)^3 G^{-1}_k\right)_\mathfrak{U}\\
&=|G_k S\,{\overline{G_k S}\,}^{-1}|_\mathfrak{U}\;
     |\overline{G_k S}\,T\,{\overline{G_k ST}\,}^{-1}|_\mathfrak{U}\;
     |\overline{G_k ST}\,S\,{\overline{G_k STS}\,}^{-1}|_\mathfrak{U}\\
&\quad\cdot
     |\overline{G_k STS}\,T\,{\overline{G_k STST}\,}^{-1}|_\mathfrak{U}\;
     |\overline{G_k STST}\,S\,{\overline{G_k STSTS}\,}^{-1}|_\mathfrak{U}\;
   |\overline{G_k STSTS}\,T\,{\overline{G_k STSTST}\,}^{-1}|_\mathfrak{U}\\
&=|G_k S\,{\overline{G_k S}\,}^{-1}|_\mathfrak{U}\;
    |\overline{G_k S}\,T\,{\overline{G_k ST}\,}^{-1}|_\mathfrak{U}\;
    |\overline{G_k ST}\,S\,{\overline{G_k STS}\,}^{-1}|_\mathfrak{U}\\
&\quad\cdot
     |\overline{G_k STS}\,T\,{\overline{G_k T^{-1}S}\,}^{-1}|_\mathfrak{U}\;
  |\overline{G_k T^{-1}S}\,S\,{\overline{G_k T^{-1}}\,}^{-1}|_\mathfrak{U}\;
    |\overline{G_k T^{-1}}\,T\,G^{-1}_k|_\mathfrak{U}.
\end{align*}
In order to determine the factors $U_{G,T}$ more precisely we must determine the residue classes $\overline{G_k S}$, $\overline{G_k STS}$,
$\overline{G_k T^{-1}}$. We have
\begin{center}
$\overline{G_k S}=G_{k^*}$ 
\; for \; $k\ne 0$, \; $\overline{G_0 S}=E$,\\
$\overline{G_k STS}=\overline{G_k T^{-1}ST^{-1}}=G_{(k-1)^*-1}$ 
\; for\;  $k\ne 1$, \; $=E$\;  for\; $k=1$,\\
$\overline{G_k T^{-1}}=G_{k-1}$\;  for\; $k\ge 1$, 
\;$=G_{p-1}$\; for\; $k=0$.
\end{center}
Thus the $U_{G,T}$ which do not belong to $\mathfrak{u}_1$ appear
in the first place for $k=0$ and 1, in the second place for $k=1$, in the
last place for $k=0$. Hence we completely determine the relations for
$G_0$ and $G_1$ first. One sees that
\[
G_1 (ST)^3 G^{-1}_1=(ST)^3
\]
in the free group on $S$ and $T$, and hence also
\[
\left(G_1(ST)^3 G^{-1}_1\right)_\mathfrak{U}
=\left((ST)^3\right)_\mathfrak{U}
\]
in the free group on the $U$, and that
\[
G_0 (ST)^3 G^{-1}_0=S(ST)^3 S^{-1}=S^2 T(ST)^3 T^{-1}S^{-2}
\]
in the free group on $S$ and $T$ and hence
\[
\left(G_0 (ST)^3 G^{-1}_0=S(ST)^3 S^{-1}\right)_\mathfrak{U}=
\left(S^2 T\right)_\mathfrak{U}
\left((ST)^3\right)_\mathfrak{U} 
\left(T^{-1}S^{-2}\right)_\mathfrak{U}
\]
in the free group on the $U$, and consequently these relations can
be omitted, as consequence relations of 
$\left((ST)^3\right)_\mathfrak{U}$.

There remain the relations free of the $U_{G,T}$:
\begin{align*}
&\quad\left(G_k(ST)^3 G^{-1}_k\right)_\mathfrak{U}\\
&=|G_k S\,{\overline{G_k S}\,}^{-1}|_\mathfrak{U}\;
     |\overline{G_k ST}\,S\,{\overline{G_k STS}\,}^{-1}|_\mathfrak{U}\;
     |\overline{G_k STST}\,S\,{\overline{G_k STSTS}\,}^{-1}|_\mathfrak{U}\\
&\quad(k=2,3,\ldots,p-1). \tag{3}
\end{align*}
One sees that the relations for
\[
G_k, G_{k^*+1}=\overline{G_k ST}\text{ and }
G_{(k^*+1)^*+1}=\overline{G_k(ST)^2}
\]
are convertible into each other by cyclic interchange, and one concludes
as in Section 3.9 that two relations (3) are always convertible into each
other by cyclic interchange when they contain the same generators.
Thus by omitting superfluous relations one can reach the stage
where each generator in (3) appears in only one relation. It is clear,
just as in Section 3.9, that the three generators that appear in a
relation are either all equal or all different.

For the process $k'=k^*+1$ this says
\[
k'''=((k^*+1)^*+1)^*+1=k
\]
and, if two of the numbers
\[
k,\quad k'=k^*+1,\quad k''=(k^*+1)^*+1
\]
are equal, then they are all equal.

\section{The relations $U^{u_{G,S}}_{G,S}$ of the group $\mathfrak{U}_p$}

The reduction of the presentation of $\mathfrak{U}_p$, by elimination 
of suitable $U_{G,T}$ until only relations of the form
\begin{equation}
U^{u_{G,S}}_{G,S}=1
\end{equation}
remain, will not be carried out in general. Here we shall work under
the assumption that only relations of the form (1) appear, i.e., that
$k=k^*$ and $k^*+1=k$. For the number of relations (1) that
remain after omission of superfluous ones, see the work cited in
Section 3.10. There one finds an elimination process for
the $U_{G,S}$.

If $x^*=x$ then
\[
x^2\equiv -1\text{ (mod $p$),}
\]
and conversely, if $y=y^*+1$ then
\[
yy^*=y(y-1)\equiv -1\text{ (mod $p$).}
\]
It must then be that
\begin{equation}
y^2-y+1\equiv 0\text{ (mod $p$)} \tag{2}
\end{equation}
or
\begin{equation}
(2y-1)^2\equiv -3\text{ (mod $p$)}. \tag{3}
\end{equation}
Conversely, each solution of the congruence (2) follows from that of
(3) and hence
\[
y=y^*+1.
\]

Now the theory of quadratic forms\footnote{Specifically, the
quadratic reciprocity theorem. (Translator's note.)} shows that the
congruence
\begin{equation}
x^2\equiv a\text{ (mod $p$)} \tag{4}
\end{equation}
(for $p$ an odd prime number) has either no solutions or two when
$a$ is not divisible by $p$, and $a$ is called a quadratic residue when
the congruence (4) is solvable, otherwise a nonresidue. We let
\[
\left(\frac{a}{p}\right)
\]
be the symbol which equals $+1$ or $-1$ according as $a$ is a
quadratic residue or nonresidue (mod $p$).

Now $\left(\frac{-1}{p}\right)=+1$ when $p=4n+1$, for any integer
$n$, otherwise it equals $-1$, and $\left(\frac{-3}{p}\right)=+1$
when $p=3n+1$, otherwise it equals $-1$.

Because of this, there are either exactly two $U_{G,S}$ with
$U^{2}_{G,S}=1$, or else none, according as $p=4n+1$ or not,
and there are two $U_{G,S}$ with $U^{3}_{G,S}=1$ according as
$p=3n+1$ or not. For $p=2$, (4) has the solution $x\equiv 1$,
for $p=3$, (3) has the solution $y\equiv -1$ (mod 3).

We give a few more numerical examples.

Firstly, we can always eliminate $U_{G_0,S}$ by Section 3.10, (1),
and $U_{G_{p-1},T}$ with the help of Section 3.10, (2).

In the case $p=2$ only two generators remain,
\[
U_{E,T},\quad U_{G_1,S},
\]
and the relation
\[
U^{2}_{G_1,S}\equiv 1.
\]

For $p=3$ the generators remaining are
\[
U_{E,T},\quad U_{G_i,S}\quad (i=1,2) 
\]
with the relations
\[
U_{G_1,S} U_{G_2,S}\equiv 1,\quad U^{3}_{G_i,S}\equiv 1.
\]
We can eliminate $U_{G_2,S}$ and obtain $U_{G_1,S}$ with
$U^{3}_{G_1,S}\equiv 1$.

In the case $p=5$ there remain initially 
\begin{align*}
&U_{E,T},\quad U_{G_i,S}\quad (i=1,2,3,4)\\
U_{G_1,S} &U_{G_4,S}\equiv 1,\quad
U^{2}_{G_2,S}\equiv 1,\quad U^{2}_{G_3,S}\equiv 1,\\
&U_{G_1,S}U_{G_3,S}U_{G_4,S}\equiv 1.
\end{align*}
We eliminate $U_{G_1,T}$ and $U_{G_4,T}$ and obtain
\[
U_{E,T},\quad U_{G_i,S};\qquad U^{2}_{G_i,S}\equiv 1\quad (i=2,3).
\]

For $p=7$ we get
\[
U_{E,T},\quad U_{G_i,S};\qquad U^{3}_{G_i,S}\equiv 1\quad (i=3,5).
\]

For $p=11$ we get the free group with the generators
\[
U_{E,T},\quad U_{G_i,S} \quad (i=4,6).
\]

\section{Commutator groups}

In the factor group modulo the commutator subgroup it can always
be decided whether two power products belong to the same residue
class, hence a system of representatives for these residue classes,
and thereby generators of the commutator subgroup, may always
be given.

We now focus on the case where $\mathfrak{F}$ is determined by two 
generators $S_1$, $S_2$ and a relation $R(S_1,S_2)$, and the factor
group $\mathfrak{A}$ of $\mathfrak{F}$ by the commutator
subgroup $\mathfrak{K}$ is infinite cyclic.

Let
\[
R(S)=\prod^{m}_{i=1} S^{s_{1i}}_1 S^{s_{2i}}_2
\]
and
\[
s_k=\sum^{m}_{i=1} s_{ki}\quad(k=1,2),
\]
so
\[
R'(S)=S^{s_1}_1 S^{s_2}_2 \equiv 1
\]
is the defining relation of $\mathfrak{A}$. And, since $\mathfrak{A}$
is an infinite cyclic group, the greatest common divisor of $s_1$ and
$s_2$, $(s_1,s_2)$, must be 1. If $s_2=0$, then $s_1=\pm 1$. If $s_1$
and $s_2$ are both $\ne 0$ then new generators $s'_1$ and $s'_2$
may always be introduced, and the relation $R(S)$ converted to
\[
R(S')=\prod^{m'}_{i=1} S'^{s'_{1i}}_1 S'^{s'_{2i}}_2,
\]
so that
\[
s'_1=\sum^{m'}_{i=1}s'_{1i}=0.
\]
This follows by induction from the following fact: if
\[
|s_1|>|s_2|\quad\text{and}\quad s_1=ns_2+r_1\quad (0\le r_1<|s_2|)
\]
and if the generators $S'_i$ are introduced by
\[
S'_2=S_2 S^{n}_1,\quad S'_1=S_1,
\]
then
\[
R(S')=\prod^{m'}_{i=1} S'^{s_{1i}}_1 (S'_2 S'^{-n}_1)^{s_{2i}}
\]
and therefore
\[
s'_1=s_1+ns_2=r_1.
\]

We now suppose that 
\[
s'_1=0,\quad s'_2=1
\]
and write $S$ and $K$ for $S'_1$ and $S'_2$ respectively, in order to
indicate that $S'_2=K$ belongs to the commutator subgroup
$\mathfrak{K}$. The elements
\[
G_s=S^s\quad(s=0,\pm 1, \pm 2,\ldots)
\]
then represent all the residue classes $\mathfrak{K}F$. They satisfy
the condition ($\Sigma$). An element
\[
\prod^{l}_{i=1}S^{s_i}K^{k_i}\quad
\text{belongs to the residue class $\mathfrak{K}G_s$ with}\quad
s=\sum^{m}_{i=1} s_i.
\]
So generators of $\mathfrak{K}$ are the elements
\[
U_{G_s,S}=S^{s} S {\overline{S^{s+1}}\,}^{-1},
\]
which belong to the class $\mathfrak{u}_1$, and
\[
U_{G_s,K}=S^{s}K S^{-s}=K_s\quad(s=0,\pm 1,\pm 2,\ldots),
\]
which belong to the class $\mathfrak{u}_2$.

Now let
\begin{equation}
R(S,K)=\prod^{m}_{i=1} S^{r_i} K^{\varepsilon_i} S^{-r_i}. \tag{1}
\end{equation}
Then
\[
G_s R G^{-1}_s=
\prod^{m}_{i=1} S^{r_i+s} K^{\varepsilon_i} S^{-r_i-s}, \tag{2}
\]
and hence all relations
\begin{equation}
(G_s R G^{-1}_s)_\mathfrak{U} \tag{3}
\end{equation}
result from $(R)_\mathfrak{U}$ when $K_i$ is replaced by $K_{i+s}$.

The automorphism induced in $\mathfrak{K}$ by transformation by
$S$ is determined by
\[
SK_s S^{-1}=K_{s+1}.
\]
The collection of elements from the residue classes
$\mathfrak{K}S^{gn}$ ($g$ fixed; $n=0,\pm 1,\pm 2,\ldots$) is an
invariant subgroup $\mathfrak{K}_g$ of $\mathfrak{F}$.
The elements
\[
G_i=S^i \quad (i=0,1,\ldots,g-1)
\]
form a complete system of representatives for the residue classes
$\mathfrak{K}_g F$. Consequently, the elements
\[
U'_{G_i,S},\quad U'_{G_i,K}\quad (i=0,1,\ldots,g-1)
\]
form a system of generators for $\mathfrak{K}_g$. The
\[
U'_{G_i,S}\quad (i=0,1,\ldots,g-2)
\]
belong to the class $\mathfrak{u}_1$, while
\[
U'_{G_{g-1},S}=S^g, \quad U'_{G_i,K}=S^i K S^{-i}\quad
(i=0,1,\ldots,g-1)
\]
belong to the class $\mathfrak{u}_2$. The defining relations of
$\mathfrak{K}_g$ are
\[
R_i=(S^i R S^{-i})_\mathfrak{U}\quad(i=0,1,\ldots,g-1).
\]
They result from the relations (2) when one sets
\[
U_{G_k,S}=U'^{l}_{G_{g-1},S} U'_{G_i,S} U'^{-l}_{G_{g-1},S}
\]
with
\[
lg+i=k\quad\text{for}\quad (0\le i <g).
\]
If we make the elements of the commutator subgroup commute
and set $SKS^{-1}=K^x$, then we obtain a group with  operator
determined by the generator $K$ and the relation 
\begin{equation}
K^{f(x)}\equiv 1, \tag{4}
\end{equation}
which follows from (1). This is because the relations resulting from
(2) are
\[
K^{x^i f(x)}
\]
and hence they are consequence relations of (4) in the sense of
Section 2.14.

\section{The Freiheitssatz (the freeness theorem)}

By considering invariant subgroups with infinite cyclic factor 
groups one can solve the word problem for groups with one
defining relation. This is achieved by the following, so-called
Freiheitssatz,\footnote{It is usual to use the German name for 
this theorem. (Translator's note.)}\footnote{\textsc{W. Magnus}, 
J. f\"ur reine und angew. Math. \textbf{163}, (1930), 3.} formulated by
\textsc{Dehn}, which states:

\emph{If $S_i$ ($i=1,2,\ldots,n$) are generators of a group 
$\mathfrak{G}$ with one defining relation $R(S_i)$ which is a short
word in the sense of Section 2.4 and which properly contains $S_n$,
then the group generated by the $S_i$ ($i=1,2,\ldots,n-1$) is free.}

The theorem holds if $R$ does not contain the generators
$S_i$ ($i=1,2,\ldots,n-1$), because $\mathfrak{G}$ is then the free
product of the free group on the $S_i$ ($i=1,2,\ldots,n-1$) and the
group generated by $S_n$ with $R(S_n)\equiv 1$.

We can then carry out the proof by complete induction on the
length $l(R)$, which is the sum of the absolute values of the 
exponents of all $S_i$ in $R$, and assume the theorem proved 
for $l-1$, because for $l=1$ we must have $R=S^{\pm 1}_n$. 
We can also assume that all $S_i$ really appear in $R$.

By allowing factors to commute, $R$ may be brought into the form
\[
S^{r_1}_1 S^{r_2}_2\cdots S^{r_n}_n.
\]
We first make the assumption \emph{that $r_n$ is zero}. With
addition of the relations $S_i\equiv 1$ ($i=1,2,\ldots,n-1$),
$\mathfrak{G}$ yields the free group with one generator $S_n$,
and if we set
\[
S^i_n S_k S^{-i}_n=S_{ki}
\]
and denote the relations resulting from $S^i_n R S^{-i}_n$ by
introducing the $S_{lm}$ by $R_i(S_{lm})$, then $R_i$ results from
$R_{i-1}$ when each $S_{lm}$ is replaced by $S_{l,m+1}$, and the
power products $R_i(S_{lm})$ thus have the same length $l'$, and
in fact $l'<l$ because in $R_i$ no element corresponds to the 
factors $S^k_n$. 

By $\mathfrak{G}_{ik}$ we mean the group 
presented on those generators
that appear in $R_i,R_{i+1}$, $\ldots$, $R_{i+k-1}$, and the $S_{lm}$
for which $m_1<m<m_2$ and $S_{lm_1}$ and $S_{lm_2}$ both 
appear in $R_i, R_{i+1},\ldots$ or $R_{i+k-1}$, and with the 
defining relations $R_i,R_{i+1},\ldots,R_{i+k-1}$.

By $\mathfrak{U}_i$ we mean the group with generators common
to $\mathfrak{G}_{i1}$ and $\mathfrak{G}_{i+1,1}$. Obviously
$\mathfrak{U}_i$ is a proper subgroup of $\mathfrak{G}_{i1}$ 
and $\mathfrak{G}_{i+1,1}$, and hence it is a free group by the
induction hypothesis. It follows that $\mathfrak{G}_{i,k+1}$ is the
free product of $\mathfrak{G}_{i,k}$ and $\mathfrak{S}_{i+k,1}$
with the amalgamated subgroup $\mathfrak{U}_{i+k}$
(Section 2.7). Finally, let $\mathfrak{B}_{ika}$ be those subgroups
of $\mathfrak{G}_{ik}$ generated by the elements $S_{la}$
($a$ fixed; $l=1,2,\ldots, n-1$) appearing among the generators
of $\mathfrak{G}_{ik}$. The groups $\mathfrak{B}_{i1a}$ are free
groups with the free generators $S_{la}$, for generators appearing in
$R_i$ cannot have the same second index, since it was assumed that
$S_n$ did not appear in $R$. Thus it follows from the induction 
hypothesis that the $\mathfrak{B}_{i1a}$ are free groups. However,
$\mathfrak{B}_{i,k+1,a}$ is the free product of $\mathfrak{B}_{ika}$
and $\mathfrak{B}_{i+k,1a}$ and hence likewise a free group on the
free generators $S_{la}$ (Section 2.7).

Now we can derive an absurdity from the assumption that a
relation $R'(S)$ between the $S_i$ ($i=1,2,\ldots,n-1$) follows
from $R$, because $R'$ would yield a relation between the $S_{la}$.

The considerations are similar \emph{when $r_k$ is nonzero} and
one of the $r_i$, say $r_1$, is zero. By addition of relations $S_i=1$
($S_i=2,3,\ldots,n$) we get the free cyclic group with generator
$S_1$. Now let
\[
S^i_1 S_k S^{-i}_1=S_{ki}
\]
and let $R_i$ be the relation resulting from $S^i_1 R S^{-i}_1$
by introduction of the $S_{lm}$. We retain the definition of the
$\mathfrak{G}_{ik}$ and $\mathfrak{U}_i$, and understand
$\mathfrak{S}_{ik}$ to be the subgroup of $\mathfrak{G}_{ik}$
generated by the $S_{lm}$ from $\mathfrak{G}_{ik}$ with
$l\ne n$. $\mathfrak{G}_{ik}$ is free, because one $S_{nk}$
properly appears in $R_i$, and it again follows from the properties
of free products with amalgamated subgroups that all the
$\mathfrak{S}_{ik}$ are free groups with free generators $S_{lm}$
($l\ne n$). If there were a relation $R'$ between the $S_i$
($i=1,2,\ldots,n-1$) alone, then this would have a relation
between the $S_{lm}$ ($l\ne n$) as a consequence.

If all $r_i\ne 0$ and $n=2$ then a relation for $S_1$ alone cannot
follow from $R$. Namely, if 
\[
\prod L_i R^{\varepsilon_i}_i L^{-1}_i=S^a_1
\]
then it follows, when we allow the factors on the left side to commute
and set $\sum \varepsilon_i=\varepsilon$, that the left side 
becomes $S^{\varepsilon r_1}_1 S^{\varepsilon r_2}_2$.
Consequently, we must have $\varepsilon=0$ and hence also
$a=0$. 

Finally, if all $r_i\ne 0$ and $n\ge 3$ we extend the group
$\mathfrak{G}$ by a generator $T_1$ with the relation
\[
T^{r_2}_1 S^{-1}_1=1.
\]
If the group with the generators $T_1,S_2,\ldots,S_{n-1}$ is free,
so also is the group generated by $T^{r_2}_1=S_1,S_2,\ldots,S_{n-1}$,
being a subgroup of a free group (Section 3.9), and in fact
$S_1,S_2,\ldots,S_{n-1}$ are free generators.

In place of $S_2$ we introduce the new generator $T_2$ into $R$ by
\[
S_2=T^{-r_1}_1 T_2.
\]
Suppose $R(S)$ is converted into $\overline{R}(T,S)$ by elimination
of $S_1$ and $S_2$. If we make the $T$ commute with the $S$ then
the element $T_1$ disappears. Further, if we construct the
invariant subgroup generated by
\[
T^i_1 T_2 T^{-i}_1=T_{2i},\quad T^i_1 S_k T^{-i}_1=S_{ki}
\]
and express
\[
T^i_1 \overline{R}(T,S) T^{-i}_1 = \overline{R}_i
\]
in terms of the $T_{2i},S_{ki}$, then this results in a word of length
$l'$ less than that of $R$. For no factors in $\overline{R}_i$ correspond
to the $T^k_1$ from $\overline{R}$, and the sum of the absolute
values of the remaining terms is equal to the sum of the exponent
values of the factors of $R$ different from $S_1$. Thus the group
generated by $T_1,T_2,S_3,\ldots,S_{n-1}$ is a free group with 
these free generators.

Among the consequences of the Freiheitssatz, one deserves
particular attention:\footnote{For proof see \textsc{Magnus} loc. cit.}

\emph{If $\mathfrak{G}$ and $\mathfrak{G}'$ are two groups with
generators $S_i, S'_i$ respectively, and defining relations $R(S),R'(S')$
respectively, and if the correspondence $\boldsymbol{I}(S_i)=S'_i$
is an isomorphism between $\mathfrak{G}$ and $\mathfrak{G}'$,
then $\boldsymbol{I}(R(S))$ is an element which results from $R'(S')$
by transformation in the free group of the $S$.} For more on the
solution of the word problem in groups with one defining relation we
refer to \textsc{W. Magnus} (Math. Ann. \textbf{106}, 295).

\section{Determination of automorphisms}

A general procedure for determining the automorphism group from
the generators and defining relations is not known. We shall collect
the most important results. If $S_i$ ($i=1,2,\ldots,n$) are  
generators of a group then an automorphism $\boldsymbol{A}$ is
determined when the elements $\boldsymbol{A}(S_i)=S'_i$ are
known as power products of the $S_i$. This is because any power 
product of the $S_i$ then goes to the element resulting from
formal replacement of the $S_i$ by $S'_i$.

Given a \emph{free commutative group} with $n$ free generators
$S_i$ ($i=1,2,\ldots,n$) and an automorphism
\[
\boldsymbol{A}(S_i)=S^{a_{i1}}_1 S^{a_{i2}}_2\cdots S^{a_{in}}_n,
\]
the determinant of the $a_{ik}$ must equal $\pm 1$; this is because the
elements $S'_i=\boldsymbol{A}(S_i)$ must again be free generators
of the group and the equations
\[
\prod_i S'^{x_i}_i=\prod_k S^{b_k}_k,
\]
respectively
\[
\sum a_{ik} x_k=b_k,
\]
must therefore be satisfied. Conversely, any substitution with
determinant $\pm 1$ defines an automorphism.

Now suppose we have a \emph{free group $\mathfrak{S}$ with free
generators} $S_i$ ($i=1,2,\ldots,n$). A closed presentation of the
automorphism group must be very intricate, but it is worth asking
about its generators and defining relations.

The automorphisms of $\mathfrak{S}$ include all permutations of
the $S_i$,
\begin{equation}
\boldsymbol{A}(S_i)=S_{k_i}\quad (i=1,2,\ldots,n). \tag{1}
\end{equation}
Furthermore, the elements
\begin{equation}
\boldsymbol{A}(S_1)=S^{-1}_1,\quad \boldsymbol{A}(S_i)=S_i
\quad (i\ne 1) \tag{2}
\end{equation}
and the elements
\begin{equation}
\boldsymbol{A}(S_1)=S_1 S_2,\quad\boldsymbol{A}(S_i)=S_i
\quad (i\ne 1) \tag{3}
\end{equation}
are again systems of free generators, because one can represent
the $S_i$ in terms of the $\boldsymbol{A}(S_i)$. 
\textsc{Neilsen}\footnote{\textsc{J. Nielsen}, Math. Ann. \textbf{79}
(1919), 269, and \textbf{91} (1924), 169.} showed that these
operations generate the automorphism group of $\mathfrak{S}$.

The automorphisms $\boldsymbol{A}(S_i)=S'_i$ obviously induce
certain automorphisms in the factor group $\mathfrak{F}$ by the
commutator subgroup. One obtains them by bringing the power
products for the $S'_i$ into the form
\[
S^{a_{i1}}_1 S^{a_{i2}}_2\cdots S^{a_{in}}_n
\]
by transposition of factors. Since $\mathfrak{F}$ is a free
\textsc{Abelian} group, the determinant of the $a_{ik}$ must equal
$\pm 1$.

Conversely, each automorphism of a free \textsc{Abelian} group 
$\mathfrak{F}$ is induced by an automorphism of $\mathfrak{S}$.
Namely, if one considers the $S_i$ in (1), (2), and (3) as generators
of $\mathfrak{F}$, then one sees that all automorphisms of
$\mathfrak{F}$ may be composed from the automorphisms (1), (2),
and (3).

A beautiful result of \textsc{Nielsen}\footnote{\textsc{J. Nielsen},
Math. Ann. \textbf{78} (1918), 385.} sharpens this connection in
the case of the free group $\mathfrak{S}_2$ on two generators.

If $\boldsymbol{A}_1$ and $\boldsymbol{A}_2$ are two
automorphisms of $\mathfrak{S}_2$ which induce the same
automorphism in the factor group by the commutator subgroup,
then there is an inner automorphism $\boldsymbol{I}$ of
$\mathfrak{S}_2$ such that $\boldsymbol{I}$ composed with
$\boldsymbol{A}_1$ yields $\boldsymbol{A}_2$, i.e.,
$\boldsymbol{A}_2=\boldsymbol{I}\boldsymbol{A}_1$.

The results for \emph{free products of finite cyclic groups} are
much simpler. We consider only groups\footnote{\textsc{O. Schreier},
Hamb. Abhandl. \textbf{3} (1924), 167.} with two generators $S_1,S_2$
and the relations
\begin{equation}
S^{a_1}_1\equiv S^{a_2}_2\equiv 1,\quad a_1>a_2. \tag{4}
\end{equation}

The automorphisms are then given by the equations
\[
\boldsymbol{A}(S_1)=LS^{r_1}_1 L^{-1},\quad
\boldsymbol{A}(S_2)=LS^{r_2}_2 L^{-1}
\]
where $r_i$ and $a_i$ are relatively prime to each other.

It is clear, first of all, that that the given transformations are
automorphisms. Conversely, in automorphisms 
$\boldsymbol{A}(S_i)=S'_i$ the element $S'_i$ must have order $a_i$,
from which it follows by Section 2.6 that the $S'_i$ must be transforms
of $S^{l_1}_1$, $S^{l_2}_2$ respectively. But the two $S'_i$ cannot be
transforms of powers of the same element $S_b$. Because power
products of the $S'_i$ would then all belong to residue classes,
modulo the commutator subgroup, represented by $S^n_b$
($n=0,\pm 1,\ldots$). Thus we must have
\[
\boldsymbol{A}(S_1)=LS^{r_1}_1 L^{-1},\quad
\boldsymbol{A}(S_2)=MS^{r_2}_2 M^{-1}.
\] 
Then 
\[
\overline{\boldsymbol{A}}(S_1)=M^{-1}LS^{r_1}_1 L^{-1}M,\quad
\overline{\boldsymbol{A}}(S_2)=S^{r_2}_2
\]
and
\[
\boldsymbol{A}^*(S_1)=S^l_2 M^{-1}LS_1L^{-1}MS^{-l}_2,\quad
\boldsymbol{A}^*(S_2)=S_2
\]
are also automorphisms. The element 
$S^l_2 M^{-1}LS_1L^{-1}MS^{-l}_2$ may now be reduced by the
method of Section 2.6, and $l$ chosen, so that the reduced form
begins and ends with $S^{\pm m}_1$. Now if $\boldsymbol{A}^*(S_1)$
were different from $S_1$, say equal to $NS_1 N^{-1}$, then each power
product of the $\boldsymbol{A}^*(S_i)$ which contained $S_1$ at all
would also contain the product $N$, which is impossible.

Closely related to the above are the automorphisms of the 
\emph{group $\mathfrak{G}$ with generators $S_1,S_2$ defined 
by}\footnote{In the special case $a_1=3$, $a_2=2$ the
automorphisms were found by \text{M. Dehn}, Math. Ann. \textbf{75}
(1914), 402. In this case the group is the group of the trefoil knot, and \textsc{Dehn} used its automorphisms to show that the left trefoil
knot is not deformable into the right trefoil knot. (Translator's note.)}
\[
S^{a_1}_1 S^{a_2}_2\equiv 1, \quad a_1>a_2.
\]
As we saw in Section 2.6, the infinite cyclic subgroup generated by
$S^{a_1}=D$ is the center $\mathfrak{Z}$, and the factor group by
$\mathfrak{Z}$ is the group defined by (4). An automorphism 
$\boldsymbol{A}$ of $\mathfrak{G}$ carries the center into
itself and must induce automorphisms in $\mathfrak{Z}$ on the one
hand and in the factor group $\mathfrak{G}/\mathfrak{Z}$ on the
other. The only automorphisms of $\mathfrak{Z}$ are
$\boldsymbol{A}(D)=D$ and $\boldsymbol{A}(D)=D^{-1}$. Hence
\begin{equation}
(\boldsymbol{A}(S_1))^{a_1}=S^{\varepsilon a_1}_1=D^\varepsilon,
\quad
(\boldsymbol{A}(S_2))^{a_2}=S^{\varepsilon a_2}_2=D^{-\varepsilon}
\quad (\varepsilon=\pm 1); \tag{5}
\end{equation}
because of the second condition an automorphism of $\mathfrak{G}$
must have the form
\[
\boldsymbol{A}(S_i)=MS^{r_i}_i M^{-1} D^{s_i}
                               =MS^{r_i+a_i s'_i}_i M^{-1}
\]
where
\[
s'_1=s_1,\quad s'_2=-s_2.
\]
Then because of (5) we must have
\[
(r_i+a_i s'_i)a_i=\varepsilon a_i,
\]
so
\[
r_i+a_i s'_i=\varepsilon.
\]
Thus we have: the automorphisms of the group $\mathfrak{G}$ are
\[
\boldsymbol{A}(S_i)=MS^\varepsilon_i M^{-1}\quad
(\varepsilon=\pm 1; i=1,2).
\]
\chapter{Line Segment Complexes}

\section{The concept of a line segment complex}

We now turn to the simplest objects of combinatorial topology,
the line segment complexes. In particular, we have to consider
paths in line segment complexes and coverings of line segment
complexes, two concepts which facilitate the connection between
line segment complexes and groups.

By a \emph{line segment complex}\footnote{The reader may
prefer the shorter term \emph{graph} or \emph{1-complex}.
However, since Reidemeister uses not only the term ``Streckencomplex''
but also its parts ``Strecken'' (line segments) and 
``Complex'' (complex) I have decided to keep them all. It may be
wordy, but at least the words fit together better than ``graph,'' ``edge,''
and ``complex.'' (Translator's note.)} 
$\mathfrak{C}$ we mean a finite or denumerably 
infinite collection of points and line segments. The relations 
``is the initial point of,'' ``is the final point of,'' and ``is equal but
oppositely directed'' between these objects may be explained with
the help of the following declarations.

A.1. \emph{If $s$ is a segment then there is always a point $p_1$
which is the initial point of $s$, and a point $p_2$ which is the
final point of $s$.}

A.2. \emph{If $s$ is a segment then there is a single oppositely
directed segment equal to $s$, $s'=s^{-1}$. The oppositely
directed segment equal to $s^{-1}$, namely $(s^{-1})^{-1}$, is $s$.}

A.3. \emph{If $p_1$ and $p_2$ are initial and final points of $s$,
then $p_2$ and $p_1$, respectively, are the initial and final points
of $s^{-1}$.}

If $p$ is an initial or final point of $s$ then $p$ is called a boundary
point of $s$, and we say that $p$ bounds $s$. If the initial and final
points $p_1$ and $p_2$ of $s$ are different then the point pair
$p_1,p_2$ is called the boundary of $s$. If $p_1=p_2=p$ then $s$ is
called a \emph{singular} segment and $p$ is called the boundary of
$s$. If all points and segments of
a complex $\mathfrak{C}'$ also belong to the complex $\mathfrak{C}$,
and if the boundary relations in $\mathfrak{C}'$ also hold in
$\mathfrak{C}$, then $\mathfrak{C}'$ is called a \emph{subcomplex}
of $\mathfrak{C}$.

Now let the pairs of equal but oppositely directed line segments be 
given a fixed numbering, with one segment of the pair 
associated with the symbol $s_i=s^1_i$ and the other with
$s^{-1}_i$. Finitely many line segments, given in a certain order,
\begin{equation}
w=S^{\varepsilon_1}_{\alpha_1}S^{\varepsilon_2}_{\alpha_2}\cdots
     S^{\varepsilon_m}_{\alpha_m}\quad(\varepsilon_i=\pm 1)
\tag{1}
\end{equation}
constitute a \emph{path} when
the final point of $S^{\varepsilon_i}_{\alpha_i}$ is the initial point
of $S^{\varepsilon_{i+1}}_{\alpha_{i+1}}$.

Thus $w$ ``traverses'' the segments $S^{\pm 1}_{\alpha_i}$.
A path (1) is called ``reduced'' when it is never the case that
\[
\alpha_i=\alpha_{i+1}\quad\text{and}\quad
\varepsilon_i+\varepsilon_{i+1}=0
\]
hold simultaneously. Striking out or inserting two such segments
is called reduction or extension of a path, respectively. A path is 
called open or closed according as the initial point $p$ of
$S^{\varepsilon_1}_1$ and the final point $p'$ of
$S^{\varepsilon_n}_n$ coincide or not. An open path ``leads''
from $p$ to $p'$ or ``connects'' $p$ with $p'$. A path (1) is
called \emph{simple} if the endpoints of the segments
$S^{\varepsilon_i}_i$ are all different from each other. The
\emph{length of a path} is the number $m$ of segments which
it contains. If (1) is a path then
\[
S^{\varepsilon_k}_{\alpha_k}S^{\varepsilon_{k+1}}_{\alpha_{k+1}}
\cdots S^{\varepsilon_{k+l}}_{\alpha_{k+l}}\quad
(1\le k, k+l\le m)
\]
is also a path, which will be called a subpath of (1). The path
\[
w^{-1}=S^{-\varepsilon_m}_{\alpha_m}\cdots
     S^{-\varepsilon_2}_{\alpha_2}S^{-\varepsilon_1}_{\alpha_1}
\]
is called the oppositely directed path equal to (1). If a path is
not simple, then it contains a closed subpath. One easily shows:
if there is a path connecting different points $p$ and $p'$ then
there is also a simple path connecting $p$ and $p'$. It
consists of the original with the closed subpaths deleted.

By a \emph{tour} we mean a path which traverses each segment
only once. If $w$ in (1) is a tour then
\[
S^{\varepsilon_i}_{\alpha_i}\ne S^{\varepsilon_k}_{\alpha_k}
\quad\text{and}\quad
S^{\varepsilon_i}_{\alpha_i}\ne S^{-\varepsilon_k}_{\alpha_k}
\quad
\text{when }i\ne k.
\]

A complex is called \emph{connected} if any two points are
connected by a path of finitely many segments. A complex
which is not connected consists of a finite or denumerable
number of disjoint, connected subcomplexes, which may be
called its \emph{components}.

Examples of complexes may be constructed so easily from,
e.g., euclidean segments that we shall not describe them further.

\section{Orders of points. Regular complexes}

A few questions that come up with finite complexes are the 
following. Let $a_0$ be the number of points of $\mathfrak{C}$,
$2a_1$ the number of directed segments, where equal but
oppositely directed segments are counted separately. If $k$ 
segments begin at a point, then $k$ is called the 
\emph{order\footnote{Today, $k$ is more commonly called the
\emph{degree} of the point, but in this book the word ``degree''
is used only when the order is constant. (Translator's note.)}
of the point}. Let $r_k$ be the number of points of order $k$.
Then
\[
\sum k r_k=2a_1,
\]
where the sum is taken over all points. Further, let
\[
\sum r_{2l+1}=r^{(1)}
\]
be the number of points of odd order. Then
\begin{align*}
\sum k r_k &= \sum 2l r_{2l} + \sum (2l+1) r_{2l+1}\\
                  &= 2\sum l r_{2l} + 2\sum l r_{2l+1} + r^{(1)},
\end{align*}
so $r^{(1)}=2r$ is always an even number.

If the same number, $k$, of segments begin at each point then
the complex is called \emph{regular of degree $k$}. This concept
is also meaningful for complexes with infinitely many elements.
If $k$ is odd and the number of points is finite, then $a_0$ must
be even because the number of segments $2a_1$ is in this case
equal to $a_0 k$. The regular complexes with $k=1$ are
decomposable into $\frac{1}{2} a_0$ components, each
consisting of a nonsingular segment and two points. The connected
regular complexes with $k=2$ consist either of one point, a
singular segment, and its inverse; or secondly of a finite number
of points $p_1,p_2,\ldots, p_{a_0}$ and $a_0$ segments
$s_1,s_2,\ldots, s_{a_0}$ together with their inverses. With
suitable orientation and numbering, $s_i$ begins at $p_i$ and
ends at $p_{i+1}$ (for $i=1,2,\ldots,a_0-1$) or $p_1$
(for $i=a_0$). The third case is infinitely many points
$p_0,p_1,p_{-1},p_2,p_{-2},\ldots$ and segments
$s_0,s_1, s_{-1},s_2,s_{-2},\ldots$; with suitable orientation 
and numbering, $s_i$ begins at $p_i$ and ends at $p_{i+1}$.
If $k>2$ then it is difficult to obtain a general view of the
regular complexes.

\section{The K\"onigsberg bridge problem}

One of the first topological questions to be asked was whether a
connected complex could be traversed in a single tour. The
so-called K\"onigsberg bridge problem is a problem requiring a
tour of a particular complex, and \textsc{Euler}\footnote{Petrop.
Comm. \textbf{8} (1741), 128.} showed that it is impossible by proving
the following general theorem: \emph{let $\mathfrak{C}$ be connected
and with only finitely many points and segments. If all points of the
complex are of even order then $\mathfrak{C}$ may be traversed by
a closed tour. If, on the other hand, $\mathfrak{C}$ contains $2r$ 
points of odd order, then there are $r$ paths}
\[
w_1,\quad w_2,\quad \ldots,\quad, w_r,
\]
\emph{and no fewer, which together traverse each segment of
$\mathfrak{C}$ exactly once.}

It is clear that there can be no fewer than $r$ paths of the kind required,
because only the initial and final point of the tour can be of odd order.
If $\mathfrak{C}$ can be traversed in $m$ tours than
\[
2r\le 2m.
\]
Now to show that exactly $r$ tours exist first look at the case where all
points are of even order. If $w$ is any tour of $\mathfrak{C}$ which
begins at $p_1$ and ends at $p_2$, then if $p_1\ne p_2$ there is 
an odd number of segments which begin at $p_2$. There is
consequently a segment emanating from $p_2$ which is not traversed
by $w$, and thus $w$ can be prolonged. If $p_1=p_2$ then $w$ is 
closed and one possibility is that there is a point $p$ on $w$ which
still bounds a segment not traversed by $w$. We then let $w'$ be the
closed path that results from a cyclic interchange of $w$ and begins at
$p$, and we can prolong $w'$ from $p$. Otherwise, $w$ traverses all
segments emanating from points through which it passes, so that the
points and segments of $w$ form a component of $\mathfrak{C}$,
which must be $\mathfrak{C}$ itself, since $\mathfrak{C}$ is
connected.

It follows that there must be a tour containing all segments of
$\mathfrak{C}$, and this tour must be closed, otherwise odd
numbers of segments would emanate from its initial and final points.

Now if $\mathfrak{C}$ is a complex with $r\ne 0$ we construct a
complex $\mathfrak{C}'$ with $r'=0$ by adding $r$ nonsingular
segments which connect the points of odd order with each other in pairs.
Then $\mathfrak{C}'$ may be traversed by a tour, and if we remove from 
it the segments not belonging to $\mathfrak{C}$ this tour falls into
exactly $r$ subtours.

\section{Trees}

\emph{A complex is called a tree if it is connected and no reduced
closed paths may be constructed from its segments.}

It follows that all segments of a tree are nonsingular and that all
reduced paths in a tree are open and simple. Two points of a tree
may be connected by only one reduced path. Namely, if $w$ and
$w'$ are two such paths, then $ww'^{-1}$ is a closed path. The
latter may be shortened by cancelling neighboring equal but
oppositely directed segments until finally all segments are
removed. Hence $w=w'$. It also follows from this fact that if
\[
S^{\varepsilon_1}_{\alpha_1}
S^{\varepsilon_2}_{\alpha_2}\cdots
S^{\varepsilon_i}_{\alpha_i}
S^{\varepsilon_{i+1}}_{\alpha_{i+1}}\cdots
S^{\varepsilon_m}_{\alpha_m}
\]
is a reduced path which connects the point $p_1$ to the point $p_2$,
and if the final point of $S^{\varepsilon_i}_{\alpha_i}$ is the point
$p_3$, then
\[
S^{\varepsilon_1}_{\alpha_1}
S^{\varepsilon_2}_{\alpha_2}\cdots
S^{\varepsilon_i}_{\alpha_i}
\]
is the reduced path that connects $p_1$ to $p_3$.

If a complex $\mathfrak{C}$ is connected, and if each point $p$ of
$\mathfrak{C}$ can be connected with any point $p'$ by only a single 
reduced path, then $\mathfrak{C}$ is a tree. If $p$ is an arbitrary point 
of the tree, and if $w_i$ are the simple paths from $p$ to all other
points $p_i$ of the tree, then any segment $s_k$ appears in all
$w_i$ with the same exponent $\varepsilon$. Supposing it 
appeared in $w_{i1}$ with $\varepsilon=+1$ and in $w_{i2}$ with
$\varepsilon=-1$, there would be a subpath $w'_{i1}$ of $w_{i1}$
ending with $s_k$, and a subpath $w'_{i2}$ of $w_{i2}$ ending
at the endpoint of $s_k$, but not with $s_k$ itself. Then 
$w'_{i1}w'^{-1}_{i2}$ would be a closed reduced path, which is a
contradiction. Thus one can label the segments in such a way that
all segments in the $w_i$ have $\varepsilon=1$ and the path $w_i$
to $p_i$ ends with $s_i$. Then $p_i$ is the final point of only one
segment, with which $w_i$ ends. The point $p$ is the initial point
of all segments which it bounds. Thus, if a finite tree
$\mathfrak{B}$ contains
$a_0$ points, then $\mathfrak{B}$ contains
\begin{equation}
a_1=a_0-1 \tag{1}
\end{equation}
pairs of oppositely directed equal segments.

After these preliminary remarks we demonstrate the important
theorem: \emph{if $\mathfrak{C}$ is a connected complex containing 
at least two distinct points $p_1$ and $p_2$, then there is a
subcomplex $\mathfrak{B}$ of $\mathfrak{C}$ which contains all
points of $\mathfrak{C}$ and is a tree.}\footnote{Such a subcomplex
$\mathfrak{B}$ of $\mathfrak{C}$ is today called a \emph{spanning
tree} of $\mathfrak{C}$. In future I will often use phrases
such as ``Let $\mathfrak{B}$ be a spanning tree of $\mathfrak{C}$''
in place of Reidemeister's phrases such as ``Let $\mathfrak{B}$ be
a tree which is a subcomplex of $\mathfrak{C}$ containing all points of 
$\mathfrak{C}$.'' (Translator's note.)}

We first suppose that $\mathfrak{C}$ contains only finitely many
points,
\[
p_1,\quad p_2,\quad \cdots,\quad p_{a_0},
\]
and prove the theorem under this hypothesis by complete induction.
In the base case there is a tree $\mathfrak{B}_2$ which contains
$p_1$ and $p_2$, because there is a simple line segment
connecting $p_1$ to $p_2$ and the complex consisting of the two
points and the pair of oppositely directed line segments between
them is a tree because it contains no simple closed path.

Now if $\mathfrak{B}_i$ is a tree which contains the points
\[
p_1,\quad p_2,\quad \cdots,\quad p_i,
\]
then either $p_{i+1}$ already occurs in this tree, in which case
$\mathfrak{B}_{i+1}=\mathfrak{B}_i$, or else
\[
S^{\varepsilon_1}_{\alpha_1}
S^{\varepsilon_2}_{\alpha_2}\cdots
S^{\varepsilon_k}_{\alpha_k}
S^{\varepsilon_{k+1}}_{\alpha_{k+1}}\cdots
S^{\varepsilon_m}_{\alpha_m}
\]
is a simple path which connects $p_1$ with $p_{i+1}$. There is then
a value $k\le m$ such that the subpath
\[
S^{\varepsilon_k}_{\alpha_k}
S^{\varepsilon_{k+1}}_{\alpha_{k+1}}\cdots
S^{\varepsilon_m}_{\alpha_m}
\]
has no point in common with $\mathfrak{B}_i$ apart from its
initial point. We now construct the complex $\mathfrak{B}_{i+1}$
consisting of the points and segments of $\mathfrak{B}_i$
together with the segments
\[
S^{\pm 1}_{\alpha_k},\quad
S^{\pm 1}_{\alpha_{k+1}},\cdots\quad
S^{\pm 1}_{\alpha_m}
\]
and their boundary points. Then $\mathfrak{B}_{i+1}$ is a tree,
and $\mathfrak{B}_{a_0}$ is the tree sought.

We now suppose that $\mathfrak{C}$ contains denumerably many
points $p_1,p_2,\ldots$ . By applying the process for finite
$\mathfrak{C}$ here we obtain an infinite sequence of subcomplexes
$\mathfrak{B}_i$ which are trees consisting of finitely many elements
and which contain the points $p_k$ for $k<i$. By $\mathfrak{B}$
we mean the complex containing all the points and segments
occurring in the $\mathfrak{B}_i$. These are all the points of
$\mathfrak{C}$. Furthermore, $\mathfrak{B}$ is a tree. Namely, if
there were a closed reduced path in $\mathfrak{B}$ then it would 
also be in the $\mathfrak{B}_k$ containing all of the (finitely many)
segments in this path, and this is a contradiction because
$\mathfrak{B}_k$ is a tree.

\section{The connectivity number}

If $\mathfrak{C}$ is a tree, then the complex $\mathfrak{B}$ in
$\mathfrak{C}$ constructed above is identical with $\mathfrak{C}$
itself. For if a segment, $s$ say, did not appear in $\mathfrak{B}$
then the complex $\mathfrak{C}$ would not be a tree. If $s$ is a
singular segment this is clear. If $s$ is not singular, let $w$ be the 
simple path in $\mathfrak{B}$ running from the final point of $s$ to the
initial point of $s$. Then $sw$ is a simple closed path and so
$\mathfrak{C}$ is not a tree. If $\mathfrak{C}$ is not a tree, then
there are either singular segments in $\mathfrak{C}$ or else different
subcomplexes which are trees containing all the points of $\mathfrak{C}$.

Let $\mathfrak{B}$ be such a tree and let $s$ be a nonsingular
segment of $\mathfrak{C}$ not occurring in $\mathfrak{B}$. We now
construct the complex $\mathfrak{C}'$ which results from
$\mathfrak{B}$ by the addition of $s$. Let $w$ be the path in
$\mathfrak{B}$ which leads from the final point $p_2$ to the initial
point $p_1$ of $s$. Then $sw$ is a simple path, and indeed up to
cyclic interchange and reversal of direction it is the only simple closed
path in $\mathfrak{C}'$. This is because, if $w'$ is a simple closed path, 
then it cannot run entirely in $\mathfrak{B}$ and hence must contain $s$
or $s^{-1}$; say, containing $s$ and beginning with $p_1$. With the
removal of $s$ we obtain from $w'$ a simple path $w^*$ in
$\mathfrak{B}$ running from $p_2$ to $p_1$, so $w^*=w$ and
$w'=sw$. Now if $s'$ is any segment of $\mathfrak{B}$ contained
in $w$, then
\[
w=w_1 s' w_2.
\]
Let $\mathfrak{B}'$ be the complex which results from $\mathfrak{C}$
by elimination of $s'$. $\mathfrak{B}'$ contains all the points of
$\mathfrak{C}$ and $\mathfrak{B}'$ is a tree different from
$\mathfrak{B}$. A tree $\mathfrak{B}'$ which results from
$\mathfrak{B}$ in this way is called a \emph{neighbor} of
$\mathfrak{B}$.

\emph{If $\mathfrak{B}$ and $\mathfrak{B}^*$ are two spanning trees 
of $\mathfrak{C}$, and if there are exactly $2k$ segments appearing in
$\mathfrak{B}^*$ but not in $\mathfrak{B}$, then 
$\mathfrak{B}=\mathfrak{B}^*$ if $k=0$, $\mathfrak{B}$ and
$\mathfrak{B}^*$ are neighbors if $k=1$, and there is a chain of
neighboring trees}
\[
\mathfrak{B},\quad
\mathfrak{B}^{(1)},\quad
\mathfrak{B}^{(2)},\quad\ldots,\quad
\mathfrak{B}^{(k)}
\]
\emph{with $\mathfrak{B}^{(k)}=\mathfrak{B}^*$ in case $k>1$.}

If $k=0$, then all segments of $\mathfrak{B}^*$ belong to
$\mathfrak{B}$. For if $\mathfrak{B}$ contained a segment which did
not appear in $\mathfrak{B}^*$, $\mathfrak{B}^*$ would not be a tree,
because $\mathfrak{B}^*$ contains all points of $\mathfrak{C}$ by
hypothesis and addition of a further segment to $\mathfrak{B}^*$
gives a complex which is not a tree.

If $k\ge 1$, let $s_1$ be a segment occuring in $\mathfrak{B}^*$
but not in $\mathfrak{B}$. The path $w$ in $\mathfrak{B}$ from the
final point of $s_1$ to the initial point then itself contains a segment
$s'_1$ which does not appear in $\mathfrak{B}^*$, otherwise
$\mathfrak{B}^*$ would contain the simple closed path $s_1 w$. We
now construct $\mathfrak{B}^{(1)}$ from $\mathfrak{B}$ by
adding $s_1$ and removing $s'_1$. Then $\mathfrak{B}^*$ contains
only $2(k-1)$ segments which do not appear in $\mathfrak{B}^{(1)}$.
Thus by iteration of the process we come to a tree $\mathfrak{B}^{(k)}$
which contains all the segments that $\mathfrak{B}^*$ does.

\emph{By the connectivity number of a connected complex 
$\mathfrak{C}$ we mean the number of segment pairs of $\mathfrak{C}$
which do not appear in a spanning tree $\mathfrak{B}$ 
of $\mathfrak{C}$.} The connectivity number of a tree is therefore 0.
If the connectivity number is $a>0$ we first have to show that the $a$
associated with complex is unique. We suppose that $\mathfrak{C}$
contains $2a$ more segments than the tree $\mathfrak{B}$ and $2a'$
more than the tree $\mathfrak{B}'$, where $\mathfrak{B}$ and
$\mathfrak{B}'$ each contain all the points of $\mathfrak{C}$ and
$a$ is finite. Now if $\mathfrak{B}'$ contains $2k$ segments which do
not belong to $\mathfrak{B}$, then $k\le a$. If $k=0$ then
$\mathfrak{B}$ is identical with $\mathfrak{B}'$, so $a=a'$. If $k=1$
then $\mathfrak{B}$ and $\mathfrak{B}'$ are neighbors and, as
one easily sees, $a=a'$. If $k>1$ there is a chain of neighboring trees
$\mathfrak{B},
\mathfrak{B}^{(1)},
\mathfrak{B}^{(2)},\ldots,
\mathfrak{B}^{(k)}=\mathfrak{B}'$ beginning with $\mathfrak{B}$
and ending with $\mathfrak{B}'$. Since $a^{(i)}=a^{(i+1)}$
always holds for the number $a^{(i)}$ of segments of $\mathfrak{C}$
not belonging to $\mathfrak{B}^{(i)}$, we have $a=a'$.

One can also give an easy indirect proof that: if $\mathfrak{C}$ 
contains infinitely
many more segments than one spanning tree $\mathfrak{B}$ 
of $\mathfrak{C}$, then $\mathfrak{C}$ contains infinitely many more
segments than any such tree.

For complexes containing only finitely elements, $a$ may be easily 
computed. \emph{If $a_0$ is the number of points, and $2a_1$ the
number of segments, then the connectivity number}
\[
a=-a_0+a_1+1.
\]
For a tree, this follows from (1) of Section 4.4. Another proof for trees
with finitely many segments is the following: it is correct for the tree
with a single pair of segments; here $a_0=2$. Supposing that the
theorem is true for $a_1=k$, we prove it for $a_1=k+1$. Let
$\mathfrak{B}$ be a tree with $a_0$ points and $2a_1=2(k+1)$
segments, with $p$ an arbitrary point and $w_i$ the unique simple
path from $p$ to $p_i$. All of these paths have finite length, hence 
there is among them one of greatest length, $w_m$, which leads to
$p_m$. But then $p_m$ can only bound the segment $s_m$ with
which $w_m$ ends, and $s_m$ can appear in no reduced path connecting 
$p$ to the other $p_i$. Then if we remove $p_m$ and $s_m$ from
$\mathfrak{B}$ we again obtain a connected complex, and in fact a
tree with $a_0-1$ points and $2(a_1-1)=2k$ segments, so
\[
0=(a_0-1)-(a_1-1)-1=a_0-a_1-1,
\]
as was to be proved.

If now $\mathfrak{C}$ is an arbitrary complex with $a_0$ points and
$2a_1$ segments, and if $\mathfrak{B}$ is a spanning tree 
of $\mathfrak{C}$, hence with $a_0$ points, then $\mathfrak{B}$ contains
exactly $2(a_0-1)$ segments, and if $a$ is the connectivity number of
$\mathfrak{C}$ then $\mathfrak{C}$ contains
\[
a_1=a+a_0-1
\]
segment pairs. This is the equation claimed.

\section{The fundamental group of a line segment complex}

Let $\mathfrak{C}$ be a connected line segment complex, and let
$p_0$ be one of its points. The closed paths in this complex
emanating from $p_0$ determine a group when we admit the empty
path consisting of $p_0$ alone and make the
following definition: two such paths are called equivalent when they
may be converted into each other by expansion and reduction. One
concludes easily, as we did in Section 2.2, that this relation is symmetrical
and transitive. Each path is equivalent to a reduced path when we also
admit the empty path consisting of $p_0$ alone. One concludes further,
as in Section 2.3, that a path is equivalent to a single reduced path, so
that only one reduced path appears in each class of equivalent paths.
The class of closed paths emanating from $p_0$ in which $w$ appears
will be denoted by $[w]$. By the product $[w_1][w_2]$ of two classes
we mean the class $[w_1 w_2]$, and one concludes as in Section 2.2
that the classes form a group under this multiplication, called the
\emph{fundamental group}\footnote{Reidemeister calls it the
\emph{Wegegruppe} (``path group,'' not unreasonably), but I have 
decided to use the term ``fundamental group,'' since it is now used 
universally. However, Reidemeister's notation is easier to understand 
knowing that``Weg'' is the German word for path. This explains the
notation $w$ for paths, and (from Section 6.1 onwards) the notation 
$\mathfrak{W}$ (fraktur W) for the fundamental group. (Translator's 
note.)} of $\mathfrak{C}$ with basepoint $p_0$. The class containing 
the empty path plays the role of the identity, and $[w^{-1}]$ is the class
inverse to $[w]$. The \emph{groupoid} (Section 1.15) of classes of
arbitrary equivalent paths may be defined similarly. The identities of the
latter correspond to the points of $\mathfrak{C}$.

If $\mathfrak{C}$ is a tree, then the fundamental group consists only
of the identity. If $\mathfrak{C}$ is not a tree then there is a reduced 
closed path, and hence also one which begins at $p_0$. The fundamental
group is therefore more than just the identity. In general it may be
shown that \emph{if $a$ is the connectivity number of $\mathfrak{C}$,
then the fundamental group is a free group with $a$ free generators.}
To show this, let $\mathfrak{B}$ be a spanning tree of $\mathfrak{C}$. 
Let
\[
s_1,\quad s_2,\quad \ldots,\quad s_k
\]
be the segments which, together with their inverses, do not appear in 
$\mathfrak{B}$. Let $p_{i,1}$ be the initial point, and $p_{i,2}$ the
final point, of the segment $s_i$. Also let $w_{i,1}$ and $w_{i,2}$
be simple paths in $\mathfrak{B}$ which lead from $p_0$ to
$p_{i,1}$ and $p_{i,2}$ respectively. We now construct the closed
path
\[
w_{i,1} s_i w^{-1}_{i,2}
\]
beginning at $p_0$, set
\[
[w_{i,1} s_i w^{-1}_{i,2}]=S_i,
\]
and assert that the $S_i$ are generators of the fundamental group.
To express an element $[w]$ in terms of the $S_i$, we take the
reduced path $w'$ contained in $[w]$. If $w'$ is not empty, then $w'$
runs through finitely many of the segments $s_i$ in succession, say
$m$ of them. If $s_{\alpha_1}$ is the first of these segments to appear
in $w'$, then $w'$ first leaves the tree $\mathfrak{B}$ at the point
$p_{\alpha_1,1}$ or $p_{\alpha_1,2}$, according as $s_{\alpha_1}$
is traversed in the positive or negative sense. If we assume the former,
then $w'$ begins with $w_{\alpha_1,1}s_{\alpha_1}$, say
\[
w'=w_{\alpha_1,1}s_{\alpha_1}w''.
\]
Then
\[
S^{-1}_{\alpha_1}[w']=[w_{\alpha_1,2}w'']
\]
contains a path which traverses only $m-1$ segments $s_i$ in succession,
and the same holds also for the reduced path in this class. It then
follows by induction that if
\[
s^{\varepsilon_1}_{\alpha_1},\quad
s^{\varepsilon_2}_{\alpha_2},\quad
\ldots,\quad
s^{\varepsilon_m}_{\alpha_m}
\]
are the segments $s^{\pm 1}_i$ which $w'$ traverses successively in
the positive direction, then in the class
\[
S^{-\varepsilon_m}_{\alpha_m}
S^{-\varepsilon_{m-1}}_{\alpha_{m-1}}\cdots
S^{-\varepsilon_1}_{\alpha_1}[w']
\]
there is a path which traverses none of the segments $s_i$. Thus the
path is entirely contained in $\mathfrak{B}$ and may be converted 
to the empty path by reduction.

We still have to see how the different generating systems of the
fundamental group depend on the different choices of the tree
$\mathfrak{B}$. Let $\mathfrak{B}$ and $\mathfrak{B}'$ be two
neighboring trees. Let $s'=s_1$ be the segment that is in
$\mathfrak{B}'$ but not in $\mathfrak{B}$, and $s$ the segment
that is in $\mathfrak{B}$ but not in $\mathfrak{B}'$. Let
\[
s_{\alpha_1},\quad s_{\alpha_2},\quad \cdots
\]
be those segments $s_i$ ($i=1,2,\ldots$) whose endpoints both
remain connected to $p_0$ after removal of $s$ from $\mathfrak{B}$,
let
\[
s_{\beta_1},\quad s_{\beta_2},\quad \cdots\quad
(s_{\gamma_1},\quad s_{\gamma_2},\quad \cdots)
\]
be those for which this holds for the initial point (respectively, final
point) but not for the final point (respectively, initial point), and let
\[
s_{\delta_1},\quad s_{\delta_2},\quad \cdots
\]
be those for which it holds for neither the initial or final point.
Then the simple paths from $p_0$ in $\mathfrak{B}'$ to
$p_{\alpha_i,1},p_{\alpha_i,2},p_{\beta_i,1},p_{\beta_i,2}$ are
identical with those in $\mathfrak{B}$; those from $p_0$ to
$p_{\beta_i,2},p_{\gamma_i,1},p_{\delta_i,1},p_{\delta_i,2}$,
on the other hand, always pass through $s'$, and indeed all in the positive 
direction (cf. Section 4.4) with a suitable orientation of $s'=s_1$.
Now let
\[
S'_i=[w'_{i1} s_i w'^{-1}_{i2}]\quad (i>1)
\]
be the generators corresponding to the paths $s_i$ relative to
$\mathfrak{B}'$, let $w'_1$ be the path in $\mathfrak{B}'$ from $p_0$
to the initial point of $s$, and $w'_2$ that to the final point of $s$, so
\[
S'_1=[w'_1 s w'^{-1}_2]
\]
is the missing generator relative to $\mathfrak{B}'$. Then one sees that
either $w'_1$ or $w'_2$, but not both, passes through the segment $s$.
It then follows that
\begin{align*}
S'_1=S^{\pm 1}_1,\quad S'_{\alpha_i}&=S_{\alpha_i},\quad
S'_{\beta_i}=S_{\beta_i}S^{-1}_1,\\
S'_{\gamma_i}=S_1 S_{\gamma_i},&\quad
S'_{\delta_i}=S_1 S_{\delta_i}S^{-1}_1.
\end{align*}
By Section 3.14, this mapping from the $S_i$ to the $S'_i$ is an
automorphism of the free group generated by the $S_i$.

\section{Coverings of complexes}

There is a relation between complexes which is very similar
to homomorphism of groups. We introduce it with the help of the
following definitions.

A complex $\mathfrak{C}$ covers a complex $\mathfrak{C}^*$ if 
each point $p$ and each segment $s$ of $\mathfrak{C}$ is
associated with a point $\boldsymbol{A}(p)=p^*$ and with a segment $\boldsymbol{A}(s)=s^*$,
respectively, of $\mathfrak{C}^*$ in the following way:
\begin{enumerate}
\item[A.1.]
Each $p^*$ corresponds to at least one $p$.
\item[A.2.]
If $p$ is the initial point of $s^{\varepsilon}_i$ then $\boldsymbol{A}(p)$ is the
initial point of $\boldsymbol{A}(s^{\varepsilon}_i)$.
\item[A.3.]
If $s^{\varepsilon_1}_{\alpha_1},s^{\varepsilon_2}_{\alpha_2},
\ldots,s^{\varepsilon_n}_{\alpha_n}$ are the segments with 
initial point $p$, then
\[
\boldsymbol{A}(s^{\varepsilon_1}_{\alpha_1}),\quad
\boldsymbol{A}(s^{\varepsilon_2}_{\alpha_2}),\quad
\ldots,\quad
\boldsymbol{A}(s^{\varepsilon_n}_{\alpha_n})
\]
are all different, and in fact they are the segments with initial point
$\boldsymbol{A}(p)$.
\item[A.4.]
$\boldsymbol{A}(s^{-1})=\boldsymbol{A}(s)^{-1}$.
\end{enumerate}
$\boldsymbol{A}$ is called a mapping of $\mathfrak{C}$ onto 
$\mathfrak{C}^*$,
or a covering of $\mathfrak{C}^*$ by $\mathfrak{C}$, or a
homomorphism\footnote{Reidemeister calls it an ``isomorphism,''
using the old terminology which called all homomorphisms
``isomorphisms,'' and distinguished our isomorphisms as
``one-to-one isomorphisms.'' Of course, the old terminology is 
confusing for modern readers, so I have used the words 
``homomorphism''  and ``isomorphism'' as we do today. 
(Translator's note.)}
from $\mathfrak{C}$ to $\mathfrak{C}^*$.
$\mathfrak{C}$ is called homomorphic to $\mathfrak{C}^*$ if
there is a homomorphism $\boldsymbol{A}$ with 
$\boldsymbol{A}(\mathfrak{C})=\mathfrak{C}^*$.

If each element of $\mathfrak{C}^*$ corresponds to only a single
element of $\mathfrak{C}$, then $\mathfrak{C}$ is said to be
\emph{isomorphic} to $\mathfrak{C}^*$. This relation is reflexive and
symmetric. If
\[
w=s^{\varepsilon_1}_{\alpha_1}s^{\varepsilon_2}_{\alpha_2}
\cdots s^{\varepsilon_n}_{\alpha_n}
\]
is a path, then the sequence of segments
\[
\boldsymbol{A}(s^{\varepsilon_1}_{\alpha_1})
\boldsymbol{A}(s^{\varepsilon_2}_{\alpha_2})
\cdots
\boldsymbol{A}(s^{\varepsilon_n}_{\alpha_n})
\]
is likewise a path, which may be denoted by $\boldsymbol{A}(w)$. If $w$ is closed,
then $\boldsymbol{A}(w)$ is also closed. Thus if $\boldsymbol{A}(w)$ is open, $w$ is also open.
If $\mathfrak{C}'$ is any subcomplex of $\mathfrak{C}$, then by
$\boldsymbol{A}(\mathfrak{C}')$ we mean the corresponding points and segments
of $\mathfrak{C}^*$ under the mapping $\boldsymbol{A}$. They likewise constitute
a complex $\mathfrak{C}'^*$.

It follows from A.2 and A.4 that if $p$ is the final point of
$s^{\varepsilon_i}_i$ then $\boldsymbol{A}(p)$ is the final point of
$\boldsymbol{A}(s^{\varepsilon_i}_i)$.

A.3 can also be expressed as follows. The segments which have $p$ as
initial point are mapped one-to-one onto the segments which have $\boldsymbol{A}(p)$
as initial point. It should be stressed that $s$ and $s^{-1}$ are different
segments. Thus it can happen, e.g., that $p$ bounds two nonsingular
segments and $\boldsymbol{A}(p)$ bounds one singular segment.

It follows from A.1 and A.3 that each segment $s^*$ corresponds to at 
least one segment $s$ for which $\boldsymbol{A}(s)=s^*$.

\emph{The homomorphism relation is transitive}, or more precisely:
if $\mathfrak{C}, \mathfrak{C}^*, \mathfrak{C}^{**}$ are three
complexes, the elements of which are denoted by
\[
p,s,\quad
p^*,s^*,\quad
p^{**}, s^{**}
\]
and if 
\[
\boldsymbol{A}_1(p)=p^*,\quad \boldsymbol{A}_1(s)=s^*
\]
is a covering of $\mathfrak{C}^*$ by $\mathfrak{C}$, and
\[
\boldsymbol{A}_2(p^*)=p^{**},\quad \boldsymbol{A}_2(s^*)=s^{**}
\]
is a covering of $\mathfrak{C}^{**}$ by $\mathfrak{C}^*$, and if we set
\[
\boldsymbol{A}_3(p)=
\boldsymbol{A}_2(\boldsymbol{A}_1(p))=p^{**},\quad
\boldsymbol{A}_3(s)=
\boldsymbol{A}_2(\boldsymbol{A}_1(s))=s^{**},
\]
then $\boldsymbol{A}_3$ is a covering of $\mathfrak{C}^{**}$ by $\mathfrak{C}$.

Namely, each $p^{**}$ in $\mathfrak{C}^{**}$ corresponds to some
$p$, because there is a $p^*$ with $\boldsymbol{A}_2(p^*)=p^{**}$ and in turn a
$p$ with $\boldsymbol{A}_1(p)=p^*$. If $p$ is the initial point of $s^\varepsilon$,
then $\boldsymbol{A}_3(p)$ is the initial point of 
$\boldsymbol{A}_3(s^\varepsilon)$, because
$\boldsymbol{A}_i(p)$ is the initial point of $\boldsymbol{A}_1(s^\varepsilon)$ and
$\boldsymbol{A}_2(\boldsymbol{A}_1(p))=\boldsymbol{A}_3(p)$ is the initial point of 
$\boldsymbol{A}_2(\boldsymbol{A}_1(s^\varepsilon))=
\boldsymbol{A}_3(s^\varepsilon)$. If $\boldsymbol{A}_3(p)=p^{**}$
then the segments beginning with $p$ are mapped one-to-one onto
those beginning with $p^{**}$, because $\boldsymbol{A}_1$ maps these segments
one-to-one onto those beginning with $\boldsymbol{A}_1(p)$, and $\boldsymbol{A}_2$ maps the
latter one-to-one onto those beginning with $p^{**}$. Finally,
\[
\boldsymbol{A}_3(s^{-1})=
\boldsymbol{A}_2(\boldsymbol{A}_1(s^{-1}))=
\boldsymbol{A}_2(\boldsymbol{A}_1(s)^{-1})=
(\boldsymbol{A}_3(s))^{-1}.
\]

If $\mathfrak{C}$ is connected, then so too is the covered complex
$\mathfrak{C}^*$, because if $p^*_1,p^*_2$ are any two points of
$\mathfrak{C}^*$ then there are two points $p_1,p_2$ with
\[
\boldsymbol{A}(p_i)=p^*_i \quad (i=1,2),
\]
and there is a path $w$ in $\mathfrak{C}$ from $p_1$ to $p_2$.
But then $\boldsymbol{A}(w)$ is a path from $p^*_1$ to $p^*_2$.

\section{Paths and coverings}

It also follows from A.3 that

\emph{If $p_1$ and $p'_1$ are two points for which}
\[
\boldsymbol{A}(p_1)=\boldsymbol{A}(p'_1)
\]
\emph{and if $w$ is a path which begins at $p_1$, then there is a 
well-defined path $w'$ which begins at $p'_1$ and for which}
\[
\boldsymbol{A}(w')=\boldsymbol{A}(w).
\]

Namely, if $s^\varepsilon$ is the initial segment of $w$ then there is 
exactly one segment $s'^{\varepsilon'}$ beginning at $p_1$ for which
$\boldsymbol{A}(s'^{\varepsilon'})=\boldsymbol{A}(s^\varepsilon)$, and if $p_2$ is the final point
of $s$, $p'_2$ that of $s'$, then $\boldsymbol{A}(p_2)=
\boldsymbol{A}(p'_2)$. The general
theorem follows by induction. If a path $w^*$ beginning at $p^*$
is reduced, $\boldsymbol{A}(p)=p^*$ and if $w$ is the path beginning at $p$ with
$\boldsymbol{A}(w)=w^*$, then by A.3 $w$ is also reduced. If $w^*$ is simple,
then $w$ is also simple, for $w$ is certainly reduced and each proper subpath
of $w$ is open, because each proper subpath of $w^*$ is open.

If $\mathfrak{B}^*$ is a tree contained in $\mathfrak{C}^*$, if $p^*$
is a point of $\mathfrak{B}^*$, if the $w^*_i$ are the simple paths in
$\mathfrak{B}^*$ from $p^*$ to the points $p^*_i$, and if also 
$\boldsymbol{A}(p)=p^*$ and the $w_i$ are the paths emanating from $p$ with
$\boldsymbol{A}(w_i)=w^*_i$, then the segments traversed by the paths $w_i$
form a complex $\mathfrak{B}$ which is likewise a tree. Namely,
$\boldsymbol{A}(\mathfrak{B})=\mathfrak{B}^*$ and this relation between the
points and segments of $\mathfrak{B}$ and $\mathfrak{B}^*$ is a
bijection. Firstly, all of $p,p_i$ are different, because all of 
$p^*,p^*_i$ are different, so $\boldsymbol{A}(p_i)=p^*_i$ is one-to-one with
respect to the points. Further, if $s^*$ is a segment of 
$\mathfrak{B}^*$ which begins at $p^*_a$, suppose $s_i$ 
($i=1,2$) are two segments of $\mathfrak{B}$ for which
\[
\boldsymbol{A}(s_1)=\boldsymbol{A}(s_2)=s^*.
\]
Then $s_1$ and $s_2$ must begin at $p_a$, and hence by A.3 we
certainly cannot have $\boldsymbol{A}(s_1)=
\boldsymbol{A}(s_2)$. Thus in fact $\mathfrak{B}$
and $\mathfrak{B}^*$ are isomorphic to each other and therefore
$\mathfrak{B}$ is a tree.

\emph{If $w^*$ is a simple closed path which begins and ends at
$p^*$, ${w^*}^k$ is the same path traversed $k$ times, and if}
\[
\boldsymbol{A}(w_k)={w^*}^k,
\]
\emph{then it may happen that the $w_k$ for $k<a$ are open paths but
$w_a$ is a closed path. In this event $w_a$ is a simple closed path.}

Namely, if
\[
w_a=w'_1 w'_2 w'_3 \cdots w'_a\quad\text{with}\quad
A(w'_i)=w^*\quad (i=1,2,\ldots,a)
\]
then the $w'_i$ are certainly simple paths. The initial points $p^{(i)}$ of
the $w'_i$, and only they, lie over the same point $p^*$ and by hypothesis
they are all different from the first, $p^{(1)}$. But they are also all
different from each other; for if $p^{(i)}=p^{(l)}$ for $i<l$ then
$p^{(1)}=p^{(l-i+1)}$, because the path over ${w^*}^{-i+1}$
emanating from $p^{(i)}$ and the path over ${w^*}^{-i+1}$
emanating from $p^{(l)}$ are subpaths of $w^{-1}_a$ which end at
$p^{(1)}$ and $p^{(l-i)}$ respectively, and they are otherwise 
identical because $p^{(i)}=p^{(l)}$. Supposing further that $w_a$
passes through the point $p'$ twice, and in fact once in $w'_i$ and the second 
time in $w'_l$, then the subpath of $w'_i$ from $p'$ to $p^{(i+1)}$
and the subpath of $w'_l$ from $p'$ to $p^{(l+1)}$ must lie over the 
same subpath ${w^*}'$ of $w^*$; for, since $w^*$ is simple, $w^*$ passes 
through the point $\boldsymbol{A}(p')$ only once, and ${w^*}'$ is determined
thereby. Thus $p^{(i+1)}=p^{(l+1)}$, contrary to what was 
previously shown. Likewise one concludes: if all paths $w_k$
emanating from $p$ and lying over ${w^*}^k$ are open, then they
are also simple.

\section{Simplicity of a covering}

With the help of the results on the trees $\mathfrak{B}$ in
$\mathfrak{C}$ for which $\boldsymbol{A}(\mathfrak{B})=\mathfrak{B}^*$ is also a
tree one can easily obtain a deeper understanding of the way 
$\mathfrak{C}$ covers a complex $\mathfrak{C}^*$. First we show:

\emph{If $\mathfrak{C}^*$ is a tree, $\mathfrak{C}$ is connected, and
$\boldsymbol{A}(\mathfrak{C})=\mathfrak{C}^*$ is a covering, then $\mathfrak{C}$
is isomorphic to $\mathfrak{C}^*$.}  If $\mathfrak{C}$ is not connected,
then $\mathfrak{C}$ separates into finitely or denumerably many trees
$\mathfrak{B}_1, \mathfrak{B}_2,\ldots$ which are isomorphic to
$\mathfrak{C}^*$.

For if $p$ is any point with $\boldsymbol{A}(p)=p^*$ then, by Section 4.8, there is a
tree $\mathfrak{B}_p$ which is a subcomplex of $\mathfrak{C}$, contains
$p$, and for which $\boldsymbol{A}(p)=\mathfrak{B}_p$. Now if $p'$ is any point of
$\mathfrak{B}_p$ and $s$ is any segment of $\mathfrak{C}$ which begins
at $p'$, then $s$ belongs to $\mathfrak{B}_p$. For $\boldsymbol{A}(s)=s^*$ and
$\boldsymbol{A}(p')={p^*}'$ belong to $\mathfrak{C}^*$, so in $\mathfrak{B}_p$
there is a segment $s'$ which begins at $p'$ and for which
$\boldsymbol{A}(s')=s^*$. But in $\mathfrak{C}$ there is only one segment which
begins at $p'$ and for which $\boldsymbol{A}(s')=s^*$. Thus $s=s'$. $\mathfrak{B}_p$
is therefore identical with $\mathfrak{C}$ in the case where $\mathfrak{C}$
is connected, and otherwise it is a component of $\mathfrak{C}$.

Now to the general case, where $\mathfrak{C}^*$ is connected but not a
tree. In this case let $\mathfrak{B}^*$ be a subcomplex of 
$\mathfrak{C}^*$ which is a tree and contains all the points of
$\mathfrak{C}^*$. If $\boldsymbol{A}(p)=p^*$ then by $\mathfrak{B}_p$ we 
mean the tree in $\mathfrak{C}$ containing $p$ and for which
$\boldsymbol{A}(\mathfrak{B}_p)=\mathfrak{B}^*$. We show:

\emph{If $p$ and $p'$ are two different points of $\mathfrak{C}$ for
which $\boldsymbol{A}(p)=\boldsymbol{A}(p')$, then the trees associated with them, 
$\mathfrak{B}_p$ and $\mathfrak{B}_{p'}$, are disjoint.}

Namely, suppose that $\mathfrak{B}_p$ and $\mathfrak{B}_{p'}$ had
the point $p''$ in common. Then there would be a simple path $w$,
lying wholly in $\mathfrak{B}_p$, running from $p$ to $p''$, and a
simple path $w'$, lying wholly in $\mathfrak{B}_{p'}$, running from $p'$
to $p''$. The paths $w$ and $w'$ are uniquely determined and
$\boldsymbol{A}(w)=
\boldsymbol{A}(w')$, because $\boldsymbol{A}(w)$ and $\boldsymbol{A}(w')$ 
run from $\boldsymbol{A}(p)=\boldsymbol{A}(p')$ to
$\boldsymbol{A}(p'')$. Consequently, $\boldsymbol{A}(w^{-1})=
\boldsymbol{A}(w'^{-1})$ too, and since these 
two paths begin at the same point, $w^{-1}=w'^{-1}$ and $p=p'$,
contrary to hypothesis. It follows further that:

\emph{If $p$ is any point of $\mathfrak{C}$ and if 
$p^{(1)}, p^{(2)},\ldots $ are all the points of $\mathfrak{C}$ for which
$\boldsymbol{A}(p^{(i)})=\boldsymbol{A}(p)$, and}
\[
\mathfrak{B}_{p^{(1)}},\quad
\mathfrak{B}_{p^{(2)}},\quad \ldots
\]
\emph{are the trees associated with the points $p^{(i)}$, then an
arbitrary point of $\mathfrak{C}$ appears in exactly one of the trees
$\mathfrak{B}_{p^{(i)}}$.}

We need only show that each point $p'$ of $\mathfrak{C}$ appears in
a $\mathfrak{B}_{p^{(i)}}$. Let $\boldsymbol{A}(p')={p^*}'$ and let ${w^*}'$ be
the simple path in $\mathfrak{B}^*$ from $p^*$ to ${p^*}'$, and $w$
the uniquely determined path ending at $p'$ for which $\boldsymbol{A}(w)={w^*}'$.
The latter  begins at a point $p^{(i)}$ over $p$ and runs wholly in 
$\mathfrak{B}_{p^{(i)}}$, thus $p'$ belongs to $\mathfrak{B}_{p^{(i)}}$.

Thus if a point $p^*$ lies under $k$ different, or denumerably many,
different points of $\mathfrak{C}$, then every point ${p^*}'$ lies under
$k$ or denumerably many different points of $\mathfrak{C}$. Likewise,
each directed segment of $\mathfrak{C}^*$ lies under $k$ or
denumerably many segments of $\mathfrak{C}$. Correspondingly, the
homomorphism $\boldsymbol{A}$ may be called \emph{$k$-to-1} or
\emph{infinite-to-one}.

\section{Coverings and permutations}

The segments of $\mathfrak{C}^*$ are divided into two classes by the
choice of tree $\mathfrak{B}^*$: the $s^*$ that belong to the tree
and the $\overline{s}^*$ that do not. Correspondingly, the
segments of $\mathfrak{C}$ are also divided into two classes: the $s$
that belong to a $\mathfrak{B}_{p^{(i)}}$, and the
$\overline{s}$ that belong to no $\mathfrak{B}_{p^{(i)}}$.

Let $\overline{s}^*$ be any segment not belonging to
$\mathfrak{B}^*$ and let
\[
\overline{s}^{(1)},\quad
\overline{s}^{(2)},\quad\ldots
\]
be the segments of $\mathfrak{C}$ for which 
$\boldsymbol{A}(\overline{s}^{(i)})=s^*$; let $p^*_1$ and $p^*_2$ be the
initial and final points of $\overline{s}^*$, and let $p^{(i)}_1$ and
$p^{(i)}_2$ be the initial and final points of $\overline{s}^{(i)}$.
A $\mathfrak{B}_{p^{(i)}}$ then contains exactly one of the final
points $p^{(n_i)}_2$ of these segments. The numbering may be arranged 
so that $\overline{s}^{(i)}$ begins in $\mathfrak{B}_{p^{(i)}}$ and
thus ends in $\mathfrak{B}_{p^{(n_i)}}$, and $p^{(i)}_1$ and
$p^{(i)}_2$ lie in $\mathfrak{B}_{p^{(i)}}$. Then the correspondence
\[
\pi=\left(\begin{array}{ccc}
                     1 & 2 & \cdots\\
                     n_1 & n_2 & \cdots
             \end{array}\right)
\]
is a permutation. One such permutation corresponds to each segment
$\overline{s}^*$. The initial and final points $p^{(i)}_1$ 
and $p^{(n_i)}_2$ of all $\overline{s}^{(i)}$ are determined by $\pi$
and $p^*_1$ and $p^*_2$.

Now we can immediately give the totality of $k$-to-one coverings
of a complex $\mathfrak{C}^*$. We construct a tree $\mathfrak{B}^*$
and $k$ trees isomorphic to $\mathfrak{B}^*$,
\[
\mathfrak{B}_{p^{(1)}},\quad
\mathfrak{B}_{p^{(2)}},\quad\ldots,\quad
\mathfrak{B}_{p^{(k)}}.
\]
Each segment $\overline{s}^*$ of $\mathfrak{C}^*$ which does not
belong to $\mathfrak{B}^*$ is associated with $k$ segments
$\overline{s}^{(1)},\overline{s}^{(2)},\ldots,\overline{s}^{(k)}$
and an arbitrary permutation of the numbers $1,2,\ldots,k$
\[
\pi=\left(\begin{array}{cccc}
                     1 & 2 & \cdots & k\\
                     n_1 & n_2 & \cdots & n_k
             \end{array}\right).
\]
If $p^{(i)}_1,p^{(i)}_2$ are the points in $\mathfrak{B}_{p^{(i)}}$ over
$p^*_1,p^*_2$, the initial and final points of $\overline{s}^*$, then
$\overline{s}^{(i)}$ begins with $p^{(i)}_1$ and ends with
$p^{(n_i)}_2$.

\section{Fundamental domains}

If we add to $\mathfrak{B}_{p^{(i)}}$ all segments $\overline{s}$
beginning at a point of $\mathfrak{B}_{p^{(i)}}$, without taking
their final points, and call the domain constructed in this way
$\mathfrak{F}_i$, then, for each point and segment of $\mathfrak{C}^*$,
$\mathfrak{F}_i$ contains exactly one element for which
\[
\boldsymbol{A}(p)=p^*,\quad 
\boldsymbol{A}(s)=s^*,\quad 
\boldsymbol{A}(\overline{s})=\overline{s}^*.
\]
For this reason, $\mathfrak{F}_i$ is called a \emph{fundamental domain}
of the covering $\boldsymbol{A}$. The $\mathfrak{F}_i$ are not complexes as long as
there are segments connecting different $\mathfrak{B}_{p^{(i)}}$,
because they do not contain final points of the latter or their equal but
oppositely directed segments. For each segment $\overline{s}$ of
$\mathfrak{F}_i$ there is a second, $\overline{s}'$, for which
$\boldsymbol{A}(\overline{s}')=\boldsymbol{A}(\overline{s})^{-1}$. 
We can now omit one of each 
such pair, say the $\overline{s}'$, and take $\overline{s}^{-1}$ in place
of it. The domain $\mathfrak{F}'_i$ resulting from $\mathfrak{F}_i$ in
this way is again a fundamental domain.

If $\mathfrak{C}^*$ contains only a single point and $2r$ singular 
segments which begin and end at $p^*$, say ${s^*_1}^{\pm 1},
{s^*_2}^{\pm 1},\ldots, {s^*_r}^{\pm 1}$, then the trees 
$\mathfrak{B}_{p^{(i)}}$ also consist of single points $p^{(i)}$
lying over $p^*$, and $2r$ segments emanating from each point
$p^{(i)}$ of $\mathfrak{C}$. $\mathfrak{C}$ is thus a regular complex
of degree $2r$.

If one collects all segments of $\mathfrak{C}$ which lie over
$s^*_i,{s^*_i}^{-1}$, together with their boundary points, into a
complex $\mathfrak{U}_i$, then the $\mathfrak{U}_i$ are again
regular complexes of degree 2, which together contain all the points
of $\mathfrak{C}$. Each segment of $\mathfrak{C}$ appears in
exactly one $\mathfrak{U}_i$. We express this state of affairs as
follows: $\mathfrak{C}$ may be \emph{decomposed} into regular
complexes $\mathfrak{U}_i$ of degree 2. One sees immediately that
the converse also holds:

\emph{If $\mathfrak{C}$ is any regular complex of degree $2r$ and if
$\mathfrak{U}_i$ ($i=1,2,\ldots,r$) are subcomplexes of $\mathfrak{C}$
which constitute a decomposition in the above sense, then there is a
mapping of $\mathfrak{C}$ onto a complex $\mathfrak{C}^*$ with a single
point and $r$ singular segments.}

In order to define the mapping $\boldsymbol{A}(s)=s^*$, let $s_k$ be any segment of
$\mathfrak{C}$ which appears in $\mathfrak{U}_i$. $\mathfrak{U}_i$
need not be connected; let $\mathfrak{U}_{i,k}$ be the component of
$\mathfrak{U}_i$ in which $s_k$ appears. We can traverse all segments of
$\mathfrak{U}_{i,k}$ in a closed path $w$ and orientate the segments
of $\mathfrak{U}_{i,k}$ in such a way that those occurring in $w$ all
have positive exponent. We then set $\boldsymbol{A}(s)=s^*_i$ for all segments
$s$ of $\mathfrak{U}_{i,k}$, and similarly for all $\mathfrak{U}_{i,k}$.

\section{Regular complexes of even order}

The following theorem\footnote{\textsc{Julius Petersen}, Acta Math.
\textbf{15}, (1891), 193--220.} is interesting in this connection:
\emph{every regular complex of even degree $2r$ may be decomposed
into $r$ regular complexes $\mathfrak{U}_i$ of degree 2.}
Consequently, every finite regular complex   may be
regarded as covering a complex $\mathfrak{C}^*$ with a single point
and $r$ singular segments.

For $2r=2$ there is nothing to prove. For $2r=4$ we argue as follows:
let $\mathfrak{C}'$ be a component; then $\mathfrak{C}'$ is also a
regular complex of degree 4 and hence it may be traversed by a tour
$w$. If $a_0$ is the number of points of $\mathfrak{C}'$ then $w$
contains $2a_0$ segments. With suitable orientation and numbering of
the segments let
\[
w=s_1 s_2 \cdots s_{2a_0}.
\]
We now construct a complex $\mathfrak{U}_1$ from all the segments
$s^{\pm 1}_{2l+1}$ and their boundary points, and a 
complex $\mathfrak{U}_2$ from the segments $s^{\pm 1}_{2l}$ and their
boundary points. Each segment $s$ of $\mathfrak{C}'$ appears in
either $\mathfrak{U}_1$ or $\mathfrak{U}_2$. Since $w$ traverses
each point $p$ of $\mathfrak{C}'$, $p$ is itself the initial point of a
segment $s^{\pm 1}_{2l+1}$ and the final point of a segment
$s^{\pm 1}_{2l}$. Thus the complexes $\mathfrak{U}_1$ and
$\mathfrak{U}_2$ each contain all the points of $\mathfrak{C}'$, and
since $w$ passes each point $p$ exactly twice, each point $p$ appears
exactly twice in the series of initial and final points of the segments
\[
s_1,\quad s_3,\quad\ldots,\quad s_{2a_0-1}.
\]
That is, beginning at each point $p$ there are exactly two segments from
$\mathfrak{U}_1$ and two from $\mathfrak{U}_2$.

In order to settle the general case, we must still prove a few lemmas.
First some definitions!

If $\mathfrak{C}$ and $\mathfrak{C}'$ are two regular complexes with
$a_0$ points and degree $2r$, if there is a subcomplex $\mathfrak{C}_n$
of $\mathfrak{C}$ which contains $2l$ segments and is homomorphic to
a subcomplex $\mathfrak{C}'_n$ of $\mathfrak{C}'$, and if there is no
subcomplex of $\mathfrak{C}$ having more than $2l$ segments which
is homomorphic to a subcomplex $\mathfrak{C}'$, then $2a_0 r-2l$ is
called the \emph{distance} between $\mathfrak{C}$ and $\mathfrak{C}'$.
Thus $l$ is zero only when $\mathfrak{C}$ contains only regular
segments and $\mathfrak{C}'$ only singular segments.

For isomorphic complexes the distance is zero and, conversely, if the
distance is zero then the complexes are isomorphic.

If $s_1$ and $s_2$ are two different and not just just oppositely
directed segments of $\mathfrak{C}$, if $p_{i1}$ are the boundary points
of $s_1$ and $p_{i2}$ are the boundary points of $s_2$, and if the
complex $\mathfrak{C}'$ results from $\mathfrak{C}$ by replacing the
segments $s_1$ and $s_2$ by two other segments $s'_1$ and $s'_2$,
where $s'_1$ is bounded by $p_{11},p_{12}$ and $s'_2$ by
$p_{21},p_{22}$, then $\mathfrak{C}$ and $\mathfrak{C}'$ are
called \emph{neighboring}. The distance between neighboring
complexes is at most 4.

\section{Modifications of regular complexes}

We now assert the theorem:

\emph{If $\mathfrak{C}$ and $\mathfrak{C}'$ are any two regular complexes
of $a_0$ points and degree $2r$ then there are two chains of complexes,
also of this kind,}
\[
\mathfrak{C}_1,\mathfrak{C}_2,\ldots,\mathfrak{C}_k;\quad
\mathfrak{C}'_1,\mathfrak{C}'_2,\ldots,\mathfrak{C}'_l,
\]
\emph{such that $\mathfrak{C}=\mathfrak{C}_1$, 
$\mathfrak{C}'=\mathfrak{C}'_1$, $\mathfrak{C}'_l$ is isomorphic to
$\mathfrak{C}_k$, and $\mathfrak{C}_i,\mathfrak{C}_{i+1}$, likewise
$\mathfrak{C}'_i,\mathfrak{C}'_{i+1}$, are neighbors.}

To prove this we first consider two complexes $\mathfrak{C}$ and
$\mathfrak{C}'$ which possess no homomorphic subcomplexes.
$\mathfrak{C}'$ then consists purely of singular segments, so we can
replace $\mathfrak{C}'$ by a neighboring complex $\mathfrak{C}''$ with
two regular segments, whose distance from $\mathfrak{C}$ is therefore
less.

Let the distance between $\mathfrak{C}$ and $\mathfrak{C}'$ be
$2a_0 r-2l>0$, and let $\mathfrak{C}_n$ and $\mathfrak{C}'_n$
respectively be the largest isomorphic subcomplexes of $\mathfrak{C}$
and $\mathfrak{C}'$.

If all points of $\mathfrak{C}_n$ have order $2r$ then $\mathfrak{C}_n$
consists of certain components of $\mathfrak{C}$, and likewise
$\mathfrak{C}'_n$ consists of certain components of $\mathfrak{C}'$.
Let $\mathfrak{C}=\mathfrak{C}_n+\overline{\mathfrak{C}}_n$ and
$\mathfrak{C}'=\mathfrak{C}'_n+\overline{\mathfrak{C}}'_n$.
Then we must have $\overline{\mathfrak{C}}_n$ consisting purely of
regular segments and $\overline{\mathfrak{C}}'_n$ consisting purely of
singular segments, otherwise $\mathfrak{C}_n$ and $\mathfrak{C}'_n$
would not be the largest homomorphic subcomplexes. Then if one
replaces  $\overline{\mathfrak{C}}'_n$ by a neighboring complex
$\overline{\mathfrak{C}}''_n$ containing two regular segments, and sets
$\mathfrak{C}''=\mathfrak{C}'_n+\overline{\mathfrak{C}}''_n$, then the
distance between $\mathfrak{C}''$ and $\mathfrak{C}$ is less than that
between $\mathfrak{C}'$ and $\mathfrak{C}$.

Now let $p_1$ be any point of $\mathfrak{C}_n$ which has order less than
$2r$ in $\mathfrak{C}_n$. Further, let $s$ be a segment emanating from
$p_1$ which belongs to $\mathfrak{C}$ but not to $\mathfrak{C}_n$,
and which ends at $p_2$. Let $\boldsymbol{A}(p_1)=p'_1$ be the point of
$\mathfrak{C}'_n$ that corresponds to $p_1$ under the homomorphism
$\boldsymbol{A}$ from $\mathfrak{C}_n$ to $\mathfrak{C}'_n$. Then $p_1$ 
has the same order in $\mathfrak{C}'_n$ as $p_1$ has in $\mathfrak{C}_n$, and
thus there is a segment $s'$ in $\mathfrak{C}'$ which emanates from $p'_1$
and does not belong to $\mathfrak{C}'_n$.

1. Now if $p_2$ belongs to $\mathfrak{C}_n$ and if $p'_2=
\boldsymbol{A}(p_2)$ then
the segments emanating from $p'_1$ that do not belong to
$\mathfrak{C}'_n$ certainly cannot end at $p'_2$, otherwise the
subcomplex of $2l+2$ segments consisting of the elements of
$\mathfrak{C}_n$ together with $s$ would be homomorphic to the
subcomplex consisting of the elements of $\mathfrak{C}'_n$ together
with a segment $s'$. Thus $s'$ ends at $p'_3\ne p'_2$. Likewise, $p_2$
has an order less than $2r$ in $\mathfrak{C}_n$, hence $p'_2$ does
similarly, and so there is a segment $s''$ emanating from $p'_2$ which
does not belong to $\mathfrak{C}'_n$ and ends at $p'_4\ne p'_1$.

We now construct a new complex $\overline{\mathfrak{C}}'$ 
from $\mathfrak{C}'$ by introducing a segment $\overline{s}'$
between $p'_1$ and $p'_2$ and a segment $\overline{s}''$ between
$p'_3$ and $p'_4$. $\overline{\mathfrak{C}}'$ is again regular, and
a neighbor of $\mathfrak{C}'$. The distance between 
$\overline{\mathfrak{C}}'$ and $\mathfrak{C}$ is smaller by at least 4.

2. If $p_2$ does not belong to $\mathfrak{C}_n$, then the segments
emanating from $p'_1$ must all end in $\mathfrak{C}'$, otherwise we
again could give a subcomplex  of $\mathfrak{C}$, of $2l+4$ segments,
homomorphic to a subcomplex of $\mathfrak{C}'$. $\mathfrak{C}'_n$
then satisfies the hypothesis that we made about $\mathfrak{C}_n$ in 1,
and we construct a complex $\overline{\mathfrak{C}}$
neighboring $\mathfrak{C}$ by the same process, where
$\overline{\mathfrak{C}}$ has a smaller distance from $\mathfrak{C}'$ than
$\mathfrak{C}$. The assertion follows by iterated application of this
process.

\section{Invariance of the decomposition}

The following theorem holds:

\emph{If $\mathfrak{C}$ and $\mathfrak{C}'$ are neighboring regular
complexes of degree $2r>4$ and if $\mathfrak{C}$ may be decomposed
into $r$ regular complexes}
\[
\mathfrak{U}_1,\quad
\mathfrak{U}_2,\quad\ldots,\quad
\mathfrak{U}_r
\]
\emph{of degree 2, then $\mathfrak{C}'$ may also be decomposed into
$r$ such complexes}
\[
\mathfrak{U}'_1,\quad
\mathfrak{U}'_2,\quad\ldots,\quad
\mathfrak{U}'_r.
\]

When $\mathfrak{C}$ and $\mathfrak{C}'$ are isomorphic this is clear.
In any case, $\mathfrak{C}$ and $\mathfrak{C}'$ contain a common
subcomplex $\overline{\mathfrak{C}}$ of $2ra_0-4$ segments. Let
$s_1,s_2; s^{-1}_1,s^{-1}_2$ be the segments of $\mathfrak{C}$
that do not belong to $\overline{\mathfrak{C}}$, and
$s'_1,s'_2; {s'}^{-1}_1,{s'}^{-1}_2$ the segments of $\mathfrak{C}'$
that do not belong to $\overline{\mathfrak{C}}$.
$\overline{\mathfrak{C}}$ contains all points of $\mathfrak{C}$,
because $2r>4$. Let $p_{k,i}$ ($k=1,2$) be the boundary points of
$s_i$ ($i=1,2$). Then, with suitable numbering, $p_{1,1}$ and
$p_{1,2}$ bound the segment $s'_1$ and hence $p_{2,1}$ and
$p_{2,2}$ bound the segment $s'_2$.

We now distinguish two cases:

1. The segments $s_1$ and $s_2$ may belong to the same
subcomplex $\mathfrak{U}_1$ of $\mathfrak{C}$. Then all segments
and points of the $\mathfrak{U}_i$ ($i\ge 2$) belong to
$\overline{\mathfrak{C}}$. We now define
\[
\mathfrak{U}'_i=\mathfrak{U}_i\quad (i=2,\ldots,r)
\]
and let $\mathfrak{U}'_1$ be the complex that results from 
$\mathfrak{C}'$ by leaving out the $\mathfrak{U}'_i$. Then
$\mathfrak{U}'_1$ is a regular complex of degree 2 which contains all
points of $\mathfrak{C}'$, because the complex of all the elements of
the $\mathfrak{U}'_i$ ($i\ge 2$) is a regular complex of degree $2r-2$
containing all points of $\mathfrak{C}'$.

2. The segments $s_1$ and $s_2$ may belong to different subcomplexes
$\mathfrak{U}_1$ and $\mathfrak{U}_2$. The all $\mathfrak{U}_i$
with $i\ge 3$ belong to $\overline{\mathfrak{C}}$ and we set
$\mathfrak{U}'_i=\mathfrak{U}_i$ ($i\ge 3$). Now the segments of
$\mathfrak{C}'$ that belong to none of the $\mathfrak{U}_i$
($i\ge 3$) constitute a regular complex of degree 4, and by Section 4.12
this may be decomposed into two regular complexes $\mathfrak{U}'_1$
and $\mathfrak{U}'_2$ of degree 2. Then the
\[
\mathfrak{U}'_i\quad (i=1,2,\ldots,r)
\]
constitute a decomposition of $\mathfrak{C}'$ into $r$ regular complexes
of degree 2.

Since there are certainly regular complexes of degree $2r$ with $a_0$
points which may be decomposed into $r$ complexes of degree 2, it
follows with the help of Section 4.13 that all regular complexes of
degree $2r$ may be decomposed into $r$ complexes of degree 2.

\section{Regular complexes of degree three}

We now single out a special class of $k$-fold coverings of
$\mathfrak{C}^*_2$, the complex of a single point with two singular
segments: the permutation associated with the segment $s^*_2$ must
leave no element fixed and must yield the identity when applied twice in
succession. Then, if $s_1$ and $s_2$ are the two segments that
emanate from $p^{(i)}$ and lie over $s^*_2$ and ${s^*_2}^{-1}$
respectively, $s_1$ and $s_2$ must both end at the same point
$p^{(l)}$. The $p^{(i)}$ may therefore be grouped in pairs, which
necessarily have degree $k$. The numbering of points may be
arranged so that the segments over $s^*_2$ lead from $p^{(i)}$
to $p^{(i+1)}$. We now construct a new complex $\mathfrak{C}_{2k}$
in which we replace each pair of segments over ${s^*_2}^{\pm 1}$
with the same initial point $p^{(i)}$ and final point $p^{(i+1)}$ by a
single segment $s'_i$. $\mathfrak{C}_{2k}$ contains $k$ points and
$3k$ segments. Three segments emanate from each point, so
$\mathfrak{C}_{2k}$ is a regular complex of degree 3. If one collects
all segments lying over $s^*_1, {s^*_1}^{-1}$ into a complex
$\mathfrak{U}_1$, and all the remaining segments of $\mathfrak{C}_{2k}$
into a complex $\mathfrak{U}_2$, then one sees that $\mathfrak{U}_1$
is a regular complex of degree 2, $\mathfrak{U}_2$ is a regular
complex of degree 1, and the two complexes constitute a decomposition
of $\mathfrak{C}_{2k}$. The corresponding covering complex
$\mathfrak{C}^*_2$ may be easily recovered from $\mathfrak{C}_{2k}$ 
and the given decomposition.

It is now natural to ask whether all regular complexes of degree 3 may be
decomposed into a regular complex of degree two and one of degree one. 
This is not the case, as is shown by the example of the complex
with points $p_i$ ($i=0,1,2,3$) and the regular segments $s_i$ with the
boundary points $p_0,p_i$ as well the singular segments $s'_i$ with the
boundary points $p_i$ ($i=1,2,3$).

In general it may be proved that an indecomposable regular complex of
degree 3 which contains no singular segment must have at least three
``leaves.'' A leaf is a subcomplex connected to the remaining points of
$\mathfrak{C}$ by just a single segment. Thus the example given has
three leaves.\footnote{Cf. the work cited on p.~96.}

Apart from this, little is known about the decomposition of regular
complexes of odd order into subcomplexes. It may be pointed out that
there is a regular complex of degree three, containing no singular
segment, which may be decomposed into a complex of degree two
and one of degree one, but not into three complexes of  
degree one.\footnote{\textsc{J. Petersen}, L'intermed \textbf{5}, (1898), 225.}

\section{Coverings and permutation groups}

We go further into the connection between the permutations
$\Pi$ and the coverings of a complex containing only one point. The
$\Pi$ generate a permutation group $\mathfrak{P}$ consisting of
all the permutations representable as power products of the $\Pi$.
Naturally, the structure of this group is closely connected with the 
structure of the covering. 

Let $\mathfrak{C}$ be a covering of the complex $\mathfrak{C}^*_r$
with one point $p$ and $r$ singular segments $s^*_1,s^*_2,\ldots,
s^*_r$ and their inverses ${s^*_i}^{-1}$. Let $S_1,S_2,\ldots,S_r$
be free generators of a free group $\mathfrak{S}$. Now if
\[
w=
s^{\varepsilon_1}_{\alpha_1}
s^{\varepsilon_2}_{\alpha_2}\cdots
s^{\varepsilon_m}_{\alpha_m}
\]
is any path in $\mathfrak{C}$ and 
$\boldsymbol{A}(s^{\varepsilon_i}_{\alpha_i})={s^*_{\beta_i}}^{\eta_i}$, then
\[
\boldsymbol{A}(w)=
{s^*_{\beta_1}}^{\eta_1}
{s^*_{\beta_2}}^{\eta_2}\cdots
{s^*_{\beta_m}}^{\eta_m},
\]
so let
\[
W=
S^{\eta_1}_{\beta_1}
S^{\eta_2}_{\beta_2}\cdots
S^{\eta_m}_{\beta_m}
\]
be the power product from $\mathfrak{S}$ associated with $w$.
Each power product of $\mathfrak{S}$ corresponds to a well-defined
path $w$ when an initial point is given in $\mathfrak{C}$. By a
relation in the $S$ we mean a power product $R(S)$
such that all the paths in $\mathfrak{C}$ corresponding to the
$R(S)$ are closed. The collection $\mathfrak{R}$ of relations
constitutes an invariant subgroup of $\mathfrak{S}$. For $R^{-1}$
is also a relation along with $R$, and $R_1 R_2$ along with $R_1,R_2$;
thus $\mathfrak{R}$ is a subgroup. And since
\[
S^{\varepsilon_i}_i R S^{-\varepsilon_i}_i\quad (i=1,2,\ldots,r)
\]
is also a relation along with $R$, $\mathfrak{R}$ is an invariant
subgroup of $\mathfrak{S}$. We now assert that the factor group
$\mathfrak{F}=\mathfrak{S}/\mathfrak{R}$ is isomorphic to the
permutation group $\mathfrak{P}$, when we understand $\Pi \Pi'$
to be the permutation resulting from first performing $\Pi$, then
$\Pi'$. If $\Pi_i$ ($i=1,2,\ldots,r$) are the permutations
associated with the segments $s^*_i$, then we claim more
precisely that
\[
\boldsymbol{I}(
S^{\varepsilon_1}_{\alpha_1}
S^{\varepsilon_2}_{\alpha_2}\cdots
S^{\varepsilon_m}_{\alpha_m}
)
= 
\Pi^{\varepsilon_1}_{\alpha_1}
\Pi^{\varepsilon_2}_{\alpha_2}\cdots
\Pi^{\varepsilon_m}_{\alpha_m}
\]
is an isomorphism between $\mathfrak{P}$ and $\mathfrak{F}$. The
mapping is certainly a homomorphism between the free group
$\mathfrak{S}$ and the group $\mathfrak{P}$. But now the elements of 
$\mathfrak{R}$ correspond to the identity permutation and, on the other
hand, each power product corresponding to the identity
permutation also belongs to $\mathfrak{R}$. This is because such power
products correspond exactly to the closed paths. Thus the association
follows.

\section{Residue class group diagrams}

Now we suppose that the complex $\mathfrak{C}$ is connected. Then the
connection between $\mathfrak{F}$, $\mathfrak{P}$ and the structure of
$\mathfrak{C}$ may be further elucidated. Let $p_0$ be an arbitrary, but
fixed, point of $\mathfrak{C}$ and let $G$ be a power product for which
the path emanating from $p_0$ is closed. The collection $\mathfrak{G}$
of these power products $G$ obviously constitute a group, a subgroup 
of $\mathfrak{F}$. For $G^{-1}$ obviously belongs to $\mathfrak{G}$
along with $G$, and $G_1 G_2$ along with $G_1,G_2$. Now if $F_1$
and $F_2$ are any two power products for which the paths emanating
from $p_0$ lead to the same point $p'$ of $\mathfrak{C}$, then $F_1$
and $F_2$ belong to the same right-sided residue class  modulo
$\mathfrak{G}$ in $\mathfrak{F}$. For $F_1 F^{-1}_2$ belongs to
$\mathfrak{G}$, and thus $F_1=G_1 F_2$. Conversely, if $F_1$ and $F_2$
are two power products which belong to the same right-sided residue
class $\mathfrak{G}F$ in $\mathfrak{F}$, then
\[
F_1=G_1 F,\quad F_2=G_2 F,
\]
and since the paths emanating from $p_0$ that correspond to the $G$
end at $p_0$, the paths emanating from $p_0$ that correspond to $F_1$
and $F_2$ end at the same point $p'$ of $\mathfrak{C}$.

If $\mathfrak{C}$ is connected then, relative to a distinguished point
$p_0$, each point of $\mathfrak{C}$ corresponds to a certain residue
class $\mathfrak{G}F$ modulo $\mathfrak{G}$ in $\mathfrak{F}$,
and for each such residue class there is a point. If the points $p'$ and
$p''$ are connected by a segment $s'$, and if $p'$ corresponds to the
residue class $\mathfrak{G}F'$, $p''$ to the residue class
$\mathfrak{G}F''$, and if $s'$ lies over $\boldsymbol{A}(s')={s^*_i}^\varepsilon$
then,
\[
\mathfrak{G}F''=\mathfrak{G}F'S^\varepsilon_i.
\]
It is therefore natural to view $\mathfrak{S}$ as the \emph{residue class
diagram} of $\mathfrak{G}$ in $\mathfrak{F}$.

Conversely, if $\mathfrak{F}$ is any group with finitely many generators
$S_1,S_2,\ldots,S_r$ and if $\mathfrak{G}$ is a subgroup of
$\mathfrak{F}$, then a residue class group diagram $\mathfrak{C}$ of
$\mathfrak{G}$ in $\mathfrak{F}$ with the generators $S_i$ may always
be constructed. Namely, let $\mathfrak{C}^*$ be the complex with one 
point and $r$ singular segments $s^*_1,s^*_2,\ldots,s^*_r$. Each
residue class $\mathfrak{G}F$ modulo $\mathfrak{G}$ in
$\mathfrak{F}$ corresponds to a point $p$ of $\mathfrak{C}$, $p'$
corresponds to the residue class $\mathfrak{G}F'$ and $p''$ to the
residue class $\mathfrak{G}F''$, and if
\[
\mathfrak{G}F''=\mathfrak{G}F' S_i
\]
then $p'$ and $p''$ may be connected by a segment $s$ which begins
at $p'$ and ends at $p''$, and $s$ covers the segment $\boldsymbol{A}(s)=s^*_i$.
Apart from the segments given in this way, $\mathfrak{C}$ contains
only their oppositely directed segments.\footnote{\textsc{O. Schreier},
Hamb. Abhdlg, \textbf{5} (1929), 180.}

Thus there are exactly $2r$ segments emanating from each point $p$
of $\mathfrak{C}$, lying over $s^*_i$ or ${s^*_i}^{-1}$. Then if we
let $\boldsymbol{A}(p)=p^*$ one sees that $\boldsymbol{A}$ is a covering of $\mathfrak{C}^*$
by $\mathfrak{C}$.

Instead of constructing the group $\mathfrak{F}$ first, one may also directly
construct the subgroup $\mathfrak{U}$ of the free group $\mathfrak{S}$
whose elements correspond to closed paths beginning at $p_0$.  Then
to each point different from $p_0$ there corresponds a residue class
$\mathfrak{U}S$ modulo $\mathfrak{U}$ in $\mathfrak{S}$. The group
$\mathfrak{R}$ defined in Section 4.16 is the intersection of the
subgroups conjugate to $\mathfrak{U}$ in $\mathfrak{S}$; the group
$\mathfrak{U}$ contains the power products of the $S_i$ that yield
elements of $\mathfrak{G}$.

Conversely, for each subgroup $\mathfrak{U}$ of $\mathfrak{S}$ there
is a such a regular complex and a covering by it of the complex of a single 
point with $r$ singular segments. \emph{Since the paths emanating from a
point $p_0$ always constitute a free group, one sees that all subgroups of
a free group are free.} (Cf. Section 3.9.)\footnote{(Translator's note.)}

If $\mathfrak{C}^*$ is any connected complex and $\mathfrak{C}$ is
a connected covering of $\mathfrak{C}^*$ then the fundamental
domains of this covering with respect to the complexes
$\mathfrak{B}_{p^{(i)}}$ may be put in one-to-one correspondence
with the residue classes $\mathfrak{G}F$ modulo a group 
$\mathfrak{G}$ in a group $\mathfrak{F}$. One obtains a residue
class group diagram of $\mathfrak{G}$ in $\mathfrak{F}$ by
contracting $\mathfrak{B}_{p^{(i)}}$ to a point $p^{(i)}$.

\section{Regular Coverings}

A particularly important class of coverings are the regular ones, 
defined as follows:

\emph{If $w$ and $w'$ are two paths of $\mathfrak{C}$ lying
over the same path $w^*$ of $\mathfrak{C}^*$, $\boldsymbol{A}(w)=
\boldsymbol{A}(w')$, and if
$w'$ is always closed when $w$ is, then the covering given by
$\boldsymbol{A}(\mathfrak{C})=\mathfrak{C}^*$ is regular.}

The simplest regular covering is a one-to-one mapping 
$\boldsymbol{A}(\mathfrak{C})=\mathfrak{C}^*$, where each $p^*$ is thus 
associated with only one $p$, and hence each $s^*$ with only one $s$.

If $w$ is a simple open path then $w'$ must also be a simple open path
when $\boldsymbol{A}(w)=
\boldsymbol{A}(w')$.  For if $w'_1$ is a closed subpath of $w'$ then
there is a well-defined subpath $w_1$ of $w$ for which 
$\boldsymbol{A}(w_1)=
\boldsymbol{A}(w'_1)$. Then $w_1$ must also be closed, contrary ot the
assumption that $w$ is simple an open. One deduces that, if $w$ is a
simple closed path, then $w'$ is also a simple closed path when
$\boldsymbol{A}(w)=\boldsymbol{A}(w')$.

One might conjecture that the composition of regular coverings is
transitive, i.~e., the following state of affairs: if $\mathfrak{C},
\mathfrak{C}', \mathfrak{C}''$ are three complexes,
$\boldsymbol{A}_1(\mathfrak{C})=\mathfrak{C}'$ is a regular mapping of
$\mathfrak{C}$ onto $\mathfrak{C}'$, and 
$\boldsymbol{A}_2(\mathfrak{C}')=\mathfrak{C}''$ is a regular mapping of
$\mathfrak{C}'$ onto $\mathfrak{C}''$, then the mapping of
$\mathfrak{C}$ onto $\mathfrak{C}''$ given by
\[
\boldsymbol{A}_2(\boldsymbol{A}_1(\mathfrak{C}))=
\boldsymbol{A}_3(\mathfrak{C})=\mathfrak{C}''
\]
is likewise regular. However, this is not the case, as examples
easily show.

Now let $\mathfrak{C}^*$ be the complex of one point and $r$
singular segments $s^*_i$. If we construct the groups $\mathfrak{F}$ 
and $\mathfrak{G}$ as in Section 4.17 then we see that $\mathfrak{G}$
is the group consisting only of the identity element $E$ of $\mathfrak{F}$.
For if the power products $G$ correspond to the closed paths
emanating from $p$, then each $G$ corresponds to a closed path in
$\mathfrak{C}$, and $G$ belongs to $\mathfrak{R}$.

If $\mathfrak{C}$ is connected then $\mathfrak{C}$ is called the
\emph{group diagram}\footnote{\textsc{M. Dehn}, Math. Ann.,
\textbf{69} (1910), 137.} of $\mathfrak{F}$ in the generators $S_i$.
Each group $\mathfrak{F}$ with a system of generators may be
associated with a group diagram; one has only to replace the
residue classes $\mathfrak{G}F$ in the method of Section 4.17 by
$EF=F$, i.~e., by the group elements themselves. \emph{A group
diagram is a regular covering of the associated complex
$\mathfrak{C}^*$.}

\section{Iterated Coverings and Groups}

If $\mathfrak{C}, \mathfrak{C}^*$ and $\mathfrak{C}^{**}$ are three
connected complexes and $\boldsymbol{A}(\mathfrak{C})=\mathfrak{C}^*$ and
$\boldsymbol{A}'(\mathfrak{C}^*)=\mathfrak{C}^{**}$ are coverings of
$\mathfrak{C}^*$ by $\mathfrak{C}$ and $\mathfrak{C}^{**}$ by
$\mathfrak{C}^*$ respectively, and if $\boldsymbol{A}(p)=p^*$, $\boldsymbol{A}'(p^*)=p^{**}$,
then the closed paths emanating from $p^{**}$ constitute a group
$\mathfrak{S}^{**}$. Let $\mathfrak{U}^*$ be the subgroup of
paths $w^{**}$ from among those corresponding to closed paths
$w^*$ emanating from $p^*$, in $\mathfrak{C}^*$, with
$\boldsymbol{A}'(w^*)=w^{**}$, and let $\mathfrak{U}$ be the subgroup of the
paths $w^{**}$ corresponding to closed paths $w$ in $\mathfrak{C}$
emanating from $p$ with 
$\boldsymbol{A}'(\boldsymbol{A}(w))=w^{**}$. Then $\mathfrak{U}$
is a subgroup of  $\mathfrak{U}^*$, because if $w$ is closed then
$\boldsymbol{A}(w)=w^*$ is also closed.

If $L^*_1, L^*_2,\ldots, L^*_h$ is a system of representatives for the
residue classes $\mathfrak{U}L^*$ modulo $\mathfrak{U}$ in
$\mathfrak{U}^{**}$, and if 
$L^{**}_1, L^{**}_2,\ldots, L^{**}_h$ is a system of representatives
for the residue classes $\mathfrak{U}^* L^{**}$ modulo 
$\mathfrak{U}^*$ in $\mathfrak{S}^{**}$, then the
$L^*_i L^{**}_j$ ($i=1,2,\ldots,h; k=1,2,\ldots,l$)
constitute a complete system of representatives for the residue classes 
modulo $\mathfrak{U}$ in $\mathfrak{S}^{**}$. If $p_{ik}$ are the
points of $\mathfrak{C}$ for which 
$\boldsymbol{A}'(\boldsymbol{A}(p_{ik}))=p^{**}$, and if
$p^*_i$ are the points of $\mathfrak{C}^*$ for which
$\boldsymbol{A}'(p^*_i)=p^{**}$ then, if the residue class 
$\mathfrak{U}L^*_i L^{**}_k$ corresponds to $p_{ik}$ and if the 
residue class $\mathfrak{U}^* L^{**}_k$ corresponds to $p^*_k$,
we have $\boldsymbol{A}(p_{ik})=p^*_k$. From this it follows, conversely, that:
\emph{given two coverings $\boldsymbol{A}''(\mathfrak{C})=\mathfrak{C}^{**}$
and $\boldsymbol{A}'(\mathfrak{C}^*)=\mathfrak{C}^{**}$, where $\boldsymbol{A}''$ corresponds
to a subgroup $\mathfrak{U}$ of $\mathfrak{S}^{**}$ and $\boldsymbol{A}'$
corresponds to a subgroup $\mathfrak{U}^*$ of $\mathfrak{S}^{**}$
and $\mathfrak{U}$ is a subgroup of $\mathfrak{U}^*$, then there is
a further covering $\boldsymbol{A}(\mathfrak{C})=\mathfrak{C}^*$ such that
$\boldsymbol{A}'(\boldsymbol{A}(\mathfrak{C}))=
\boldsymbol{A}''(\mathfrak{C})$.}

If $\mathfrak{C}^{**}$ is a complex of one point and $r$ singular
segments $s^{**}_i$ ($i=1,2,\ldots, r$), $\mathfrak{C}^*$ is a residue class
group diagram of the group $\mathfrak{G}$ in $\mathfrak{F}$ with the
generators $S_i$ ($i=1,2,\ldots,r$), if $\mathfrak{C}$ is the group
diagram of $\mathfrak{F}$ in the generators $S_i$ ($i=1,2,\ldots,r$), 
and if $\boldsymbol{A}'(\mathfrak{C}^*)=\mathfrak{C}^{**}$ and 
$\boldsymbol{A}''(\mathfrak{C})=\mathfrak{C}^{**}$ are the corresponding coverings,
then there is also a covering $\boldsymbol{A}(\mathfrak{C})=\mathfrak{C}^*$.
Because if the covering $\boldsymbol{A}''(\mathfrak{C})=\mathfrak{C}^{**}$ belongs
to the subgroup $\mathfrak{U}$ of the fundamental group 
$\mathfrak{S}^{**}$ of $\mathfrak{C}^{**}$, and the covering
$\boldsymbol{A}'(\mathfrak{C}^*)=\mathfrak{C}^{**}$ belongs to the subgroup
$\mathfrak{U}^*$ of the fundamental group $\mathfrak{S}^{**}$ then
$\mathfrak{U}$ consists exactly of the representations of the identity of the
group $\mathfrak{F}$ in the generators $S_i$, and hence the
$\boldsymbol{A}'(w^*)=w^{**}$ for which $[w^{**}]$ belongs to $\mathfrak{U}$
must be in the residue class group diagram of a subgroup $\mathfrak{G}$
of $\mathfrak{G}$ of $\mathfrak{F}$, as well as in $\mathfrak{C}^*$.
Analogous results hold for any connected complex $\mathfrak{C}^{**}$.

We obtain a special covering when we apply the covering construction to the
subgroup of $\mathfrak{S}^*$ consisting of the identity element alone. The
complex $\mathfrak{C}$ obtained in this way is called the \emph{universal 
covering complex} of $\mathfrak{C}^*$. if $\mathfrak{C}'$ is a connected 
complex and $\boldsymbol{A}(\mathfrak{C}')=\mathfrak{C}^*$ is a covering of
$\mathfrak{C}^*$ by $\mathfrak{C}'$ then, by the theorem proved above,
there is also a covering $\boldsymbol{A}'(\mathfrak{C})=\mathfrak{C}'$ of
$\mathfrak{C}'$ by the universal covering complex $\mathfrak{C}$ of
$\mathfrak{C}^*$.

\emph{The universal covering complex $\mathfrak{C}$ is a tree.} Namely, 
let $w$ be a closed reduced nonempty path in $\mathfrak{C}$. Then
$\boldsymbol{A}(w)=w^*$ is a closed path in $\mathfrak{C}^*$ and $[w^*]$ is the identity
element of the fundamental group $\mathfrak{S}^*$ of $\mathfrak{C}^*$,
so $[w^*]$ is not reduced and likewise $w$ is not reduced.

If $\mathfrak{C}^*$ is a complex of one point and $r$ singular segments
then the universal covering complex can be viewed as the group diagram of
the free group with generators $S_i$ ($i=1,2,\ldots,r$). This group diagram
is therefore a tree.

\section{Transformations into Itself}

If $\mathfrak{C}$ is a multiple regular covering of $\mathfrak{C}^*$, then
$\mathfrak{C}$ admits a group of mappings onto itself of the following kind:
under the mapping, corresponding elements always remain over the same
element of $\mathfrak{C}^*$. If $p^{(i)}$ and $p^{(k)}$ are two points over
$p^*$ then there is one and, when $\mathfrak{C}$ is connected, only one
transformation which carries $p^{(i)}$ to $p^{(k)}$. The transformations
constitute a group which is homomorphic to the group $\mathfrak{F}$,
respectively $\mathfrak{P}$, and in fact it is isomorphic when $\mathfrak{C}$
is connected. From now on we assume that $\mathfrak{C}$ is connected.

If $\mathfrak{C}$ is a group diagram we define our transformations as follows:
if $p_F$ is the point corresponding to the group element $F$ of $\mathfrak{F}$
and $F'$ is an arbitrary element of $\mathfrak{F}$ then the transformation
$\boldsymbol{I}_{F'}$ of the complex $\mathfrak{C}$ corresponding to $F'$
sends the point $p_F$ to $p_{FF'}$, and the segments $s$ over $s^*_i$
emanating from  $p_F$ go to the segments $s'$ over $s^*_i$ emanating from
$p_{FF"}$. The final point $p_{FS_i}$ of $s$ in fact goes to the final point
$p_{F'FS_i}$ of $s'$; $\boldsymbol{I}_{F'}$ is therefore an automorphism
of $\mathfrak{C}$. There is exactly one $\boldsymbol{I}_{F'}$ which carries
$p_{F_1}$ to $p_{F_2}$.

The $\boldsymbol{I}_F$ constitute a group isomorphic to $\mathfrak{F}$, 
because
\[
\boldsymbol{I}_{F''}(\boldsymbol{I}_{F'}(\mathfrak{C}))=
\boldsymbol{I}_{F'' F'}(\mathfrak{C}).
\]
Finally, if $\boldsymbol{I}$ is any automorphism of $\mathfrak{C}$ which 
exchanges elements over the same element of $\mathfrak{C}^*$ and which
carries $p_F$ to $p_{F'F}$ then 
$\boldsymbol{I}_{F'^{-1}}(\boldsymbol{I}(\mathfrak{C}))$ is a mapping
which leaves the point $p_F$ fixed, and hence also the segments emanating
from $p_F$, and their final points, and so on. Thus
$\boldsymbol{I}_{F'^{-1}}(\boldsymbol{I}(\mathfrak{C}))$ is the identity
map.

If $\mathfrak{C}$ is a connected covering of an arbitrary complex
$\mathfrak{C}^*$, then the subcomplexes $\mathfrak{B}_{p^{(i)}}$ may be
associated with the group elements $F_i$ of $\mathfrak{F}$. If
\[
FF_i=F_{k_i}
\]
then $F$ corresponds to the automorphism $\boldsymbol{I}_F$ sending
$\mathfrak{B}_{p^{(i)}}$ to $\mathfrak{B}_{p^{(k_i)}}$.
$\boldsymbol{I}_F$ is uniquely and consistently determined by this, as one
easily establishes using Section 4.17. One derives the properties of the group
of the $\boldsymbol{I}_F$ analogously as with group diagrams.

The converse also holds: if $\mathfrak{C}$ is a connected complex and
$\mathfrak{I}'$ is a group of automorphisms $\boldsymbol{I}$ of
$\mathfrak{C}$, and if a transformation from $\mathfrak{I}'$ which leaves
a point $p$ of $\mathfrak{C}$ fixed is the identity mapping, then there is a
connected complex $\mathfrak{C}^*$ and a regular covering $\boldsymbol{A}$
 such that
$\boldsymbol{A}(\mathfrak{C})=\mathfrak{C}^*$. The transformation group
$\mathfrak{I}$ associated with the covering is identical with 
$\mathfrak{I}'$.\footnote{With these observations, Reidemeister has come close 
to stating explicitly the 1977 theorem of Bass and Serre that a group acting
freely on a tree is free. (Translator's note.)}

For example, if $\mathfrak{C}$ is the group diagram of a group 
$\mathfrak{F}$ and $\mathfrak{V}$ is a subgroup of $\mathfrak{F}$,
then one can construct a complex $\mathfrak{C}'$ and a regular covering
$\boldsymbol{A}(\mathfrak{C})=\mathfrak{C}'$ for which $\mathfrak{V}$ is the
transformation group. If $\mathfrak{C}^*$ is the complex with a single 
point and singular segments corresponding to the generators of
$\mathfrak{F}$, and if $\boldsymbol{A}(\mathfrak{C})=\mathfrak{C}^*$ is the
covering brought about by these generators, then there is also a covering
$\boldsymbol{A}'(\mathfrak{C}')=\mathfrak{C}^*$, and in fact $\boldsymbol{A}'$ is regular or not
according as $\mathfrak{V}$ is an invariant subgroup of $\mathfrak{F}$
or not.

One can also elucidate the process for determining generators and
defining relations of subgroups by means of group diagrams. Retaining
the notation of the last section, let $\mathfrak{B}'_{p'}$ be a spanning tree 
of $\mathfrak{C}'$, and
let $\mathfrak{B}_{p^{(i)}}$ be the subcomplex of $\mathfrak{C}'$
for which $\boldsymbol{A}(\mathfrak{B}_{p^{(i)}})=\mathfrak{B}'_{p'}$. Then
corresponding to the simple paths $w'_k$ in $\mathfrak{B}'_{p'}$
emanating from $p'$, respectively the paths $w^{(i)}_k$ in
$\mathfrak{B}_{p^{(i)}}$ emanating from $p^{(i)}$, we have power
products $L_k$ in the generators of $\mathfrak{F}$ which yield a complete
system of representatives for the residue classes $\mathfrak{V}L$ and
which satisfy the condition $(\Sigma)$ of Section 3.5. One obtains from
$\mathfrak{C}$ a group diagram $\mathfrak{C}_\mathfrak{V}$ of
$\mathfrak{V}$ by shrinking the segments of each 
$\mathfrak{B}_{p^{(i)}}$ to the point $p^{(i)}$.

From this one sees again that \emph{subgroups of free groups are free}
(cf. Section 3.9).
Because, when $\mathfrak{F}$ is a free group on free
generators, then $\mathfrak{C}$ is a tree; hence
$\mathfrak{C}_\mathfrak{V}$ is also a tree.

\chapter{Surface Complexes}

\section{The Concept of a Surface Complex}

By a surface complex\footnote{Usually known today as a 2-complex.
Likewise, what Reidemeister calls a ``surface piece'' would usually be
called a 2-cell today. However, since Reidemeister does not intend to
discuss the general case of an $n$-complex, it is more faithful to his
mindset to use the ``point'', ``line segment'' and ``surface'' terminology.
(Translator's note.)} 
$\mathfrak{F}$  we mean a finite or denumerable collection of points,
line segments and surface pieces satisfying the following conditions.

A.1. \emph{The points and line segments of $\mathfrak{F}$ form a
line segment complex $\mathfrak{C}$.}

A.2. \emph{If $f$ is a surface piece of $\mathfrak{F}$ then there is a 
closed path $w$ of $\mathfrak{C}$ which circumscribes $f$ once in
the positive sense. If $w'$ is a second path which circumscribes $f$
once in the positive sense then $w'$ is convertible to $w$ by a
cyclic interchange.}

A.3. \emph{If $w$ circumscribes the surface piece $f$ positively, then
$w^{-1}$ circumscribes $f$ negatively. For each surface piece $f$ there
is an oppositely directed $f^{-1}$. The surface piece $(f^{-1})^{-1}$ 
equals $f$, and $w^{-1}$ circumscribes $f^{-1}$ positively.}

The complex consisting of the elements of $w$ is called the boundary
of $f$. An element of the boundary bounds $f$. The surface piece $f$
is called singular if its boundary path is not simple.

We are mainly concerned with special surface complexes called
two-dimensional manifolds $\mathfrak{M}$. These satisfy four further
axioms, of which we initially give only the first three.

A.4. \emph{The line segment complex $\mathfrak{C}$ is connected.}

A.5. \emph{If $s$ is a segment of $\mathfrak{C}$, then $s$ appears
in the boundary of some surface piece $f$.}

A.6. \emph{A boundary path $w$ runs through a segment $s$ at most
twice. If the boundary path $w$ of a surface piece $f$ runs through the
segment $s$ only once then there is exactly one surface piece
$f'\ne f^{\pm 1}$, the boundary of which, $w'$, runs through $s$, and 
then only once. If $w$ runs through the segment $s$ twice, then $s$
bounds only the surface piece $f$ and the oppositely oriented $f^{-1}$. }

If $s$ is a segment which belongs to the boundary of $f$ then there is
a boundary path of $f$ which begins with $s$. We denote it by
$sw$, with the understanding that $w$ may be empty. If $s$ occurs
twice in a boundary path $sw$ of $f$, then the latter is denoted by
\[
sw=sw_1 s^\varepsilon w_2\qquad (\varepsilon=\pm 1),
\]
where the $w_i$ may again be empty. There is then a second
boundary path of $f$ that begins with $s$, namely
\begin{align*}
sw'&=sw_2 s w_1\quad\text{when}\quad \varepsilon=+1,\\
sw'&=sw^{-1}_1 s^{-1} w^{-1}_2\quad\text{when}\quad \varepsilon=-1.
\end{align*}
If $s$ occurs only once in the boundary path of $f$, then there is a
second surface piece $f'$ with a boundary path of the form
$sw'$. In any case, we can speak of the two boundary paths $sw$ and
$sw'$ of $\mathfrak{M}$ which begin with $s$. We call the two
segments $s_\alpha,s_\beta$ neighbors of $s$ when $s_\alpha$ follows
$s$ in $sw$ and $s_\beta$ follows $s$ in $sw'$.

The segments $s_\alpha$ and $s_\beta$ both begin at the point where
$s$ ends. If $s_\alpha$ is a neighbor of $s$, then $s^{-1}$ is
a neighbor of $s^{-1}_\alpha$. If $s_\alpha$ is a neighbor of
$s$ and $s_\beta$ is a neighbor of $s$ and $s_\alpha\ne s_\beta$,
then $s_\alpha,s_\beta$ are the two neighbors of $s$. If $s$ is
a neighbor of $s^{-1}$, then there is a boundary path $ss^{-1}w$,
and consequently the second boundary path beginning with $s$ is
$ss^{-1}w^{-1}$. Then both segments neighboring $s$ are
equal to $s^{-1}$.

\section{Stars}

The segments emanating from a point,
\begin{equation}
s_{\alpha_1},\quad
s_{\alpha_2},\quad\ldots,\quad
s_{\alpha_m}, \tag{1}
\end{equation}
are said to form a \emph{star} when their ordering is such that
$s_{\alpha_{i-1}}$ and $s_{\alpha_{i+1}}$ are the neighbors of
$s^{-1}_{\alpha_{i}}$. The segments (1) form a \emph{closed star}
when
\[
s_{\alpha_m},\quad
s_{\alpha_1},\quad
s_{\alpha_2},\quad\ldots,\quad
s_{\alpha_m},\quad
s_{\alpha_1}
\]
is a star. Thus a closed star can consist of a single segment. The sequence
\begin{equation}
s_{\alpha_i},\quad
s_{\alpha_{i+1}},\quad\ldots,\quad
s_{\alpha_{i+l}}, \tag{2}
\end{equation}
is also a star when (1) is and $l\ge 1$; (2) is called a substar of (1). A star
in which each segment appears only once is called \emph{simple}.

\emph{If $p$ is a point from which only finitely many segments emanate, 
then there is a simple closed star composed of those segments.}

Namely, if
\[
s_{\alpha_1},\quad
s_{\alpha_2},\quad\ldots,\quad
s_{\alpha_m}
\]
is a star then we can certainly extend it by a further segment
$s_{\alpha_{m+1}}$. Since $s_{\alpha_m}$ is a neighbor of
$s^{-1}_{\alpha_{m-1}}$, $s_{\alpha_{m-1}}$ is
a neighbor of $s^{-1}_{\alpha_{m}}$; now let $s_{\alpha_{m+1}}$
be the second segment neighboring $s^{-1}_{\alpha_{m}}$.

Thus there is certainly a star in which a segment $s_1$ beginning at
$p$ appears twice, and hence also a star
\[
s_{\alpha_1},\quad
s_{\alpha_2},\quad\ldots,\quad
s_{\alpha_n},\quad s_{\alpha_{n+1}}
\]
in which $s_{\alpha_1}= s_{\alpha_{n+1}}$ and all the remaining
$ s_{\alpha_{i}}$ are different from each other and from $ s_{\alpha_{1}}$.
Now, either $n=1$, so that $s_{\alpha_1}=s_{\alpha_2}$ and both the
segments neighboring $s^{-1}_{\alpha_1}$ are $s_{\alpha_1}$ itself,
so $s_{\alpha_1}$ alone is a closed star; or else $n>1$, in which case
\[
s_{\alpha_1},\quad
s_{\alpha_2},\quad\ldots,\quad
s_{\alpha_n}
\]
is a closed star, because
\[
s_{\alpha_n},\quad
s_{\alpha_1},\quad
s_{\alpha_2},\quad\ldots,\quad
s_{\alpha_n},\quad s_{\alpha_{1}}
\]
is a star because $s_{\alpha_1}$ is a neighbor of $s_{\alpha_2}$
and to $s_{\alpha_n}$, and $s_{\alpha_2}$ is different from
$s_{\alpha_n}$.

A proper substar of a simple star (1) is not closed, otherwise
$s_{\alpha_i},s_{\alpha_{i+1}},\ldots,s_{\alpha_{i+l}}$ would be a
closed substar of (1) and this would mean 
$s_{\alpha_{i-1}}=s_{\alpha_{i+l}}$ when $i>1$,
$s_{\alpha_{i}}=s_{\alpha_{i+l+1}}$ when $i=1$.

Moreover, one sees that two stars consisting of segments beginning at
$p$, which coincide in their first and last elements and contain the 
same elements, are identical.

Given two stars with the same number of elements,
say $m$, such that the first and second terms of the first coincide with
the last and second last, respectively, of the second, the $k$th term of the 
first is identical with the $(m-k)$th term of the second. It then follows:

\emph{If a segment $s$ appears in a simple closed star $(1)$
consisting of segments which begin at $p$, then all simple closed stars
of this kind result from $(1)$ or from}
\[
s_{\alpha_m},\quad
s_{\alpha_{m-1}},\quad,\ldots,\quad
s_{\alpha_2},\quad
s_{\alpha_1}
\]
\emph{by cyclic interchanges, and all closed stars result from a $k$-fold
repetition of a simple closed star.} On the other hand, each segment
appears in some closed star.

\section{Manifolds}

Our last axiom for manifolds can now be stated:

A.7. \emph{If $s_1$ and $s_2$ are two segments which begin at $p$
then there is a star of segments beginning at $p$,}
\[
s_{\alpha_1},\quad
s_{\alpha_2},\quad\ldots,\quad
s_{\alpha_m},
\]
\emph{with $s_{\alpha_1}=s_1$ and $s_{\alpha_m}=s_2$.}

It follows from this and the theorems of the previous section that
there is no closed star around $p$ which does not contain all
the segments beginning at $p$.

If only finitely many segments begin at $p$ then there is a simple
closed star of all these segments. Then we can say that the segments
beginning at $p$ are cyclically ordered in a certain way. If infinitely
many segments begin at $p$ then all stars of these segments are
open and consequently simple. The segments are ``linearly''
ordered.

A few more conclusions about the boundary paths of surface pieces
may be drawn from axiom A.7: if there is a boundary path $w$ which
is not reduced, and hence equals $w_1 ss^{-1} w_2$, and if $p$ is the
final point of $s$, then only a single segment emanates from $p$,
for $s^{-1}$ constitutes a closed star. It follows, further, that such a
segment $s$ cannot be singular. Otherwise the segments $s$ and
$s^{-1}$ would emanate from the final point $p$ of $s$, and $s$
and $s^{-1}$ together would constitute a closed star. Similarly,
one sees that a subpath
\[
sk_1 k_2 \cdots k_r s^{-1}\quad\text{where}\quad
k_i=s_{i1}s_{i2}s^{-1}_{i1}s^{-1}_{i2}
\]
and the $s,s_{ik}$ are singular segments with boundary point $p$
cannot appear in a boundary path. For 
\[
s, v_1, v_2, \ldots, v_r\quad\text{with}\quad
v_i=s^{-1}_{i1},s_{i2},s_{i1},s^{-1}_{i2}
\]
constitutes a star, and in fact a simple closed star, which does not contain 
the segment $s^{-1}$ emanating from $p$.

If $\mathfrak{M}$ is a manifold and if the boundary paths of the
surface pieces $f_1,f_2,\ldots,f_k$ are simple paths which together 
include each point at most once, then the manifold that results from
$\mathfrak{M}$ by leaving out the surface pieces $f_1,f_2,\ldots,f_k$
may be called a \emph{manifold with boundary}.

Examples of manifolds are the canonical normal forms given in
Section 5.11, polyhedra in Euclidean space, \textsc{Riemann}
surfaces, or the subdivided surfaces described in the  introduction.

\section{An Auxiliary Manifold}

The most important consequence of axiom A.7 in Section 5.3 is the
existence of the dual manifold $\mathfrak{D}$ \emph{for each
manifold $\mathfrak{M}$ in which only finitely many segments
emanate from each point $p$.}

We first derive an auxiliary manifold $\mathfrak{M}''$ from
$\mathfrak{M}$ in the following way:

Each point $p$ of $\mathfrak{M}$ corresponds to a certain point
$\Delta_1(p)=p''$ of $\mathfrak{M}''$. Each segment $s$ of
$\mathfrak{M}$ also corresponds to a point $\Delta_1(s)=q''$ of
$\mathfrak{M}''$ and in fact
\[
\Delta_1(s)=\Delta_1(s^{-1}).
\]
Finally, each surface piece $f$ of $\mathfrak{M}$ corresponds to a
point $\Delta_1(f)=r''$ of $\mathfrak{M}''$, and we have
\[
\Delta_1(f)=\Delta_1(f^{-1}).
\]
Conversely, under this correspondence each point of $\mathfrak{M}''$
is associated with either a point, a segment pair $s^{\pm 1}$, or a
surface piece pair $f^{\pm 1}$ of $\mathfrak{M}$.

Further, each segment $s$ corresponds to a path of two segments in
$\mathfrak{M}''$,
\[
\Delta_2(s)=s''_1 s''_2.
\]
If $p_1,p_2$ are the initial and final points of $s$, then $\Delta_1(p_1)$
is the initial point of $s''_1$, $\Delta_1(p_2)$ is the final point of $s''_2$,
and $\Delta_1(s)$ is the final point of $s''_1$ and the initial point of
$s''_2$. In addition,
\[
\Delta_2(s^{-1})=(\Delta_2(s))^{-1}={s''_2}^{-1} {s''_1}^{-1}.
\]
Moreover, two segments $t''_1,t''_2$ emanate from each point
$\Delta_1(s)$, ending at those points $\Delta_1(f_1)$ and $\Delta_1(f_2)$
corresponding to the regions $f_1$ and $f_2$ bounded by $s$. Each
segment of $\mathfrak{M}''$ begins or ends at a point $\Delta_1(s)$.
The points $\Delta_1(s)$ are of order four, the points $\Delta_1(f)$
are of order $k$ when a simple boundary path of $f$ has $k$
elements.

Finally we come to the surface pieces of $\mathfrak{M}''$. These are
all quadrilaterals, and each $k$-gon of $\mathfrak{M}$ corresponds to
exactly $k$ quadrilaterals in $\mathfrak{M}''$. If
\[
w=s_{\alpha_1}s_{\alpha_2}\cdots s_{\alpha_k}
\]
is a simple boundary path of $f$, let
\begin{align*}
\Delta_2(w)&=
\Delta_2(s_{\alpha_1})\Delta_2(s_{\alpha_2})\cdots\Delta_2(s_{\alpha_k})\\
&=
\left({s''_{\alpha_1}}^{(1)} {s''_{\alpha_1}}^{(2)}\right)
\left({s''_{\alpha_2}}^{(1)} {s''_{\alpha_2}}^{(2)}\right)\cdots
\left({s''_{\alpha_k}}^{(1)} {s''_{\alpha_k}}^{(2)}\right)
\end{align*}
and let $t''_{\alpha_1},t''_{\alpha_2},\ldots,t''_{\alpha_k}$ respectively
be the segments that lead from $\Delta_1(f)$ to the
$\Delta_1(s_{\alpha_i})$. Then there are exactly $k$ quadrilaterals
$f''_{\alpha_1},f''_{\alpha_2},\ldots,f''_{\alpha_k}$ and in fact let
\begin{center}
$t''_{\alpha_1}{s''_{\alpha_1}}^{(2)}
{s''_{\alpha_2}}^{(1)}{t''_{\alpha_2}}^{-1}$
be a boundary path of $f''_{\alpha_1}$,

$t''_{\alpha_i}{s''_{\alpha_i}}^{(2)}
{s''_{\alpha_{i+1}}}^{(1)}{t''_{\alpha_{i+1}}}^{-1}$
be a boundary path of $f''_{\alpha_i}$,
\end{center}
and finally let
\begin{center}
$t''_{\alpha_k}{s''_{\alpha_k}}^{(2)}
{s''_{\alpha_{1}}}^{(1)}{t''_{\alpha_{1}}}^{-1}$
be a boundary path of $f''_{\alpha_k}$.
\end{center}

\emph{The surface complex $\mathfrak{M}''$ so described is a
manifold.} For the one-dimensional complex contained in it is
connected, because any points $\Delta_1(f)$ and $\Delta_1(s)$
are connected to points $\Delta_1(p)$, and each point $\Delta_1(p)$
is connected to all the others of the same kind. Each segment of a
path $\Delta_2(s)$ is run through at most once in a boundary path,
and each such segment bounds two different surface pieces. The
same holds for the segments $t''$ bounded by a point $\Delta_1(f)$
corresponding to a $k$-gon $f$ ($k>1$). A segment $t''$ bounded by
a $\Delta_1(f)$ corresponding to 1-gon is run through twice by the
boundary path of the single surface piece in $\mathfrak{M}''$
corresponding to $f$. Finally, axiom A.7 is also satisfied by the three
kinds of points $\Delta_1(p),\Delta_1(s),\Delta_1(f)$.

\section{Dual Manifolds}

For each point $p$ of $\mathfrak{M}$ we now collect together all the 
surface pieces $f''$ of $\mathfrak{M}''$, together with the segments and
points of their boundaries, into a surface complex $\mathfrak{P}$
whose boundary contains the point $p''=\Delta_1(p)$.
Each segment emanating from $\Delta_1(p)$
appears in the boundary of exactly two surface pieces of $\mathfrak{P}$.
Let
\[
f''_{\beta_1},\quad
f''_{\beta_2},\quad\ldots,\quad
f''_{\beta_m}
\]
be all the distinct surface pieces of $\mathfrak{P}$, ordered in such a
way that $f''_{\beta_i}$ and $f''_{\beta_{i+1}}$, $f''_{\beta_m}$
and $f''_{\beta_1}$, have exactly one boundary segment in common
if $m>1$, and oriented in such a way that the common boundary segment
of $f''_{\beta_i}$ and $f''_{\beta_{i+1}}$ is traversed in the opposite
directions for positive orientations of $f''_{\beta_i}$ and $f''_{\beta_{i+1}}$. 
Now let ${t''_{\beta_i}}^{(1)}{t''_{\beta_i}}^{(2)}$ be the piece of the
boundary path of $f''_{\beta_i}$ from a segment $t$. Then 
${t''_{\beta_i}}^{(1)}$ begins at a $\Delta_1(s)$ while ${t''_{\beta_i}}^{(2)}$
ends at such a point; ${t''_{\beta_i}}^{(1)}$ ends, and 
${t''_{\beta_i}}^{(2)}$ begins, at a $\Delta_1(f)$. The path
\begin{equation}
w=
{t''_{\beta_1}}^{(2)}
{t''_{\beta_2}}^{(1)}
{t''_{\beta_2}}^{(2)}
{t''_{\beta_3}}^{(1)}\cdots
{t''_{\beta_m}}^{(2)}
{t''_{\beta_1}}^{(1)}
\tag{1}
\end{equation}
is a closed path we will call the boundary path of $\mathfrak{P}$.

We now construct the complex $\mathfrak{D}$ \emph{dual to}
$\mathfrak{M}$ by associating with each path $t^{(2)} t^{(1)}$ a
segment $\Delta(t^{(2)}t^{(1)})=t'$ of $\mathfrak{D}$,  with
each complex $\mathfrak{P}$ a pair of surface pieces 
$\Delta'(\mathfrak{P})={f'}^{\pm 1}$ of $\mathfrak{D}$, and with each 
point $\Delta_1(f)$ of $\mathfrak{M}''$ a unique point
$\Delta'(\Delta_1(f))=p'$ of $\mathfrak{D}$. If  $t^{(2)} t^{(1)}$
begins and ends at $p''_1$ and $p''_2$ respectively then
$\Delta'(t^{(2)} t^{(1)})$ begins and ends at $\Delta'(p''_1)$ and
$\Delta'(p''_2)$ respectively. Let
\[
\Delta'((t^{(2)} t^{(1)})^{-1})=(\Delta'(t^{(2)} t^{(1)}))^{-1},
\]
and if (1) is the boundary path of $\mathfrak{P}$ let
\[
\Delta'(w)=
\Delta'({t''_{\beta_1}}^{(2)}{t''_{\beta_2}}^{(1)})
\cdots
\Delta'({t''_{\beta_m}}^{(2)}{t''_{\beta_1}}^{(1)})
\]
be a boundary path of $\mathfrak{P}$.

One easily sees that: \emph{the dual complex $\mathfrak{D}$ is a
manifold}. For the segments $s'$ are traversed by the boundary paths in
the way required by A.6, and the stars of the $s'$ come directly from the 
stars of the segments $t''$ around the points $\Delta_1(f)$ of
$\mathfrak{M}''$.

The structure of $\mathfrak{D}$ may be described as follows: each
point $p$ of $\mathfrak{M}$ corresponds to a unique complex
$\mathfrak{P}$ and hence also to the surface piece pair
$\Delta(p)={f'}^{\pm 1}$ of $\mathfrak{D}$; each surface piece of
$\mathfrak{M}$ corresponds to exactly one point $\Delta(f)=p$ of
$\mathfrak{D}$ and
\[
\Delta(f^{-1})=\Delta(f).
\]
Each segment $s$ of $\mathfrak{M}$ corresponds firstly to a point
$\Delta_1(s)$, and through this point goes a path $t^{(1)}t^{(2)}$
and its inverse. The resulting correspondence between the segment
pairs $s^{\pm 1}$ and the segment pairs 
$\Delta(s^{\pm 1})={s'}^{\pm 1}$ is likewise one-to-one. 
Further, the boundary relations change as follows. If $p$ bounds the
segment $s$ or $s^{-1}$, then $\Delta(p)$ is bounded by the segments
$\Delta(s^{\pm 1})$. If $s^{\pm 1}$ bounds a surface piece $f$, then
one of the segments ${s'}^{\pm 1}$ is bounded by $\Delta(f)$,
where ${s'}^{\pm 1}=\Delta(s^{\pm 1})$.

\section{Dual Line Segment Complexes}

As a result of the second relation one can construct the line segment
complex $\mathfrak{D}_1$ of $\mathfrak{D}$ directly from $\mathfrak{M}$.
If we imagine the segment pairs of $\mathfrak{D}$ numbered, and the
segments of a pair denoted by $s'_i,{s'_i}^{-1}$, then it follows from the
first condition that the boundary path of $\Delta(p)$ is also determined,
provided this boundary path is a simple path and the segments
traversed by a boundary path can be combined into a path in only one way.
One can improve the statement of the structural relation between
$\mathfrak{M}$ and $\mathfrak{D}$ by remarking that if
\[
s_{\alpha_1},\quad
s_{\alpha_2},\quad\ldots,\quad
s_{\alpha_m}
\]
is a simple closed star around $p$ then there is a simple boundary path
around $\Delta(p)$,
\[
{s'_{\alpha_1}}^{\varepsilon_1}
{s'_{\alpha_2}}^{\varepsilon_2}\cdots
{s'_{\alpha_m}}^{\varepsilon_m}\qquad(\varepsilon_1=\pm 1)
\]
such that
\[
\Delta(s^{\pm 1}_{\alpha_i})={s'_{\alpha_i}}^{\pm 1}.
\]

In this way the boundary paths of surface pieces of $\mathfrak{D}$ 
are determined when no segment of $\mathfrak{D}$ is singular. Namely,
$\varepsilon_i$ is determined by the condition that (1) is a closed path.

The construction of the dual line segment complex is meaningful for all
complexes satisfying just the axioms A.1 to A.6.

$\mathfrak{D}_1$ is called the line segment complex dual to $\mathfrak{F}$
when each surface piece of $\mathfrak{F}$ is asssociated with a point
$\Delta(f)=p'$ of $\mathfrak{D}_1$, $\Delta(f)=\Delta(f^{-1})$, and
each $p'$ also corresponds to a pair $f^{\pm 1}$, and the pairs $s^{\pm 1}$
of $\mathfrak{F}$ are mapped one-to-one onto the pairs ${s'}^{\pm 1}$
of $\mathfrak{D}_1$ by
\[
\Delta(s^{\pm 1})={s'}^{\pm 1}
\]
in such a way that the boundary relations correspond dually. This definition
is consistent, because a segment $s$ of $\mathfrak{F}$ either bounds a
surface piece doubly, or else it bounds two different surface pieces singly.
If axiom A.7 is also satisfied, then $\mathfrak{D}$ is connected. Namely,
each star in $\mathfrak{F}$ corresponds to a closed path in 
$\mathfrak{D}_1$. If the star $v$ contains a segment equal, though perhaps
oppositely directed, to a segment of a star $v'$, then these correspond to paths
in $\mathfrak{D}_1$ running through the same segment. If $f_1$ and $f_2$
are two surface pieces of $\mathfrak{F}$ which contain in their boundary
points $p_1,p_2$ bounding the same segment $s$, then one can connect
$\Delta(f_1)$ and $\Delta(f_2)$ in $\mathfrak{D}_1$. For if $s_i$ is a
segment belonging to the boundary path of $f_i$ and bounded by $p_i$
then there is a star $v_i$ of the segments around $p_i$ which contains $s$
and $s_i$ ($i=1,2$). But in $\mathfrak{D}_1$ the $v_i$ correspond to
two paths $w_i$ which run through both the segments $\Delta(s^{\pm 1})$,
and consequently they permit $\Delta(f_1)$ to be connected to
$\Delta(f_2)$ by a path in $\mathfrak{D}_1$. It then follows by induction
that $\mathfrak{D}_1$ is connected; since $\mathfrak{F}$ is connected,
any two surface pieces $f_1$ and $f_2$ contain points $p_1$ and $p_2$
respectively which may be connected by a simple path.

For a manifold $\mathfrak{F}$ this means that any two surface pieces of
$\mathfrak{F}$ may be embedded in a chain of surface pieces in which
any two successive members meet along a common segment.

In order to complete the dual surface complex $\mathfrak{D}$ from
$\mathfrak{D}_1$ we cannot omit the hypothesis that each point of
$\mathfrak{F}$ bounds only finitely many segments.

\section{Elementary Transformations}

The most important properties of a complex, e.g. those required for the
existence of coverings, are preserved by certain simple alterations of the
complexes, the so-called elementary transformations. It is therefore
useful to divide the complexes into classes, in which a given complex is
grouped with all those resulting from it by elementary transformations.

If $\mathfrak{F}$ is any surface complex containing the segment $s$
beginning at $p_1$ and ending at $p_2$, then by an \emph{elementary
extension of the first kind} we mean the construction of a complex
$\mathfrak{F}'$ which contains all the elements of $\mathfrak{F}$
apart from $s^{\pm 1}$, and in place of $s^{\pm 1}$ it has segments
$s^{\pm 1}_1,s^{\pm 1}_2$ and an additional point $p'$. The segment
$s_i$ begins at $p_i$ ($i=1,2$) and ends at $p'$. In all boundary paths
of surface pieces which contain $s$, $s$ is replaced by $s_1 s^{-1}_2$
and $s^{-1}$ by $s_2 s^{-1}_1$.

By an \emph{elementary reduction of the first kind} we mean the inverse 
process, which carries $\mathfrak{F}'$ to $\mathfrak{F}$. $\mathfrak{F}$
contains all points of $\mathfrak{F}'$ except a $p'$ which bounds exactly
two segments; the latter are replaced by a single segment. 

By an \emph{elementary extension of the second kind} we mean the 
following construction of a complex $\mathfrak{F}'$ from $\mathfrak{F}$:
let $f$ be a surface piece of $\mathfrak{F}$ with boundary path
\[
w=w_1 w_2,
\]
where $w_1$ begins at $p_1$ and ends at $p_2$. The new complex
$\mathfrak{F}'$ is constructed by taking all elements of $\mathfrak{F}$
apart from $f^{\pm 1}$, and in its place taking two surface pairs
$f^{\pm 1}_1, f^{\pm 1}_2$ and a segment $s'$ which begins at $p_1$
and ends at $p_2$, together with $s'^{-1}$; $w_1 s'^{-1}$ is the boundary 
path of $f_1$ and $s'w_2$ is the boundary path of $f_2$. We can have 
$p_1=p_2$ and $w_1$ or $w_2$ may also be empty.

By an \emph{elementary reduction of the second kind} we mean the inverse
process, converting $\mathfrak{F}'$ to $\mathfrak{F}$. Thus $\mathfrak{F}$
contains all the segments of $\mathfrak{F}'$ except $s'^{\pm 1}$, and all the
surface pieces except $f^{\pm 1}_1,f^{\pm 1}_2$, which are replaced by
$f^{\pm 1}$. The segment $s'$ appears exactly once in the boundary paths
of $f_1$ and $f_2$, and in no other boundary path of a surface piece in
$\mathfrak{F}'$.

Reductions and extensions which involve different elements of the original
complex can be interchanged with each other or thought of as occurring
simultaneously. E.g. if $s_1,s_2,\ldots,s_n$ are $n$ segments of a complex
$\mathfrak{F}$ and if $n$ new points $p_1,p_2,\ldots,p_n$ are
introduced successively so as to replace the segments in question by $2n$
segments $s_{1i},s_{2i}$ ($i=1,2\ldots,n$), converting $\mathfrak{F}$
to a complex $\mathfrak{F}'$ in which $p_i$ bounds just $s_{1i}$ and
$s_{2i}$, then these modifications can be performed in any order, and they
can be thought of simultaneously as a single modification of general type.
Such modifications, composed of interchangeable reductions and extensions
of a complex $\mathfrak{F}$, will be called \emph{elementary transformations}.
They can also be composed of infinitely many such modifications.

Two complexes $\mathfrak{F}$ and $\mathfrak{F}'$ are called elementarily
related if there is a chain of complexes
\[
\mathfrak{F}=\mathfrak{F}_1,\quad
\mathfrak{F}_2,\quad\ldots,\quad
\mathfrak{F}_r=\mathfrak{F}',
\]
beginning with $\mathfrak{F}$ and ending with $\mathfrak{F}'$ and in
which $\mathfrak{F}_{i+1}$ results from $\mathfrak{F}_i$ by an
elementary transformation. If $\mathfrak{F}$ is elementarily related to
$\mathfrak{F}'$ then so is $\mathfrak{F}'$ to $\mathfrak{F}$, and if
$\mathfrak{F}$ is elementarily related to $\mathfrak{F}'$, and
$\mathfrak{F}'$ to $\mathfrak{F}''$, then $\mathfrak{F}$ is also
elementarily related to $\mathfrak{F}''$.

Each complex $\mathfrak{F}$ is elementarily related to one with no 
singular elements. First, one can replace all singular segments by
non-singular ones by elementary extensions of the first kind. If
$w=w_1 w_2$ is a boundary path of $f$ in the resulting complex then
$w$ certainly runs through two different points $p_1$ and $p_2$, and
indeed we may take $w_1$ running from $p_1$ to $p_2$. We add a 
segment $s$ which goes from $p_1$ to $p_2$ and replace $f$ by
$f_1$ and $f_2$ with boundary paths $w_1 s^{-1}$ and $sw_2$
respectively. We replace $s$ by $s_1 s_2$ by adding a new point $q$.
Now
\[
w_1=s_{\alpha_1} s_{\alpha_2}\cdots s_{\alpha_k}
\]
runs successively through the points
\[
p_1,\quad
p_{\alpha_1},\quad
p_{\alpha_2},\quad\ldots,\quad
p_{\alpha_k}=p_2
\]
and
\[
w_2=s_{\alpha_{k+1}} s_{\alpha_{k+2}}\cdots s_{\alpha_m}
\]
through the points
\[
p_2=p_{\alpha_{k+1}},\quad\ldots,\quad
p_{\alpha_m}=p_1,
\]
so let
\[
t_{\alpha_1},\quad\ldots,\quad t_{\alpha_{k-1}},\quad
t_{\alpha_k}=s_2,\quad
t_{\alpha_{k+1}},\quad\ldots,\quad
t_{\alpha_{m-1}},\quad
t_{\alpha_m}=s^{-1}_1
\]
be $m$ different segments which begin at $q$ and end at $p_{\alpha_i}$.
Then $f_i$ is replaced by the $m$ triangles $f_{\alpha_i}$ with
boundary path $t_{\alpha_i} s_{\alpha_i} t^{-1}_{\alpha_{i+1}}$, so
the $f_{\alpha_i}$ are not singular. In this way we can eliminate each
singular surface piece.

\section{Elementary Relatedness of Manifolds}

If $\mathfrak{M}$ is a manifold and $\mathfrak{M}'$ is elementarily
related to $\mathfrak{M}$, then $\mathfrak{M}'$ is also a manifold.

If $\mathfrak{M}$ is a manifold of finitely many elements, and
$\mathfrak{D}$ is its dual, then $\mathfrak{M}$ is elementarily related
to $\mathfrak{D}$. This is because the passage from $\mathfrak{M}$ to
$\mathfrak{M}''$ may be accomplished by elementary extensions,
and that from $\mathfrak{M}''$ to $\mathfrak{D}$ by elementary
reductions.

For a manifold $\mathfrak{M}$ the following modifications may be
accomplished by elementary transformations: 

If $p_1$ and $p_2$ are two points of $\mathfrak{M}$ connected by 
a segment $s$, and if $f_1$ and $f_2$ are two distinct surface pieces
with simple boundaries $sw_1$ and $s^{-1}w_2$ containing the
segment $s$, then let $\mathfrak{M}^*$ be the complex which results 
from $\mathfrak{M}$ by omitting the segment $s^{\pm 1}$,
replacing $p_1$ and $p_2$ by a point $p^*$, and replacing the
segments $s_i$ in $\mathfrak{M}$ which begin at $p_1$ or $p_2$ by
segments $s^*_i$ which begin at $p^*$ and end at the same points.
The boundary paths in $\mathfrak{M}^*$ result from those of
$\mathfrak{M}$ by replacing segments $s_i$ by the corresponding
$s^*_i$ and eliminating occurrences of the segment $s$.
$\mathfrak{M}^*$ is then a manifold elementarily equivalent to
$\mathfrak{M}$. Namely, if $\mathfrak{D}$ and $\mathfrak{D}^*$ are
the manifolds dual to $\mathfrak{M}$ and $\mathfrak{M}^*$
respectively then $\mathfrak{D}^*$ results from $\mathfrak{D}$ by
an elementary reduction of the second kind. We will call these
transformations \emph{reductions of the third kind}.

They can always be carried out when a segment other than $s$
emanates from one of the two points $p_i$. The transformations of 
the inverse type, in which one point is replaced by two, may be simply 
called \emph{extensions of the third kind}.

If $f_1$ and $f_2$ are two distinct surface pieces of a manifold, if 
$w_1 s^{-1}$ is a boundary path of $f_1$, $sw_2$ is a boundary path
of $f_2$, and if either $w_1$ or $w_2$ is non-empty, then $f_1$ and
$f_2$ may be replaced by the surface piece $f$ with boundary path
$w_1 w_2$.

If the boundary path of $f_1$ and $f_2$ is equal to $s$, then
$\mathfrak{M}$ consists only of these two surface pieces and this segment,
together with their inverses and a single point $p$. This is so because
only one segment $s$ emanates from $p$, so $p$ cannot be connected
to any other point of the one-dimensional complex contained in 
$\mathfrak{M}$. But, since this complex is connected, it consists only of
$s,s^{-1}$ and $p$. It follows that:

\emph{Each manifold $\mathfrak{M}$ of finitely many elements is
elementarily related to a manifold $\mathfrak{M}'$ which contains only 
one surface piece pair, or else to a manifold which contains two
surface piece pairs and a singular segment, together with their inverses.}

Namely, suppose there was a manifold with $k>1$ surface pieces in which 
the surface piece $f_1$ cannot be combined with any other $f_2$ by an
elementary reduction. Then a boundary path of $f_1$ must run twice
through all segments it contains. But then the surface piece $f_1$ must
correspond to a point $\Delta(f_1)$ in the dual complex which bounds 
only singular segments. Thus if there were another surface piece 
$f_2\ne f_1,f^{-1}_1$ the dual complex would not be connected. It
then follows also:

\emph{A manifold $\mathfrak{M}$ of finitely many elements is
elementarily related to a manifold of two surface piece pairs, a line
segment pair, and a point; or to a manifold of a surface piece
pair, two points and two line segment pairs; or to a manifold which contains
only a surface piece pair $f^{\pm 1}$ and a point $p$.}

Namely, let $\mathfrak{M}'$ be a manifold elementarily equivalent to
$\mathfrak{M}$ which contains only one surface piece pair, and let
$\mathfrak{D}'$ be the manifold dual to $\mathfrak{M}'$. Then
$\mathfrak{D}'$ contains only one point $p'$, and furthermore
$\mathfrak{D}'$ may be converted, by elementary reductions of the
second kind, into a manifold $\mathfrak{D}^*$ which contains only
the point $p'$ and either a single surface piece or two surface pieces
and one line segment. Since $\mathfrak{M}$ and $\mathfrak{M}'$,
$\mathfrak{M}'$ and $\mathfrak{D}'$, $\mathfrak{D}'$ and
$\mathfrak{D}^*$ are elementarily related, so are $\mathfrak{M}$
and $\mathfrak{D}^*$.

The latter may be obtained by transformations of
manifolds that unify two points. Such transformations can
be carried out except when only a single segment emanates from each 
point. But then the manifold consists of only two points and a pair of
regular segments, since the line segment subcomplex is connected.
The dual of such a manifold consists of two surface piece pairs, a pair
of singular segments, and a point.

\section{Reduction to Normal Form}

One can generally bring a manifold $\mathfrak{M}$ into normal form
in various ways by elementary reductions. However, these possibilities 
are easy to survey. Let $\mathfrak{C}$ be the one-dimensional line
segment complex contained in $\mathfrak{M}$, $\mathfrak{C}'$ the
one-dimensional line segment complex manifold of the manifold
$\mathfrak{D}=\Delta(\mathfrak{M})$ dual to $\mathfrak{M}$,
$\mathfrak{B}$ any spanning tree of
$\mathfrak{C}$, and $\mathfrak{B}'$ such a tree in $\mathfrak{C}'$.
We call $\mathfrak{B}$ and $\mathfrak{B}'$ \emph{compatible} if
there is no segment of $\mathfrak{B}'$ which corresponds to a segment
of $\mathfrak{B}$ under the dual mapping. It is not claimed initially that
compatible trees $\mathfrak{B}$ and $\mathfrak{B}'$ must exist.
However, we show:

\emph{Each pair of compatible trees $\mathfrak{B},\mathfrak{B}'$
corresponds to a chain of reductions of the second and third kind, unique
up to the order of the reductions,
which brings the manifold into normal form. Conversely, each chain of 
such reductions corresponds to a pair of compatible trees $\mathfrak{B}$,
$\mathfrak{B}'$.}

Namely, if $s$ is a segment in $\mathfrak{B}$, then $s$ can be removed
by a reduction of the third kind when $\mathfrak{M}$ is not already in
normal form. For if $s$ is a regular segment and no other segments
emanate from its endpoints then $\mathfrak{M}$ is in fact in normal
form. Further, if $s^*$ is a segment corresponding to a segment
$s'$ of $\mathfrak{B}'$ under the dual mapping then $s^*$ may be 
removed by a reduction of the second kind, or else $\mathfrak{M}$ is
a normal form because two different surface pieces meet along $s^*$.
Now if $\overline{\mathfrak{M}}$ is a manifold resulting from
$\mathfrak{M}$ by a reduction of the third kind, let
$\overline{\mathfrak{B}}$ be the complex resulting from
$\mathfrak{B}$ by elimination of the segment $s$ and the point $p_2$,
with all the segments that emanated from $p_2$ now emanating from
$p_1$. $\overline{\mathfrak{B}'}$ is identical with $\mathfrak{B}'$.
Then $\overline{\mathfrak{B}}, \overline{\mathfrak{B}'}$ again
constitute a pair of compatible trees. One reaches an analogous
conclusion for reductions of the second kind.

Conversely, given a chain of reductions which bring $\mathfrak{M}$
into a normal form with one point and one surface piece, the collection of 
segments removed by reductions of the third kind form a spanning tree
$\mathfrak{B}$ of $\mathfrak{M}$, and the
segments corresponding to those removed by reductions of the second
kind in the dual complex form a tree $\mathfrak{B}'$, compatible with
$\mathfrak{B}$, which contains all points of the dual complex. This is
certainly true for those manifolds $\mathfrak{M}$ which may be brought 
into normal form by a single reduction or extension. Now if $\mathfrak{B}$
is any tree which contains all points of the manifold $\mathfrak{M}$,
and if $\overline{\mathfrak{M}}$ goes into $\mathfrak{M}$ by an
extension of the third kind, in which the regular segment $s$ is added
and the point $p$ is divided into the points $p_1$ and $p_2$, then the
complex $\overline{\mathfrak{B}}$ of all the points of
$\overline{\mathfrak{M}}$, the segments in $\mathfrak{B}$, and those
corresponding to them in $\overline{\mathfrak{M}}$, and the newly added
segment $s$, is again a tree. If $\mathfrak{B},\mathfrak{B}'$ were a
pair of compatible trees associated with the manifold $\mathfrak{M}$
via a chain of reductions to the normal form, then the trees
$\overline{\mathfrak{B}}$ and $\overline{\mathfrak{B}'}=\mathfrak{B}'$
are a pair of compatible trees of $\overline{\mathfrak{M}}$, and in fact
those associated with the chain of reductions of 
$\overline{\mathfrak{M}}$ to its normal form. One reaches this conclusion
analogously when $\overline{\mathfrak{B}}$ results from $\mathfrak{M}$
by an extension of the second kind.

If the normal form has two points and a surface piece pair, then the complex of
segments eliminated by reductions of the third kind consists of two
components, $\mathfrak{B}_1$ and $\mathfrak{B}_2$, which are both
trees. On the other hand, the line segment complex of the dual manifold
$\mathfrak{B}'$ is a tree. We can construct a  tree compatible with 
$\mathfrak{B}'$ from $\mathfrak{B}_1$ and $\mathfrak{B}_2$ by
adding one of the segments contained in the normal form to
$\mathfrak{B}_1$ and $\mathfrak{B}_2$.

\section{Neighboring Normal Forms}

In order to survey all classes of elementarily equivalent manifolds of 
finitely many elements we need only classify the manifolds 
containing a single point and a surface piece pair. We show that each such
manifold may be brought into one of the normal forms given in
Section 5.11. This takes place through a series of steps we first
describe as follows:

In $\mathfrak{M}$ is a manifold and
\[
w=w_1 s w_2 w_3 s w_4,\quad\text{respectively,}\quad
w=w_1 s w_2 w_3 s^{-1} w_4
\]
is a boundary path of a surface piece $f$, let $\mathfrak{M}'$ be the
manifold consisting of a surface piece $f'$ with boundary path
\begin{equation}
w'=s^* w_3 w^{-1}_1 s^* w^{-1}_2 w_4, \tag{1}
\end{equation}
respectively,
\begin{equation}
w'=s^* w_3 w^{-1}_1 {s^*}^{-1} w^{-1}_2 w_4, \tag{2}
\end{equation}
a manifold in which the segment $s$ is replaced by $s^*$.
$\mathfrak{M}'$ is elementarily related to $\mathfrak{M}$, and
$\mathfrak{M}'$ \emph{is called a neighbor of $\mathfrak{M}$.}

Namely, if
\[
w=w_1 s w_2 w_3 s w_4,
\]
let $\mathfrak{M}''$ be the manifold of two surface pieces $f_1$
and $f_2$ with respective boundary paths
\[
w_1 s w_2 {s^*}^{-1}\quad\text{and}\quad s^* w_3 s w_4.
\]
Then $s w_2  {s^*}^{-1} w_1$ or $w^{-1}_1 s^* w^{-1}_2 s^{-1}$ 
is also a boundary path of $f_2$. We now construct $\mathfrak{M}'$ in 
which we replace $f_1$ and $f_2$ by the surface piece $f'$ with boundary 
path
\[
w^{-1}_1  s^* w^{-1}_2 w_4 s^* w_3.
\]
The latter is a cyclic interchange of (1) and, since $\mathfrak{M}$ and
$\mathfrak{M}''$, $\mathfrak{M}''$ and $\mathfrak{M}'$ are
elementarily related, the assertion is proved in this case.

If
\[
w=w_1 s w_2 w_3 s^{-1} w_4,
\]
let $\mathfrak{M}''$ be the manifold of two surface pieces $f_1$ and $f_2$
with the respective boundary paths
\[
w_1 s w_2 {s^*}^{-1}\quad\text{and}\quad
s^* w_3 s^{-1} w_4.
\]
Then $w_2 {s^*}^{-1} w_1 s$ is also a boundary path of $f_1$ and
$s^{-1}w_4 s^* w_3$ is one of $f_2$. We now construct the manifold
$\mathfrak{M}'$ in which $f_1$ and $f_2$ are replaced by $f'$ with the
boundary path
\[
w_2 {s^*}^{-1} w_1 w_4 s^* w_3,
\]
which is a cyclic interchange of (2). The assertion now follows.

Any one of the modifications given can be replaced by a chain of simple 
steps in which one subdivides a surface piece in such a way that at least
one triangle is produced.

\section{Canonical Normal Forms}

We proceed to our goal in stages.\footnote{For this section see also:
\textsc{F. Levi}, Geometrischen Konfigurationen, Leipzig, 1929.}
\medskip

Theorem 1. \emph{If $w$ is a boundary path of $f$ and $s$ is a segment
traversed twice in the same direction by $w$, then there is a manifold
$\mathfrak{M}^*$, elementarily equivalent to $\mathfrak{M}$, with the
boundary path}
\[
w^*={s^*_1}^2 {s^*_2}^2 \cdots {s^*_m}^2 w^*_1.
\]
\emph{Here $w^*_1$ is empty or it traverses each remaining segment
once in each direction}: $s^2$ denotes the path $ss$.

If 
\[
w=s w_2  s w_4
\]
is a boundary path of $f$, then by Section 5.10 the manifold $\mathfrak{M}'$,
in which $f'$ has the boundary path
\[
{s^*_1}^2 w^{-1}_2 w_4,
\]
is a neighbor of $\mathfrak{M}$. Now let
\[
{s^*_1}^2 {s^*_2}^2 \cdots {s^*_i}^2 w_1 s_{i+1} w_2 s_{i+1} w_3
\]
be the boundary path of a manifold $\mathfrak{M}^*$
elementarily related to $\mathfrak{M}$. Then
\[
w_1 s_{i+1} w_2 s_{i+1} w_3{s^*_1}^2 {s^*_2}^2 \cdots {s^*_i}^2
\]
is also a boundary path and hence (with
$w_3{s^*_1}^2 {s^*_2}^2 \cdots {s^*_i}^2$ playing the
role of $w_4$ in Section 5.10)
\[
s^*_{i+1} w_1 w^{-1}_2 s^*_{i+1} w_3 {s^*_1}^2 {s^*_2}^2 
\cdots {s^*_i}^2
\]
is the boundary path of a manifold $\overline{\mathfrak{M}^*}$
neighboring $\mathfrak{M}^*$. Letting $w_2 w^{-1}_1$ play the role 
of $w_4$ in Section 5.10 we obtain 
\[
{s^*_{i+1}}^2 w_1 w^{-1}_2  w_3 {s^*_1}^2 {s^*_2}^2 
\cdots {s^*_i}^2,
\]
or,
\[
{s^*_1}^2 {s^*_2}^2 \cdots {s^*_i}^2 {s^*_{i+1}}^2 w'
\]
as the boundary path of a manifold elementarily related to
$\mathfrak{M}$. The assertion follows by induction.

To further normalize the boundary path we prove:
\medskip

Theorem 2. \emph{Let $w=w_1 w_2 w_3$ be a boundary path in a
manifold $\mathfrak{M}$, where $w_2$ runs through any segment in both
directions if it runs through it at all, and let}
\[
w_2=w_{21} k_1 w_{22} k_2 \cdots w_{2r} k_r w_{2,r+1},
\]
\emph{where}
\[
k_i=s_{i1} s_{i2} s^{-1}_{i1} s^{-1}_{i2}.
\]
\emph{Then there are two segment pairs, $s_1,s^{-1}_1$ 
and $s_2,s^{-1}_2$, among the $w_{2i}$ which mutually separate each
other or else the $w_{2i}$ are empty.} That is, with suitable
numbering $w_2$ runs first through $s_1$, then $s_2$, then $s_1$
again, then $s_2$ again.

We consider the subpaths of $w_2$ of the form $sw' s^{-1}$, where
$s$ is traversed by the $w_{2i}$. There is a certain shortest path of this
kind; let it be $s_1 w' s^{-1}_1$. Then $w'$ is certainly not empty,
for the $s$ are singular segments and by Section 5.3 a path of singular 
segments is reduced. But $w'$ also cannot consist only of the elements
$k_i$, because
\[
s k_i k_{i+1} \cdots k_{i+l} s^{-1}
\]
cannot appear in a boundary path by Section 5.3. If $w'$ passes through
$s_2$, where $s_2$ is also traversed by the $w_{2i}$, then $w'$ runs
through the segment $s_2$ only once. Otherwise, the subpath 
$s_2 w' s^{-1}_2$ of $w_2$ would be shorter than $s_1 w' s^{-1}_1$.
Thus the line segment pairs $s_1,s^{-1}_1$ and $s_2,s^{-1}_2$
separate each other in $w_2$. 
\medskip

Theorem 3. \emph{If the paths $w_{2i}$ in (1) are not empty then there
is a manifold neighboring $\mathfrak{M}$ with a boundary path}
\[
w'=w_1 w'_2 w_3,
\]
\emph{where $w'_2$ contains one more element of the form
$k=s_1 s_2 s^{-1}_1 s^{-1}_2$ then $w$.}

The path $w_2$ contains two line segment pairs,  $s_1,s^{-1}_1$ 
and $s_2,s^{-1}_2$, which mutually separate.

We therefore set
\[
w_2 = 
w'_{21} s_1 w'_{22} s_2 w'_{23} s^{-1}_1 w'_{24} s^{-1}_2 w'_{25}.
\]
Now we go to a neighboring manifold by replacing $s_2$ by $s^*_2$
and construct the new boundary path
\[
w_1 w'_{21} s_1 s^*_2 s^{-1}_1 w'_{24} w'_{23} {s^*_2}^{-1}
w'_{22} w'_{25} w_3,
\]
then go to a further neighboring manifold in which we replace $s_1$ by
$s^*_1$ and construct the boundary path
\[
w_1 w'_{21} w'_{23} w'_{24} s^*_1 s^*_2 {s^*_1}^{-1} {s^*_2}^{-1}
w'_{22} w'_{25} w_3.
\]
Here we have
\[
s^*_1 s^*_2 {s^*_1}^{-1} {s^*_2}^{-1}=k^*
\]
as a new element of the required form, while none of the paths $k_i$ 
in $w_2$ has been destroyed. It follows immediately from
the second and third theorems that:
\medskip

Theorem 4. \emph{If a boundary path $w$ of a manifold runs through
each segment in both directions, then there is an elementarily related
manifold having the boundary path}
\begin{equation}
w'=k_1 k_2 \cdots k_g\quad\text{\emph{with}}\quad
k_i=s_{1i} s_{2i} s^{-1}_{1i} s^{-1}_{2i}\quad
(i=1,2,\ldots,g). \tag{2}
\end{equation}

If a boundary path $w$ of a manifold runs through a segment twice in 
the same direction then there is an elementarily related manifold whose
boundary path is either 
\[
w'=s^2_1 s^2_2\cdots s^2_g
\]
or 
\begin{equation}
w'=s^2_1 s^2_2\cdots s^2_m k_1 k_2 \cdots k_n \tag{3}
\end{equation}
where
\[
k_i=s_{1i} s_{2i} s^{-1}_{1i} s^{-1}_{2i}\quad (i=1,2,\ldots,n).
\]

We show finally:
\medskip

Theorem 5. \emph{If a boundary path of a manifold runs through 
a segment twice in
the same direction, then there is an elementarily related manifold with
boundary path}
\[
w'=s^2_1 s^2_2 \cdots s^2_g.
\]

To prove this we show that the ``commutator'' $k_1$ in (3) may be
replaced by two ``squares'' $s^2$. Let
\[
w''=s^2_m s_{11} s_{12} s^{-1}_{11} s^{-1}_{12} w_1 
\]
be a boundary path which results from (3) by cyclic interchange.
We go over to a neighboring manifold in which $s_m$ is replaced 
by $s^*_m$ and take
\[
\overline{w}''=
s^*_m s^{-1}_{12} s^{-1}_{11} s^*_m s^{-1}_{11} s^{-1}_{12} w_1
\]
as boundary path; then replace $s_{11}$ by $s^*_{11}$ and take
\[
\overline{\overline{w}}''=
s^*_m s^{-1}_{12} {s^*_{11}}^{-1} {s^*_{11}}^{-1} 
{s^*_m}^{-1} s^{-1}_{12} w_1
\]
as boundary path, and finally replace $s_{12}$ by $s^*_{12}$ and take
\[
w^*=s^*_m s^*_m s^*_{11} s^*_{11} s^*_{12} s^*_{12} w_1
\]
as boundary path. The manifolds with the boundary paths (2) and (4),
together with the manifold of a surface piece pair, two points, and a
line segment pair, may be called ``canonical normal forms.''

\section{Normal Forms with Retention of a Segment}

To refine the results above we investigate what normal forms are obtainable
when a segment $s$ of the original boundary path is retained throughout,
while the remaining segments may be replaced as before.

Firstly, the form
\[
w=s^2_1 s^2_2 \cdots s^2_m w_1 s w_2 s^\varepsilon w_3\qquad
(\varepsilon=\pm 1)
\]
may be attained by introduction of suitable segments $s_1,s_2,\ldots,s_m$,
where $w_1,w_2,w_3$ together run through all their segments in both
directions. By introduction of
\[
s'_i=s_i w_1,\quad s''_i =w^{-1}_i s'_i
\qquad (i=1,2,\ldots,m)
\]
this may be replaced by
\[
w'={s''_1}^2 {s''_2}^2 \cdots {s''_m}^2 s w_2 s^\varepsilon w_3 w_1.
\]
Hence, with a new notation, we can start with the form
\begin{equation}
w=s^2_1 s^2_2 \cdots s^2_m s w_1 s^\varepsilon w_2
\qquad(\varepsilon=\pm 1). \tag{1}
\end{equation}

1. Now either $w_1$ already traverses one of its segments in both directions,
in which case the same holds for $w_2$, and by introduction of new
segments $w_1$ and $w_2$ may each be replaced by a path of
commutators. Then, by Section 5.3, $\varepsilon$ must be $+1$.
If $w_2$ begins with $s_\alpha s_\beta s^{-1}_\alpha s^{-1}_\beta$
then the subpath $s s_\alpha s_\beta s^{-1}_\alpha s^{-1}_\beta$ may
be replaced successively by
\begin{equation}
s'_\alpha s_\beta {s'_\alpha}^{-1} s s^{-1}_\beta,\quad
s'_\alpha s'_\beta s {s'_\alpha}^{-1} {s'_\beta}^{-1},\quad
s''_\alpha s s'_\beta  {s''_\alpha}^{-1} {s'_\beta}^{-1},\quad
s''_\alpha s''_\beta  {s''_\alpha}^{-1} {s''_\beta}^{-1} s
\tag{2}
\end{equation}

By applying such steps we convert $w$ into
\[
s^2_1 s^2_2 \cdots s^2_m w'_1 w'_2 s^2,
\]
and then, if $m\ne 0$, into a path of squares, as in Section 5.11:
\[
s^2_1 s^2_2 \cdots s^2_m  s^2.
\]

2. Otherwise, there is a segment $t$ which is traversed once by $w_1$
and once, in the opposite sense, by $w_2$. Then we can take
\[
w_1 = w_{11} t w_{12}=t'
\]
as a new segment and obtain
\begin{equation}
w = s^2_1 s^2_2 \cdots s^2_m st' s^\varepsilon w'_1 {t'}^{-1} w'_2,
\tag{3}
\end{equation}
where each commutator that appeared in $w_1$ and $w_2$ also appears
in $w'_1$ and $w'_2$. Now either $w'_1$ and $w'_2$ consist only of
commutators, or $w'_1 w'_2$ is empty, or there is a mutually separating
pair of segments in
\[
t' s^\varepsilon w'_1 {t'}^{-1} w'_2
\]
which do not belong to the commutators. But then, using the modifications
in Section 5.10, this segment pair may be replaced by a new one which is
a commutator. As a result,
\[
t' s^\varepsilon w'_1 {t'}^{-1} w_2\quad\text{becomes}\quad
w''_1 s^\varepsilon w''_2.
\]

By repeating these modifications we obtain the form (1) or (3) where 
$w_i,w'_i$ respectively consist only of commutators and, using the
deformations in (2), we can convert this into
\[
w = s^2_1 s^2_2 \cdots s^2_m w_3 st s^\varepsilon t^{-1},
\]
where $w_3$ consists only of commutators. Hence when $m\ne 0$ we
can again convert it into squares.

Further, we can also arrive at $\varepsilon=+1$ in the case $m=0$.
For
\[
st s^{-1} t^{-1} s^2_1 w'
\]
may  be converted into
\[
st  s'_1 ts s'_1 w'
\]
by introduction of $s^{-1}t^{-1}s_1=s'_1$. It follows that:

\emph{If the boundary path runs through each segment in both directions,
then the second normal form may be attained with retention of an
arbitrary segment.}

\emph{If the boundary path does not run through each segment in 
both directions, there are four cases to distinguish:}

Let the boundary be $w$ and the retained segment $s$.

1. The path $w$ runs through $s$ twice in the same direction
($\varepsilon=+1$). We have 
\[
\text{a)}\quad w'=k_1 k_2 \cdots k_r s^2
\qquad\text{or}\qquad
\text{b)}\quad w'=k_1 k_2 \cdots k_r st s^{-1}t^{-1}.
\]

2. There is a segment other than $s$ which $w$ runs through twice in the
same direction. We have
\[
\text{a)}\quad w'=s^2_1 s^2_2 \cdots s^2_n s^2
\qquad\text{or}\qquad
\text{b)}\quad w'=s^2_1 s^2_2 \cdots s^2_n st s^{-1}t^{-1}.
\]

We will see in Sections 6.6 and 6.7 that these cases cannot be reduced to
each other.

\section{Orientability. Characteristic}

In order to conclude the classification of manifolds with finitely many
elements, we must determine whether normal forms can be elementarily related
to each other. This is not the case. The proof is based on the invariance
of two properties of a manifold under elementary transformations: 
\emph{orientability} and  \emph{characteristic}.\footnote{Usually
called the (negative) \emph{Euler} characteristic today. (Translator's note.)}

Let $f^{\pm 1}_1,f^{\pm 1}_2,\ldots, f^{\pm 1}_n$ be the
different surface pieces of a manifold. Let $w^\varepsilon_i$ be a
be a positive boundary path of $f^\varepsilon_i$ ($\varepsilon=\pm 1$).
Now if the surface pieces
\[
f^{\varepsilon_1}_1,\quad 
f^{\varepsilon_2}_2,\quad\cdots,\quad
f^{\varepsilon_n}_n
\]
may be chosen in such a way that the boundary paths $w^{\varepsilon_i}_i$
traverse each segment $s$ ($\varepsilon_i=\pm 1$) once in the positive
sense and once in the negative sense, the manifold is called orientable.

If 
\[
f^{\varepsilon_1}_1,\quad 
f^{\varepsilon_2}_2,\quad\cdots,\quad
f^{\varepsilon_n}_n
\]
is a choice which satisfies our condition, then so too is
\[
f^{-\varepsilon_1}_1,\quad 
f^{-\varepsilon_2}_2,\quad\cdots,\quad
f^{-\varepsilon_n}_n.
\]
However, there is no choice apart from these two. Because, if the
exponent $\eta_1$ is chosen in $f^{\eta_1}_1$ then the exponent for
all surface pieces meeting it is already determined. Then one deduces
the assertion with the help of Section 5.6.

The \emph{invariance of orientability} under elementary transformations
is easy to see. If $\mathfrak{M}'$ results from $\mathfrak{M}$ by an
elementary transformation of the first kind then we retain the choice of surface 
pieces for $\mathfrak{M}'$. If $\mathfrak{M}'$ results from $\mathfrak{M}$
by an elementary transformation of the second kind, and if $f_i$ is
replaced by $f_{i1},f_{i2}$ as a result, let
\[
w_i=w_{i1} w_{i2}
\]
and let $w_{i1}{s'}^{-1}$ be a positive boundary path of $f_{i1}$
and $s' w_{i2}$ a positive boundary path of $f_{i2}$. We now replace
\[
f^{\varepsilon_i}_i\quad\text{by}\quad f^{\varepsilon_i}_{i1},
f^{\varepsilon_i}_{i1}
\]
in the choice of surface pieces and note that the segments appearing
in the $w_i$, as well as the new segment $s'$, are traversed as the
rule requires.

One sees similarly: if $\mathfrak{M}'$ is orientable, so is $\mathfrak{M}$.
Namely, if $f^{\varepsilon_{i1}}_{i1},f^{\varepsilon_{i2}}_{i2}$
are two surface pieces meeting along $s'$ which appear in a choice of
surface pieces orienting $\mathfrak{M}'$, then let $w_{i1}{s'}^{-1}$
be a boundary path of $f^{\varepsilon_{i1}}_{i1}$ and let
$s' w_{i2}$ be a boundary path of $f^{\varepsilon_{i2}}_{i2}$.
By hypothesis these boundary paths run through $s'$ in opposite
directions. We now replace
\[
f^{\varepsilon_{i1}}_{i1},f^{\varepsilon_{i2}}_{i2}\quad
\text{by}\quad
f^{\varepsilon_i}_i,
\]
with the boundary path $w_{i1} w_{i2}$.

Conversely, it follows immediately that if $\mathfrak{M}$ is a
non-orientable manifold and $\mathfrak{M}'$ results from
$\mathfrak{M}$ by an elementary transformation then
$\mathfrak{M}'$ is also non-orientable. For if $\mathfrak{M}'$
were orientable, $\mathfrak{M}$ would be too.

As soon as a cycle of segments beginning at a point contains more than
two segments one can distinguish between a positive cycle 
$s_{\alpha_1},s_{\alpha_2},\ldots,s_{\alpha_r}$
and a negative one
$s_{\alpha_r},s_{\alpha_{r-1}},\ldots,s_{\alpha_1}$. For orientable
manifolds one can also determine positive paths in such a way that the
following holds: if
\[
\ldots,s_{\alpha_1},s,s_{\alpha_2},\ldots\quad\text{and}\quad
\ldots,s_{\beta_2},s^{-1},s_{\beta_1},\ldots
\]
are the positive cycles in which the segment $s$ appears and if
$sw_1,sw_2$ are the two simple boundary paths beginning with $s$,
then with suitable numbering $w_i$ begins with $s^{-1}_{\alpha_i}$
and ends with $s_{\beta_i}$. If positive cycles cannot be defined in
this way then the surface is not orientable. However, if they can be,
it does not follow that the surface must be orientable, as one sees
from the surface with one cycle $s,s^{-1}$ and the boundary path
$s^2$.

By the \emph{characteristic} of a surface complex $\mathfrak{F}$
of finitely elements---$a_0$ points, $2a_1$ line segments, and $2a_2$
surface pieces---we mean the number
\[
c=-a_0+a_1-a_2.
\]

If $\mathfrak{F}'$ is a complex resulting from $\mathfrak{F}$
by an elementary extension of the first kind, then the number of
points $a'_0$ equals $a_0+1$, the number $2a'_1$ of segments
equals $2a_1+2$, the number $2a'_2$ of surface pieces equals
$2a_2$ and hence
\[
-a_0+a_1-a_2 = -a'_0+a'_1-a'_2.
\]
One concludes similarly for elementary reductions of the first kind.
If $\mathfrak{F}'$ results from $\mathfrak{F}$ by an elementary
extension of the second kind, then the numbers $a'_0,a'_1,a'_2$
for points, lines, and surfaces are
\[
a'_0=a_0,\quad
2a'_1=2a_1+2,\quad
2a'_2=2a_2+2,
\]
so the invariance of $c$ again follows. One concludes similarly for
elementary reductions of the second kind.

If we apply these two concepts to the normal forms of manifolds
we see that the manifolds of Section 5.11 (2) are orientable, and the
manifolds of Section 5.11 (4) are not. In the manifold of one surface
piece, one line segment, and two points the chacteristic has value $-2$
and in the manifolds of Section 5.11 (2) the value $2g-2$ ($g$ is
called the \emph{genus} of these manifolds),  and in the manifolds
of Section 5.11 (4) the value $g$. Thus different normal forms are
not elementarily related to each other. Further:

\emph{Two manifolds are elementarily related if they are both
orientable or both non-orientable, and have the same characteristic.}

The manifolds of Section 5.11 (2) can be visualized as spheres with
$g$ handles, the manifolds of Section 5.11 (4) as spheres with $g$
projective planes inserted.
\chapter{Groups and Surface Complexes}

\section{The Fundamental Group of a Surface Complex}

We now address ourselves to the connection between groups and surface
complexes which is realized on the one hand by paths in the complexes,
and on the other by coverings of complexes, quite analogously as with line
segment complexes. Let $\mathfrak{F}$ be a connected surface complex
and let $\mathfrak{C}$ be the connected line segment complex contained
in it.

Let $\mathfrak{W}_{\mathfrak{C}}$ be the group of closed paths in
$\mathfrak{C}$ that originate at $p_0$. With the help of boundary paths
in $\mathfrak{F}$ we now define an invariant subgroup $\mathfrak{R}$
of $\mathfrak{W}_{\mathfrak{C}}$ and a system of defining relations of
a group 
$\mathfrak{W}_{\mathfrak{C}}/\mathfrak{R}=\mathfrak{W}_\mathfrak{F}$, 
which we will call
the group of closed paths of $\mathfrak{F}$ beginning at $p_0$.

Let $\mathfrak{B}$ be a tree of segments of $\mathfrak{F}$ which 
contains all the points of $\mathfrak{F}$ and let $S_1$, $S_2$, $\ldots$,
$S_n$ be the system of generators of the fundamental 
group\footnote{Recall from Section 4.6 that Reidemeister calls the
fundamental group the ``path group,'' and that he uses $\mathfrak{W}$
(fraktur W) to denote this group because the German word for path is
``Weg.'' (Translator's note.)} 
of 
$\mathfrak{C}$ with initial point $p_0$. Let $f_i$ be an arbitrary surface 
piece of $\mathfrak{F}$, and $w_i$ its boundary path that begins and
ends at the point $p_k$. Further, let $w'_k$ be the simple path from $p_0$
to $p_k$ in the tree $\mathfrak{B}$, so that $w'_k w_i w'^{-1}_k$ is a 
closed path beginning at $p_0$ and if
\[
[w'_k w_i w'^{-1}_k]=R_i(S_l)
\]
is the power product in the $S$ associated with this path, then $R_i$ is
called the \emph{defining relation associated with the surface piece $f_i$}.
The collection of relations obtained in this way is the system of
defining relations for $\mathfrak{W}_\mathfrak{F}$. The product $R_i$
is determined up to cyclic interchange and cancellation of formal inverses
by $f_i$, for the segments in $w_i$ alone determine the power product
$R_i$ since $w'_k$ occurs in $\mathfrak{B}$.

Our problem is to show that \emph{the group determined in this way is
independent of the choice of $\mathfrak{B}$}. So let $\mathfrak{B}'$ be
a tree neighboring $\mathfrak{B}$ in $\mathfrak{C}$, let
\[
S'_1,\; S'_2,\; \ldots,\, S'_n
\]
be the generator system associated with $\mathfrak{B}'$ for the fundamental
group of $\mathfrak{C}$ with initial point $p_0$, and let $R'_i(S'_k)$ be a
relation associated with the boundary path $w_i$ of the surface piece $f_i$.
But we know from Section 4.6 that the $S'_k$ may be expressed in terms of the
$S_k$. If we now replace $R'_i(S'_k)$ by a power product 
$\overline{R}_i(S_k)$ in the $S_k$, then $\overline{R}_i(S_k)$ results from
$R_i(S_k)$ by elementary manipulations in the free group of the $S$. Thus
all consequence relations of the $R'_i(S'_k)$ likewise result from consequence
relations of the $R_i(S_k)$. One derives the converse similarly. Thus the
invariant subgroup $\mathfrak{R}$ is independent of the choice of tree
$\mathfrak{B}$ and hence so is 
$\mathfrak{W}_\mathfrak{F}=\mathfrak{W}_\mathfrak{C}/\mathfrak{R}$.

\emph{The groups of paths beginning at $p_0$ and $p'_0$ are isomorphic}.
Namely, if we take the same tree to define the generators of both groups
and the same boundary paths $w_i$ to define the defining relations, then
the defining relations are exactly the same. Thus we speak of the fundamental
group of $\mathfrak{F}$ to cover all these groups.

\section{Invariance of the Fundamental Group under Elementary 
Transformations}

\emph{The fundamental groups of elementarily related surface complexes
are isomorphic}.

Let $\mathfrak{F}'$ be a surface complex resulting from $\mathfrak{F}'$
by an elementary extension of the first kind, in which the segment $s$ of
$\mathfrak{F}$ with initial point $p_i$ and final point $p_k$ is replaced by
the two segments $s'_1,s'_2$ with $p'$ the new point, at which $s'_1$ ends
and $s'_2$ begins. Now if $s$ is a nonsingular segment there is a tree
$\mathfrak{B}$ which contains $s$ and all the points of $\mathfrak{C}$.
If we replace the segment $s$ in $\mathfrak{B}$ by $s'_1,s'_2,p$ then 
$\mathfrak{B}$ yields a spanning tree $\mathfrak{B}'$ of
$\mathfrak{F}'$ and one sees that the definitions of generators and defining 
relations for the paths emanating from $p_0\ne p'$ in $\mathfrak{F}$ and
$\mathfrak{F}'$ are formally identical when we base them on $\mathfrak{B}$
and $\mathfrak{B}'$ respectively. If $s=s_1$ is a singular segment, so
$p_i=p_k$, and $\mathfrak{B}$ is any spanning tree of $\mathfrak{F}$, let 
$\mathfrak{B}'$ be the tree that consists of $\mathfrak{B}$
with the addition of $s'_1$ and $p'$. If $s_1,s_2,\ldots,s_n$ are the segments
of $\mathfrak{F}$ that do not belong to $\mathfrak{B}$, then
$s'_2,s_2,\ldots,s_n$ are the segments of $\mathfrak{F}'$ that do not belong to
$\mathfrak{B}'$, and the generators corresponding to the segments
$s_2,s_3,\ldots, s_n$ are identical for $\mathfrak{F}$ and $\mathfrak{F}'$.
If
\[
S_1=[w s_1 w^{-1}]
\]
is the generator corresponding to the segment $s_1$ in $\mathfrak{F}$ then
\[
S'_1=[w s'_1 s'_2 w^{-1}]
\]
is a boundary path that does not contain $s_1$, corresponding to the one that 
does not contain $s'_1$, and thus it remains unaltered. If a boundary path in
$\mathfrak{F}$ runs through the segment $s_1$ exactly $k$ times, then the 
corresponding path in $\mathfrak{F}'$ runs through the path $s'_1 s'_2$
exactly $k$ times, and the power product corresponding to this path in
$\mathfrak{F}'$ therefore results from that in the $S_k$ when $S_1$ is
replaced by $S'_1$.

Now let $\mathfrak{F}'$ be a surface complex resulting from $\mathfrak{F}$
by an elementary extension of the second kind, in which the segment $s'_0$,
which goes from $p_i$ to $p_k$, is inserted and the surface piece $f$ with
boundary path $w_1 w_2$ is replaced by two surface pieces with the boundary
paths $w_1 s'^{-1}_0, s'_0 w_2$. If $\mathfrak{B}$ is any spanning tree of 
$\mathfrak{F}$, then $\mathfrak{B}$ has the same property in relation to
$\mathfrak{F}'$. The generators $S_1,S_2,\ldots,S_n$ in $\mathfrak{F}$ are
also generators in $\mathfrak{F}'$, and in addition there is a generator $S'_0$
corresponding to the segment $s'_0$. The relations of $\mathfrak{F}$, except
that for the boundary path of $f$, remain the same. If $R(S_k)$ is the power
product corresponding to the boundary path $w_1 w_2$ of $f$, and if
$R'_1(S_k,S'_0), R'_2(S_k,S'_0)$ are those corresponding to the paths
$w_1 s'^{-1}_0,s'_0 w_2$ respectively, then $R'_1 R'_2$ is convertible into
$R(S_k)$ by elementary reductions in the free group of the $S'_0,S_k$.
Thus the fundamental groups associated with $\mathfrak{F}$ and
$\mathfrak{F}'$ are isomorphic by the theorem of \textsc{Tietze} in
Section 2.10.

A special isomorphism $\boldsymbol{I}$ between the two fundamental groups
$\mathfrak{W}_\mathfrak{F}$ and $\mathfrak{W}_{\mathfrak{F}'}$ with
respect to a basepoint $p$ common to both complexes $\mathfrak{F}$ and
$\mathfrak{F}'$ is obtained by the following argument. If $w$ is a path that
appears in $\mathfrak{F}'$ as well as in $\mathfrak{F}$, let
\begin{equation}
\boldsymbol{I}([w])=[w]'. \tag{1}
\end{equation}

If $\mathfrak{F}$ and $\mathfrak{F}^{(n)}$ are two complexes convertible
into each other by elementary transformations, and if
\[
\mathfrak{F}^{(1)}=\mathfrak{F},\quad
\mathfrak{F}^{(2)},\quad
\ldots,\quad
\mathfrak{F}^{(n)}
\]
is a chain of complexes in which any two in succession are related by an
elementary extension or reduction, and if the point $p$ appears in all the
complexes $\mathfrak{F}^{(i)}$, and $\boldsymbol{I}_k$ is the 
isomorphism established by (1) between the fundamental groups with
basepoint $p$ in $\mathfrak{F}^{(k)}$ and $\mathfrak{F}^{(k+1)}$,
\[
\boldsymbol{I}_k([w]^{(k)})=[w]^{(k+1)},
\]
then
\[
\boldsymbol{I}([w])=\boldsymbol{I}_{n-1}
                               \left(\boldsymbol{I}_{n-2}\cdots
                                \boldsymbol{I}_1([w]^{(1)})\right)
                               =[w]^{(n)}
\]
is called \emph{the isomorphism effected by the chain} $\mathfrak{F}^{(i)}$
between the fundamental groups $\mathfrak{W}_{\mathfrak{F}^{(1)}}$ and
$\mathfrak{W}_{\mathfrak{F}^{(n)}}$. If there are two different chains that
connect $\mathfrak{F}^{(1)}$ and $\mathfrak{F}^{(n)}$, and if they effect
isomorphisms $\boldsymbol{I}$ and $\boldsymbol{I}'$, then
$\boldsymbol{I}\boldsymbol{I}'^{-1}$ is an automorphism of the
fundamental group $\mathfrak{W}_{\mathfrak{F}^{(1)}}$.

\section{Homotopy and Homology}

With the help of the fundamental group, we can classify all the closed paths
of a complex independently of the choice of basepoint. Namely, if $w$ is a closed
path beginning at $p_0$, $p_1$ is any other point, and $w_1$ is a path from
$p_1$ to $p_0$, then $w_1 w w^{-1}_1$ is a closed path beginning at $p_1$
and corresponding to the element $[w_1 w w^{-1}_1]$ in the fundamental
group with basepoint $p_1$. If $w'_1$ is any other path from $p_1$ to $p_0$,
then $[w'_1 w w'^{-1}_1]$ is a transform of $[w_1 w w^{-1}_1]$, namely
\begin{equation}
[w'_1 w w'^{-1}_1]=[w'_1 w^{-1}_1][w_1 w w^{-1}_1][w'_1w^{-1}_1]^{-1}.
\tag{1}
\end{equation}
Thus a closed path $w$ corresponds to a class of transformed elements in the
fundamental group that may be denoted by $\{w\}$. Naturally, different
closed paths can belong to the same class of transformed 
elements.\footnote{This class is what we today call a \emph{conjugacy class}.
In fact, Reidemeister calls them conjugate elements in the next paragraph
but one.
(Translator's note.)} If $w$ and $w^{*}$ are two such elements, then
\begin{equation}
\{w\}=\{w^{*}\} \tag{2}
\end{equation}
in the fundamental group with basepoint $p_1$, and this also holds for
any other basepoint $p_2$. For the fundamental groups are in fact 
isomorphically related to each other by the paths. Thus the relation
$\{w\}=\{w^{*}\} $ depends only on the paths $w$ and $w^{*}$
themselves. We call paths that satisfy (1) or (2) \emph{homotopic}
to each other. This relation is transitive. Obviously, $w$ is
homotopic to any path resulting from itself by cyclic interchange.
Further, if $w=w_1 w_2 w_3$ and if $w_2 w'^{-1}_2$ is a simple
boundary path of a surface piece, then $w$ is homotopic to
$w_1 w'_2 w_3$. Conversely, one can also use these two theorems to
define homotopy.

A path $w$ is called null homotopic when $\{w\}$ is the identity. Thus
the problem of deciding whether a path is null homotopic is the same as
the word problem, while deciding whether two paths are homotopic is the
same as the transformation problem, for the fundamental group.
(Cf. Sections 1.14 and 1.15.)

Corresponding to a class $\{w\}$ of conjugate elements there is a 
well-defined element of the factor group $\mathfrak{W}/\mathfrak{K}$ 
of the fundamental group by the commutator group, which may be
denoted by $\langle w \rangle$. Two paths $w$ and $w^{*}$ are 
called \emph{homologous} 
to each other when $\langle w \rangle=\langle w^{*} \rangle$.
Homology is transitive. We say that a path runs through the segment
$s_i$ exactly $k=m-n$ times if it runs through $m$ times in the positive 
sense and $n$ times in the negative sense, so two paths are certainly
homologous if they traverse the same segments equally often. A path is
called null homologous when $\langle w \rangle$ is the identity element,
i.e. when $\{ w\}$ consists of elements of the commutator group. By
Section 2.13, homology of curves is always decidable.

\section{Simple Paths on Manifolds}

The concepts of homotopy and homology have special interest in the case
of simple closed paths on manifolds. Namely, the homotopy and homology
properties of such paths depend on whether and how the manifold may
be decomposed. We say that \emph{the path $w$ separates the manifold}
$\mathfrak{M}$ when each tree $\mathfrak{B}'$ in the dual complex
$\mathfrak{D}$ contains at least one segment corresponding to a
segment of $w$ under the dual mapping. If on the other hand there is
such a tree containing no such segment then we say that $w$ does not
separate the manifold.

Now let $w$ be a simple path of the manifold $\mathfrak{M}$, let
$\mathfrak{B}_1$ be a tree consisting of segments of $w$ that includes
all points traversed by $w$, and let $\mathfrak{B}$ be a tree that contains
the segments of $\mathfrak{B}_1$ and all the points of $\mathfrak{M}$.
If one removes the segments of $\mathfrak{B}$ by reductions of the
third kind, then $\mathfrak{M}$ is converted into an equivalent manifold
$\mathfrak{M}^{*}$
containing only a single point, in which $w$ corresponds to a single,
and of course singular, segment $s^{*}$.

Now either $w$ does not separate the manifold $\mathfrak{M}$ and
hence $s^{*}$ does not separate the manifold $\mathfrak{M}^{*}$,
in which case in the manifold dual to $\mathfrak{M}^{*}$ there is a
tree ${\mathfrak{B}^{*}}'$ containing all points but not the segment
corresponding to $s^{*}$, and so $\mathfrak{M}^{*}$ may be
converted into an equivalent manifold of a single point and a surface
piece with $s^{*}$ in its boundary.

Or else $w$ separates $\mathfrak{M}$ and $s^{*}$ separates
$\mathfrak{M}^{*}$ and no such tree ${\mathfrak{B}^{*}}'$ exists.
Then $\mathfrak{M}^{*}$ may be converted into a manifold
$\mathfrak{M}^{**}$ of two surface pieces $f^{**}_1,f^{**}_2$
which meet only along $s^{*}$, for otherwise the two surface pieces
could be replaced by a single one without removing $s^{*}$ and thus
there would be a tree ${\mathfrak{B}^{*}}'$ in the manifold dual to
$\mathfrak{M}^{*}$ containing all points but not the segment
corresponding to $s^{*}$. If we denote two boundary paths in
$\mathfrak{M}^{**}$ beginning with $s^{*}$ by
${s^{*}}^{-1}r^{**}_1, {s^{*}}^{-1}r^{**}_2$ then, using the
manipulations described in Section 5.9, the $r^{**}_i$ may be
converted into one of the normal forms 
$s^{2}_1 s^{2}_2 \cdots s^{2}_l$
or $k_1 k_2 \cdots k_l$ where 
$k_i=s_{i1} s_{i2} s^{-1}_{i1} s^{-1}_{i2}$;
one of these paths may also be empty.

A few simple theorems follow from this.

\emph{If $\mathfrak{M}$ is an orientable manifold and $w$ is a simple 
path that separates $\mathfrak{M}$, then $w$ is null homologous and
conversely: if $w$ is a simple path and null homologous then $w$
separates the manifold.}

\emph{If $\mathfrak{M}$ is not orientable and $w$ is a simple path
that separates $\mathfrak{M}$, then either $w$ is null homologous or
there is path $w'$ for which}
\[
\langle w \rangle = \langle w' w' \rangle = \langle w' \rangle^2
\]
\emph{and conversely: if $w$ is a simple path and either null
homologous or homologous to a twice-traversed path $w' w'$ then
$w$ separates the manifold $\mathfrak{M}$}.

For when $w$ separates we have
\[
\langle w \rangle = \langle s^{*} \rangle = \langle r^{**}_1 \rangle
\]
and when $w$ does not separate we have
\[
\langle w \rangle = \langle s^{*} \rangle,
\]
where $s^{*}$ is a segment that appears in an equivalent normal form.
Hence $\langle s^{*} \rangle$ is neither the identity element nor the
square of another group element.

We will show later that the elements $[r^{**}_1]$ and $[r^{**}_2]$ 
can only be the identity when one of these paths, say $r^{**}_1$, is
empty. It then follows that

\emph{A simple path is null homotopic if and only if it may be converted
by reduction of the manifold into a singular segment which is the complete
boundary of a surface piece}.

\section{Intersection Numbers}

The membership of paths in well-defined homology classes has an 
important geometric
consequence: the appearance of intersection 
points.\footnote{\textsc{H. Poincar\'e}, Rendic. d. Palermo \textbf{13}, 
314, (1899).} In this connection we will confine ourselves to orientable 
surfaces. Then we can associate an index with each intersection point
$p$ of two paths $w_1$ and $w_2$ by the following rule: let
\[
w_i=w_{i1}s^{-1}_{i1}s_{i2} w_{i2}
\]
where the $s_{ik}$ are the segments beginning at the intersection $p$.
When the pairs of segments $s_{i1},s_{i2}$ ($i=1,2$) do not mutually 
separate in the cycle of segments beginning at $p$, and thus the paths
$w_1$ and $w_2$ do not cross at $p$, then $p$ receives the \emph{index
zero}; otherwise \emph{index $+1$ or $-1$} according as
$s_{12}, s_{22}, s_{11}, s_{21}$ appear in that order in the positive or
negative cycle of segments beginning at $p$. We define the index for a
common initial or final point of closed paths analogously as for
intersection points.

Now if $w_1,w_2$ are \emph{closed paths} that run through no common
segments and pass each intersection point only once, the sum of the
indices of all those points is called the \emph{intersection number}
\[
N(w_1,w_2)
\]
of the two paths $w_1$ and $w_2$.

Obviously
\[
N(w_1,w_2)=-N(w_2,w_1),\quad
N(w^{-1}_1,w_2)=-N(w_1,w_2).
\]
One easily establishes the following additional property of the
intersection number: if $w_1=w_{11} w_{12}$, where the $w_{1i}$ are
themselves closed paths, then
\[
N(w_1,w_2)=N(w_{11},w_2)+N(w_{12},w_2).
\]
If $w'_1$ is a path homotopic to $w_1$ which has no segment in common
with $w_2$ then
\[
N(w'_1,w_2)=N(w_1,w_2).
\]

This is proved most conveniently by subdividing the manifold into triangles
and following the effect of deformation over a triangle.

In particular, if $w_1$ and $w'_1$ are two closed paths which begin at the 
same point $p$, and if $[w'_1]=[w_1]$ in the fundamental group with
basepoint $p$, then it is permissible to speak of the intersection of the
classes $[w_1]$ and $[w_2]$, since it is consistent to set
\[
N([w_1],[w_2])=N(w_1,w_2)
\]
and we again have
\[
N([w_{11}][w_{12}],[w_2])=N([w_{11}],[w_2])+N([w_{12}],[w_2]).
\]
If we now fix $[w_2]$ we see that the intersection numbers relative to
$[w_2]$ form a commutative group homomorphic to the fundamental
group; a path $k$ for which $[k]$ belongs to the commutator group
always has the intersection number $N(k,w_2)=0$, or, expressed in
another way: \emph{the intersection number $N(w_1,w_2)$ depends
only on the homology class $\langle w_1 \rangle$}. We can therefore
define a function $N(\langle w_1 \rangle, w_2)$. And because
\[
N(w_1,w_2)=-N(w_2,w_1)
\]
the same holds for the second argument. In order to complete the 
determination of
\[
N(\langle w_1 \rangle, \langle w_2 \rangle)
\]
we set up a system of generators that corresponds to the canonical form
of the manifold from Section 5.11 (2), and note that
\[
N(S_i,S_k)=N(T_i,T_k)=N(S_i,T_k)=0\qquad (i \ne k)
\]
and, with suitable ordering of generators,
\[
N(S_i,T_i)=1.
\]

Moreover, since for two paths $w_k$ of the homology classes
$\prod_i S^{a_{ki}}_i T^{b_{ki}}_i$ ($k=1,2$) we have
\begin{align*}
N(w_1,w_2)&=N\left(\prod_i S^{a_{ki}}_i T^{b_{ki}}_i,w_2\right)\\
&=\sum a_{1i}N(S_i,w_2)+\sum b_{1i} N(T_i,w_2)\\
&=\sum a_{1i} a_{2k} N(S_i,S_k)+\sum a_{1i} b_{2k} N(S_i,T_k)\\
&+\sum b_{1i} a_{2k} N(T_i,S_k)+\sum b_{1i}b_{2k} N(T_i,T_k),
\end{align*}
$N(w_1,w_2)$ yields the bilinear form
\[
\sum (a_{1i}b_{2i}-a_{2i}b_{1i}).
\]
If $S'_i$ and $T'_i$ are another system of generators for the fundamental
group, likewise associated with the paths of a second normal form, and if
$\prod_i S'^{a'_{ki}}_i T'^{b'_{ki}}_i$ are the homology of the paths
$w_k$, then we must have
\[
\sum (a'_{1i}b'_{2i}-a'_{2i}b'_{1i}) = \sum (a_{1i}b_{2i}-a_{2i}b_{1i}).
\]

Finally, here is one more application of the formula: if $S_i,T_i, S'_i, T'_i$
retain the same meaning in the factor group by the commutator group,
\[
S'_i=\prod S^{\alpha_{ik}}_k \prod T^{\beta_{ik}}_k,\quad
T'_i=\prod S^{\gamma_{ik}}_k \prod T^{\delta_{ik}}_k,
\]
then
\[
N(S'_i,T_i)=\alpha_{ii}\quad\text{and}\quad
N(S_i,T'_i)=\delta_{ii},
\]
and hence
\[
s=\sum \alpha_{ii} + \sum \delta_{ii},
\]
the trace of the matrix by which the passage from $S_i,T_i$ to $S'_i,T'_i$
takes place in the factor group, is equal to the intersection number
\[
\sum N(S'_i,T_i)+\sum N(S_i,T'_i).
\]

The investigation of the intersection points of closed paths may be greatly
refined by dividing the intersection points into classes which depend only 
on the homotopy classes of the paths.\footnote{\textsc{R. Baer}, J. f\"ur
reine und angew. Math. \textbf{156} (1927) 231.}

\section{One-sided and Two-sided Paths}

In order to describe the properties of simple paths somewhat more precisely
we introduce the concepts of ``one-sided'' and ``two-sided'' paths.

Let $w$ be any simple closed path of the manifold $\mathfrak{M}$, which passes through the points $p_0,p_1,\ldots,p_n=p_0$ and the segments
$s_1,s_2,\ldots,s_n$ in that order. New segments
\[ 
s_{i1},\quad s_{i2} \qquad (i=1,2,\ldots,n)
\]
may now be introduced and $\mathfrak{M}$ thereby extended
to the manifolds $\mathfrak{M}_{i1},\mathfrak{M}_{i2}$ in succession. 
The segments
$s_{i1},s_{i2}$ begin at $p_{i-1}$ and end at $p_i$, and if $s_i w_{i1}$ and 
if $s_i w_{i2}$ are the two boundary paths in $\mathfrak{M}_{i-1,2}$
beginning with $s_i$ then new surface pieces may always be introduced with
the boundary paths $s_i s^{-1}_{i1}$ and $s_i s^{-1}_{i2}$. It is
assumed that the paths $w_{i1},w_{i2}$ do not run through the segment
$s_i$. The path  $w$ is embedded in 2-gons in the new manifold
$\mathfrak{M}_{n,2}$. The distribution of the second indices 1, 2 in the
new manifold can be arranged in such a way that the ordering
\[
s_{i1}, s_i, s_{i2}, \ldots, s^{-1}_{i+1,2}, s^{-1}_{i+1}, s^{-1}_{i+1,1},
\ldots \quad (i=1,2,\ldots,n)
\]
is obtained in the stars of the points $p_i$. But then the star of the point 
$p_0$ is completely determined and either has the ordering
\[
s_{n1},s_n, s_{n2},\ldots,s^{-1}_{12},s^{-1}_1, s^{-1}_{11},\ldots
\]
or else the ordering
\[
s_{n2},s_n,s_{n1},\ldots,s^{-1}_{12},s^{-1}_1, s^{-1}_{11},\ldots .
\]

In the first case we call the path \emph{two-sided}, in the second case
\emph{one-sided}.

Because if we construct the path
\[
w_k=s_{1k} s_{2k}\cdots s_{nk}\qquad (k=1,2)
\]
then in the first case we always remain on the same side of $w$, intuitively
speaking, while in the second case we go from one side of $w$ to the other.

\emph{If a manifold possesses a one-sided path, then it is not orientable}.
For otherwise the stars of the points $p_0,p_1,\ldots, p_n$ could not be
oriented as we did in Section 5.13. \emph{On the other hand, one-sided
and two-sided paths may always be given in nonorientable manifolds}.

We show that a singular segment that appears in a boundary path of the
form $s^2 w$ is one-sided, whereas a segment that appears in the form
$stst^{-1}$ is two-sided. Namely, we construct new manifolds by
insertion of the segment $s'$ and in the first case take $s' s^{-1}$ and 
$s' s w$ as new boundary paths; in the second case take $s's^{-1}$
and $sts't^{-1}$ as new boundary paths. Then in the first case the cycle
of segments beginning at the point of the normal form has the substar
$s^{-1},s'^{-1},s$, so $s'$ goes from one side of $s$ to the other; in
the second case it has the substar $s^{-1},s'^{-1},t^{-1},s',s$, so $s'$
remains on the same side of $s$. Thus we have obtained a geometric
interpretation of the appearance of a segment in boundary paths of the
two forms $s^2 w$ and $stst^{-1}w$, and \emph{it is shown that the
forms} 1a \emph{and} 1b, \emph{as well as the forms} 2a \emph{and}
2b, \emph{in Section} 5.12 \emph{are not reducible to each other}.

One can define one-sided and two-sided quite analogously for non-simple
paths and show that homotopic curves are always of the same type.
Separating curves are obviously two-sided.

\section{Simple Strips}

The manifolds with the paths $w,w_1,w_2$, constructed in the previous 
section, may be subjected to the following extensions of the third kind:
in place of each point $p_i$ we take three points, with new segments
$t_{i1}$ connecting the boundary points $p_{i1}$ and $p_i$, and $t_{i2}$
connecting the boundary points  $p_{i2}$ and $p_i$. If
\[
s^{-1}_{i1},s^{-1}_i, s^{-1}_{i2},s_{\alpha_{i1}},s_{\alpha_{i2}},\ldots,
s_{\alpha_{ir}},s_{i+1,2},s_{i+1},s_{i+1,1}, s_{\beta_{i1}},
s_{\beta_{i2}},\ldots,s_{\beta_{il}}
\]
is the star of segments through $p_i$, then
\begin{align*}
&s^{-1}_{i1},s_{\beta_{il}},\ldots,s_{\beta_{i1}},s_{i+1,1},t_{i1}
\quad\text{becomes the star of }p_{i1},\\
&s^{-1}_{i2},s_{\alpha_{i1}},\ldots,s_{\alpha_{ir}},s_{i+1,2},t_{i2}
\quad\text{becomes the star of }p_{i2},\\
&s^{-1}_i, t^{-1}_{i1}, s_{i+1}, t_{i2}
\quad\text{becomes the star of }p_i.
\end{align*}
This converts $w_1$ and $w_2$ into two simple paths which terminate in
the case of a two-sided path $w$ and which result in a single simple 
closed path $w_1 w_2$ in the case of a one-side path $w$. The 2-gons
around $p_i$ have been converted into 4-gons. The complex of these
quadrilaterals and their boundary elements is a manifold with boundary
which may be called a \emph{strip}. In the case of a one-sided path we
get a strip convertible into a manifold which is a projective plane by
insertion of  a surface piece with the boundary path $w_1 w_2$. The strip
is a \textsc{M\"obius} band. In the case of a two-sided path we get a strip
convertible into a sphere by the addition of two surface pieces with the
boundary paths $w_i$---a cylindrical band.

Conversely, if we remove from $\mathfrak{M}$ the points $p_0,\ldots,
p_{n-1}$, the segments $s_1,\ldots,s_n$, $t_{i1}$, $t_{i2}$, and the 
quadrilaterals
of the strip then, in the case where $w$ does not separate $\mathfrak{M}$,
$\mathfrak{M}$ is converted into a manifold with boundary, 
$\mathfrak{R}$, which can be converted into a manifold $\mathfrak{M}'$ 
without boundary by the addition of one or two surface pieces with the
boundary paths $w_1 w_2$ and $w_i$ respectively. Otherwise, i.e. when
$w$ separates $\mathfrak{M}$, $\mathfrak{M}$ breaks into two manifolds
with boundary, $\mathfrak{R}_1$ and $\mathfrak{R}_2$.

In the case of non-separating cuts $w$ on non-orientable manifolds one
can further distinguish between \emph{two types of non-separating cuts},
according as \emph{the manifold resulting from $\mathfrak{M}$ from the
cut along $w$ is orientable or not}. One is easily convinced that this
property enables us to distinguish between simple paths which may be
embedded in a boundary path of the types 1a, b of Section 5.12, and
those for which this is not possible and which consequently reduce to the
types 2a, b of Section 5.12. For when one cuts a manifold with a boundary
path of type 1a,b along $s$ the result is an orientable manifold, whereas
a boundary path of type 2a,b  results in a non-orientable manifold.
Consequently, \emph{we have now shown that the four types of boundary 
path realized by non-orientable manifolds in Section 5.12 are not
reducible to each other}.

\section{Normal Forms and Fundamental Groups}

The preceding section showed that the group elements $[w]$ corresponding
to simple indecomposable curves $w$ constitute a special class. It is natural
to ask how all the elements of this class are determined. As Section 6.4
shows, this question is closely related to another. 
\emph{Let $w_1,w_2,\ldots, w_g$ be closed curves emanating from a point
$p$, which pass into $g$ singular segments of a normal polygon under
suitable reduction of the manifold. The $[w_i]$ ($i=1,2,\ldots,g$) are a special
$g$-tuple of group elements}, and the question is how to obtain them. They
certainly constitute a system of generators. We can formulate our question
precisely as: \emph{if $[w_i]$ and $[w'_i]$ are two such $g$-tuples, how
may the $[w'_i]$ be expressed in terms of the 
$[w_i]$?} \footnote{\textsc{K. Reidemeister}, J. f\"ur reine und angew. Math.
Henselschrift (1932).}

Now by Section 5.9 the reduction of a manifold to normal form is established
by giving two compatible trees $\mathfrak{B},\mathfrak{B}'$ which contain all
points of $\mathfrak{M}$ and its dual manifold respectively. So let
$\mathfrak{B},\mathfrak{B}'$ be a pair for which the $[w_i]$ go into the singular
segments of the normal form and let 
$\overline{\mathfrak{B}},\overline{\mathfrak{B}}'$ be a pair for which this
happens to the $[w'_i]$. The group elements $[w_i]^{\pm 1}$ and 
$[w'_i]^{\pm 1}$ are determined up to numbering these tree pairs, and for 
that reason we first investigate the meaning of our question with regard to
these tree pairs.

We can embed $\mathfrak{B}$ and $\overline{\mathfrak{B}}$ in a 
\emph{chain of neighboring trees} in $\mathfrak{M}$
\[
\mathfrak{B}=\mathfrak{B}_1,\mathfrak{B}_2,\ldots,
\mathfrak{B}_n=\overline{\mathfrak{B}}
\]
which contain all the points of $\mathfrak{M}$. For each $\mathfrak{B}_i$
there is therefore a compatible $\mathfrak{B}'_i$ which contains all the points of
$\mathfrak{M}'$. One also sees immediately that two trees
$\mathfrak{B}'_{i1}$ and $\mathfrak{B}'_{ik}$ compatible with $\mathfrak{B}_i$
may be embedded in a chain of compatible trees
$\mathfrak{B}'_{i2},\ldots,\mathfrak{B}'_{i,k-1}$. Finally, if
$\mathfrak{B}_i$ and $\mathfrak{B}_{i+1}$ are two neighboring trees then
there is either a tree $\mathfrak{B}'_i$ compatible with both $\mathfrak{B}_i$
and $\mathfrak{B}_{i+1}$, or else two neighboring trees $\mathfrak{B}'_i$
and $\mathfrak{B}'_{i+1}$ compatible with $\mathfrak{B}_i$ and
$\mathfrak{B}_{i+1}$ respectively. Namely, if $s_{i+1}$ is the segment
that appears in $\mathfrak{B}_{i+1}$ but not in $\mathfrak{B}_i$, and if
$w_i$ is a simple path which connects the boundary points of $s_{i+1}$ in 
$\mathfrak{B}_i$, and hence passes through the segment $s_i$ which
appears in $\mathfrak{B}_i$ but not in $\mathfrak{B}_{i+1}$, then either the 
path $s_{i+1}w_i$ separates $\mathfrak{M}$ or not. In the first case there is a 
tree compatible with $\mathfrak{B}'_i$ as well as with $\mathfrak{B}_i$. In the
second case each tree $\mathfrak{B}'_i$ compatible with $\mathfrak{B}_i$
contains the segment $s'_{i+1}$ corresponding to $s_{i+1}$. Now if one 
replaces $s'_{i+1}$ by the segment $s'_i$ dual to $s_i$ then a tree
$\mathfrak{B}'_{i+1}$ neighboring $\mathfrak{B}'_i$ and compatible
with $\mathfrak{B}_{i+1}$ results from $\mathfrak{B}'_i$.

\section{Manifolds with Two Surface Pieces and Two Points}

We now focus on the case where there are two neighboring trees
$\mathfrak{B}_i,\mathfrak{B}_{i+1}$ for which there is no tree
$\mathfrak{B}'$ compatible with both. In this case \emph{the manifold
may be altered by elementary transformations in such a way that there is
a tree $\mathfrak{B}_{i,i+1}$ in the new manifold $\mathfrak{M}^*$,
neighboring  both $\mathfrak{B}_i$ and $\mathfrak{B}_{i+1}$, and
in the dual manifold of $\mathfrak{M}^*$ there are trees
${\mathfrak{B}^{*}_i}'$ and ${\mathfrak{B}^{*}_{i+1}}'$
compatible with $\mathfrak{B}_i$ and $\mathfrak{B}_{i,i+1}$ on the
one hand, and with $\mathfrak{B}_{i,i+1}$ and $\mathfrak{B}_{i+1}$
on the other}.

If $\mathfrak{B}'_i$ and $\mathfrak{B}'_{i+1}$ are two neighboring
trees compatible with $\mathfrak{B}_i$ and $\mathfrak{B}_{i+1}$
respectively, then we remove by reduction all segments which appear in
$\mathfrak{B}_i$ as well as in $\mathfrak{B}_{i+1}$, and all which 
correspond to those appearing in $\mathfrak{B}'_i$ as well as in
$\mathfrak{B}'_{i+1}$. In this way $\mathfrak{M}$ is converted into a
manifold with two surface pieces $f_1,f_2$ and two points $p_1,p_2$
and trees consisting of single segments. The one appearing in
$\mathfrak{B}_i$ is called $s_i$ and the one appearing in 
$\mathfrak{B}_{i+1}$ is called $s_{i+1}$. Then the segment in
$\mathfrak{B}'_i$ is the dual to $s_{i+1}$ and the segment in
$\mathfrak{B}'_{i+1}$ is the dual to $s_i$. If $r_1,r_2$ are simple
boundary paths of $f_1,f_2$ then
\[
r_1=s_i r_{11} s^{\varepsilon_1}_{i+1} r_{12},\quad
r_2=s_i r_{21} s^{\varepsilon_2}_{i+1} r_{22}.
\]
Here $r_{11},r_{12}$ together pass through each segment either twice
or not at all. The same holds for $r_{21}$ and $r_{22}$. Otherwise the
closed path $s_i s^{-1}_{i+1}$ ($s_i$ and $s_{i+1}$ each run from $p_1$
to $p_2$) would not separate the manifold. 

Now, on the one hand, all segments apart from $s_i$ and $s_{i+1}$
could be singular. Then either the segment $s_\alpha$ following $s_i$
in $r_1$ is singular or the last segment $s_\beta$ of $r_1$ is singular
(otherwise we would have $s_\alpha=s_\beta=s^{\varepsilon_1}_{i+1}$).
If, say, $s_\alpha$ is singular it begins at $p_1$ and we divide $f_1$ into
$f_{11}$ and $f_{12}$ by introducing the segment $t$ with boundary
points $p_1$ and $p_2$. And when $r_1=s_i s_\alpha r_3$ we take
$r_{11}=tr_3$ as the boundary path of $f_{11}$ and 
$r_{12}=s_i s_\alpha t^{-1}$ as the boundary path of $f_{12}$. Then
the path $ts^{-1}_{i+1}$ does not separate the manifold. For $f_{11}$
and $f_{12}$ meet along $s_\alpha$ and $f_{11}$ and $f_2$ meet
along $s_{i}$. The tree consisting of the single segment $t$ therefore
satisfies the required conditions. We proceed quite analogously when
$s_\alpha=s^{\varepsilon_1}_{i+1}$ and $s_\beta$ is singular.

If, on the other hand, there is another regular segment $s_\gamma$, 
apart from $s_i$ and $s_{i+1}$, that goes from $p_1$ to $p_2$, then 
the tree consisting of the segment $s_\gamma$ satisfies our conditions.
For the paths $s_\gamma s^{-1}_i$ and $s_\gamma s^{-1}_{i+1}$ 
do not separate the manifold.

\section{Elementary Relatedness and Isomorphism}

In order to be able to precisely express the results of the previous sections
we introduce, by analogy with the concept of a ``neighboring normal
polygon'' in Section 5.10 (``normal polygon'' is the same as ``normal
form of a manifold''), the concept of ``dually neighboring normal polygons'':

Two normal polygons are called \emph{dually neighboring} if they are
convertible into each other by an extension and subsequent reduction of
the third kind.

Then our result reads as follows:

\emph{Let $\mathfrak{N}$ and $\mathfrak{N}'$ be two normal polygons,
and let}
\[
\mathfrak{N}=\mathfrak{M}^{(1)},\ldots,\mathfrak{M}^{(n)},
\mathfrak{M}^{(n+1)},\ldots,\mathfrak{M}^{(m)}=\mathfrak{N}'
\]
\emph{be a chain of manifolds, where $\mathfrak{M}^{(i)}$ results from
$\mathfrak{M}^{(i-1)}$ by an elementary subdivision for $i=1,2,\ldots,n$
and by an elementary reduction for $i=n+1,\ldots, m$. Let 
$\boldsymbol{I}([w])=[w]'$ be the induced isomorphism between the groups
of $\mathfrak{N}$ and $\mathfrak{N}'$. Then there is a chain of normal
polygons}
\[
\mathfrak{N}=\mathfrak{N}^{(1)},\mathfrak{N}^{(2)},\ldots,
\mathfrak{N}^{(m)}=\mathfrak{N}'
\]
\emph{of which any two in succession are neighboring or dually neighboring, 
and the isomorphism $\boldsymbol{I}'([w])$ induced by this chain is identical
with $\boldsymbol{I}([w])$}.

Now let 
${s^{(i)}_1}^{\pm 1},{s^{(i)}_2}^{\pm 1},\ldots,{s^{(i)}_g}^{\pm 1}$
be the  $2g$ segments which appear in $\mathfrak{N}^{(i)}$, let
${S^{(i)}_1}^{\pm 1}$, ${S^{(i)}_2}^{\pm 1}$, \ldots, 
${S^{(i)}_g}^{\pm 1}$
be the corresponding group elements, let $\boldsymbol{I}_i$ be the
isomorphism between $\mathfrak{N}^{(i)}$ and $\mathfrak{N}^{(i+1)}$
determined by Section 6.2, and in fact let
\[
T^{(k)}_i=\boldsymbol{I}_i(S^{(i+1)}_k).
\]
Then the elements $T^{(k)}_i$ may be expressed in terms of the $S^{(i)}_l$
in such a way that these power products are a system of free generators for
the free group determined by the $S^{(i)}_l$. If $\mathfrak{N}^{(i)}$ and
$\mathfrak{N}^{(i+1)}$ are neighbors this follows from the formulas 
given in Section 6.2. If they are dually neighboring one notes that the 
connection number of the line segment complex is preserved when the
transition complex has two points and one surface piece. The assertion then
follows from the facts on exchange of generators proved in Section 4.6.

But then the analogous result follows for arbitrary elementarily related 
normal forms. If $\boldsymbol{I}$ is the isomorphism between the
fundamental groups of $\mathfrak{N}^{(1)}$ and $\mathfrak{N}^{(n)}$
determined by $T^{(1)}_k=\boldsymbol{I}(S^{(n)}_k)$ then the
$T^{(1)}_k$ may be expressed as power products in the $S^{(1)}_k$ in
such a way that these products constitute a system of free generators of the
free group determined by the $S^{(1)}_k$.

There remains the problem of making the system of generators $T^{(1)}_k$
obtained in this way more easily visualizable. For this purpose we refine the
previous result by the following theorem.

\emph{If $\mathfrak{N}$ and $\mathfrak{N}'$ are two dually neighboring
normal polygons then there is a chain}
\[
\mathfrak{N}^{(1)}=\mathfrak{N},\mathfrak{N}^{(2)},\ldots,
\mathfrak{N}^{(n)}=\mathfrak{N}'
\]
\emph{of neighboring normal polygons which establishes the same 
isomorphism between the groups of $\mathfrak{N}$ and
$\mathfrak{N}'$ as that induced by the original transformation}.

Certainly there is a chain $\mathfrak{N}^{(1)}=\mathfrak{N},
\mathfrak{N}^{(2)},\ldots,\mathfrak{N}^{(n)}=\mathfrak{N}'$
of dually neighboring normal polygons in which the passage from
$\mathfrak{N}^{(i)}$ to $\mathfrak{N}^{(i+1)}$ is effected by a
complex in which the newly introduced, and then eliminated, point
bounds exactly three line segments. The dual situation has already 
been dealt with in Section 5.10. I claim that $\mathfrak{N}^{(i)}$
and $\mathfrak{N}^{(i+1)}$ are then neighboring also.

To prove this we consider the complex of a single surface piece $f$
and two points $p_1,p_2$ that accomplishes the transition. Let $s_i$,
$s^{-1}_{i+1}$, and $s_1$ be the three segments that go from $p_1$ 
to $p_2$; $\mathfrak{N}^{(i)}$ and $\mathfrak{N}^{(i+1)}$ result
from reduction of $s_i$ and $s_{i+1}$ respectively. The boundary path of
$f$ runs over $p_2$ three times and contains the subpaths
$s_is^{-1}_1$, $s_1 s_{i+1}$ or $s^{-1}_{i+1} s^{-1}_1$, 
$s_i s_{i+1}$ or $s^{-1}_{i+1} s^{-1}_i$. Say
\[
r=s_i s^{-1}_1 r_1 s_1 s_{i+1} r_2 s_i s_{i+1} r_3.
\]
Then the boundary path remaining in $\mathfrak{N}^{(i)}$ is
\[
r^{(i)}=s^{-1}_1 r_1 s_1 s_{i+1} r_2 s_{i+1} r_3,
\]
and in $\mathfrak{N}^{i+1}$
\[
r^{(i+1)}=s_i s^{-1}_1 r_1 s_1  r_2 s_i  r_3.
\]
If we now construct the triangle $s_i s^{-1}_1 t$ and eliminate $s_1$, 
then the boundary path
\[
t^{-1} r_1 t s_i r_2 s_i r_3
\]
results from $\mathfrak{N}^{(i+1)}$. But apart from notation
($t\sim s$, $s_{i+1}\sim s_i$) this is identical with the boundary path of
$\mathfrak{N}^{(i)}$.

Further, in both transformations of $\mathfrak{N}^{(i)}$ into
$\mathfrak{N}^{(i+1)}$ the same isomorphism between their groups
takes place.

One can describe the connection between the two transformations
intuitively by saying that in one case one lets the endpoint of the segment
$s_1$ slide along the segment $s_i$, respectively $s_{i+1}$, while in the
second case one constructs a triangle from the initial and final  positions
and the segment $s_i$, respectively $s_{i+1}$, across which the transition
from the initial to final position is carried out.

The passage between two neighboring normal polygons can finally be carried
out by a chain of extensions and reductions in which each newly introduced, 
and then eliminated, surface piece is a triangle. If we call these alterations
triangle transformations then we have the theorem:

\emph{If $\mathfrak{N}$ and $\mathfrak{N}'$ are two elementarily
related normal polygons and if $\boldsymbol{I}$ is the isomorphism
induced between the groups of $\mathfrak{N}$ and $\mathfrak{N}'$
by the elementary transformation of $\mathfrak{N}$ into
$\mathfrak{N}'$, then the same isomorphism may be effected by a
chain}
\[
\mathfrak{N}=\mathfrak{N}^{(1)},\mathfrak{N}^{(2)},
\ldots,\mathfrak{N}^{(n)}=\mathfrak{N}',
\]
\emph{where $\mathfrak{N}^{(i)}$ results from $\mathfrak{N}^{(i-1)}$
by a triangle transformation}.

\section{Some Problems}

Using the results of the preceding section we may derive a method for
presenting generators for the groups of automorphisms of the
fundamental group which result from elementary transformations and
mappings of isomorphic manifolds (cf. Section 6.12). The exchange of
generators in the fundamental group due to elementary modifications of
the normal form yields all automorphisms immediately in the case of the
torus, because the torus has only the canonical form as normal form;
however this is not the case for surfaces of higher genus, because they always
have different normal forms which are elementarily related but not
isomorphic to to each other. Further, the exchanges of generators
induced by elementary transformations in this case constitute a groupoid
(Section 1.15), the identities of which correspond to the different normal
forms of the surface. By the final result of Section 6.10 one can take the
transformations
\[
S'_i=S_i\quad (i\ne a),\qquad
S'_a=S_a S^{\pm 1}_b\quad\text {or}\quad S'_a=S^{\pm 1}_b S_a
\]
as generators of this groupoid, corresponding to the modifications which
convert one normal form into another.

Incidentally, one sees that the automorphism group of the torus determined
by elementary transformations is identical with the automorphism group of
the free commutative group on two generators. The analogous theorem has
been proved by \textsc{Dehn} and 
\textsc{Nielsen}\footnote{\textsc{J. Nielsen}, Acta. Math. \textbf{50},
(1927), 191.} for the groups of the remaining orientable 
manifolds.\footnote{Reidemeister here describes the relation between
homeomorphic mappings of a surface and automorphisms of its fundamental 
group rather loosely. To be more
precise, the surface mappings should be taken modulo isotopy and the 
automorphisms should be taken modulo inner automorphisms. (Note that the
fundamental group of the torus is abelian, so its inner automorphism
group is trivial.) Then the Dehn-Nielsen-Baer theorem states that the
group of homeomorphisms modulo isotopy is isomorphic to the group of
automorphisms of the fundamental group modulo inner automorphisms.
(Translator's note.)}

The automorphism group for $p=2$ has been determined by
\textsc{Baer}\footnote{\textsc{R. Baer}, J. f\"ur reine und angew. Math.
\textbf{160} (1928), 1.} in a way different from that sketched here.

Conversely, one can ask to what extent the elementary transformations
are characterized by the induced isomorphisms of the fundamental group.
To make this question precise we remark that: by \textsc{Tietze} we can
subdivide two elementarily related manifolds so that two isomorphic
subcomplexes result. More precisely, the following holds: if
$\mathfrak{M}$ and $\mathfrak{M}'$ are two manifolds resulting from
each other by a chain of extensions and reductions, then one can arrange
the transition in such a way that first only extensions appear, and then only
reductions, while the isomorphisms of the fundamental group induced by
both chains are the same.\footnote{\textsc{H. Tietze} Mon. f. Math. u.
Phys. Jahrg. 19, p.1 and \textsc{E. Biltz}, Math. Zeitschr. \textbf{18} (1923), 1.}
Now if $\mathfrak{M}$ and $\mathfrak{M}'$ are canonical normal forms
($\mathfrak{M}$ and $\mathfrak{M}'$ are then isomorphic), and
$\mathfrak{M}^{*}$ is the manifold in the chain that contains the most
elements, $w_i$ and $w'_i$ ($i=1,2,\ldots,q$) are the paths in
$\mathfrak{M}^{*}$ corresponding to the segments $s_i$ and $s'_i$ of
$\mathfrak{M}$ and $\mathfrak{M}'$, and if the induced automorphism
is induced is inner, what is the relation between the paths $w_i$ and $w'_i$?

One would conjecture a theorem corresponding to one of 
\textsc{Baer}\footnote{\textsc{R. Baer}, J. f\"ur reine und angew. Math.
\textbf{159} (1928), 101.} concerning continuous manifolds, that the $w_i$
may be deformed into $w'_i$ in $\mathfrak{M}^{*}$ or a subdivision of
$\mathfrak{M}^{*}$ in such a way that in each intermediate position
$\overline{w}_i$ there is always a system of simple paths meeting in
only one point. In connection with this there is the still unproved
combinatorial theorem: simple homotopic paths are always 
``isotopic,''\footnote{\textsc{R. Baer}, J. f\"ur reine und angew. Math.
\textbf{159} (1928), 101.} i.e., with suitable subdivision of the initial
manifold they may always be transformed into one another in such a way
that the intermediate positions are also simple paths.

On the basis of these results one can see that the elementary transformations
by which one manifold is converted to another, isomorphic to it, and which
therefore induce automorphisms in the fundamental group may be
classified quite analogously to the continuous mappings of continuous
manifolds.

One may compare these questions with the determination of simple paths
on a manifold by \textsc{Dehn} and \textsc{Baer}, which, as stated at the
beginning of Section 6.8, is connected with the determination of the
automorphism group.

\section{Coverings of Surface Complexes}

Coverings of surface complexes may be defined similarly to those for line
segment complexes, and the basic theorems are likewise obtained as they
were earlier.

If $\mathfrak{F}$ and $\mathfrak{F}^{*}$ are two surface complexes we
say that $\mathfrak{F}$ covers $\mathfrak{F}^{*}$ when there is a mapping
$\boldsymbol{A}(\mathfrak{F})=\mathfrak{F}^{*}$ of the points, line
segments, and surface pieces of $\mathfrak{F}$ onto those of 
$\mathfrak{F}^{*}$ satisfying the following conditions.

A.1. \emph{If $\mathfrak{C}$ is the line segment complex contained in
$\mathfrak{F}$ and $\mathfrak{C}^{*}$ is that contained in 
$\mathfrak{F}^{*}$, then the mapping 
$\boldsymbol{A}(\mathfrak{C})=\mathfrak{C}^{*}$ is a covering of $\mathfrak{C}^{*}$ by $\mathfrak{C}$, which satisfies the conditions of 
Section 4.17}.

A.2. \emph{For each surface piece $f$ of $\mathfrak{F}$ there is a
well-defined surface piece $\boldsymbol{A}(f)=f^{*}$, and
$\boldsymbol{A}(f^{-1})=(\boldsymbol{A}(f))^{-1}$}.

For the sake of convenience, before we formulate the last condition we will
suppose that in $\mathfrak{F}^{*}$ each simple positive boundary path
$w$ of a surface piece $f^{*}$ bounds only this one surface piece
positively.

A.3. \emph{If $r$ is a simple positive boundary path of the surface piece
$f$ of $\mathfrak{F}$ then $\boldsymbol{A}(r)=r^{*}$ is a simple
positive boundary path of $\boldsymbol{A}(f)$. If $w'$ is any path of
$\mathfrak{F}$ for which $\boldsymbol{A}(w')=r^{*}$, where $r^{*}$
is a boundary path of $f^{*}$, then there is also a surface piece $f'$ in
$\mathfrak{F}$ of which $w'$ is a simple positive boundary path and for
which $\boldsymbol{A}(f')=f^{*}$}.

$\mathfrak{F}$ is called homomorphic to $\mathfrak{F}^{*}$ when there
is a covering $\boldsymbol{A}(\mathfrak{F})=\mathfrak{F}^{*}$.
Homomorphism is transitive. If 
$\mathfrak{F},\mathfrak{F}^{*},\mathfrak{F}^{**}$ are three surface
complexes and if $\boldsymbol{A}(\mathfrak{F})=\mathfrak{F}^{*}$
and $\boldsymbol{A}'(\mathfrak{F}^{*})=\mathfrak{F}^{**}$ are
coverings of $\mathfrak{F}^{*}$ and $\mathfrak{F}^{**}$ by
$\mathfrak{F}^{*}$, then
\[
\boldsymbol{A}''(\mathfrak{F})=
\boldsymbol{A}'(\boldsymbol{A}(\mathfrak{F}))=\mathfrak{F}^{**}
\]
is the mapping of $\mathfrak{F}$ onto $\mathfrak{F}^{**}$ obtained
via $\boldsymbol{A}$ and $\boldsymbol{A}'$. Then the line segment
complex of $\mathfrak{F}$ covers that of $\mathfrak{F}^{**}$, each
$f^{**}$ corresponds to an $f$, and 
$\boldsymbol{A}''(f^{-1})=(\boldsymbol{A}''(f))^{-1}$. Finally, if
$w^{**}$ is any simple boundary path of a simple surface piece $f^{**}$,
and $\boldsymbol{A}''(w)=\boldsymbol{A}'(w^*)=w^{**}$ and
$\boldsymbol{A}(w)=w^*$, then $w^*$ is a simple positive boundary
path of the surface piece $f^*$ with $\boldsymbol{A}(f)=f^*$ and $w$
is a simple positive boundary path of the surface piece $f$ with
$\boldsymbol{A}(f)=f^*$ and $\boldsymbol{A}''(f)=f^{**}$.

\section{Coverings of Line Segment and Surface Complexes}

One can immediately survey those coverings of a line segment complex
$\mathfrak{C}^*$ by $\mathfrak{C}$ that are extendible to a covering
of the surface complex $\mathfrak{F}^*$ by $\mathfrak{F}$. 
\emph{It is necessary and sufficient that in the covering}
\[
\boldsymbol{A}(\mathfrak{C})=\mathfrak{C}^*
\]
\emph{all paths $w$ that lie over simple boundary paths of $\mathfrak{C}^*$
are closed}. It is necessary by A.3 in Section 6.1 because the boundary path
of a surface piece is always closed; it is sufficient because a line segment
complex $\mathfrak{C}$ which covers $\mathfrak{C}^*$ in the way described
may be immediately extended to a surface complex $\mathfrak{F}$ which
covers $\mathfrak{F}^*$. Namely, if
${f^{*}_1}^{\pm 1}, {f^{*}_2}^{\pm 1},\ldots $ are the surface pieces of
$\mathfrak{F}^*$, ${r^{*}_i}^{\pm 1}$ is a simple positive boundary path
of ${f^{*}_i}^{\pm 1}$ ($i=1,2,\ldots$), and if $r_{i1},r_{i2},\ldots$ are
the paths in $\mathfrak{C}$ over $r^*_i$, then we add the surface pieces
$f^{\pm 1}_{i1},f^{\pm 1}_{i2}\ldots$ ($i=1,2,\ldots$) to $\mathfrak{C}$
and define $r_{ik}$ to be a simple positive boundary path of $f_{ik}$. The
complex $\mathfrak{F}$ consisting of $\mathfrak{C}$ and the 
$f^{\pm 1}_{ik}$ then covers $\mathfrak{F}^*$.

For if $r$ is any path which lies over a simple boundary path $r^*$,
$\boldsymbol{A}(r)=r^*$, then $r$ is also a simple boundary path. 
This is because,
by Section 4.8, $r$ results from cyclic interchange in one of the paths
$r^{\pm 1}_{ik}$ over ${r^*_i}^{\pm 1}$. Thus $r$ bounds at least one
surface piece $f$.

It is possible that $r$ bounds several surface pieces $f_1,f_2,\ldots,f_n$;
namely, when there are $n$ different cyclic interchanges
$r^{(1)},r^{(2)},\ldots,r^{(n)}$ of $r$ for which
$\boldsymbol{A}(r^{(i)})={r^*}^{\varepsilon_i}$ ($\varepsilon_i=\pm 1$).
In this case it can happen, e.g., that $r$ bounds exactly two different
surface pieces $f_{11}$ and $f_{12}$ for which 
$\boldsymbol{A}(f_{11})=\boldsymbol{A}(f_{12})$; namely, when there
are two cyclic interchanges $r^{(11)}$ and $r^{(12)}$ of $r$, for which
$\boldsymbol{A}(r^{(11)})=\boldsymbol{A}(r^{(12)})$.
$\boldsymbol{A}(r^{(11)})={r^{(11)}}^*$ must then be carried into itself
by a cyclic transformation; i.e., we must have ${r^{(11)}}^*=w^* w^*$.

We will replace the condition just found for extendibility of a covering of a
line segment complex $\mathfrak{C}^*$ to one of a surface complex
$\mathfrak{F}^*$ by a condition on the permutations associated with
each covering of a line segment complex $\mathfrak{C}^*$ by
$\mathfrak{C}$ (Section 4.10). Let $\mathfrak{B}^*$ be any spanning tree 
of $\mathfrak{C}^*$, let 
${s^*_1}^{\pm 1},{s^*_2}^{\pm 1},\ldots$ be the segments of 
$\mathfrak{C}^*$ that do not appear in $\mathfrak{B}^*$, let
$\pi_1,\pi_2,\ldots$ be the corresponding permutations, and let
$S_1,S_2,\ldots$ be the corresponding generators of the fundamental
group of $\mathfrak{C}^*$. Further, let $r^*_1,r^*_2,\ldots$ be all the
simple boundary paths of surface pieces of $\mathfrak{C}^*$. Then 
the permutations corresponding to these paths $r^*_i$ must all be the identity
permutation. Thus, by Section 4.16, the $\pi_i$ must satisfy the relations
of the fundamental group of the surface complex. Thus the group
generated by the $\pi_i$ is a representation of the fundamental group 
$\mathfrak{W}$ of $\mathfrak{F}^*$. Conversely, each such
permutation group may be associated with a line segment complex
$\mathfrak{C}$ with $\boldsymbol{A}(\mathfrak{C})=\mathfrak{C}^*$
and hence also with a surface complex $\mathfrak{F}$ with
$\boldsymbol{A}(\mathfrak{F})=\mathfrak{F}^*$.

\section{The Fundamental Group of the Covering Complex}

Under the hypotheses of the previous section it is clear from Section 4.17
that the group of the closed paths in $\mathfrak{C}$ beginning at $p$ is
homomorphic to a subgroup of the group of closed paths in
$\mathfrak{C}^*$ beginning at $\boldsymbol{A}(p)=p^*$. We now
also assert:

\emph{The group $\mathfrak{W}$ of classes of closed paths in the surface 
complex $\mathfrak{F}$ beginning at $p$ is isomorphic to a subgroup 
$\mathfrak{U}$ of the group $\mathfrak{W}^*$ of closed paths in the
surface complex $\mathfrak{F}^*$ beginning at 
$\boldsymbol{A}(p)=p^*$}.

If $w$ is a closed path beginning at $p$ and $\boldsymbol{A}(w)=w^*$,
then $[w],[w^*]$ respectively are the elements of the groups
$\mathfrak{W},\mathfrak{W}^*$ to which $w,w^*$ respectively belong,
then we define the mapping $\boldsymbol{I}$ by
\[
\boldsymbol{I}([w])=[w^*]=[\boldsymbol{A}(w)]
\]
and claim that $\boldsymbol{I}$ maps $\mathfrak{W}$ one-to-one onto
a subgroup $\mathfrak{U}$ of $\mathfrak{W}^*$. Since
\[
[w^{-1}]=[w]^{-1},\quad 
[\boldsymbol{A}(w^{-1})]=[\boldsymbol{A}(w)]^{-1},
\]
and since also
\[
[w_1][w_2]=[w_1 w_2]\quad\text{and}\quad
[\boldsymbol{A}(w_1)][\boldsymbol{A}(w_2)]=[\boldsymbol{A}(w_1 w_2)],
\]
we have only to prove one-to-oneness; i.e., to show that if $w_1$ and $w_2$
are two closed paths beginning at $p$ and $[w_1]=[w_2]$ then
$[\boldsymbol{A}(w_1)]=[\boldsymbol{A}(w_2)]$; and conversely, if
$w^*_1, w^*_2$ are two closed paths beginning at $p^*$ and if
$w_1,w_2$ are two beginning at $p$, likewise closed, for which
\[
\boldsymbol{A}(w_i)=w^*_i\quad (i=1,2)
\]
and if $[w^*_1]=[w^*_2]$, then $[w_1]=[w_2]$ also.

But that comes down to showing: if $[w]$ is the identity element of
$\mathfrak{W}$ then $[\boldsymbol{A}(w)]$ is the identity element of
$\mathfrak{U}$, and conversely, if $[\boldsymbol{A}(w)]$ is the
identity element of $\mathfrak{U}$ then $[w]$ is the identity element
of $\mathfrak{W}$. Now if $[w]$ is the identity element then $w$ may
be altered by elementary extensions and reductions until it consists of
subpaths of the form $w' r w'^{-1}$, which likewise begin and end at $p$
and for which each $r$ is a simple boundary path in $\mathfrak{F}$.
But then
\[
\boldsymbol{A}(w' r w'^{_1})=
\boldsymbol{A}(w')\boldsymbol{A}(r)\boldsymbol{A}(w'^{-1})
\]
and $\boldsymbol{A}(r)$ is likewise a simple boundary path in
$\mathfrak{F}^*$, so $[\boldsymbol{A}(w' r w'^{-1})]$ is the
identity element of $\mathfrak{W}^*$. The first part of the claim follows
from this. And, since $r$ is a boundary path in $\mathfrak{F}$ 
when $\boldsymbol{A}(r)=r^*$ is a boundary path in $\mathfrak{F}^*$,
the converse part of the claim also follows.

If $\mathfrak{F}$ is connected then each point $p_i$ for which
$\boldsymbol{A}(p_i)=p^*$ may be connected to $p=p_0$. The
collection of paths from $p$ to $p_i$ then corresponds to a residue class
$\mathfrak{U}G$ in the group $\mathfrak{W}^*$ of the complex
$\mathfrak{F}^*$. The corresponding permutations stand in the same
relation to the residue classes as was described in Section 4.17. Likewise 
in analogy with Section 4.17, a covering of $\mathfrak{F}^*$ by 
$\mathfrak{F}$ may be constructed for each subgroup $\mathfrak{U}$
of $\mathfrak{W}^*$.

One can imagine applying the process for determining generators and
relations of subgroups to these surface complexes. A system of paths $w_i$
from $p_0$ to the $p_i$ yields a complete system of representatives $G_i$
for the residue classes $\mathfrak{U}G$ modulo $\mathfrak{U}$ in
$\mathfrak{W}^*$. The $w_i$ constitute a tree when the $G_i$ satisfy
the \textsc{Schreier} condition $(\Sigma$) of Section 3.6. The boundary 
paths of the surface pieces $f_i$ lying over the same surface piece $f^*$ 
yield the
relations $G_i R G^{-1}_i$ when $R$ is the relation corresponding to $f^*$.
The free group determined by the generators $U_{G,S}$ of Section 3.7  is
the fundamental group of the line segment complex $\mathfrak{C}$
associated with $p_0$. Now one has everything needed to follow the proofs
of the theorems of Section 3.7 step by step for surface complexes.

The invariance of the surface complex group under elementary 
transformations corresponds here to the theorem: if $\mathfrak{F}^*$ and
${\mathfrak{F}^*}'$ are two elementarily related surface complexes, and if
$\mathfrak{F}$ covers $\mathfrak{F}^*$, then there is a well-defined
surface complex $\mathfrak{F}'$, elementarily related to $\mathfrak{F}$,
which covers ${\mathfrak{F}^*}'$.

\section{Regular Coverings}

If we take  $\mathfrak{U}$ to be an invariant subgroup of $\mathfrak{W}^*$
then $\mathfrak{C}$ covers the complex $\mathfrak{C}^*$ regularly, and
conversely. In this case one sees that the transformations of $\mathfrak{C}$
constructed in Section 4.20 may be extended to \emph{transformations of
the surface complex $\mathfrak{F}$ into itself}. Namely, if $r$ is a simple 
boundary path in $\mathfrak{F}$ and $\boldsymbol{A}(r)=r^*$ is a
simple boundary path in $\mathfrak{F}^*$ that we associate uniquely with
the surface piece $f^*$, then there is exactly one surface piece $f$, simply
bounded by $r$, for which $\boldsymbol{A}(f)=f^*$; and if $r$ goes to
$\boldsymbol{I}(r)=r'$ by a mapping of $\mathfrak{C}$ into itself, then
$\boldsymbol{A}(r)=\boldsymbol{A}(r')$ and hence $r'$ is also a
boundary path of a well-defined surface piece $f'$ for which
$\boldsymbol{A}(f')=f^*$. If $\boldsymbol{I}(f)=f'$, then
$\boldsymbol{I}(f^{-1})={f'}^{-1}$, so 
$\boldsymbol{I}(\mathfrak{F})=\mathfrak{F}^*$ is in fact a mapping of
$\mathfrak{F}$ onto itself that satisfies conditions A.1 to A.3 of
Section 6.12. Each surface piece $f$ goes to a surface piece $f'$ for which
$\boldsymbol{A}(f)=\boldsymbol{A}(f')$ under all these mappings, but
not necessarily in the same way.

According to Section 6.13 it can happen that a boundary path $r$ in
$\mathfrak{F}$ bounds two different surface pieces $f_1$ and $f_2$
with $\boldsymbol{A}(f_1)=\boldsymbol{A}(f_2)$. Then $r$ goes to itself
under a transformation $\boldsymbol{I}(\mathfrak{F})=\mathfrak{F}$.
An example is afforded by the complex $\mathfrak{F}^*$ consisting of
a single surface piece, one point, and a segment with $w=s^* s^*$ as
simple boundary path. The latter is covered by a complex of two surface
pieces, two points, and two segments with the simple boundary path
$s_1 s_2$ for both $f_1$ and $f_2$.

The connected complex $\mathfrak{F}$ that covers $\mathfrak{F}^*$
and corresponds as described to the subgroup $\mathfrak{U}$ of
$\mathfrak{W}^*$ consisting only of the identity element is called the
\emph{universal covering complex}. If $\mathfrak{F}'$ is an arbitrary
connected covering complex of $\mathfrak{F}^*$ then there is always
a covering $\boldsymbol{A}(\mathfrak{F})=\mathfrak{F}'$ of
$\mathfrak{F}'$ by the universal covering complex. If $w$ is a closed
path in $\mathfrak{F}$ and $\boldsymbol{A}(w)=w^*$ is the path 
corresponding to it, then $[w^*]$ is the identity element of the group
$\mathfrak{W}^*$ of the complex $\mathfrak{F}^*$. If
$\mathfrak{F}^*$ contains only one point, then the line segment
complex contained in $\mathfrak{F}$ is the group diagram of
$\mathfrak{W}^*$ in the generators $[s^*_i]$ that correspond to the
segments $s^*_i$ of $\mathfrak{F}^*$.

\section{Coverings of Manifolds}

We now assume that the covered complex $\mathfrak{F}^*$ is a
manifold. We can apply the considerations of the previous sections when
$\mathfrak{M}^*$ contains no simple boundary path $r^*$ bounding two 
surface pieces $f^*_1$ and $f^*_2\ne {f^*_1}^{\pm 1}$. In this
exceptional case $\mathfrak{M}^*$ is just a sphere. The path $r^*$
runs through each of its segments $s^*$ only once. Either there is only one
such segment, and $\mathfrak{M}^*$ is in fact a sphere, or there are various 
such segments and $r^*=s^* {r^*}'$. Then we combine $f^*_1$ and
$f^*_2$ into a surface piece $f^*$ with the boundary path 
${r^*}'{{r^*}'}^{-1}$. Then one easily sees, by Section 5.3, that ${r^*}'$
is a simple path and hence $\mathfrak{M}^*$ is a sphere. Since a covering of 
the sphere by the sphere is the identity---because the fundamental group is
the identity---we can now assume that each simple boundary path $r^*$ of
$\mathfrak{F}^*$ bounds only one surface piece ${f^*}^{\pm 1}$. The
most important theorem, which we will prove, reads:

\emph{If $\mathfrak{F}$ is  a connected covering of a manifold
$\mathfrak{M}^*$,}
\[
\boldsymbol{A}(\mathfrak{F})=\mathfrak{M}^*,
\]
\emph{then $\mathfrak{F}$ is itself a manifold}.

First we have to show that a segment $s$ of $\mathfrak{F}$ either appears
twice in the boundary path of some surface piece and in no other boundary,
or else $s$ appears once in the boundary of exactly two surface pieces. Let
$\boldsymbol{A}(s)=s^*$. If $s^* r^*_1$ and $s^*r^*_2$ are the two 
boundary paths of $\mathfrak{M}^*$ beginning with $s^*$ and if
\begin{equation}
r^*_1\ne r^*_2, \tag{1}
\end{equation}
then there are also two different paths $sr_i$ with
\[
\boldsymbol{A}(sr_i)=s^* r^*_i\qquad (i=1,2)
\]
and each simple boundary path $sr$ that begins with $s$ is either $sr_1$ or
$sr_2$. Now either $s^* r^*_i$ runs through the segment $s^*$ only once, 
and consequently bounds two different surface pieces $f^*_1$ and $f^*_2$, 
so that $sr_i$ likewise runs through the segment $s$ only once and $s$
bounds two different surface pieces $f_{11}$ and $f_{12}$ with
$\boldsymbol{A}(f_{1i})=f^*_i$; or else $s^* r^*_1$ runs through $s^*$
twice, so that $s^* r^*_i$ are boundary paths of the same surface piece
$f^*$. Then either $sr_1$ runs through $s$ only once; in which case
$sr_1$ and $sr_2$ are distinct paths, not convertible to each other by
cyclic interchange or reversal, and hence they bound distinct surface pieces
$f_1$ and $f_2$ with $\boldsymbol{A}(f_i)=f^*$ ($i=1,2$); or else
$sr_1$ runs through $s$ twice, in which case $sr_2$ is a cyclic
interchange of $(sr_1)^\varepsilon$ ($\varepsilon=\pm 1$) and the $sr_i$
are therefore simple boundary paths of the same surface piece $f$.

In each of these cases only a single surface piece is spanned by the boundary
path $sr_i$. For $s^* r^*_i$ certainly cannot be put in the form
${w^*}^k$. Only $k=2$ comes into question, but then
\[
s^* r^*_i=s^* r_{11} s^* r_{11},
\]
which is excluded by $r^*_1\ne r^*_2$. So if $r^*_1\ne r^*_2$ the assertion
is proved.

Now let
\begin{equation}
r^*_1=r^*_2. \tag{2}
\end{equation}
Then $s^* r^*_1$ must run through the segment $s^*$ twice, so that either
\[
s^* r^*_1=s^* w^*_1 {s^*}^{-1} w^*_2
\]
and hence
\[
s^* r^*_2=s^* {w^*_1}^{-1} {s^*}^{-1} {w^*_2}^{-1},
\]
so that
\[
w^*_i={w^*_i}^{-1}
\]
is the empty path, which means
\[
s^* r^*_1=s^* {s^*}^{-1}
\]
and $\mathfrak{M}^*$ is the sphere, or else
\begin{align*}
s^*r^*_1&=s^* w^*_1 s^* w^*_2\\
s^*r^*_2&=s^* w^*_2 s^* w^*_1
\end{align*}
and $w^*_1=w^*_2=w^*$. Here $s^* w^*$ runs through each segment
of $\mathfrak{M}^*$ exactly once. Now either $sr_1=sr_2$ runs through $s$,
and hence each other segment, twice, in which case $\mathfrak{F}$ is
identical with $\mathfrak{M}^*$; or else $sr_1=sr_2$ runs through $s$ only
once, in which case there are exactly two segments $s_{i1},s_{i2}$ over 
each segment $s^*_i$ of $\mathfrak{M}^*$ and two surface pieces $f_1$
and $f_2$ over the one surface piece $f^*$ of $\mathfrak{M}^*$. Thus the
assertion is also proved in this exceptional case.

It is easy to see that condition A.7 of Section 5.3 is satisfied; the theorem 
above then follows immediately.

\chapter{Branched Coverings}

\section{The Concept of a Branched Covering}

In the case of manifolds it is of particular interest to introduce a new
type of covering, the \emph{branched covering}, in contrast to which
those previously considered may be called \emph{unbranched
coverings}. The theory of branched coverings originates from a strictly
combinatorial treatment of \textsc{Riemann} surfaces and planar
discontinuous groups.

If $\mathfrak{C}$ and $\mathfrak{C}^*$ are two surface complexes,
we say that $\mathfrak{C}$ is a branched covering of $\mathfrak{C}^*$
when there is a mapping $\boldsymbol{A}(\mathfrak{C})=\mathfrak{C}^*$
of the points, line segments, and surface pieces of $\mathfrak{C}$ onto
those of $\mathfrak{C}^*$ satisfying the conditions A.1 and A.2 of
Section 6.12 and, instead of A.3, satisfying the following.

A.31. \emph{If $r$ is a simple positive boundary path of a surface piece
$f$ of $\mathfrak{C}$, let $\boldsymbol{A}(r)=w^*$ be a $k$-tuply
bounding path $r^*$ of the surface piece $\boldsymbol{A}(f)=f^*$, so
$w^*={r^*}^k$. If $r=w_1 w_2 \cdots w_k$ and
$\boldsymbol{A}(w_i)=r^*$ then $w_i\ne w_l$ for $i\ne l$}.

A.32. \emph{If $w'$ is any path of $\mathfrak{C}$ for which
$\boldsymbol{A}(w')=r^*$, where $r^*$ is a boundary path of a surface 
piece $f^*$ then there is exactly one surface piece $f'$ in $\mathfrak{C}$
in the boundary $r'$ of which $w'$ appears. Let
$\boldsymbol{A}(r')$ be a $k$-tuply bounding path of $f^*$, so 
$\boldsymbol{A}(r')={r^*}^k$. Then $k-1$ is called the branching number
of $f'$}.

We are again assuming that in $\mathfrak{C}^*$ each simple positive
boundary path of a surface piece $f^*$ bounds this surface piece only 
positively. If the surface complex $\mathfrak{C}$ covers the surface
complex $\boldsymbol{A}'(\mathfrak{C})=\mathfrak{C}^*$ and
$\mathfrak{C}^*$ covers the complex  
$\boldsymbol{A}(\mathfrak{C}^*)=\mathfrak{C}^{**}$ then the
mapping of the elements of $\mathfrak{C}$ onto those of
$\mathfrak{C}^{**}$ defined by
\[
\boldsymbol{A}''(\mathfrak{C})=
\boldsymbol{A}'(\boldsymbol{A}(\mathfrak{C}))=\mathfrak{C}^{**}
\]
is also a covering, and it is branched or unbranched according as one of
the mappings $\boldsymbol{A}$ or $\boldsymbol{A}'$ is branched or
neither are.

If the covering of a surface complex $\mathfrak{M}^*$ by a surface
complex $\mathfrak{C}$ is of finite order $o$, then there is
a simple relation between the branching number $k-1$ and the order.
\emph{If $f^*$ is a surface piece of $\mathfrak{C}^*$, if}
\[
f^{\pm 1}_1,\quad
f^{\pm 1}_2,\quad
\ldots,\quad
f^{\pm 1}_r
\]
\emph{are the surface pieces over ${f^*}^{\pm 1}$, if
$\boldsymbol{A}(f_i)=f^*$, and if $k_i-1$ is the branching number of
$f_i$, then}
\begin{equation}
o=\sum^{r}_{i=1} k_i. \tag{1}
\end{equation}
For if the boundary path $r^*$ of $f^*$ begins with $s^*$ then by
condition A.31 there are exactly $k_i$ subpaths in the boundary path
$w_i$ of $f_i$ that lie over $r^*$, and consequently $k_i$ different
segments
\[
s_{il}\qquad (l=1,2,\ldots,k_i)
\]
that lie over $s^*$. Further, $s_{il}$ and $s_{jm}$ must be different
segments over $s^*$ when $i\ne j$. Otherwise the boundary paths of
$f_i$ and $f_j$ would be identical, and then by A.32 we should have
$f_i=f_j$. Thus there are at least
\[
\sum^{r}_{i=1} k_i
\]
different segments over $s^*$. On the other hand, since each segment
over $s^*$ must be contained in a boundary path of an $f_i$ (again by
A.32), (1) is in fact satisfied. As a result, a surface piece $f^*$ is
covered by at most $o$ surface pieces $f_i$.

If $\mathfrak{C}_1,\mathfrak{C}^*_1$ are the line segment complexes
of $\mathfrak{C},\mathfrak{C}^*$ respectively then the branched
covering $\boldsymbol{A}(\mathfrak{C})=\mathfrak{C}^*$ is already
determined by the mapping
$\boldsymbol{A}(\mathfrak{C}_1)=\mathfrak{C}^*_1$.
For the boundary paths in $\mathfrak{C}$ are characterized as those that,
on the one hand, lie over boundary paths of $\mathfrak{C}^*$, and on
the other hand are closed.

It follows from this remark that one can define \emph{regular branched
coverings} just as for unbranched coverings, as those for which the
covering of the associated line segment complex is regular, and one may further
define mappings of regular coverings into themselves which exchange
elements of $\mathfrak{C}$ lying over the same element of
$\mathfrak{C}^*$. On the other hand, the close connection between
the fundamental groups of $\mathfrak{C}$ and $\mathfrak{C}^*$ is
destroyed by branched coverings.

In the case of regular coverings the branching number of surface pieces
that lie over the same $f^*$ are all equal. In this case, if $f$ lies over
$f^*$ with branching number $k-1\ne 0$ then there is a transformation
of $\mathfrak{C}$ into itself which displaces the boundary path of $f$
into itself cyclically and carries $f$ into itself. \emph{Among the regular
coverings the branched ones are thus characterized as those that admit
transformations with fixed elements (namely, fixed surface pieces)}.

Interest in these coverings started because of the theorem that
\emph{branched coverings of manifolds are again manifolds}. For when
each segment $s^*$ appears in only two boundary paths, $s^* w^*_1$
and $s^* w^*_2$, then the same is true for each segment $s$ lying
over $s^*$.

\section{Self-transformations and Automorphisms}

There are three kinds of group to consider with a regular covering
$\boldsymbol{A}(\mathfrak{C})=\mathfrak{C}^*$: The fundamental 
group $\mathfrak{W}$ of the complex $\mathfrak{C}$, the
fundamental group $\mathfrak{W}^*$ of the complex 
$\mathfrak{C}^*$, and the group $\mathfrak{T}$ of mappings of
$\mathfrak{C}$ into itself that corresponds to the regular covering.
We have previously investigated the connection between $\mathfrak{W}$
and $\mathfrak{W}^*$ on the one hand, and between $\mathfrak{T}$
and $\mathfrak{W}^*$ on the other, in the case of unbranched
coverings, and we will now apply ourselves to the connection between
$\mathfrak{T}$ and $\mathfrak{W}$ for arbitrary regular coverings.

Let $p_0,p_1,p_2,\ldots$ be the points of $\mathfrak{C}$ over the same
point $p^*$ of $\mathfrak{C}^*$, let $\mathfrak{W}_i$ be the group
of closed  paths of $\mathfrak{C}$ beginning at $p_i$, and let $T_i$
be the mapping of $\mathfrak{C}$ onto itself that carries $p_0$ to
$p_i$. Then the mapping $T_i$ induces an isomorphism
$\boldsymbol{I}_i$ between the groups $\mathfrak{W}_0$ and
$\mathfrak{W}_i$. But a class of isomorphisms between the groups
$\mathfrak{W}_0$ and $\mathfrak{W}_i$ is already known. Namely, if
$w$ is a closed path beginning at $p_i$, and $h_i$ is a path from $p_0$
to $p_i$, then $h_i w h^{-1}_i$ is a closed path beginning at $p_0$
and the mapping
\begin{equation}
[w]\rightarrow [h_i w h^{-1}_i] \tag{1}
\end{equation}
is such an isomorphism $\boldsymbol{I}^*_i$. Thus the transformations
$\boldsymbol{I}_i{\boldsymbol{I}^*_i}^{-1}$  are automorphisms of
the group $\mathfrak{W}_0=\mathfrak{W}$, and in fact they form the
residue class of automorphisms determined by some representative
$A_i$ of the inner automorphisms of $\mathfrak{W}$. The collection of
these residue classes of automorphisms constitutes a group
$\mathfrak{A}$. When none of the transformations $T$ other than the
identity induces an inner automorphism, and if $\mathfrak{I}$ is the
group of inner automorphisms of $\mathfrak{W}_0$, which by
Section 1.12 is an invariant subgroup of $\mathfrak{A}$, then the
factor group $\mathfrak{A}/\mathfrak{I}$ is homomorphic to
$\mathfrak{T}$, and indeed isomorphic.

\section{Principal Group of a Regular Covering}

We can throw more light on the connection between $\mathfrak{W}$
and $\mathfrak{T}$ by constructing the universal covering complex
$\mathfrak{K}$ of $\mathfrak{C}$. \emph{Let
$\boldsymbol{A}'(\mathfrak{K})=\mathfrak{C}$ be the mapping of
$\mathfrak{K}$ onto $\mathfrak{C}$ and let 
$\boldsymbol{A}''(\mathfrak{K})=\mathfrak{C}^*$ be the mapping
of $\mathfrak{K}$ onto $\mathfrak{C}^*$ composed from
$\boldsymbol{A}'$ and $\boldsymbol{A}$. Then $\boldsymbol{A}''$
is a regular covering of $\mathfrak{C}^*$ by $\mathfrak{K}$}.
$\boldsymbol{A}''$ is first of all a covering by transitivity.
Regularity means: if $w$ is a closed path of $\mathfrak{K}$ and
$\overline{w}$ is a path of $\mathfrak{K}$ for which
$\boldsymbol{A}''(w)=\boldsymbol{A}''(\overline{w})$ then
$\overline{w}$ is also closed. We now construct
$\boldsymbol{A}'(w)$ and $\boldsymbol{A}'(\overline{w})$ and
distinguish the two cases
\begin{equation}
\boldsymbol{A}'(w)=\boldsymbol{A}'(\overline{w})
\tag{1}
\end{equation}
and
\begin{equation}
\boldsymbol{A}'(w)\ne \boldsymbol{A}'(\overline{w}).
\tag{2}
\end{equation}
In the first case $\boldsymbol{A}'(w)$ must be a contractible path on
$\mathfrak{C}$ by definition of the universal covering complex in
Section 6.15, and because $w$ is closed it follows that $\overline{w}$
lies over the same contractible path and hence must also be closed.
In the second case we use the fact that
$\boldsymbol{A}(\boldsymbol{A}'(w))=
\boldsymbol{A}(\boldsymbol{A}'(\overline{w}))$, so that
$\boldsymbol{A}'(\overline{w})$ results from $\boldsymbol{A}'(w)$
by a transformation in $\mathfrak{T}$. Thus if
$\boldsymbol{A}'(w)$ bounds a surface piece of $\mathfrak{C}$
the same must hold for $\boldsymbol{A}'(\overline{w})$, and
consequently $\overline{w}$ is also closed in $\mathfrak{K}$.

If $\mathfrak{C}^*$ is a complex that contains only a single
point---and one can convert any connected complex into this form
without altering its group, and think of $\mathfrak{C}$ analogously
altered---then by Section 4.18 the line segment complex
$\mathfrak{K}_1$ contained in $\mathfrak{K}$ is a group diagram,
because the covering 
$\boldsymbol{A}''(\mathfrak{K}_1)=\mathfrak{C}^*_1$ is
regular. The generators and defining relations of the group
$\mathfrak{V}$ determined by this group diagram may be read off
from those of $\mathfrak{W}^*$. If $s^*_1,s^*_2,\ldots$ are the
singular segments of $\mathfrak{C}^*$ and if
$r^*_1(s^*),r^*_2(s^*),\ldots,r^*_n(s^*)$ are the simple 
boundary paths of the surface pieces $f^*_1,f^*_2,\ldots,f^*_n$
of $\mathfrak{C}^*$, and if also the surface piece $f^*_i$ is
covered $(k_i-1)$-tuply by the branched covering, then the generators
$S_1,S_2,\ldots$ of $\mathfrak{V}$ correspond uniquely to the
$s^*_1,s^*_2,\ldots$ and one obtains the defining relations by
replacing the $s^*_i$ by $S_i$ in the power products
$r^*_i(s^*)^{k_i}$, because the boundary paths of $\mathfrak{K}$
and $\mathfrak{C}$ lie over the paths $r^*_i(s^*)^{k_i}$.
\emph{$\mathfrak{V}$ is called the principal group of the regular
covering $\boldsymbol{A}(\mathfrak{C})=\mathfrak{C}^*$}.

\section{Structure of the Principal Group}

The structure of the principal group $\mathfrak{V}$ is connected
with $\mathfrak{W}_0$ and the automorphisms of $\mathfrak{W}_0$
from $\mathfrak{A}$. The paths $w$ of the group diagram
$\mathfrak{K}_1$ that begin and end at points lying over $p_0$ of
$\mathfrak{C}$ determine a \emph{subgroup $\mathfrak{U}$ of
$\mathfrak{V}$ which is isomorphic to $\mathfrak{W}_0$}. One 
obtains the power products of $\mathfrak{U}$ by replacing the
$s^*_i$ by $S_i$ in $\boldsymbol{A}''(w)=w^*(s^*)$.

If one moves any point $p$ in $\mathfrak{K}_1$ by an element of
$\mathfrak{U}$, then the path corresponding to this element ends
at a point $p'$ over the same point in $\mathfrak{C}$ as $p$,
$\boldsymbol{A}'(p)=\boldsymbol{A}'(p')$, because the
covering $\boldsymbol{A}(\mathfrak{C})=\mathfrak{C}^*$ is
regular. Thus the group $\mathfrak{U}$ does not depend on $p_0$.
For this reason, \emph{$\mathfrak{U}$ is an invariant subgroup of
$\mathfrak{V}$}. Namely, if $U$ is any word of $\mathfrak{U}$
then, when carried from any point, $U$ determines a path with a
closed image in $\mathfrak{C}$, $SUS^{-1}$ determines a path
whose image in $\mathfrak{C}$ is a path $sw\overline{s}^{-1}$.
Here $w$ is the closed path corresponding to $U$ and hence 
$s=\overline{s}$, because both segments end at the same point and lie 
over the same segment $s^*$. Thus $sws^{-1}$ is again closed, and so
$SUS^{-1}$ belongs to $\mathfrak{U}$.

At the same time one sees that \emph{the automorphism of
$\mathfrak{U}$ effected by}
\[
S\mathfrak{U}S^{-1}=\mathfrak{U}
\]
\emph{is associated with a definite automorphism of the group given
above}. For if $p$ is any point of $\mathfrak{C}$ and if $w$ and
$w'=sws^{-1}$ are the paths corresponding to the elements $U$ and
$SUS^{-1}$, then the subpath $w$ of $w'$ either begins and ends at 
a point $p'$ different from $p$, or else begins and ends at $p$. In the
first case we can use $s$ in place of $h$ in the formula (1) of Section
7.2 and write the automorphism effected by
$S\mathfrak{U}S^{-1}=\mathfrak{U}$ as
\[
[w]\rightarrow [sws^{-1}].
\]
In the second case $s$ is a closed path of $\mathfrak{C}$ and thus 
belongs to the subgroup $\mathfrak{U}$, so 
$S\mathfrak{U}S^{-1}=\mathfrak{U}$ is an inner automorphism
of $\mathfrak{U}$.

By Section 4.20 the group diagram $\mathfrak{K}_1$ admits a
simply transitive group of transformations $\overline{\mathfrak{V}}$
isomorphic to the group of paths $\mathfrak{V}$. The subgroup
$\overline{\mathfrak{U}}$ of $\overline{\mathfrak{V}}$
corresponding to $\mathfrak{U}$
exchanges those points that lie over the same point of $\mathfrak{C}$.
These mappings of $\mathfrak{K}_1$ onto itself may be extended to
mappings of $\mathfrak{K}$ onto itself.

A mapping $T$ of $\mathfrak{C}$ which carries the point $p_0$ to 
the point $p_i$ corresponds to a class of mappings of the covering 
$\mathfrak{K}$; namely, all those mappings that carry a point over
$p_0$ to a point over $p_i$. The collection of these transformations
forms a residue class modulo $\overline{\mathfrak{U}}$.
\emph{Thus the group $\mathfrak{T}$ is isomorphic to the factor
group $\overline{\mathfrak{V}}/\overline{\mathfrak{U}}$ or
$\mathfrak{V}/\mathfrak{U}$ respectively}.

Likewise it follows that $\mathfrak{V}$ is homomorphic to the group
$\mathfrak{A}$ of automorphisms of $\mathfrak{W}_0$ given above.
Namely, if $V$ is any element of $\mathfrak{V}$, then this element is
associated with a definite automorphism from $\mathfrak{A}$ by
\begin{equation}
VUV^{-1}=U'. \tag{1}
\end{equation}
This relation effects a homomorphism of $\mathfrak{V}$ onto
$\mathfrak{A}$, which is one-to-one when there are no two elements 
of $\mathfrak{V}$ that induce the same automorphism in $\mathfrak{U}$;
i.e., when
\[
V_0 U V^{-1}_0=U \tag{2}
\]
being true for all $U$ in $\mathfrak{U}$ implies $V_0=1$. In this case
the mapping of $\mathfrak{C}$ onto itself corresponding to $V_0$ is
the identity. When none of the transformations $T$ induces an inner 
automorphism in $\mathfrak{U}$, so that no automorphism (1) is an
inner automorphism of $\mathfrak{U}$ unless $V$ belongs to
$\mathfrak{U}$, then $V_0$ in (2) belongs to $\mathfrak{U}$ and
indeed, by (2), to the center of $\mathfrak{U}$. In summary we can say:

\emph{The group $\mathfrak{A}$ of automorphisms and the group
$\mathfrak{V}$ are isomorphic when}

\emph{1. each transformation $T$ in $\mathfrak{T}$ different from the
identity induces an automorphism which is not an inner automorphism,
and}

\emph{2. the center of $\mathfrak{U}$ consists only of the identity
element}.

\section{Group Diagrams and Manifolds}

Each group diagram of a finite group can be converted in various ways
to a regular covering of an oriented 
manifold.\footnote{cf. \textsc{E. Steinitz}: ``Polyeder und
Raumeinteilungen,''  Math. Enzykl. III, A, B 12 \S 48, 49.} This is done
by cyclically ordering the segments emanating from a point in the
group diagram in a certain way. Let $s_{ik}$ be the segments
emanating from the point $p_k$, let $S^{\pm 1}_{i}$
($i=1,2,\ldots,r$) be the generators corresponding to the segments
$s_{ik}$, and let $\boldsymbol{A}(s_{ik})=S_i$ ($i=1,2,\ldots,r$).
Further, let
\[
s_{\alpha_1,k},\quad
s_{\alpha_2,k},\quad
\ldots,\quad
s_{\alpha_{2r},k}
\]
be a star corresponding to one of the points $p_k$ ($k=1,2,\ldots,o$)
that contains each of the segments emanating from $p_k$ exactly
once. \emph{If we now construct the star}
\begin{equation}
\boldsymbol{A}(s_{\alpha_1,k}),\quad
\boldsymbol{A}(s_{\alpha_2,k}),\quad
\ldots,\quad
\boldsymbol{A}(s_{\alpha_{2r},k}) \tag{1}
\end{equation}
\emph{then the same ``star'' of the $S$ results for all $k$}.

With the help of these stars we define a class of paths that we convert
into simple positive boundary paths of our manifold, and which we will
already call boundary paths. \emph{The closed path}
\[
s_{\beta_1 \gamma_1} s_{\beta_2 \gamma_2}
\cdots s_{\beta_m \gamma_m}
\]
\emph{is called a positive boundary path $r(s)$ when}
\[
s^{-1}_{\beta_i,\gamma_i},\quad
s_{\beta_{i+1},\gamma_{i+1}}\quad (i=1,2,\ldots,m-1)\quad
\text{and}\quad
s^{-1}_{\beta_m,\gamma_m},\quad s_{\beta_1,\gamma_1}
\]
\emph{appear consecutively in the order given in the star of the point at 
which they begin. The path $r(s)$ is called a simple boundary path when 
$r(s)$ contains no subpath which is likewise a boundary path according to
the definition}.

Each segment may be embedded in at least one simple positive 
boundary path which runs through it positively, and two different
such paths are convertible into each other by cyclic interchange.
Likewise, each segment may be embedded in at least one simple
positive boundary path which runs through it negatively, and two 
different such paths are convertible into each other by cyclic
interchange. Each boundary path runs through the same segment
only once in the same direction, but it can run through a
segment in both directions. Thus the boundary paths determined in
the above manner have the properties laid down in Section 5.13 for
boundary paths of oriented manifolds.

If 
\[
r_i(s)\qquad (i=1,2,\ldots,q)
\]
is a system of boundary paths from which all the simple positive
boundary paths may be obtained by cyclic interchange of the $s$,
then we introduce $q$ surface pieces $f_1,f_2,\ldots,f_q$ with the
boundary paths $r_1,r_2,\ldots,r_q$. We then have an oriented
manifold $\mathfrak{M}$.

The closed paths of the group diagram then fall into two classes:
those that are contractible to a point in $\mathfrak{M}$, and those
that are not. If $w$ is contractible in $\mathfrak{M}$ and if
$\overline{w}$ is another path for which
$\boldsymbol{A}(w)=\boldsymbol{A}(\overline{w})$, then
$\overline{w}$ is also contractible in $\mathfrak{M}$. For if $w$
is a simple boundary path and
$\boldsymbol{A}(w)=\boldsymbol{A}(\overline{w})$ then
$\overline{w}$ is also a boundary path because of the condition on 
the stars of the points $p_k$.

\emph{It follows from this} that \emph{mappings of the group
diagram may be extended to mappings of the manifold} 
$\mathfrak{M}$. Of course, such mappings do not in general
exchange the surface pieces transitively among themselves.

If one constructs the manifold $\overline{\mathfrak{M}}$ dual to
$\mathfrak{M}$ then $\overline{\mathfrak{M}}$ likewise admits a
group of mappings, under which the surface pieces corresponding
to the points of $\mathfrak{M}$, or to those of the group diagram,
respectively, are exchanged with each other.

The boundary paths $r_i(s)$ yield a class of relations $R_i(S)$
when one replaces the $s$ in $r_i(s)$ by $\boldsymbol{A}(s)$;
at most $o$ of these relations are identical with each other, when
$o$, as above, is the order of the group. If
\[
R_i(S)\qquad (i=1,2,\ldots,e)
\]
are the different relations then, for each $R_i$, there is a maximal
exponent $k_i$ such that
\[
R_i={R^*_i}^{k_i}.
\]
\emph{Using the $R^*_i$ we can construct a manifold
$\mathfrak{M}^*$ which is regularly covered by the manifold
$\mathfrak{M}$}. Namely, we associate each generator $S^{\pm 1}_i$
with a singular segment ${s^*_i}^{\pm 1}$ beginning and ending
at $p^*$, and as simple positive boundary paths of the surface
pieces $f^*_i$ we take the $r^*_i(s^*)$ ($i=1,2\ldots,e$) which result
when one replaces the $S_l$ in $R^*_i(S)$ by $s^*_l$. The star of
segments $s^*$ at $p^*$ then results from the cycle (1) by 
exchanging the $S$ for the corresponding $s^*$, and in fact
$\mathfrak{M}$ covers $\mathfrak{M}^*$ regularly.

If we construct the universal covering complex $\mathfrak{K}$ of
$\mathfrak{M}$ then the line segment complex $\mathfrak{K}_1$
of $\mathfrak{K}$ forms the group diagram of the principal group
$\mathfrak{V}$ of this covering with the generators
\[
S_l\qquad (l=1,2,\ldots,r)
\]
and relations
\[
R_i={R^*_i}^{k_i}\qquad (i=1,2,\ldots,e).
\]

$\mathfrak{K}$ is called a \emph{planar complex} and 
$\mathfrak{K}_1$, with the cyclic ordering of segments emanating
from a point, is called a \emph{planar group diagram of the first
kind, and $\mathfrak{V}$ is called a planar group of the first kind}.

\section{Point-type Branching}

Under certain conditions one can also define a manifold from the group
diagram of an infinite group with finitely many generators, in which one
favors a particular cycle of generators; namely, when they give a proper
closed boundary path. However, this need not be the case. It can happen,
when one wants to embed a segment in a boundary path, that this path
does not close. Nevertheless, one obtains a well-defined class of indefinitely
extendible paths from the initial segment, or else an infinite open path, 
which may be called an \emph{infinite boundary path}.

One sees from this that a group diagram with infinite boundary paths
can be regarded as the one-dimensional dual complex $\mathfrak{D}_1$
of a manifold $\overline{\mathfrak{M}}$, obtained when one allows
each point $p_i$ of the group diagram to correspond with a surface piece
$\overline{f}_i$ and each class of cyclically related boundary paths, as
well as each infinite boundary path, with a point $\overline{p}_k$ of
$\overline{\mathfrak{M}}$. In this way we get a complex that admits a
group of transformations under which the surface pieces $\overline{f}_i$
are exchanged simply transitively with each other. Such complexes were
originally employed for the representation of arbitrary groups with
finitely many generators.\footnote{\textsc{W. v. Dyck}, Math. Ann. 
\textbf{20} (1882).}

Analogously, one can construct the universal covering complex
$\overline{\mathfrak{K}}$ of $\overline{\mathfrak{M}}$.
$\overline{\mathfrak{K}}$ is called a covering complex of the second kind;
the one-dimensional complex $\mathfrak{K}_1$ dual to $\mathfrak{K}$,
with the cyclic ordering of segments at a point, is called a \emph{planar
group diagram of the second kind}; and the associated group is called a
\emph{planar group of the second kind}. If one spans each finite
boundary path determined by the cyclic order of the segments at the points
by a surface piece, then each simple closed path of $\mathfrak{K}_1$
bounds a surface piece.

We return to the case where the manifold $\mathfrak{M}$ and its dual
$\overline{\mathfrak{M}}$ may both be constructed from the group 
diagram.  $\mathfrak{M}$ then covers $\mathfrak{M}^*$ and in fact
$\boldsymbol{A}(\mathfrak{M})=\mathfrak{M}^*$. One can then
define a manifold $\overline{\mathfrak{M}}^*$ dual to 
$\mathfrak{M}^*$ and a covering
$\boldsymbol{B}(\overline{\mathfrak{M}})=\overline{\mathfrak{M}}^*$
effected by the dual mapping. If the covering is unbranched then this is
also true of $\boldsymbol{B}$, and conversely. If $\boldsymbol{A}$
is branched, however, then $\boldsymbol{B}$ is a new type of covering:
\emph{the branching here takes place at points}.

Such a covering cannot be constructed for general surface complexes,
because stars do not correspond to points there. On the other hand, for
manifolds it is very easy to give interesting examples of point-type branching,
e.g., in the construction of \textsc{Riemann} surfaces in function
theory and the fundamental domains of discontinuous groups.
Moreover, the two types of covering are so closely connected that it does 
not matter which definition is taken as the starting point.

\section{Elementarily Related Coverings}

Single covering complexes are of less interest than classes of complexes
related to each other by elementary transformations. If the manifold
$\mathfrak{M}$ covers the manifold $\mathfrak{M}^*$ and if
$\boldsymbol{A}(\mathfrak{M})=\mathfrak{M}^*$ is unbranched,
and if $\mathfrak{M}^*$ is convertible to ${\mathfrak{M}^*}'$ by
elementary transformations, then $\mathfrak{M}$ may also be
converted to an $\mathfrak{M}'$ by elementary transformations and
a covering $\boldsymbol{A}'(\mathfrak{M})={\mathfrak{M}^*}'$
defined which agrees with $\boldsymbol{A}$ on those segments that
$\mathfrak{M},\mathfrak{M}'$ as well as 
$\mathfrak{M}^*,{\mathfrak{M}^*}'$ have in common. We have
already shown this in Section 6.14. The result is similar for branched
coverings. \emph{If $\mathfrak{M}^*$ is altered by an elementary
extension then $\mathfrak{M}$ may be analogously altered, and a 
new branched covering defined for the new manifolds. On the other
hand, the reductions possible in $\mathfrak{M}^*$ cannot always be
matched by alterations of $\mathfrak{M}$}. Namely, if
\[
\boldsymbol{A}(f_i)=f^*_i\qquad (i=1,2)
\]
and if $f_i$ covers the surface piece $f^*_i$, branched in both cases,
and if $f^*_1$ and $f^*_2$ are both bounded by the segment $s^*$,
then segments lying over $s^*$ cannot be eliminated, because
several such appear in the boundaries of $f_1$ and $f_2$. On the other
hand, if the covering, e.g., of $f^*_2$ by $f_2$ is unbranched, then one
can successively fuse the surface pieces lying over $f^*_2$ with those
lying over $f^*_1$ by elementary reductions. Reductions of the first
and third kind may always be carried out. An analogue holds for
reductions of the first kind in the case of coverings with point-type
branching.

If $\boldsymbol{A}(\mathfrak{M})=\mathfrak{M}^*$ and 
$\boldsymbol{A}'(\mathfrak{M}')={\mathfrak{M}^*}'$ are two covering
complexes related to each other by elementary transformations in the
way described, then they may be called \emph{elementarily related}
for short.

\section{Normal Forms of Coverings}

After these preliminary remarks we can now bring the branched coverings
of the complex $\mathfrak{M}^*$ into normal form and prove the
following theorem: \emph{If the covering 
$\boldsymbol{A}(\mathfrak{M})=\mathfrak{M}^*$ is branched along
the surface pieces $f^*_1,f^*_2,\ldots,f^*_n$ with boundary paths
$r^*_1,r^*_2,\ldots,r^*_n$ and orders $k_i-1$ ($k_i\le k_{i+1}$)
then there is an equivalent covering 
$\boldsymbol{A}(\mathfrak{N})=\mathfrak{N}^*$ in which
$\mathfrak{N}^*$ contains one point and $n$ surface pieces $f^*_i$
with the boundary paths} 
\[
r^*_i(s^*)\quad (i=1,2,\ldots,n-1),\quad
r^*_n(s^*)=s^*_1 s^*_2 \cdots s^*_{n-1} r^*(s^*)
\]
\emph{where $r^*(s^*)$ is the boundary path of a normal polygon}.  

First of all, using reductions of the third kind, all points of
$\mathfrak{M}^*$ may be coalesced into a single one, and the surface
pieces covered without branching may be eliminated by reductions of 
the second kind. Now if the boundary path of $f^*_1$ in the resulting 
manifold still runs e.g. through different segments we divide $f^*_1$
into $f^*_{11}$ with the boundary path $r^*_{11}=t^*$ and
$f^*_{12}$ with the boundary path $r^*_{12}={t^*}^{-1}r^*_1$.
If $f_1$ lies over $f^*_1$ and if the boundary path $r_1$ of $f_1$
equals $r_{11}r_{12}\cdots r_{1k}$ with
$\boldsymbol{A}(r_{1i})=r^*_1$, then $f_1$ is converted into $f_{11}$
with the boundary path $t_1 t_2 \cdots t_k$ and $f_{1i}$ with the
boundary path $t^{-1}_i r_{1i}$ ($i=1,2,\ldots,k$), $f^*_{1i}$ is
$k$-tuply branched and $f^*_{1i}$ is covered without branching.
Therefore, $f^*_{1i}$ may be merged with ${f^*_n}'$ by elementary
reduction. By iteration of these steps we reach a form  in which $n-1$
surface pieces $f^*_i$ appear, bounded only by a segment $t^*_i$.
The boundary $r^*_n$ of the last surface piece runs through each of
these $t^*_i$ only once, and through the remaining segments twice.

If also $r^*_n={t^*_\alpha}^{-1}s^* r^*_{n1}$ or
$r^*_n={t^*_\alpha}^{-1}{t^*_\beta}^{-1} r^*_{n1}$ then in four
steps we can always reach the boundary path
${r'_n}^*=s^*{{t'_\alpha}^*}^{-1} r^*_{n1}$ or
${\overline{t}^*_\beta}^{-1}{\overline{t}^*_\alpha}^{-1}r^*_{n1}$ 
respectively. In the first case we subdivide $f^*_n$ by a segment $u^*$ 
into $f^*_{n1}$ with the boundary path ${t^*_{\alpha}}^{-1}s^*u^*$
and $f^*_{n2}$ with boundary path ${u^*}^{-1} r^*_{n1}$; 
$f^*_{n1}$ is then covered unbranched and $f^*_{n2}$ is covered
branched, and $f^*_{n1}$ may then be merged with $f^*_\alpha$
along $t^*_\alpha$ to form ${f'_\alpha}^*$ with the boundary path
$s^* u^*$. Then we subdivide ${f'_\alpha}^*$ by ${t'_\alpha}^*$
into a surface piece ${f''_\alpha}^*$ with boundary path ${t'_\alpha}^*$
and ${f'^*_{n1}}$ with boundary path $s^* u^* {{t'_\alpha}^*}^{-1}$
and merge ${f'^*_{n1}}$ and ${f'^*_{n2}}$ along $u^*$, since they
are covered unbranched. In the second case one proceeds analogously.

By successive applications of these operations one can obtain $f^*_n$
with the boundary path
\[
t^*_1 t^*_2 \cdots t^*_{n-1} r^*(s^*).
\]
Here $r^*(s^*)$ is the boundary of a normal polygon, which one can
reduce to the canonical normal form by Section 5.11.

In the case of point-type branched coverings one obtains normal forms
with a surface piece and $n$ branch points of order $k_i=1$. The
branch points $p_i$ of order $k_i-1$ ($k_i\le k_{i+1}$)
($i=1,2\ldots,n-1$) each bound a regular segment $t^*_i$ which ends at $p_n$; the boundary path of the surface is then
\[
t^*_1 {t^*_1}^{-1} t^*_2 {t^*_2}^{-1} \cdots
t^*_{n-1}{t^*_{n-1}}^{-1} r^*(s^*),
\]
where $r^*(s^*)$ is the boundary of a normal polygon of the second
kind.

\section{Principal Groups in Normal Form}

The preceding section gives the following result for principal groups:

\emph{In normal form, the principal group $\mathfrak{V}$ has the 
relations}
\begin{align*}
T^{k_i}_i&=1\qquad (i=1,2,\ldots,n-1)\\
R_n={R^*_n}^{k_n}&=(T_1 T_2 \cdots T_{n-1} R(S))^{k_n}=1.
\tag{1}
\end{align*}
\emph{Here $R(S)$ results from $r^*(s^*)$ by replacing the $s^*$
by $S$. Thus the structure of the group $\mathfrak{V}$ is determined
by the branching numbers and the genus of the covering complex
$\mathfrak{M}$. In the case of a point-type branched covering, if at
least one point has infinite branching order then only relations of the
form}
\begin{equation}
T^{k_i}_i=1 \tag{2}
\end{equation}
\emph{remain. Then the group is the free product of cyclic groups}.

While it is possible to immediately derive the genus of the covered
manifold from the structure of the group in the case of surface-type
branched coverings (one constructs the factor group by the commutator
group and sets all elements of finite order equal to the identity,
obtaining a free \textsc{Abelian} group with $2g$ generators, where
$g$ is the genus sought), for point-type branching this is in general not
the case. For example, the group diagram of a free group with two generators
may be realized in the plane in such a way that the manifold covered
has genus 0 or 1.

The number of branching surfaces or branch points, respectively, and
the branching numbers may be divided into two classes with the help 
of the following theorem.

Let $T_i,S_k$ be the generators of the group $\mathfrak{V}$ in the
normal form, which satisfy the relations (1), (2) respectively.
\emph{Then if $Q$ is any element of finite order, so that}
\[
Q^q=1,\qquad (q\ne 0,1)
\]
\emph{then}
\[
Q=LT^a_i L^{-1}\quad\text{or}\quad Q=L{R^*_n}^a L^{-1}.
\]
\emph{On the other hand, no $T_i$ is the transform of an element
$T^a_l$ with $i\ne l$ or of an element ${R^*_n}^a$}. Assuming
this theorem is correct, then the $k_i$ are obviously determined by
the group structure. The proof of the theorem follows from the
solution of the word problem, which is carried out for planar groups
of the first kind with the restriction $k_i\ge 5$ in Sections 7.15 to
7.17, and for planar groups of the second kind in Section 2.6.

Using the relation to the covered complex, one can very easily
elucidate the structure of the diagrams of planar groups in normal form. 
In the case of
groups of the first kind, each segment bounds exactly two surface pieces.
Surrounding each point, in the normal form there are $n-1$ surface pieces
which correspond to the $n-1$ relations $T^{k_i}_i$ ($i=1,2,\ldots,n-1$), 
the latter being connected by a surface piece corresponding to the relation 
$R_n$. Between $T^{k_1}_1$ and $T^{k_n}_n$ there are $2g$
surface pieces corresponding to the relation $R_n$. The $n+2g-1$ surface 
pieces around each point corresponding to the relation $R_n$ correspond 
to the $n+2g-1$ different relations resulting from $R_n$ by cyclic
permutation. Thus one can construct the group diagrams of our groups
by geometric rules. 

\section{Properties of Elementarily Related Coverings}

If the coverings $\boldsymbol{A}(\mathfrak{M})=\mathfrak{M}^*$
and $\boldsymbol{A}'(\mathfrak{M}')={\mathfrak{M}^*}'$ are equivalent,
then the associated groups $\mathfrak{T,A,V,U}$ and
$\mathfrak{T',A',V',U'}$ are isomorphic. More precisely, there is an 
isomorphism between $\mathfrak{V}$ and $\mathfrak{V}'$ which maps the
subgroup $\mathfrak{U}$ of $\mathfrak{V}$ onto the subgroup
$\mathfrak{U}'$ of $\mathfrak{V}'$. This isomorphism alone suffices to
guarantee the equivalence. However, the equivalence of arbitrary
coverings may be reduced to the equivalence of planar coverings.

If $\boldsymbol{A}(\mathfrak{K})=\mathfrak{N}^*$ and
$\boldsymbol{A}'(\mathfrak{K}')={\mathfrak{N}^*}'$ are two planar
coverings connected by a chain of elementary transformations, then this
determines an isomorphism between their groups $\mathfrak{V}$ and
$\mathfrak{V}'$. If $T_i,S_i$ and $T'_i,S'_i$ are the generators of
$\mathfrak{V}$ and $\mathfrak{V}'$ respectively, and if
$\boldsymbol{I}(\mathfrak{V}')=\mathfrak{V}$ is the given isomorphism
then
\[
\boldsymbol{I}(T'_i)\quad\text{and}\quad \boldsymbol{I}(S'_i)
\]
are certain power products in the $T_i,S_i$. If the two complexes
$\mathfrak{K}$ and $\mathfrak{K}'$ are isomorphic as well, then the
relations between the $T_i,S_i$ are exactly the same as those between
$T'_i,S'_i$ and the transformation
$\overline{T}_i=\boldsymbol{I}(T_i), 
\overline{S}_i=\boldsymbol{I}(S_i)$
is consequently an automorphism of the group $\mathfrak{V}$. The 
collection of automorphisms effected by elementary transformations in
this way constitutes a group $\mathfrak{G}$.

\emph{Now if $\boldsymbol{A}(\mathfrak{M})=\mathfrak{M}^*$ and
$\boldsymbol{A}'(\mathfrak{M}')={\mathfrak{M}^*}'$ are two
coverings in normal form, and $\mathfrak{V,U,V',U'}$ are the
associated groups, then for the two coverings to be elementarily
related it is necessary and sufficient that under the mapping}
\[
S_i\rightarrow S'_i,\quad T_i\rightarrow T'_i
\]
\emph{of $\mathfrak{V}$ onto $\mathfrak{V}'$, $\mathfrak{U}$ is also
mapped onto $\mathfrak{U}'$, or else that there is an automorphism in
$\mathfrak{G}$ that maps $\mathfrak{V}'$ onto itself and carries
$\mathfrak{U}'$ to $\mathfrak{U}''$, where $\mathfrak{U}$ goes
to $\mathfrak{U}''$ under the mapping $S_i\rightarrow S'_i$,
$T_i\rightarrow T'_i$}.

Given two coverings for which the group $\mathfrak{V}$ is isomorphically
related to $\mathfrak{V}'$ by means of 
$\boldsymbol{I}_1(\mathfrak{V})=\mathfrak{V}'$, $\mathfrak{U}$ goes to
$\mathfrak{U}'$, and if also an isomorphism 
$\boldsymbol{I}_2(\mathfrak{V})=\mathfrak{V}'$ is effected between
$\mathfrak{V}$ and $\mathfrak{V}'$ by elementary transformations, under
which the subgroup $\mathfrak{U}$ goes to the subgroup $\mathfrak{U}''$,
then $\mathfrak{U}''$ results naturally from $\mathfrak{U}'$ by an
automorphism of $\mathfrak{V}'$, 
$\boldsymbol{A}(\mathfrak{V}')=\mathfrak{V}'$. However, this is not to 
say that $\boldsymbol{A}$ belongs to $\mathfrak{G}$. The question
whether any isomorphism $\boldsymbol{I}_1$ of $\mathfrak{V}$ onto
$\mathfrak{V}'$ guarantees the elementary relatedness of the associated
coverings comes down to the question of whether the automorphisms of
$\mathfrak{V}$ effected by elementary transformations are all the
automorphisms of $\mathfrak{V}$ or not. For unbranched coverings
$\mathfrak{G}$ is the full automorphism group (cf. Section 6.1),
otherwise this is only a tentative conjecture.

\section{Subgroups of Planar Groups}

Let $\mathfrak{K}$ be a planar covering complex with
$\boldsymbol{A}(\mathfrak{K})=\mathfrak{M}^*$,
\emph{let $\mathfrak{V}$ be the group of transformations of the complex
$\mathfrak{K}$ onto itself and let $\mathfrak{U}$ be any subgroup of
$\mathfrak{V}$. Then there is a complex $\mathfrak{M}$ and a regular 
covering $\boldsymbol{A}'(\mathfrak{K})=\mathfrak{M}$ such that two
elements of $\mathfrak{K}$ lie over the same element of $\mathfrak{M}$
if and only if they can be carried into each other by a transformation from
$\mathfrak{U}$}.

We can first construct the one-dimensional complex $\mathfrak{M}_1$
of $\mathfrak{K}$
and a covering $\boldsymbol{A}'(\mathfrak{K}_1)=\mathfrak{M}_1$
by Section 4.20, and then carry over the cyclic ordering of segments
through a point from $\mathfrak{K}_1$ to $\mathfrak{M}_1$, and hence
construct $\mathfrak{M}$ with the help of the boundary paths now
determined in $\mathfrak{M}_1$. Then there is also a covering
$\boldsymbol{A}''(\mathfrak{M})=\mathfrak{M}^*$. This likewise
follows as in Section 4.20. $\boldsymbol{A}''$ is a regular covering or not
according as $\mathfrak{U}$ is an invariant subgroup of $\mathfrak{V}$
or not. If $\boldsymbol{A}''$ is an unbranched covering, then so is
$\boldsymbol{A}'$; if $\boldsymbol{A}''$ is a branched covering then
$\boldsymbol{A}'$ can be branched or unbranched. Only in the latter
case is $\mathfrak{U}$ isomorphic to the fundamental group of 
$\mathfrak{M}$.

However, $\mathfrak{U}$ always possesses a planar group diagram. In
order to construct such a group diagram we construct a complete system of
fundamental domains of the complex $\mathfrak{K}_1$ modulo the group
$\mathfrak{U}$, as in Section 4.17. Let $\mathfrak{B}_1$ be a tree, the
points of which cannot be carried to each other by transformations in
$\mathfrak{U}$, but such that each point of $\mathfrak{K}$ may be carried
to a point of $\mathfrak{B}_1$ by such a transformation. Let
$\mathfrak{B}_2,\mathfrak{B}_3,\ldots$ be the trees that result from 
$\mathfrak{B}_1$ by transformations in $\mathfrak{U}$. Further, let
$\mathfrak{B}$ be the tree of $\mathfrak{M}$ for which
$\boldsymbol{A}'(\mathfrak{B}_1)=\mathfrak{B}$. If we now contract the 
trees $\mathfrak{B},\mathfrak{B}_i$ to single points $p,p_i$ we get two new 
complexes $\mathfrak{M}',\mathfrak{K}'$ and a covering
$\boldsymbol{A}'(\mathfrak{K}')=\mathfrak{M}'$, and $\mathfrak{K}'_1$
yields a group diagram of $\mathfrak{U}$; for the regularly covered
manifold $\mathfrak{M}'$ now contains only one point. $\mathfrak{U}$
is obviously the principal group of the covering
$\boldsymbol{A}'(\mathfrak{K})=\mathfrak{M}'$. It follows further that
subgroups of planar groups of the first (second) kind are themselves
planar groups of the first (second) kind.

\section{Branching Numbers of Subgroups}

One can obtain a simple criterion for deciding whether the paths in
$\mathfrak{K}'$ corresponding to the relations ${R^*_i}^{k_i}(S)$
in $\mathfrak{V}$ are branched or unbranched in the covering
$\boldsymbol{A}'(\mathfrak{K}')=\mathfrak{M}'$ from the behavior
of the residue classes of the subgroup $\mathfrak{U}$ of
$\mathfrak{V}$. Let
\[
L_1,\quad
L_2,\quad
\ldots,\quad L_n
\]
be  a complete system of representatives of the residue classes
$\mathfrak{U}L$ in $\mathfrak{V}$. The tree $\mathfrak{B}$ of
$\mathfrak{M}$ and the $L_i$ can be chosen in such a way that the
simple paths in the tree emanating from a point $p$ correspond uniquely 
to the power products $L_1,L_2,\ldots,L_n$ in the generators of
$\mathfrak{V}$. The same holds for the simple paths in the trees
$\mathfrak{B}_i$ emanating from the points $p_i$. If $V$ is any
element of $\mathfrak{V}$, then by $\overline{V}$ we understand the
representative of the residue class $\mathfrak{U}V$, as in Section 3.1.

Now let ${R^*}^k$ be one of the relations ${R^*_i}^{k_i}$.
Obviously
\[
\overline{VR^k}=\overline{V}
\]
because $R^k\equiv 1$. If $\overline{VR^l}=\overline{V}$ for a particular
$V$, then also
\[
\overline{VR^{2l}}=\overline{\overline{VR^l}R^l}
=\overline{\overline{V}R^l}=\overline{VR^l}=\overline{V},
\]
and hence in general $\overline{VR^{ml}}=\overline{V}$. It follows
that for each $V$ there is a smallest $l$ that divides $k$ and for which
$\overline{VR^l}=\overline{V}$. \emph{If $k=lm$ then we may put}
\[
\overline{V}R^k \overline{V}^{\,-1}
=(\overline{V}R^l\overline{VR^l}^{\,-1} 
         \overline{VR^l}R^l \overline{VR^{2l}}^{\,-1}\cdots
         \overline{VR^{(k-1)l}}^{-1}R^l \overline{V}^{\,-1})
=(\overline{V}R^l\overline{V}^{\,-1})^m.
\]
\emph{Then the path corresponding to this power product in
$\mathfrak{K}$ or $\mathfrak{K}'$ covers a path in $\mathfrak{M}'$
with branching order $m-1$}. it is clear that the order of branching is 
at least $m-1$. To see that it is not greater we must establish that the
power product that results from
\[
\overline{V}R^l\overline{V}^{-1},
\]
when we express it in terms of of the generators of $\mathfrak{U}$,
cannot be written as a formal power of $(R'')^a$. But now, in the case 
where $\mathfrak{V}$ has at least two defining relations in the normal
form (1) of Section 7.9, we can set $R=SR'$, where the generator $S$
no longer appears in $R'$. Then
\[
\overline{V}R^l \overline{V}^{\,-1}
=\prod^{l-1}_{i=0}\overline{VR^i}S\overline{VR^i S}^{\,-1}
                              \overline{VR^i}SR'\overline{VR^{i+1}}^{\,-1}
\]
and the
\[
\overline{VR^i}S\overline{VR^i S}^{\,-1}\qquad(i=0,1,\ldots,l-1)
\]
are formally different generators of $\mathfrak{U}$ by the definition
of $l$. One argues quite similarly when $\mathfrak{V}$ has only one
defining relation in normal form.

One sees at the same time that the branching numbers belonging to a
subgroup are divisors of the branching numbers of the group itself. If
\[
\overline{V{R^*_i}^l}\ne \overline{V}\qquad (l=1,2,\ldots,k_i-1)
\]
for all elements $V$ and all $R^*_i$, then the covering associated
with the subgroup $\mathfrak{U}$ is unbranched and so
$\mathfrak{U}$ is isomorphic to the fundamental group of the
manifold $\mathfrak{M}'$.

The paths $\overline{V}_1 R^k \overline{V}^{\,-1}_1$ and
$\overline{V}_2 R^k \overline{V}^{\,-1}_2$ can be covered with 
branching of different orders. However, if $\mathfrak{U}$ is an
invariant subgroup then the order is the same for all $V$. Namely,
if $\overline{V_1 R^l}=\overline{V}$ then $\overline{R^l}\equiv 1$ 
and hence in general $\overline{VR^l}=\overline{V}$.

Analogous theorems hold for coverings of group diagrams of the 
second kind, including one for subgroups of planar groups of the 
second kind. Thus subgroups of a free product of cyclic groups are
again free products of cyclic groups, as was proved purely group
theoretically in Section 3.9.  The subgroups for which all closed
boundary paths are covered without branching are free groups. The
criteria for unbranched covering of a path are the same.

\section{Automorphisms of Groups of Manifolds} 

In Section 7.2 we referred to the close connection between the
automorphisms induced by mappings $T$ of a regular covering
complex onto itself and the structure of the group $\mathfrak{V}$
belonging to the covering. One can use this connection to investigate
the properties of an individual automorphism. We will give a simple
example, confining ourselves to investigating automorphisms in the
case of a cyclic mapping group $\mathfrak{T}$, since each automorphism
realized by transformation of a regular covering complex onto itself
will be among those induced in cyclic mapping groups.

Accordingly, let $\mathfrak{V}$ be a group in the normal form
\[
T^{k_i}_i\equiv 1,\quad (T_1 T_2 \cdots 
T_{n-1} R(S))^{k_n}=R_n\equiv 1,
\]
let $\mathfrak{U}$ be an invariant subgroup of $\mathfrak{V}$ of
index $q$, isomorphic to the fundamental group of a manifold, and
let
\[
V^i\qquad (i=0,1,\ldots,q-1)
\]
be a complete set of representatives of $\mathfrak{U}$. Then if
\[
\overline{T_i}=V^{\alpha_i}\qquad 0\le\alpha_i <q
\]
in the notation of Section 3.1 we must have $\alpha_i\ne 0$,
otherwise the path corresponding to $T^{k_i}_i$ would be covered
branched. If the greatest common divisor of $\alpha_i$ and $q$ is
\[
(\alpha_i,q)=\delta_i\quad\text{and}\quad \alpha_i=\beta_i \delta_i,
\quad q=\gamma_i\delta_i,
\]
then $\overline{T^l}=\overline{T^{l'}}$ if and only if
\[
l\equiv l'\quad(\text{mod } \gamma_i).
\]
Since $\overline{T^{k_i}_i}=1$, $k_i\equiv 0$ (mod $\gamma_i$),
otherwise we would have $k_i\ne \gamma_i$ and the path 
corresponding to $T^{k_i}_i$ would again be covered branched.
The analogous result holds for $R_n$. Thus all $k_i$ are divisors of $q$.

We now focus on the case where $q$ is a prime number and so all
$k_i=q$. We can then choose $T^i_1$ ($i=0,1,\ldots,q-1$) as
representatives of the residue classes. If we denote the generators 
$T_i,S_i$ of the subgroup in the process by
\[
T_{ik}\qquad (i=1,2\ldots,n-1;k=0,1,\ldots,q-1)
\]
and $S_{ik}$, then the relations of the second kind say that
\[
T_{1k}\equiv 1\qquad (k=0,1,\ldots,q-2)
\]
and it follows from $T^q_1\equiv 1$ that $T_{1,q-1}\equiv 1$. Any of the $q$
relations that follow from $T^q_1\equiv 1$ and $R_n\equiv 1$ are cyclic 
interchanges of each other. The former reads
\begin{equation}
T_{i l_1} T_{i l_i} \cdots T_{i l_q}\equiv 1 \tag{1}
\end{equation}
where
\[
l_1,\quad l_2,\quad \ldots,\quad l_q
\]
is a permutation of the numbers $0,1,\ldots,q-1$. The latter relation
contains each $T_{ik}$ exactly once and each $S_{ik}$ once with
exponent $+1$ and once with exponent $-1$. Under the automorphism
in $\mathfrak{U}$ effected by $T_1$
\[
T_{ik},S_{ik} \text{ with }k<q-1\text{ go to }
T'_{ik}=T_{i,k+1}, S'_{ik}=S_{i,k+1}
\]
and 
\[
S_{i,q-1}, T_{i,q-1} \text{ go to } S'_{i,q-1}=S_{i0}, T'_{i,q-1}=T_{i0}.
\]

We now ask about the form of the automorphism induced in the factor
group $\mathfrak{F}$ of $\mathfrak{U}$ by the commutator group by
the automorphism effected by transformation $T_1$. $\mathfrak{F}$
is a free commutative group. The relations can now in fact be solved,
since (1) makes it possible to express $T_{i,q-1}$ in terms of the 
remaining $T_{ik}$:
\[
T_{i,q-1}=T^{-1}_{i0} T^{-1}_{i1}\cdots T^{-1}_{i,q-2}.
\]
On the other hand the last relation, resulting from $R_n$, becomes a
consequence relation since the $S_{ik}$ cancel out. Accordingly, we
take $T_{ik},S_{ik}$ with $k<q-1$ as generators of $\mathfrak{F}$
and construct the matrix that expresses the passage from
$T_{ik},S_{ik}$ to $T'_{ik},S'_{ik}$. The latter has, on its main 
diagonal, $-1$ in the rows corresponding to the $T'_{i,q-2}$ and 0 
elsewhere. \emph{The trace of the matrix is therefore equal to
$-n+1$}. This is a special case of a general theorem proved by other
methods by \textsc{Nielsen}.\footnote{\textsc{J. Nielsen}, Acta Math.
\textbf{50} (1927) 189.}

One sees that \emph{the automorphisms in $\mathfrak{F}$ induced
by transformation with $T$ are uniquely determined by the numbers
$q$, $n$, and $g$}. One can ask for a complete characterization of
the automorphisms induced in $\mathfrak{U}$ itself. For this purpose
one has to consider on the one hand the different automorphisms of
finite order $q$ which the elements of $\mathfrak{V}$ induce in
$\mathfrak{U}$ by transformation. For in general $V^q$ brings
about an inner automorphism of $\mathfrak{U}$ which is not the
identity. It is not hard to see that those elements $V$ for which
$V\mathfrak{U}V^{-1}=\mathfrak{U}'$ is of finite order are
exactly the elements of $\mathfrak{V}$ of finite order, and that the
finite order elements of $\mathfrak{V}$ are exactly the elements
\[
LT^l_i L^{-1},\quad LR^l_n L^{-1}.
\]
On the other hand, one has to investigate the totality of invariant
subgroups $\mathfrak{U}$ of index $q$ in $\mathfrak{V}$, and their
equivalence under the automorphism group $\mathfrak{G}$ of
$\mathfrak{V}$ determined by elementary transformations.

\section{The Word Problem for Planar Groups}

In conclusion we will consider the word problem for planar group 
diagrams.\footnote{Cf. \textsc{M. Dehn}, Math. Ann. \textbf{72}
(1912) 413.} To the extent that one can search for each word as a
path in the group diagram, and construct the group diagram in the
plane by laying down successive polygons, one can regard the word
problem as already solved. However, this remark alone gives no
insight into how one can decide the word problem, since one does not
yet know which words are equal to the identity. However, the word
problem can also be solved in this strict sense for many of these
groups.

The fundamental group of the sphere is the identity, that of the 
projective plane is the cyclic group of order 2, that of the torus is the
free \textsc{Abelian} group with two generators. Since the groups
corresponding to planar group diagrams of the second kind are free
products of cyclic groups there are no difficulties in all these cases.

We now take the fundamental group $\mathfrak{W}$ of an orientable
manifold of genus $g>1$ and a system of $2g$ generators $S_i$ and
the relations $R(S)$ that correspond to boundary paths of a normal
form. In $R$ there are two generators $S_1,S_2$ that mutually
separate each other. If we also set
\[
S_i\equiv 1\qquad (i\ne 1)
\]
then a cyclic factor group results. Thus there is an invariant subgroup
$\mathfrak{I}$ of $\mathfrak{W}$, the residue classes of which are
represented by
\[
S^k_1\qquad(k=0,\pm 1,\ldots).
\]
We claim that \emph{$\mathfrak{I}$ is a free group with infinitely many
generators}. By applying the process we obtain generators
\[
S_{ik}=S^k S_i S^{-k}
\]
and $S_{1k}=1$ as relations of the second kind. If we express
$S^l R S^{-l}$ in terms of the $S_{ik}$ as prescribed then an
$S_{2,a+l}$ and an $S_{2,a+l+1}$ always appear in the latter only
once, with exponents $\pm 1$. Thus one can solve all relations by
again retaining $S_{20}$ as generator and eliminating all the remaining
$S_{2k}$. The result is a free group, with generators
\[
S_{20},S_{ik}\qquad (i=3,4,\ldots,2g;k=0,\pm 1,\ldots).
\]

Since each element of $\mathfrak{W}$ may be brought into the form
$S^l_1 S$, where $S$ belongs to $\mathfrak{I}$, and since the word
problem for $\mathfrak{I}$ may be solved as in Section 2.3, the word 
problem for $\mathfrak{W}$ is also solved.

One can solve the word problem for planar groups of the first kind by
similar considerations, whereby one seeks an invariant subgroup
isomorphic to the fundamental group of an orientable manifold.

The same ideas can be applied to handle the groups $\mathfrak{W}'$
of non-orientable manifolds. Here it is useful to begin with a normal form
with boundary path $r(s)=s^2_1 s^2_2 \cdots s^2_g$. One constructs
the invariant subgroup $\mathfrak{U}$ corresponding to the factor group
$S_1=S_i$ ($i=2,3,\ldots,g$), $S^2_1=1$. If one takes the identity 
element and $S_1$ as representatives of the residue classes modulo
$\mathfrak{U}$ and
\[
S_{ik}\qquad (i=1,2,\ldots,g;k=0,1)
\]
as generators of $\mathfrak{U}$, then $S_{10}=1$ on the basis of the
relations of the second kind. $R(S)$ implies two relations which both
contain all $S_{ik}$ ($i\ge 2$) and $S_{11}$ exactly once with
exponent $+1$. By elimination of $S_{11}$ one obtains a relation
between the $S_{ik}$ which contains each $S_{ik}$ once with exponent
$+1$ and once with exponent $-1$. As a subgroup of a manifold group,
$\mathfrak{U}$ is also the group of a manifold $\mathfrak{M}$, and
$\mathfrak{M}$ is orientable, as one sees from the form of the relations
for $\mathfrak{U}$. Consequently, the  word problem is solved in
$\mathfrak{U}$ and hence also in $\mathfrak{W}'$.

\section{Word Problems in Planar Group Diagrams}

A sharper theorem is the following one on planar group diagrams of the
first kind, the relations of which have length at least 4, and for which there
is consequently a normal form with $k_i\ge 4$ ($i=1,2,\ldots,n-1$), or 
which possess only the one relation ${R^*}^k$. Thus these groups
include the groups of orientable manifolds. Under the hypotheses just
given we have:

\emph{If $W$ is a word that represents the identity element then $W$
contains a subword $W^*$ which can be made into a defining relation
or its inverse by inserting suitable factors $S^{\pm 1}_\alpha$, 
$S^{\pm 1}_\beta$}. If the length of all defining relations is greater
than 4 then it is clear that the word problem can be solved on the basis
of this theorem.\footnote{Because extending $W^*$ to a defining
relation, then canceling this defining relation, causes a net decrease in the
length of the word. This is essentially \emph{Dehn's algorithm}, found
in the paper of \textsc{M. Dehn}, Math. Ann. \textbf{72}, (1912).
(Translator's note.)}

If $w$ is a path corresponding to a word $W$ in the group diagram,
then $w$ contains a subpath $\overline{w}$ which is a simple closed
path and hence bounds a regular surface piece. If the theorem is proved
for the words $W$ that correspond to simpler closed paths $w$, then it
holds in general. We prove the claim for these words by induction, with
the help of the following theorems.

Theorem 1. \emph{If $w$ is a simple closed path of the group diagram
which is not a simple boundary path, then two other simple closed paths
$w_1$ and $w_2$ may always be given which run through at most two 
common segments and, when $w_1,w_2$ are combined and reduced,
they yield a path which results from $w$ by cyclic interchange}.

Theorem 2. \emph{If $w$ is a simple closed path, then there are at least
four subpaths $r$ (called critical subpaths) of $w$ which may be converted
to simple boundary paths of surface pieces in the interior of $w$ by the 
insertion of two segments. If $w$ has exactly $k$ subpaths $S_1 S_2$
which appear in the boundary paths of surface pieces not in the interior of
$w$---``re-entrant vertices'' as we will say---then $w$ has at least $4+k$
critical subpaths}.

\section{Re-entrant Vertices and Critical Subpaths}

Assuming Theorem 1, we prove Theorem 2 by complete induction. It is
convenient in the following proof to imagine a Euclidean net of squares.
Theorem 2 holds for simple boundary paths and for simple paths that
contain exactly two surface pieces in their interior. This is
because boundary paths of two surface pieces which meet along a segment
have only this one segment in common. If one takes two segments out of
a path $w$ satisfying Theorem 2 then at most two critical subpaths are
destroyed, so we will suppose that $w$ is a simple boundary path. Because
when two critical subpaths have segments in common they bound the 
same surface piece $f$; $w$ is then either identical with the boundary of 
$f$ or else $w$ runs through all boundary segments of $f$ except one.
If two critical subpaths are destroyed by taking segments out, then $s$ is
an inner segment of a critical subpath.

It follows that: \emph{if $w_1$ and $w_2$ are two paths satisfying the
statement of Theorem 2 and if $w=w'_1 w'_2$ is a simple path
resulting from $w_1 w_2$ by omission of a segment, then $w$ also satisfies
our theorem}. This is clear if $w$ contains no new re-entrant vertices.
For if $k_i$ is the number of re-entrant vertices of $w_i$, then the number
of critical subpaths in $w'_i$ is at least $2+k_i$. Suppose now that, say, the
endpoint of $w'_1$ is a re-entrant vertex of $w$. Then the segment eliminated 
by reduction is certainly not at the same time an inner segment of a critical
subpath of $w_1$ and $w_2$, so $w'_1$ or $w'_2$ still contains at least
$3+k_1$ or $3+k_2$ critical subpaths. If the final point of $w'_2$ is also a
re-entrant vertex then the segment removed by reduction again appears in
no critical subpath of $w_1$ and $w_2$, or in only one from $w_1$ or
$w_2$. Thus the number of critical subpaths contained in the $w'_i$ is
either $4+k_i$ or $3+k_i$.

\emph{We now assume that two segments are eliminated from each of
$w_1,w_2$ in the reduction of $w=w'_1 w'_2$}.

1. Either $w_1$ loses a re-entrant vertex in the process, in which case at
most two critical subpaths are destroyed by omission of the two end
segments of $w_1$ and $w'_1$ has at least $2+k_1$ critical subpaths.

1a. Either $w_2$ is now a simple boundary path, so that $w=w'_1 w'_2$
has at least  $2+k_1+1$ critical subpaths and $k_1-1$ re-entrant vertices
when new re-entrant vertices do not result from from the join. But if this
is the case, fewer critical subpaths of $w_1$ can be destroyed by leaving
out the reduced segments.

1b. Otherwise $w_2$ is not a simple boundary path, so $w'_2$ contains
at least $2+k_2$ critical subpaths and the theorem remains correct in case at
most one re-entrant vertex results, because $w$ contains $4+k_1+k_2$
critical subpaths. If, on the other hand, two re-entrant vertices result, then either
there are $3+k_1$ critical subpaths still in $w'_1$ or, if that is not the case, the 
two segments in $w_2$ removed by reduction either constitute a critical
subpath or else they belong to no such path at all. Thus the theorem also
holds in this case.

2. Finally we suppose that $w_1$ and $w_2$ lose no re-entrant vertex. If 
no re-entrant vertex results, then the theorem is correct. However, one
re-entrant vertex may result and $w'_1$ may really have two critical 
subpaths fewer than $w_1$. Then at most one critical subpath is destroyed
in $w_2$. If two re-entrant vertices result and $w'_1$ contains two 
critical subpaths fewer than $w_1$ then no critical subpaths at all are
destroyed in $w_2$; and if $w'_1$ contains one critical subpath fewer than 
$w_1$ then in $w_2$ at most one is destroyed.

Thus the proof of Theorem 2 has been completed, assuming Theorem 1.

\section{Simple Paths in Planar Complexes}

We now turn to the proof of Theorem 1. \emph{If $w$ is a simple closed path}
and $w$ contains no point in its interior, then the assertion is correct.
\emph{If $w$ contains points then we have to show that there is at least one
point among these from which two segments lead to points of $w$}. Suppose
this is not the case. We construct the complex $\mathfrak{C}$ consisting of
the points in the interior of $w$ and all the segments connecting these points
with each other. The points of $\mathfrak{C}$ then have order $2m$ or
$2m-1$.

Now we construct the ``boundary'' of the complex $\mathfrak{C}$ in the
following way. If $p$ is a point of order $2m-1$ and $s$ is the segment from
$p$ to a point of $w$, then two surface pieces meet along $s$, with positive 
simple boundary paths $r_1 s$ and $s^{-1}r_2$. The $r_i$ each run through 
a second point of order $2m-1$ and we can therefore set
\[
r_i=r_{i1}r_{i2}
\]
and conclude that $r_{12}$ and $r_{21}$ only run through segments of
$\mathfrak{C}$ that begin and end at points of order $2m-1$, and no other
such points. These paths $r_{12}$ and $r_{21}$ may be called the
``auxiliary paths'' bounded by $p$. Different auxiliary paths, being subpaths
of simple boundary paths, cannot intersect. By linking up the auxiliary paths
we put together reduced closed paths $w_i$ that may be called ``boundary
paths of $\mathfrak{C}$.'' All segments of $\mathfrak{C}$ branching off a
boundary path lie on the same side of it, say the negative side. Each path
$w_i$ contains a simple closed subpath $w_{i0}$ which begins and ends at
$p_0$ say. All points lying in the interior of $w_{i0}$ belong to $\mathfrak{C}$
and are of order $2m$. Consequently, the complex $\mathfrak{C}_{i0}$
consisting of the points and line segments of $w_{i0}$ and the points, line
segments, and surface pieces in its interior has the following properties.
\begin{enumerate}
\item
All points in the interior are of order $2m$ with $m>1$.
\item
All points on the boundary different from $p_0$ are of order $2m$ or $2m-1$.
\item
The simple boundary paths of the surface pieces in $\mathfrak{C}_{i0}$ run
through at least four line segments.
\item
The complex $\mathfrak{C}_{i0}$ is topologically equivalent to a surface piece.
\end{enumerate}
One easily verifies that such a complex cannot exist.

Let $a_{01}$ be the number of boundary points distinct from $p_0$, let
$a_{02}$ be the number of interior points, and let $2a_1$ and $2a_2$ 
respectively be the numbers of line segments  and surface pieces of 
$\mathfrak{C}_{i0}$. Then on the one hand
\begin{equation}
2a_1\ge (2m-1)a_{01}+2ma_{02}+2, \tag{1}
\end{equation}
and on the other hand
\begin{equation}
2a_1\ge 4a_2+a_{01}+1. \tag{2}
\end{equation}
Also, by Section 5.3,
\begin{equation}
a_1=a_{01}+a_{02}+a_2. \tag{3}
\end{equation}
Hence, by addition of (1) and (2) and subtraction of 4 times (3) we obtain
\[
0\ge (2m-4)(a_{01}+a_{02})+3,
\]
which is a contradiction.

\section{Planar Group Diagrams and Non-Euclidean Geometry}

The fact that the solution of the word problem in the preceding section depends on
certain restrictions has its origin in various properties of the groups themselves.
The groups with the relations
\begin{align*}
S^n&=1, \tag{1}\\
S^2&=T^2=(ST)^n=1, \tag{2}\\
S^2&=T^3=(ST)^k=1\quad(k=3,4,5) \tag{3}
\end{align*}
are finite groups, otherwise [for $k>5$] infinite; (2) is called the dihedral group,
(3) are called the tetrahedral, octahedral, and icosahedral groups for $k=3,4,5$
respectively. These groups may be represented by rotations of the sphere or motions
of spherical geometry, as is well known. The groups
\begin{align*}
STS^{-1}T^{-1}&=1, \tag{4}\\
S^2&=T^3=(ST)^6=1, \tag{5}\\
S^2&=T^4=(ST)^4=1, \tag{6}\\
S^2&=T^2=U^2=(STU)^2=1 \tag{7}
\end{align*}may be represented by motions of Euclidean geometry, and the 
remaining plane groups by motions of non-Euclidean geometry.

The proof of this is best obtained in connection with point-type branched covering 
complexes. The surface pieces of the complex, which are transitively exchanged with 
each other by the group of the complex, may be realized as polygons bounded by
straight lines, and hence the whole complex may be realized by a rectilinear
polygonal net. The angles of these polygons may be derived from the branching 
numbers of the group; e.g., the fundamental domain for (4) is the parallelogram,
that for (5) an isosceles triangle with top angle 120$^\circ$, that for
(6) an isosceles right-angled triangle, and for (7) an equilateral triangle.
The corresponding groups of Euclidean motions are in the case of (4) just 
translations; the others contain rotations that leave a vertex of the fundamental domain fixed and rotate through the angle of the 
fundamental domain or a multiple of it,
and in the case of (5) and (6) also rotations through 180$^\circ$ about the 
midpoint of the base of the fundamental domain and, in case of (7), 180$^\circ$
rotations about the midpoints of the three sides of the triangle.

The elements of finite order $o$ in the group also correspond to rotations 
through an angle $2\pi/o$ in the case of representation by motions of
non-Euclidean geometry. One can use this connection to prove that all
elements of finite order in a planar group are transforms of $T_i$ or
$R^*_n$ in the notation of Section 7.9 (1).

The connection with non-Euclidean geometry, discovered by \textsc{Poincar\'e}, 
has been used in many investigations, e.g., in the important work of
\textsc{Nielsen}\footnote{\textsc{J. Nielsen}, Acta Math. \textbf{50} (1927) 189
and \textbf{53} (1929) 1; \textsc{H. Gieseking} Analytische Untersuchungen 
\"uber topologische Gruppen (Dissertation) M\"unster 1912.} 
on the mappings of surfaces and the fixed point problem. Conversely, the theory
of planar groups has also been of help in the investigation of these groups of
motions,\footnote{\textsc{E. Hecke}, Hamb. Abhdl. \textbf{8} (1930) 271 and
\textsc{H. Rademacher} Ibid. \textbf{7} (1929) 134.}
which are of great importance because of their connection with automorphic
functions and uniformization.

Finally in this connection we refer to the interesting questions raised and
answered by \textsc{Steinitz},\footnote{\textsc{E. Steinitz} ``Polyeder und
Raumeinteilung,'' Math. Enzykl. III, AD 12.} on the topological
chacterization of complexes equivalent to the sphere that may be realized as
convex polyhedra in Euclidean space.

\end{document}